\documentclass[a4paper,11pt,reqno]{amsart} 

\usepackage{hyperref}

\usepackage[top=3cm,left=2.5cm,right=2.5cm,bottom=3cm]{geometry}
\usepackage{relsize}
\usepackage{lineno}

\usepackage{amsmath,amsthm,pb-diagram,amssymb,comment}
\usepackage{amsfonts,graphicx,color}
\usepackage[inline]{enumitem} 
\usepackage{fancyhdr, multicol}
\usepackage{thmtools,cleveref}
\usepackage[normalem]{ulem}
\usepackage{stmaryrd}
\usepackage{textcomp}
\usepackage{mathtools}
\usepackage{fontawesome}
\usepackage{textcomp}

\usepackage{tikz,float}
\usepackage{mleftright}
\usepackage{mathrsfs}
\usepackage{tikz-cd}
\usepackage{dsfont}

\definecolor{bluet}{rgb}{0.0, 0.4, 0.6}

\hypersetup{
    colorlinks=true,
    linkcolor=bluet,
    filecolor=bluet,
    urlcolor=bluet,
    citecolor=bluet,
}
\definecolor{caribeangreen}{rgb}{0.0, 0.2, 0.6}

\setlist[enumerate,itemize,1]{topsep=0ex,itemsep=1ex}
\setlist[enumerate,itemize,2]{topsep=1ex,itemsep=1ex}
    \newcommand{\Cv}{\mathrm{Cv}}
    \newcommand{\Ecal}{\mathcal{E}}
    \newcommand{\NE}{\mathrm{NE}}
    \newcommand{\sqce}{\sqsubset^{\mathrm ce}}
    
    \DeclareMathOperator{\dom}{dom}

    \DeclareMathOperator{\ran}{ran}

    \newcommand{\thzfc}{\mathrm{ZFC}}

    \newcommand{\Bwf}{\mathcal{B}}

    \newcommand{\Ewf}{\mathcal{E}}

    \newcommand{\Iwf}{\mathcal{I}}
    \newcommand{\Jwf}{\mathcal{J}}
    
    \newcommand{\Mwf}{\mathcal{M}}
    \newcommand{\Nwf}{\mathcal{N}}
    
    \newcommand{\Pwf}{\mathcal{P}}
    \newcommand{\Swf}{\mathcal{S}}

    \newcommand{\bfrak}{\mathfrak{b}}
    \newcommand{\cfrak}{\mathfrak{c}}
    \newcommand{\dfrak}{\mathfrak{d}}

    \DeclareMathOperator{\pts}{\mathcal{P}}
    \newcommand{\sminus}{\smallsetminus}
    \newcommand{\menos}{\smallsetminus} 

    \newcommand{\frestr}{{\upharpoonright}}
    \DeclareMathOperator{\add}{add}
    \DeclareMathOperator{\cov}{cov}
    \DeclareMathOperator{\non}{non}
    \DeclareMathOperator{\cof}{cof}

    \newcommand{\Bor}{\mathbb{B}}
    \newcommand{\Cor}{\mathbb{C}}
    \newcommand{\Dor}{\mathbb{D}}
    \newcommand{\Eor}{\mathbb{E}}

    \newcommand{\Loc}{\mathbb{LOC}}
    
    \newcommand{\Por}{\mathbb{P}}
    \newcommand{\Pbb}{\mathbb{P}}
    \newcommand{\Qor}{\mathbb{Q}}

    \newcommand{\Qnm}{\dot{\mathbb{Q}}}

    \newcommand{\Ncal}{\mathcal{N}}
    \newcommand{\Mcal}{\mathcal{M}}

    \newcommand{\Q}{\mathds{Q}}
    \newcommand{\R}{\mathds{R}}

    \DeclareMathOperator{\cf}{cf}

    \newcommand{\la}{\langle}
    \newcommand{\ra}{\rangle}

\newcommand{\Ibf}{\mathbf{I}}
\newcommand{\fin}{\mathrm{fin}}
\newcommand{\varp}{\varepsilon}

\DeclareMathOperator{\limit}{\mathrm{lim}}

\newcommand{\mbf}{\mathbf{m}}

\newcommand{\tbf}{\mathbf{t}}

\newcommand{\Fr}{\mathrm{Fr}}

\newcommand{\Cn}{\mathrm{Cn}}
\newcommand{\Lc}{\mathrm{Lc}}

\newcommand{\Hcal}{\mathcal{H}}
\newcommand{\Scal}{\mathcal{S}}
\newcommand{\id}{\mathrm{id}}
\newcommand{\Ed}{\mathrm{Ed}}

\DeclareMathOperator{\hgt}{\mathrm{ht}}

\newcommand{\leqT}{\mathrel{\mbox{$\preceq_{\mathrm{T}}$}}}
\newcommand{\eqT}{\mathrel{\mbox{$\cong_{\mathrm{T}}$}}}

\newcommand{\baire}{{}^\omega\omega}

\newcommand{\Ce}{\mathrm{Ce}}

\newcommand{\set}[2]{\left\{#1 \colon #2\right\}}
\newcommand{\setcolon}[2]{\{#1 \colon #2\}}
\newcommand{\Seq}[2]{\la #1 \colon #2\ra}
\newcommand{\largeset}[2]{\left\{#1 \colon #2\right\}}

\newcommand{\Fn}{\mathrm{Fn}}
\newcommand{\gen}{\mathrm{gn}}

\newcommand{\cantor}{{}^\omega2}

\newcommand{\Pbf}{\mathbf{P}}

\newcommand{\nsqsubset}{\not\sqsubset}

\newcommand{\redq}[1]{{\color{red} #1}}

\newcommand{\greenunq}[1]{{\color{greenun} #1}} 
\newcommand{\redunq}[1]{{\color{redun} #1}}

\newcommand{\FAM}{\mathrm{FAM}}
\newcommand{\ZFC}{\mathrm{ZFC}}

\newcommand{\Tree}{\mathrm{Tr}}

\newcommand{\CH}{\mathrm{CH}}

\newcommand{\Int}{\mathrm{int}}

\newcommand{\Ber}{\mathrm{Ber}}

\newcommand{\p}{\mathrm{Pr}}
\newcommand{\B}{\mathrm{B}}

\newcommand{\Card}{\mathrm{Card}}
\newcommand{\Ord}{\mathrm{Ord}}

\newcommand{\E}{\mathrm{E}}
\newcommand{\Var}{\mathrm{Var}}

\newcommand{\Leb}{\mathrm{Leb}}

\newcommand{\forces}{\Vdash}

\newcommand{\Bin}{\mathrm{Bin}}

\newcommand{\Lev}{\mathrm{Lev}}

\newcommand{\At}{\mathrm{At}}

\newcommand{\Be}{\mathcal{B}\boldsymbol{e}}

\newcommand{\suc}{\mathrm{succ}}

\newcommand{\bfb}{\mathbf{b}}
\newcommand{\bfc}{\mathbf{c}}

\newcommand{\bfI}{\mathbf{I}}
\newcommand{\bfi}{\mathbf{i}}

\newcommand{\bfm}{\mathbf{m}}

\newcommand{\bfp}{\mathbf{p}}

\newcommand{\bft}{\mathbf{t}}

\newcommand{\bbB}{\mathbb{B}}
\newcommand{\bbC}{\mathbb{C}}
\newcommand{\bbD}{\mathbb{D}}
\newcommand{\bbE}{\mathbb{E}}

\newcommand{\bbL}{\mathbb{L}}
\newcommand{\bbN}{\mathbb{N}}
\newcommand{\bbP}{\mathbb{P}}
\newcommand{\bbQ}{\mathbb{Q}}
\newcommand{\bbR}{\mathbb{R}}

\newcommand{\cB}{\mathscr{B}}
\newcommand{\cC}{\mathscr{C}}

\newcommand{\cI}{\mathscr{I}}

\newcommand{\gb}{\mathfrak{b}}

\newcommand{\gc}{\mathfrak{c}}
\newcommand{\gd}{\mathfrak{d}}

\newcommand{\calT}{\mathcal{T}}
\newcommand{\calC}{\mathcal{C}}

\newcommand{\calA}{\mathcal{A}}

\newcommand{\calB}{\mathcal{B}}
\newcommand{\calP}{\mathcal{P}}
\newcommand{\calS}{\mathcal{S}}
\newcommand{\calN}{\mathcal{N}}

\newcommand{\calM}{\mathcal{M}}
\newcommand{\calG}{\mathcal{G}}
\newcommand{\calL}{\mathcal{L}}

\newcommand{\calI}{\mathcal{I}}
\newcommand{\calF}{\mathcal{F}}

\newcommand{\calY}{\mathcal{Y}}

\newcommand{\calH}{\mathcal{H}}

\newcommand{\finseqb}{{}^{ < \omega} \omega}
\newcommand{\finseq}{{}^{ < \omega} 2}
\newcommand{\finseqz}{{}^{ < \omega} Z}

\newcommand{\rest}{{\restriction}}

\newenvironment{PROOF}[2][\proofname.]
   {\begin{proof}[#1]\renewcommand{\qedsymbol}{\textsquare$_{\mathrm{#2}}$}}
   {\end{proof}}

\newcommand{\RRR}{\mathds R}
\newcommand{\QQQ}{\mathds Q}
\newcommand{\ZZZ}{\mathds Z}
\newcommand{\concat}[2]{#1{}^{\frown}#2}
\newcommand{\power}[2]{{}^{#1} \!#2} 

\renewcommand{\setminus}{\smallsetminus}
\newcommand{\conj}{\mathrel{\mbox{\scriptsize $\wedge$}}}

\DeclareMathOperator{\trunk}{trunk}

\DeclareMathOperator{\otp}{otp}

\definecolor{sub0}{RGB}{29,32,137}
\definecolor{sub1}{RGB}{1,71,157}
\definecolor{sub2}{RGB}{1,104,183}
\definecolor{sub3}{RGB}{0,160,234}
\definecolor{sug}{RGB}{0,154,68}
\definecolor{suy}{RGB}{208,219,1}

\definecolor{redun}{rgb}{0.0, 0.4, 0.6}
\definecolor{greenun}{rgb}{0.0, 0.4, 0.6}
\definecolor{dodger}{rgb}{0.0,0.5,1.0}
\definecolor{carrotorange}{rgb}{0.93, 0.57, 0.13}

\newcommand\subsetdot{\mathrel{\ooalign{$\subset$\cr
  \hidewidth\hbox{$\cdot\mkern3mu$}\cr}}} 

\newcommand{\Andres}[1]{{\color{violet}Andres says: #1}}
\newcommand{\Diego}[1]{{\color{dodger}Diego says: #1}}
\newcommand{\Miguel}[1]{{\color{carrotorange}Miguel says: #1}}

\title[A general theory of iterated forcing using fams]{A general theory of iterated forcing\\ using finitely additive measures}

\author{Miguel A.~Cardona}
\address{Einstein Institute of Mathematics,
The Hebrew University of Jerusalem. Givat Ram, Jerusalem, 9190401, Israel}
\email{\href{mailto:miguel.cardona@mail.huji.ac.il}{miguel.cardona@mail.huji.ac.il}}
\urladdr{\url{https://sites.google.com/mail.huji.ac.il/miguel-cardona-montoya/home-page}}

\author{Diego A.~Mej\'ia}
\address{Graduate School of System Informatics, Kobe University. 1-1 Rokkodai-cho, Nada-ku, Kobe, Hyogo 657-8501 Japan}
\email{\href{mailto:damejiag@people.kobe-u.ac.jp}{damejiag@people.kobe-u.ac.jp}}
\urladdr{\url{http://www.researchgate.com/profile/Diego\_Mejia2}}

\author{Andr\'es F. Uribe-Zapata}
\address{TU Wien, Faculty of Mathematics and Geoinformation, Institute of Discrete Mathematics and Geometry, Wiedner Hauptstrasse 8--10, A--1040 Vienna, Austria }
\email{\href{mailto:andres.zapata@tuwien.ac.at}{andres.zapata@tuwien.ac.at}}
\urladdr{\url{https://sites.google.com/view/andres-uribe-afuz}}


\thanks{The first author was partially supported by the Slovak Research and Development Agency under Contract No.~APVV-20-0045 and by Pavol Jozef \v{S}af\'arik University in Ko\v{s}ice at a postdoctoral position; the second author was supported by the Grant-in-Aid for Scientific Research (C) 23K03198, Japan Society for the Promotion of Science; and the third author was supported by the Austrian Science Fund (FWF): project number P33895. Also, the third author would like to sincerely thank the second author for the wonderful, patient, and dedicated advice on his master's thesis, which allowed the subsequent development of this work.} 

\date{\today}

\makeatletter
\def\l@subsection{\@tocline{2}{0pt}{2.5pc}{5pc}{}}
\makeatother

\begin{document}

\makeatletter
\def\@roman#1{\romannumeral #1}
\makeatother

\newcounter{enuAlph}
\renewcommand{\theenuAlph}{\Alph{enuAlph}}

\theoremstyle{plain}
  \newtheorem{theorem}{Theorem}[section]
  \newtheorem{corollary}[theorem]{Corollary}
  \newtheorem{lemma}[theorem]{Lemma}
  \newtheorem{mainlemma}[theorem]{Main Lemma}
  \newtheorem{mainproblem}[theorem]{Main Problem}
  \newtheorem{construction}[theorem]{Construction}
  \newtheorem{prop}[theorem]{Proposition}
  \newtheorem{clm}[theorem]{Claim}
  \newtheorem{fact}[theorem]{Fact}
  \newtheorem{exer}[theorem]{Exercise}
  \newtheorem{question}[theorem]{Question}
  \newtheorem{problem}[theorem]{Problem}
  \newtheorem{cruciallem}[theorem]{Crucial Lemma}
  \newtheorem{conjecture}[theorem]{Conjecture}
  \newtheorem{assumption}[theorem]{Assumption}
  \newtheorem{teorema}[enuAlph]{Theorem}
  \newtheorem{mainLemma}[enuAlph]{Main Lemma} 
  \newtheorem*{corolario}{Corollary}
  \newtheorem*{extthmsuc}{Extension theorem at successor steps}
  \newtheorem*{extthmlim}{Extension theorem at limit steps}
\theoremstyle{definition}
  \newtheorem{definition}[theorem]{Definition}
  \newtheorem{example}[theorem]{Example}
  \newtheorem{remark}[theorem]{Remark}
    \newtheorem{hremark}[theorem]{Historical Remark}
    \newtheorem{observation}[theorem]{Observation}
  \newtheorem{notation}[theorem]{Notation}
  \newtheorem{context}[theorem]{Context}

  \newtheorem*{defi}{Definition}
  \newtheorem*{acknowledgments}{Acknowledgments}

\numberwithin{equation}{theorem}
\renewcommand{\theequation}{\thetheorem.\arabic{equation}}

\def\sectionautorefname{Section}
\def\subsectionautorefname{Subsection}

\begin{abstract}
Based on the work of Shelah, Kellner, and T\u{a}nasie (Fund.\ Math., 166(1--2):109--
136, 2000 and Comment.\ Math.\ Univ.\ Carolin., 60(1):61--95, 2019), and the recent developments in the third author's master's thesis, we develop a general theory of iterated forcing using finitely additive measures. For this purpose, we introduce two new notions: on the one hand, we define a new linkedness property, called \emph{$\mu$-$\mathrm{FAM}$-linked} and, on the other hand, we generalize the notion of \emph{intersection number} to forcing notions, which justifies the limit steps of our iteration theory. Our theory also generalizes iterations with ultrafilters, which have played an important role in the proof of the consistency of Cicho\'n's maximum. 
    
We further show that any iteration constructed with our theory preserves strong unbounded families and what we call \emph{anti-Bendixson families}, which play a central role in preserving witnesses of $\mathrm{cov}(\mathcal{N})$ of singular size (even of countable cofinality). 
We also show that our iteration method does not increase $\mathrm{non}(\mathcal{E})$, the smallest size of a set of reals that cannot be covered by an $F_\sigma$ measure zero set. 
Finally, we apply our theory to prove a new separation of the left-hand side of Cicho\'n's diagram where $\mathrm{cov}(\mathcal{N})$ is possibly singular, even with countable cofinality.
\end{abstract}

\maketitle

{
\small
\hypersetup{linkcolor=black}
\tableofcontents
}

\section{Introduction}\label{SecIntro}

As a result of Kurt Gödel's work (see \cite{GodelI}), it is known that there are statements that cannot be proven true or false within $\ZFC$, the axiomatic system in which most modern mathematics is formalized. These statements are commonly referred to as \emph{independent} or \emph{undecidable} from $\ZFC$. Although several examples of independent statements come from different areas of mathematics, such as the \emph{Whitehead Problem} in group theory (see \cite{Shelah74}), the \emph{Borel Conjecture} in measure theory (see \cite{Laver76}), and the \emph{Kaplansky Conjecture} on Banach algebras in functional analysis (see \cite{Esterle78}, \cite{Dales79} and \cite{WoodinPhD}), the most famous example is undoubtedly the \emph{Continuum Hypothesis} ($\CH$) (see \cite{Godel39} and \cite{Cohen63}). This hypothesis states that given an infinite non-empty set $A \subseteq \mathbb{R}$ either $A$ has size $\aleph_{0}$ —the size of the set of natural numbers\footnote{In this work, the set of natural numbers starts from $0$ and is denoted by $\omega$.} $\omega$— or it has the same cardinality as the set of real numbers $\mathbb{R}$. Equivalently, in terms of cardinal numbers, this can be expressed as $\gc = \aleph_{1}$, where $\gc \coloneqq 2^{\aleph_{0}}$ is called the \emph{size of the continuum} and represents the cardinality of the set of real numbers, and $\aleph_{1}$ is the smallest cardinal number greater than $\aleph_{0}$.

The Continuum Hypothesis became one of the most important problems in the history of mathematics, to the point of being the first problem on David Hilbert's celebrated list of problems presented in 1900 (see \cite{Hilbert02}). To prove that a statement $\varphi$ is independent from $\ZFC$, it is necessary to show that both $\varphi$ and its negation are \emph{consistent} with $\ZFC$, where ``$\varphi$ consistent with $\ZFC$'' means that adding $\varphi$ to $\ZFC$ does not lead to contradictions whenever $\ZFC$ is consistent. Typically, this is done by constructing a universe —called a model— where all the axioms of $\ZFC$ plus $\varphi$ hold. In the case of the Continuum Hypothesis, its consistency with $\ZFC$ was proven in 1939 by Kurt Gödel (see \cite{Godel39}). On the other hand, in 1963, Paul Cohen introduced an innovative powerful method known as the \emph{forcing method}, which proved the consistency of the negation of $\CH$ with $\ZFC$ and settled the independence of such hypothesis (see \cite{Cohen63}). This achievement earned Cohen the Fields Medal, the highest honor in mathematics. The forcing method not only settled the Continuum Hypothesis problem but also opened new horizons in mathematics where it is possible to construct different universes where certain mathematical statements such as $\CH$ can take different truth values.

Before the proof of the independence of the Continuum Hypothesis from $\ZFC$, Georg Cantor made numerous unsuccessful attempts to prove or refute it. However, these efforts significantly contributed to enhancing our understanding of the structure of real numbers. His investigations led to the development of new branches of set theory, such as Descriptive Set Theory and Infinitary Combinatorics. The first one —the definability theory of the continuum, according to Akihiro Kanamori (see \cite{KanamoriDST})— was later fully developed by mathematicians of the caliber of Henri Lebesgue, René-Louis Baire, and Émile Borel. In an attempt to find a solution to the Continuum Hypothesis, Cantor was trying to obtain definable descriptions of subsets of the real numbers in order to classify their size. For instance, Cantor showed that if $A \subseteq \mathbb{R}$ is a Lebesgue measurable set and has positive measure, then the cardinality of $A$ is $2^{\aleph_{0}}$. This implied that a potential counterexample to the Continuum Hypothesis could not come from a set with these characteristics. 

On the other hand, as part of the Infinatory Combinatorics, we have the so-called \emph{cardinal invariants}, which are the main object of study in this work. These cardinals, also known as \emph{cardinal characteristics}, are cardinal numbers that capture combinatorial properties of infinite spaces, like the real line. For example, the answer to the question: “How many Lebesgue measure zero sets are necessary to cover the real line?” is a cardinal invariant called the \emph{covering of the null ideal}, denoted by $\cov(\calN)$. 
Inequalities related to these cardinals reflect some combinatorial properties. For instance, the fact that $\aleph_{0}$-many null sets (i.e.\ of Lebesgue measure zero) are not sufficient to cover the real line —because the countable-many union of countable sets is countable— ``translates'' to the inequality $\aleph_{0} < \cov(\calN)$. On the other hand, the fact that $\gc$-many null sets are sufficient to cover $\mathbb{R}$ —since $\bigcup_{r \in \mathbb{R}} \{ r \}$ covers the real line— implies that $\cov(\calN) \leq 2^{\aleph_{0}}$. 

Other combinatorial questions related to the structure of the real numbers relevant to this work are:

\begin{enumerate} 
    \item How many null sets are necessary for their union not to be a null set? 
    
    \item What is the smallest size a non-null set can have? 
    
    \item What is the smallest size of a family of null sets that is $\subseteq$-cofinal among the null sets? (See \autoref{o00}~\ref{o00-2}.)
\end{enumerate}

The answers to these questions are given by the following cardinals characteristics: the \emph{additivity of the null ideal}, denoted by $\add(\calN)$; the \emph{uniformity of the null ideal}, denoted by $\non(\calN)$; and the \emph{cofinality of the null ideal}, denoted by $\cof(\calN)$, respectively.

These cardinal characteristics are formalized as follows. 

\begin{definition}
Given a set $X$ and $\calI\subseteq \pts(X)$, we can define the \emph{cardinals invariants associated with $\calI$} as:

  \begin{tabular}{rr@{~}l}
      \qquad\quad & $\cov(\calI)$ & $\coloneqq \min \{ \vert \mathcal{F} \vert \colon \mathcal{F} \subseteq \mathcal{I} \conj \bigcup \mathcal{F} = X \},$ called the \emph{covering of $\calI$}.\\[1ex]
    
       & $\add(\calI)$ & $\coloneqq \min \{ \vert \mathcal{F} \vert \colon \calF \subseteq \calI \conj \bigcup \calF \notin \calI \},$ called the \emph{additivity of $\calI$}. \index{$\add(\calI)$}\\[1ex]

       & $\non(\calI)$ & $\coloneqq \min \{ \vert Y \vert \colon Y \subseteq X \conj Y \notin \calI \},$ called the \emph{uniformity of $\calI$}.\\[1ex]

       & $\cof(\calI)$ & $\coloneqq \min \{ \vert \calF \vert \colon \calF \subseteq \calI \conj \forall A\in \mathcal{I} \exists B\in \mathcal{F} \ (A \subseteq B) \},$ called the \emph{cofinality of $\calI$}.\\
    \end{tabular}
\end{definition}
    
    If $\calI$ is $\subseteq$-downwards closed,\footnote{I.e.\ $B\in \calI$ and $A\subseteq B$ implies $A\in\calI$.} $X\notin\calI$ and $X=\bigcup\calI$, then it is not hard to prove the inequalities in \autoref{f-4}. For instance, $\add(\calI) \leq \non(\calI)$ because, whenever $A \subseteq X$ is a witness of $\non(\calI)$, $\add(\calI) \leq \vert \calF \vert = \vert A \vert = \non(\calI)$, where $\calF \coloneqq \{ \{ a \} \colon a \in A \}$. 

    \begin{figure}[ht]
    \centering
    
    \begin{tikzpicture}[scale=1.6]
    
    \small{
    
    \node (azero) at (-1,1) {$\aleph_0$};
    
    \node (addI) at (1,1) {$\add(\calI)$};
    
    \node (covI) at (2,2) {$\cov(\calI)$};
    
    \node (nonI) at (2,0) {$\non(\calI)$};
    
    \node (cofI) at (4,2) {$\cof(\calI)$};
    
    \node (sizX) at (4,0) {$\gc$};
    
    \node (sizI) at (5,1) {$|\calI|$};
    
    \node (sizPX) at (7,1) {$2^{|X|}$};
    
    \draw (azero) edge[->] (addI);
    
    \draw (addI) edge[->] (covI);
    
    \draw (addI) edge[->] (nonI);
    
    \draw (covI) edge[->] (sizX);
    
    \draw (nonI) edge[->] (sizX);
    
    \draw (covI) edge[->] (cofI);
    
    \draw (nonI) edge[-,draw=white,line width=3pt] (cofI);
    \draw (nonI) edge[->] (cofI);
    
    \draw (sizX) edge[->] (sizI);
    
    \draw (cofI) edge[->] (sizI)
    
          (sizI) edge[->] (sizPX);
    }
    \end{tikzpicture}
    
    \caption{Order relationships between the cardinal invariants associated with $\calI$. An arrow indicates that $\ZFC$ proves $\leq$.}
    \label{f-4}
\end{figure}

Other cardinal characteristics are associated with preorders.

\begin{definition}\label{o00}
    A \emph{preordered set} is a pair $\la P,\leq\ra$ where $P$ is a non-empty set and $\leq$ is a reflexive and transitive relation on $P$, usually called a \emph{preorder}. 
    \begin{enumerate}[label = \rm (\arabic*)]
        \item A set $X\subseteq P$ is \emph{bounded in $\la P,\leq\ra$} if there is some $y\in P$ such that $x\leq y$ for all $x\in X$.
        \item\label{o00-2} A set $D\subseteq P$ is \emph{cofinal in $\la P,\leq\ra$} if, for any $x\in P$, there is some $y\in D$ such that $x\leq y$.
        \item The following cardinal invariants are associated with a preorder $\la P,\leq\ra$:\vspace{1ex}

         \begin{tabular}{r@{~}l}
            $\bfrak(P,\leq)$ & $\coloneqq\{ |F| \colon F \text{ is unbounded in }\la P,\leq\ra \}$, called the \emph{completeness of $\la P,\leq\ra$.}\\[1ex]

            $\cof(P,\leq)$ & $\coloneqq\{ |D| \colon D \text{ is cofinal in }\la P,\leq\ra \}$, called the \emph{cofinality of $\la P,\leq\ra$.}\\[1ex]
        \end{tabular}

        When the relation $\leq$ is understood, we just write $\bfrak(P)$ and $\cof(P)$.

        \item When $\la L,\leq\ra$ is a linear order without a maximal element, $\bfrak(L)=\cof(L)$, which is usually denoted by $\cf(L)$.

        \item An ordinal number $\alpha$ is \emph{regular} if $\cf(\alpha)=\alpha$. Otherwise, we say that the ordinal is \emph{singular}.
    \end{enumerate}
\end{definition}

Representative examples of cardinal invariants are those associated with $\calN$ and $\calM$: the collection of Lebesgue-null subsets of real numbers and the collection of meager (first-category) subsets of real numbers, respectively; and the \emph{bounding number} $\bfrak\coloneqq \bfrak({}^\omega\omega,\leq^*)$ and the \emph{dominating number} $\dfrak\coloneqq \cof({}^\omega\omega,\leq^*)$, where $\leq^*$ is the preorder on the set ${}^\omega\omega$ of functions from $\omega$ into $\omega$ defined by $x \leq^{\ast} y$ iff there exists $N < \omega$ such that, for any $n \geq N, $ $x(n) \leq y(n)$, in which case we say that $y$ \emph{dominates} $x$. 
All their possible inequalities in $\ZFC$ are illustrated in the so-called \emph{Cicho\'n's diagram} (see \autoref{f25}) in accordance with the following rule: an inequality between two cardinals can be proved in $\ZFC$ if, and only if, the inequality can be illustrated with an arrow in the diagram. 

Cicho\'n's diagram was originally settled through a series of papers \cite{Ba1984,FrCichon,Miller2,RSt,Roth1,Roth2,Trussc,Miller2,JS,BJS}. Surveys can be found in~\cite{BJ} (see also~\cite{B10} and~\cite{BartInv}).

\begin{figure}[ht]
    \centering
    
    \begin{tikzpicture}
        \small{
        
        \node (aleph1) at (-2,3) {$\aleph_{1}$};
        
        \node (addn) at (0,3){$\add(\calN)$};
        
        \node (covn) at (0,7){$\cov(\calN)$};
        
        \node (nonn) at (9,3) {$\non(\calN)$} ;
        
        \node (cfn) at (9,7) {$\cof(\calN)$} ;
        
        \node (addm) at (3,3) {$\add(\calM)$} ;
        
        \node (covm) at (6,3) {$\cov(\calM)$} ;
        
        \node (nonm) at (3,7) {$\non(\calM)$} ;
        
        \node (cfm) at (6,7) {$\cof(\calM)$};
        
        \node (b) at (3,5) {$\gb$};
        
        \node (d) at (6,5) {$\gd$};
        
        \node (c) at (11,7) {$\gc$};




        
        
        \draw (aleph1) edge[->] (addn)
        
              (addn) edge[->] (covn)
        
              (covn) edge [->] (nonm)
        
              (cfm)edge [->] (cfn)
        
              (cfn) edge[->] (c);

        \draw [dashed] (addm) edge [->]  (covm) 

        (addm) edge [->] (b)

        (nonm)edge [->] (cfm)

        (d)  edge[->] (cfm);

        \draw (addn) edge [->]  (addm)
        
           (covm) edge [->]  (nonn)
        
           (nonn) edge [->]  (cfn);
        
        \draw (b)  edge [->] (nonm);
        
        \draw (covm) edge [->] (d);
        
        \draw (b) edge [->] (d);
        
        }
    \end{tikzpicture}
    
        \caption{Cicho\'n's diagram. The arrows mean $\leq$ and the dotted arrows represent $\add(\Mwf)=\min\{\bfrak,\cov(\Mwf)\}$ and $\cof(\Mwf)=\max\{\dfrak,\non(\Mwf)\}$. It is well-known that this diagram is \emph{complete} in the sense that no other inequality can be added to the diagram.}   
    \label{f25}
\end{figure}

One of the fundamental questions in the study of cardinal invariants is, consistently, what values they can take. Generally, these cardinals lie between $\aleph_{1}$ and $\mathfrak{c}$. Therefore, if we assume the continuum hypothesis is false, several possibilities arise. The typical answer to this question is that they can take ``any regular value'', although there are notable exceptions as we will see later. Another central question regarding cardinal invariants is the so-called \emph{separation of cardinals}: given two cardinal invariants $\mathfrak{j}$ and $\mathfrak{k}$, we say that they can be \emph{separated} if there exists a model of $\ZFC$ in which either $\mathfrak{j} < \mathfrak{k}$ or $\mathfrak{k} < \mathfrak{j}$ holds. A model of $\ZFC$ where two or more cardinals can be separated simultaneously is called a \emph{constellation} or a \emph{separation}. For example, in \cite{GMS}, a constellation of Cichoń's diagram where all the cardinals on the left side of the diagram are separated was constructed. The main tool used to tackle cardinal separation problems is called \emph{iterated forcing}. This technique was developed in the works of Robert M. Solovay and Stanley Tennenbaum (see \cite{ST}), who constructed a method in which Cohen's forcing method is iterated along a transfinite sequence. This method was used to solve problems in infinite combinatorics, such as the consistency of \emph{Martin's Axiom}, and in general topology, such as the consistency of the \emph{Suslin Hypothesis}. The construction of constellations of Cichoń's diagram has been one of the most important problems in the study of the combinatorics of real numbers, and it has greatly motivated the development of increasingly sophisticated iterated forcing techniques. For example, in a notable result, Martin Goldstern, Jakob Kellner, and Saharon Shelah~\cite{GKS} developed the technique of Boolean ultrapowers to successfully construct, under large cardinal assumptions, a model of a constellation of Cicho\'n's diagram where all the non-dependent cardinals are separated. Afterwar, with the second author of this paper, they introduced forcing intersected with submodels to remove the large cardinal assumption (see \cite{GKMS}).
This constellation is known as \emph{Cichoń's maximum} and represents the most significant result up to date in the study of this diagram.

Another question in the context of cardinal invariants is related to the cofinalities. We usually have the following heuristic rule: given a cardinal invariant $\mathfrak{k}$, it is either possible to show in $\ZFC$ that $\mathfrak{k}$ is regular, or it can be shown that in the Cohen model or in random's model, $\mathfrak{k}$
takes the value of the continuum. In the latter case, $\mathfrak{k}$ is consistently singular with non-countable cofinality. However, there are some remarkable exceptions. One of them —due to Jörg Brendle— is $\mathfrak{a}$, the \emph{almost-disjointness number},\footnote{This number is defined as the smallest size of an infinite MAD family. Recall that two infinite subsets of $\omega$ are called \emph{almost-disjoint} if their intersection is finite. Also, $A \subseteq [\omega]^{\aleph_{0}}$ is an almost-disjoint family if its members are pairwise almost disjoint, and it is a \emph{MAD family} if it is maximal with respect to being an almost-disjoint family.} which may have countable cofinality. In particular $\mathfrak{a} = \aleph_{\omega}$ is consistent with $\ZFC$ (see \cite{brendle03}). Another example comes from Cicho\'n's diagram.  The cofinalities in the diagram have been extensively studied (see e.g.\ \cite{Miller1}, \cite{Cofinalities}, and \cite{Br}). By the late eighties, it was known that all the cardinals in Cicho\'n's diagram have uncountable cofinality, except one: $\cov(\calN)$ (see \cite[Sec.~5.1]{uribethesis}). Fremlin conjectured that $\cov(\calN)$ has uncountable cofinality, which was partially answered by Bartoszy{\'n}ski~\cite{bartoCov}, who showed that $\cov(\Nwf) \leq \bfrak$ implies that $\cov(\Nwf)$ has uncountable cofinality. However, this problem was open for almost 20 years until Saharon Shelah \cite{Sh00} found a negative solution in the year 2000, thus constituting the first notable example of a cardinal characteristic that does not satisfy the rule described above.

To prove that $\cov(\calN)$ may have countable cofinality, Saharon Shelah constructed a finite-support iteration using partial random forcing and sequences of \emph{finitely additive measures (fams)} on $\calP(\omega)$ to preserve a witness of $\cov(\Nwf)$ without any restriction to the cofinality of its size. 
In the first part of Shelah's iteration, Cohen reals are added to produce a family of Lebesgue measure zero sets, of the desired size, that covers the real line. Afterwards, partial random forcing is used to destroy covering families of smaller size, while preserving the covering family added by the Cohen reals. Such a family is easily preserved when its size is regular, but in the singular case, further efforts are required. For this purpose, Saharon Shelah considered an additional property of the Cohen reals, which we call the \emph{anti-Bendixson property} (see \cite[$(\ast\ast)_{\bar Q}$]{Sh00} and \autoref{a50}). Preserving the anti-Bendixson property until the end of the iteration ensures that the covering family of singular size is preserved~\cite[Lem.~2.7]{Sh00}. Finally, Shelah guarantees the preservation of the anti-Bendixson property by using the finitely additive measures constructed along the iteration~\cite[Lem.~3.3]{Sh00}. We remark that no other method to preserve anti-Bendixson families is known.

Since \cite{Sh00}, the most important subsequent paper for the development of this work is \cite{KST} where new significant contributions appeared, not only in the applications but also in the development of the method. They introduced the notion of \emph{strong fam limit for intervals} (see \cite[Def.~1.7]{KST}), which formalizes requirements to extend iterations at successor steps using fams. Moreover, they proved that the forcing $\tilde{\bbE}$, a variation of a ccc tree-creature forcing that Haim Horowitz and Saharon Shelah introduced in \cite{HSh},\footnote{And later studied in \cite{M24Anatomy}.} is suitable to iterate with fams.  Using this, they forced a constellation of Cicho\'n's diagram where the entire left-hand side of the diagram is separated and $ \gb < \cov(\calN)$, and applied Boolean ultrapowers to this forcing construction to force another constellation of Cicho\'n's maximum. 

In the master thesis of the third author \cite{uribethesis}, 
a first version of a general theory of iterated forcing using finitely additive measures was introduced (see \cite[Ch.~4]{uribethesis}). Moreover, by generalizing preservation results of anti-Bendixson families from \cite{Sh00}, it was possible to construct an iteration forcing a constellation of the left-hand side of Cicho\'n's diagram allowing $\cov(\calN)$ singular, even with countable cofinality (see \cite[Ch.~5]{uribethesis}). 

This paper expands the ideas and results of the third author's master thesis. We present a new and improved version of the forcing theory using finitely additive measures. This theory allows extensions of ultrafilters, hence generalizing the method of FS (finite support) iterations using ultrafilters from~\cite{GMS}, which is one of the main ingredients in the construction of the first model of Cicho\'n's maximum~\cite{GKS}.

To achieve this, we introduce two new notions: a property we call \emph{$(\Xi, \bar{I}, \varp)$-linkedness} (see \autoref{i2}), which generalizes an essential component of the concept of strong fam limits from \cite{KST}, and the generalization of the notion of \emph{intersection number} for forcing notions (see \autoref{s2}). The first notion is used to define the property of \emph{$\mu$-$\calY$-linkedness} for forcing notions (where $\mu$ is an infinite cardinal and $\calY$ is a class of objects with finitely additive measures, see \autoref{i10}), a generalization of the concept of strong fam limits, which is the cornerstone of our iteration theory. 
The second notion, based on Kelley's work~\cite{Kelley59}, is essential to deal with the limit steps of our iteration theory. Depending on how $\calY$ is composed, we can guarantee that the components of a $\mu$-$\calY$-linked forcing notion have a ``large'' intersection number (see \autoref{i11.01}). In this case, the result presented below, which generalizes a natural fact for Boolean algebras with a finitely additive measure (see e.g.~\cite[Sec.~4.1]{uribethesis}), is the key point in our iteration theory to deal with the limit steps.

\begin{mainLemma}[\autoref{s14}]\label{Main:inter}
     Let $\bbP$ be a forcing notion, $\delta \in[0, 1]$ and $Q \subseteq \bbP$ such that $\Int^{\bbP}(Q) \geq \delta.$ Let $\varp > 0,$ $n < \omega$ and $\bar{p} = \langle p_{i} \colon i < n \rangle \in Q^{n}.$ Define
     $$ \Sigma \coloneqq \{ \sigma \in {}^{n}2 \colon \exists q \in \bbP  \forall i < n \ [ (\sigma(p_{i}) = 0 \Rightarrow q \leq p_{i}) \conj  (\sigma(p_{i}) = 1 \Rightarrow p_{i}  \perp  q)] \}.$$ 
     Then, there exists a function $f \colon \Sigma \to [0, 1]\cap \QQQ$ such that: 
     $$\forall i < n\ \left( \sum \{ f(\sigma) \colon \sigma \in \Sigma \conj \sigma(i) = 0 \} > \delta - \varp \right) \text{ and }  \sum_{\sigma \in \Sigma} f(\sigma) = 1.$$
\end{mainLemma}

Our theory studies a suitable framework we propose to iterate with $\mu$-$\calY$-linked forcing notions (see \autoref{k28}). The iterations in this framework are called \emph{$\calC$-iterations}, where $\calC$ is a set of objects, composed of finitely additive measures, which we call a \emph{fam-iteration bedrock}. 
We succeeded in proving general extension theorems, stated in a simplified way as follows.

\begin{teorema}[Generalized extension theorem at successor steps --\autoref{k39}]
    Any $\calC$-iteration $\Por_\pi$ of length $\pi$ can be extended to a $\calC$-iteration of length $\pi+1$ after iterating at $\pi$ with a forcing notion that is $\mu$-$\calY_\pi$-linked in some suitable forcing sub-extension of $\Por_\pi$.
\end{teorema}

\begin{teorema}[Generalized extension theorem at limit steps --\autoref{k42}]
    If $\gamma$ is a limit ordinal, $\Por_\gamma$ is a finite support iteration of length $\gamma$ such that each initial segment of the iteration is a $\calC$-iteration then, under certain conditions, the whole iteration is a $\calC$-iteration.
\end{teorema}

Depending on the iterations, the previous theorems can be applied to construct iterations with ultrafilters as in \cite{GMS,GKS,Ye}.

Our general theory produces iterations that are \emph{$\kappa$-$\Fr$-Knaster}. This notion comes from a property introduced in~\cite{mejiavert}, called \emph{$\Fr$-linkedness}, that provides a general framework of posets not adding dominating reals.  J\"org Brendle and the first two authors proved that forcing notions with such a property preserve $\gb$ small~\cite{mejiavert,BCM}. Therefore, the iterations with finitely additive measures do not increase $\gb$, which comes from the following general result that we prove in~\autoref{k62} and~\autoref{a68}. 

\begin{teorema}\label{itFr}
    Let $\kappa$ be a regular uncountable cardinal. 
    Under certain conditions on the bedrock $\calC$, any $\calC$-iteration of length $\pi\geq\kappa$, where each iterand has ${<}\kappa$ many linked components, is $\kappa$-$\Fr$-linked and forces $\bfrak\leq\kappa$ and $|\pi|\leq\dfrak$.
\end{teorema}

Our iteration theory also permits preserving and forcing anti-Bendixson families, which generalizes Shelah's~\cite[Lem.~3.3]{Sh00} (the particular case of iterating restricted random forcing and small posets).

\begin{teorema}[\autoref{a64}]\label{itantiB}
    Let $\kappa$ be a regular uncountable cardinal. 
    Under certain conditions on the bedrock $\calC$, any $\calC$-iteration of length $\pi\geq\kappa$, where each iterand has ${<}\kappa$ many linked components, forces that the Cohen reals added along the iteration form a $\kappa$-anti-Bendixson family of size $|\pi|$. 
\end{teorema}

We also provide in~\autoref{a60} our reformulation of Shelah's result~\cite[Lem.~2.7]{Sh00} stating that certain iterations preserving anti-Bendixson families (obtained from Cohen reals) preserve covering families of measure zero sets regardless of the cofinality of their size.

The iterations with fams have good preservation properties for $\non(\Ecal)$, where $\Ecal$ is the ideal on $\R$ generated by the $F_\sigma$ measure zero sets.

\begin{teorema}[\autoref{a69}, see also~\ref{e090}]\label{itnonE}
    Let $\kappa$ be a regular uncountable cardinal. 
    Under certain conditions on the bedrock $\calC$, any $\calC$-iteration of length $\pi\geq\kappa$, where each iterand has ${<}\kappa$ many linked components, forces $\non(\Ecal)\leq\kappa$ and $|\pi|\leq\cov(\Ecal)$.
\end{teorema}

This result was announced by the first and second authors in international conferences~\cite{MnonE,CMUZnonE} as a generalization of Bartozsy\'nski's and Shelah's results about random forcing~\cite{BS1992,Car23}.

The original model from Shelah~\cite{Sh00} yields a constellation of Cicho\'n's diagram of the form $\add(\Nwf) = \bfrak < \cov(\Nwf) < \cfrak$ with $\cov(\Nwf)$ singular (even with countable cofinality), where $\cfrak$ is also allowed (in principle) to be singular. This and \autoref{Mthm:sepcd}  (and~\ref{Mthm:sepcd2}) below are the only known constellations, obtained by finite support iterations of ccc posets, where some cardinal of the left side of Cicho\'n's diagram is singular. On the other hand, a bit more is known for the right side. The model of Cicho\'n's maximum, as well as most constellations obtained from finite support iterations, only allows $\cfrak$ to be singular. Some improvements appear in~\cite{mejiavert} and~\cite{GKMScol} where two values on the right side are singular, and very recently the second author, with Goldstern, Kellner, and Shelah, constructed a model of Cicho\'n's maximum where all five non-dependent cardinal characteristics of the right side of the diagram are singular, see~\cite{Mejmini} (however, large cardinals are used in this proof).

Some other forcing techniques yield constellations of Cicho\'n's diagram where $\non(\Mwf)$ is singular, namely, large-products with creatures (see~\cite{FGKS}) and some models using a large measure algebra (see~\cite[Sec.~5]{GKMScol}, where $\cov(\Nwf)$ can also be obtained singular but with uncountable cofinality). However, it is unknown how to force $\non(\Mwf)$ singular with finite support iterations.

As an application, we aim to force a constellation of Cicho\'n's diagram where the cardinals on the left side are pairwise different and $\cov(\Nwf)$ is singular (even with countable cofinality). For this purpose, we use Brendle~\cite{Br} and Judah and Shelah~\cite{JS} preservation theory to understand the effect of our iterations on other cardinal characteristics of the continuum. Using this, we can conclude that, in Shelah's model, $\cov(\Mwf)=\cfrak$ and $\non(\Mwf) = \cov(\Nwf)^+$, the latter when the cofinality of $\cov(\Nwf)$ is forced to be smaller than $\bfrak$. We also separate $\add(\Nwf)$ and $\bfrak$ in our application, which yields the following main result. 

\begin{teorema}[\autoref{a95}]\label{Mthm:sepcd}    
    Let $\theta \leq \kappa$ be uncountable regular cardinals and let $\lambda$ and $\chi$ be cardinals such that $\lambda$ is $\aleph_{1}$-inaccessible, $\cof([\chi]^{<\nu})=\chi$ for $\nu\in\{\theta,\kappa,\lambda\}$, $\cf(\lambda)<\kappa<\lambda\leq \chi = \chi^{\aleph_0}$ and $\log\chi<\lambda$.  
    Then, there exists a $\kappa$-$\Fr$-Knaster $\rm{ccc}$ forcing notion that forces the constellation of~\autoref{f1300} and $\non(\Ecal)=\kappa$. 

    \begin{figure}[ht]
    \centering
    
    \begin{tikzpicture}
        \small{
        
        \node (aleph1) at (-2,3) {$\aleph_{1}$};
        
        \node (addn) at (0,3){$\add(\calN)$};
        
        \node (covn) at (0,7){$\cov(\calN)$};
        
        \node (nonn) at (9,3) {$\non(\calN)$} ;
        
        \node (cfn) at (9,7) {$\cof(\calN)$} ;
        
        \node (addm) at (3,3) {$\add(\calM)$} ;
        
        \node (covm) at (6,3) {$\cov(\calM)$} ;
        
        \node (nonm) at (3,7) {$\non(\calM)$} ;
        
        \node (cfm) at (6,7) {$\cof(\calM)$};
        
        \node (b) at (3,5) {$\gb$};
        
        \node (d) at (6,5) {$\gd$};
        
        \node (c) at (11,7) {$\gc$};
        
        
        \draw (aleph1) edge[->] (addn)
        
              (addn) edge[->] (covn)
        
              (covn) edge [->] (nonm)
        
              (cfm)edge [->] (cfn)
        
              (cfn) edge[->] (c);

        \draw (addm) edge [->]  (covm) 

        (addm) edge [->] (b)

        (nonm)edge [->] (cfm)

        (d)  edge[->] (cfm);

        \draw (addn) edge [->]  (addm)
        
           (covm) edge [->]  (nonn)
        
           (nonn) edge [->]  (cfn);
        
        \draw (b)  edge [->] (nonm);
        
        \draw (covm) edge [->] (d);
        
        \draw (b) edge [->] (d);

        \draw (-3,2) -- (12,2);

        \draw (-3, 2) -- (-3, 8);

        \draw (-3, 8) -- (12, 8);

        \draw (12, 2) -- (12, 8);
        }

        \draw[redun] (4.5, 2) -- (4.5, 8);

        \draw[redun] (-3, 5.6) -- (4.5, 5.6);

        \draw[redun] (-1.4, 2) -- (-1.4, 5.6);

        \draw[redun] (1.5, 2) -- (1.5, 5.6);

        \draw[redun] (1.5, 5.6) -- (1.5, 8);


        \filldraw[color=greenun, fill=greenun!10, very thick](2.25, 6.25) circle (0.4);   

        \node[greenun] at (2.25, 6.25) {$\lambda^{+}$};      

        \filldraw[color=greenun, fill=greenun!10, very thick](7.5,5) circle (0.4);

        \node[greenun] at (7.5, 5) {$\chi$}; 

        \filldraw[color=greenun, fill=greenun!10, very thick](2.2, 4.5) circle (0.4);

        \node[greenun] at (2.2, 4.5) {$\kappa$}; 

        \filldraw[color=greenun, fill=greenun!10, very thick](-0.7, 4.5) circle (0.4);

        \node[greenun] at (-0.7, 4.5) {$\theta$}; 

        \filldraw[color=greenun, fill=greenun!10, very thick](-1.5, 7) circle (0.4);

        \node[greenun] at (-1.5, 7) {$\lambda$}; 

        \filldraw[color=greenun, fill=greenun!10, very thick](-1.5, 7) circle (0.4);

        \node[greenun] at (-1.5, 7) {$\lambda$}; 
    \end{tikzpicture}
    
    \caption{A separation of the  left side of Cicho\'n's diagram with $\cov(\calN)$ possibly singular.}
    \label{f1300}
\end{figure}
\end{teorema}

This theorem uses the following notation:

\begin{enumerate}[label=\normalfont(\arabic*)]
    \item A cardinal $\lambda$ is \emph{$\theta$-inaccessible} if $\mu^\nu < \lambda$ for any cardinals $\mu<\lambda$ and $\nu<\theta$.
    \item $\log\chi$ is the minimal cardinal $\nu$ such that $\chi\leq 2^\nu$
\end{enumerate}

The assumption $\log\chi<\lambda$ is needed (as in Shelah's original result) to apply a famous theorem by Engelking and Kar\l owicz (see \autoref{k20}) to construct a suitable bedrock $\calC$ and a $\calC$-iteration forcing the desired constellation. However, this assumption is a limitation to force the same constellation over a model of GCH: under GCH, $\log\chi<\lambda\leq\chi$ only occurs when $\lambda=\chi$ is a successor cardinal. 

To overcome this difficulty, we use a trick from~\cite{GKS} to discard the assumption ``$\log\chi<\lambda$" when GCH holds in the ground model. More generally:

\begin{teorema}[\autoref{a96}]\label{Mthm:sepcd2}
    In \autoref{Mthm:sepcd}, the hypothesis ``\/$\log\chi <\lambda$'' can be replaced by ``either $2^{\aleph_0}<\kappa$ or $2^{<\kappa}<\lambda$''.
\end{teorema}

The trick consists of using a generic extension of $\Por^*\coloneqq\Fn_{<\mu}(\chi,2)$ (the forcing adding $\chi$-many higher Cohen reals in ${}^\mu 2$)\footnote{Without assuming $2^{<\mu}=\mu$, when $\mu\leq\chi$ is regular, this forcing is ${<}\mu$-closed, $(2^{<\mu})^+$-cc, and it forces $|{}^{<\mu}2|=\mu$ and $|\chi|\leq 2^\mu$.} for $\mu=\aleph_1$ or $\mu=\kappa$ according to each case. Engelking's and Kar\l owicz's Theorem can be used to construct a suitable bedrock $\calC$ in the $\Por^*$-extension, so the desired forcing is constructed in the ground model as a finite support iteration of ccc posets such that it is a $\calC$-iteration from the point of view of the $\Por^*$-generic extension. We borrow the fams from the $\Por^*$-extension to show that the iteration will force, in $V$, the desired constellation.

Finally, throughout this paper, the reader will observe that the development of the forcing theory we propose depends on a series of results related to the topological and analytical properties of finitely additive measures on Boolean algebras, such as those stated in \autoref{m48}, \autoref{m420}, \autoref{t85}, and \autoref{t93}, among others, found in \autoref{2}. Due to the lack of references that adequately meet our needs, the authors have decided to prepare a complementary article, exclusively dedicated to the study of finitely additive measures on Boolean algebras (see \cite{CMUP}), where these results are developed in detail. Although that work is still in progress, the reader can, for the time being, consult \cite[Ch.~3]{uribethesis}, where similar results are presented, and in many cases, the exact same ones required for this work. In a more general context, the reader can refer to \cite{BhaskaraRa}.

\section{Notation and preliminaries}\label{sec:pre}

\subsection{General notation}\ 

We denote by $\ZZZ,$ $\QQQ$, and $\RRR$ the sets of integers, rational and real numbers respectively. Notice the difference between the symbols ``$\QQQ$'' and ``$\bbQ$'', and ``$\RRR$'' and ``$\bbR$'', since $\bbQ$ and $\bbR$ will be used to denote forcing notions. If $I \subseteq \RRR$ is an interval, we define $I_{\QQQ} \coloneqq I \cap \QQQ.$ We use the symbol $\omega$ to denote the set of natural numbers, $\Ord$ denotes the proper class of ordinal numbers, and $\Card$ denotes the proper class of the cardinal numbers.  The cardinal of the real numbers $2^{\aleph_{0}}$ is denoted by $\gc.$

Let $A, B$ be sets, $\alpha \in \Ord$ and $\kappa, \lambda \in \Card.$ Define $$[A]^{\kappa} \coloneqq \{ X \subseteq A \colon \vert X \vert = \kappa\} \text{ and } [A]^{< \kappa} \coloneqq \{ X \subseteq A \colon \vert X \vert < \kappa \}.$$  We denote by ${}^{A} B$ the set of functions $f$ from $A$ into $B$ and ${}^{< \alpha} A \coloneqq \bigcup \{ {}^{\xi} A \colon \xi < \alpha  \}.$ Similarly, ${}^{\leq \alpha}A \coloneqq {}^{< \alpha +1 }A.$ For any $t \in {}^{< \alpha} A$ we define its \emph{length} by $\lg(t) \coloneqq \dom(t).$  We use the symbols ``$\langle$'' and ``$\rangle$'' to denote sequences and ``$\langle \ \rangle$'' to denote the empty sequence. When $A, B$ are non-empty and $\kappa$ is an infinite cardinal, $\Fn_{<\kappa}(A, B)$ is the set of partial functions from $A$ into $B$ with domain of size ${<}\kappa$. 

We denote by $\rm{trcl}(A)$ the \emph{transitive closure of $A$} and $\calH(\kappa) \coloneqq \{A \colon \vert \rm{trcl}(A) \vert  < \kappa \}.$ The \emph{order type} of $A$ is denoted by $\otp(A).$ If $\la A, R\ra$  is a preoder and $B \subseteq A,$ set $B^{\uparrow} \coloneqq \{ a \in A \colon \exists b \in B\ (b\mathrel{R} a) \}.$ Again, if $B \subseteq A,$  we denote by $\chi_{B}$ the characteristic function of $B$ over $A.$ 

We use $\cB$ and $\cC$ to denote Boolean algebras, and $a, b, c$, etc.\ to denote elements in Boolean algebras. Regarding the atomic structure,  for a Boolean algebra $\cB$, we define: 

\begin{enumerate}[label=\normalfont(\arabic*)]
    \item For any $b \in \cB$ and $d \in \{ 0, 1 \}$:  
    $$b^{d} \coloneqq 
    \left\{ \begin{array}{ll}
            b &   \text{if $d = 0$,} \\[1ex]
            {\sim }b &  \text{if  $d = 1$.}
        \end{array}
    \right.$$ 

    \item For $B \subseteq \cB$ and $\sigma \in \Fn(B, 2)$, $ a_{\sigma} \coloneqq \bigwedge_{b \in \dom(\sigma)} b^{\sigma(b)}. $ \index{$a_{\sigma}$}
\end{enumerate}

We denote by $\At_{\cB}$ the set of all atoms of $\cB$. In the case when $\cB$ is generated by some finite subset $B$, $\At_{\cB} = \{ a_{\sigma} \colon \sigma \in {}^{B}2 \conj a_{\sigma} \neq 0_\cB \}$, so $\At_{\cB}$ is finite, $\vert  \At_{\cB} \vert \leq 2^{\vert B \vert}$, and $|\cB|=2^{|\At_{\cB}|}.$

Now, we review the following notation about trees. Given a non-empty set $Z$, say that $T\subseteq\finseqz$ is a \emph{tree} when $\la \ \ra\in T$ and, if $s\subseteq t$ and $t\in T,$ then $s\in T$. Denote $\Lev_n(T) \coloneqq T\cap{}^{n} Z$ the \emph{$n$-th level of $T$}. For $\rho \in T$, we define $\suc_{T}(\rho) \coloneqq  \{\varrho \in T \colon   \rho \subseteq \varrho\} \cap \Lev_{\lg(\rho) + 1}(T)$. When the context is clear, we simply write  ``$\suc(\rho)$'' instead of ``$\suc_{T}(\rho)$''. The \emph{height} of the tree $T$ is defined by $\hgt(T) \coloneqq \sup\{\lg(t)+1\colon t\in T\}$.  An \emph{infinite branch} of $T$ is an element of $z \in {}^{\omega}Z$ such that, for any $n < \omega,$ $z \rest n \in T.$ The set of infinite branches of $\calT$ is denoted by $[\calT].$ Finally, if $t \in T$, then $T\wedge t \coloneqq \set{s\in T}{s \text{ is comparable with $t$ in $T$}}$.

Let $T\subseteq\finseqz$ be a tree. Say that $s\in T$ is a \emph{splitting node of $T$} if $\vert \suc_{T}(\rho) \vert > 1$, and we let $\max(T) \coloneqq \{ \rho \in T \colon \suc_{T}(\rho) = \emptyset \},$ the set of \emph{maximal nodes} of $T.$ In the case that $T$ has some splitting node, we define $\trunk(T)$ as the splitting node of shortest length. Say that $T$ is \emph{perfect} if, for any $\rho \in T$, there exists some $\eta \in T$ such that $\rho \subseteq \eta$ and $\eta$ is a splitting node in $T$. Finally, $T$ is a \emph{well-pruned} tree if, for any $\rho \in T$ with $\lg(\rho)+1<\hgt(T)$, $\suc_{T}(\rho) \neq \emptyset$.

When $X$ is a topological space, $\Bwf(X)$ denotes the $\sigma$-algebra of \emph{Borel subsets of $X$}. The \emph{Lebesgue measure} $\Leb$ in the \emph{Cantor space} $\cantor$ is the (completion) of the product measure of the uniform measure on $\{0,1\}$ (i.e.\ which assigns measure $\frac{1}{2}$ to $\{0\}$ and $\{1\}$). We denote by $\calN$ the ideal of Lebesgue measure zero subsets of $\cantor$. The values cardinal invariants $\add(\Nwf)$, $\cov(\Nwf)$, $\non(\Nwf)$ and $\cof(\Nwf)$ do not depend on the spaces $\R$ and $\cantor$. For practicality, we work in the Cantor space.

We assume the reader to be familiar with basic techniques of set theory~\cite{Kunen,Je2,Ke2} and, in particular, with forcing and its iteration theory. Recall the following forcing notions:

\begin{enumerate}[label = \normalfont(\arabic*)]
     \item  The forcing notion $\Cor_H \coloneqq \Fn(H\times \omega,2)$ (or $\Fn(H\times\omega,\omega)$, according to our convenience), ordered by reverse inclusion, is the \emph{forcing adding Cohen reals indexed by $H$}. \emph{Cohen forcing} is $\Cor\coloneqq\Cor_\omega$ (or any atomless countable forcing notion).

    \item \emph{Random forcing}, denoted by $\bbB,$ is the complete Boolean algebra $\calB({}^\omega 2)/\calN$, which is forcing equivalent to $\calB(\cantor) \menos \calN$ ordered with $\subseteq$. This forcing is usually used to increase $\cov(\Nwf)$.
\end{enumerate}

For forcing arguments, we denote the ground model by $V$. For two forcing notions $\Por$ and $\Qor$, we write $\Por\subsetdot\Qor$ when $\Por$ is a complete suborder of $\Qor$, i.e.\ the inclusion map from $\Por$ into $\Qor$ is a complete embedding. If $\Seq{\Por_\alpha}{\alpha\leq\beta}$ is a $\subsetdot$-increasing sequence of forcing notions (like an iteration) and $G$ is a $\Por_\beta$-generic over $V$, we denote, for $\alpha\leq\beta$, $G_\alpha\coloneqq \Por_\alpha\cap G$ and $V_\alpha  \coloneqq V[G_\alpha]$. When $\Por_{\alpha+1}$ is obtained by a two-step iteration $\Por_\alpha\ast\Qnm_\alpha$, $G(\alpha)$ denotes the $\Qnm[G_\alpha]$-generic
set over $V_\alpha$ such that $V_{\alpha+1}= V_\alpha[G(\alpha)]$ (i.e.\ $G_{\alpha+1} = G_\alpha\ast G(\alpha)$). We use $\Vdash_\alpha$ to denote the forcing relation for $\Por_\alpha$, and $\leq_\alpha$ to denote its preorder (although we just use $\leq$ when clear from the context). 

For a forcing notion $\Por$, define the order $\leq^\bullet$ in $\Por$ as $q\leq^\bullet p$ iff for any $r\leq q$, $r$ and $p$ are compatible, i.e., they have a common extension. Recall that $\Pbb$~is \emph{separative} iff ${\leq^\bullet}$ equals ${\le}$.

Recall the following stronger versions of the chain condition of a forcing notion.
\begin{definition}
    Let $\Pbb$ be a forcing notion and $\kappa$ an infinite cardinal.
        \begin{enumerate}[label=\rm(\arabic*)]
        \item \emph{$\Pbb$ has the $\kappa$-cc (the $\kappa$-chain condition)} if every antichain in $\Por$ has size ${<}\kappa$. \textit{$\Por$ has the ccc (the countable chain contidion)} if it has the $\aleph_{1}$-cc.
        
        \item For $n<\omega$, $B\subseteq \Por$ is $n$-\textit{linked} if, for every $F\subseteq B$ of size $\leq n$, $\exists q\in \Por\forall p\in F\ (q\leq p)$. When $n=2$ we just write \emph{linked}.
        \item $C\subseteq \Por$ is  \textit{centered} if it is $n$-linked for every $n<\omega$.
        \item $\Por$ is $\kappa$-$n$-\textit{linked} if $\Por=\bigcup_{\alpha<\kappa}P_\alpha$ where each $P_{\alpha}$ is $n$-linked. When $\kappa=\omega$, we say that $\Pbb$ is $\sigma$-$n$-\textit{linked}. In the case $n=2$, we just write \emph{$\kappa$-linked} and \emph{$\sigma$-linked}. 
        \item $\Por$ is $\kappa$-\textit{centered} if $\Por=\bigcup_{\alpha<\kappa}P_\alpha$ where each $P_{\alpha}$ is centered. When $\kappa=\omega$, we say that $\Por$ is $\sigma$-\textit{centered}.
        
        
        \item $\bbP$ is \emph{$\kappa$-$n$-Knaster} if any subset of $\Por$ of size $\kappa$ contains an $n$-linked subset of size $\kappa$. When $n=2$ we just write \emph{$\kappa$-Knaster}. 
    \end{enumerate}
\end{definition}

\subsection{Probability trees}\ 

To prove the general extension theorem for limit steps (see \autoref{k40} and \autoref{k42}), some notions and results about probability trees that the second and third authors have developed in~\cite[Chapter~2]{uribethesis} and \cite{MU} are required. Such notions and results are briefly presented in this subsection. We begin fixing basic probability notation. Recall that $ \Omega \coloneqq \la\Omega, \calA, \p\ra$  is a \emph{probability space} if $\Omega$ is a non-empty set, $\calA$ is a $\sigma$-algebra on $\Omega,$ and $\p \colon \calA \to [0, 1]$ is a measure such that $\p(\Omega) = 1.$ In this case, we say that $\p$ is a \emph{probability measure on $\Omega$}. Elements in $\calA$ are called \emph{events} and, if $E, F \in \calA,$ then $\p(E)$ is called the \emph{probability of success of $E$}. When  $\p(E \cap F) = \p(E) \cdot \p(F),$  we say that $E$ and $F$ are \emph{independent events}. Also,  we say that a function $X \colon \Omega \to \R$ is a \emph{random variable on $\Omega$} if, for any $a \in \R,$ $\{ o \in \Omega \colon X(o) \leq a \} \in \calA$ (i.e.\ $X$ is an $\calA$-measurable function). If $X$ has \emph{Bernoulli distribution with parameter $p$}, we write $X \sim \Ber(p).$ Similarly, if $X$ has  \emph{binomial distribution with parameters $n, p$}, we write $X \sim \Bin(n, p).$

\begin{example}\label{p26}
    Let $n \in \bbN $ and $p \in [0, 1].$  We define $\Omega_{n} \coloneqq \{ i \in \bbN  \colon i \leq n \},$ $\calA_{n} \coloneqq \calP(\Omega_{n})$ and  $\p_{n} \colon  \calA_{n} \to \R$ such that, for any $i \leq n,$ 
    $$ \p_n(\{i\}) \coloneqq \binom{n}{i} p^{i} (1 -p)^{n - i}.$$ It is cleat that $(\Omega_{n}, \calA_{n}, \p_{n})$ is a probability space and the identity function $\B_{n, p} \colon \Omega_{n} \to \R$ is a random variable. This corresponds to the binomial distribution.
\end{example}

Now, we define the notion of \emph{probability tree}.

\begin{definition}[{\cite[Def.~2.3.1]{uribethesis}}]\label{p70.0}
    Let $A$ be a non-empty set. Say that $T\subseteq{}^{<\omega}A$ is a \emph{probability tree} if it is a well pruned tree on ${}^{<\omega}A$   
    with associated probability spaces $\la \suc_t,\calA_t,\p^T_t\ra$
    for any $t \in T \menos \max(T)$, 
    where $[\suc(t)]^{< \aleph_{0}} \subseteq \calA_{\rho}$.
\end{definition}

From \cite[Thm.~2.3.2]{uribethesis}, we know that every probability tree $T$ with finite levels induces a probability space in each of its levels.

\begin{theorem}[{\cite[Thm.~2.3.10]{uribethesis}, \cite{MU}}]\label{p70}
    Let $n^{\ast} < \omega$ and $T = {}^{\leq n^{\ast}} 2$ be the complete binary tree of height $n^{\ast} +1$ endowed with a probability tree structure. Define the random variable $Y \colon \Lev_{n^{\ast}}(T) \to \RRR$ by $Y(t) \coloneqq \vert \{ n < n^{\ast} \colon t(n) = 0 \} \vert \text{ for any $t \in \Lev_{n^{\ast}}(T).$}$ Assume that there exists some $p \in [0, 1]$ such that, for any $t \in T \menos \max(T),$ $p \leq p_{t} \coloneqq \p^T_{t}(\concat{t}{\langle 0 \rangle}).$ Then, $ \p_{\Lev_{n^{\ast}(T)}}[Y \leq z] \leq \p_{\Omega_{n^{\ast}}}[\B_{n^{\ast},  p} \leq z] \text{ for all $z \in \RRR.$}$
\end{theorem}

\section{Finitely additive measures: extension and integration}\label{2}

In this section, we review the framework of finitely additive measures that we are going to use throughout this paper. It is recommended to refer to \cite[Chapter~3]{uribethesis} and \cite{CMUP} for details, and in a more general context to \cite{BhaskaraRa}. 

\begin{definition}\label{m4}
    Let $\cB$ be a Boolean algebra. A \emph{finitely additive measure (fam) on $\cB$} is a function $\Xi \colon \cB \to [0,\infty]$ satisfying:    
    
        \begin{enumerate}[label=\normalfont(\roman*)]
            \item \label{m4a} $\Xi (0_{\cB})=0$,
            
            \item \label{m4b} $\Xi(a\vee b)=\Xi(a)+\Xi(b)$ whenever $a,b\in\cB$ and $a \wedge b= 0_{\cB}$.
        \end{enumerate}
\end{definition}

We say that $\Xi$ is \emph{finite} if $\Xi(1_\cB)<\infty$. When $\Xi(1_\cB)=1$ we say that $\Xi$ is a \emph{finitely additive probability measure}.

We now give an example of a  finitely additive probability measure that will appear in several places throughout this work: 

\begin{example}\label{m37}\ 

    \begin{enumerate}
        \item Let $\cB$ be a Boolean algebra, $\Xi$ a finitely additive measure on it, and $b \in \cB$ with positive finite measure. We define the function 
        $\Xi_{b} \colon \cB \to [0,1]$ by $\Xi_{b}(a) \coloneqq \frac{\Xi(a \wedge b)}{\Xi(b)}$ for any $a \in \cB$. It is clear that $\Xi_{b}$ is a finitely additive probability measure. 
        
        \item\label{m37.b} Let $X$ be a non-empty set and denote by ${\rm c}^X$ its \emph{counting measure}. For a finite non-empty set $u \in \calP(X)$, we define $\Xi^u \coloneqq {\rm c}^X_u$, i.e.\ 
        for any $x \in \calP(X),$ $\Xi^{u}(x) \coloneqq \frac{\vert x \cap u \vert}{ \vert u \vert}$, which we call the \emph{uniform measure with support $u$}. 
    \end{enumerate}    
\end{example}

\begin{example}\label{m37.1}
    Assume that $K$, $L$ are non-empty sets and $\Xi$ is a fam on $\calP(K)$.  Then, any function $h \colon K \to L$ induces a natural fam $\Xi_{h}$ on $\calP(L)$ as follows. For any $A \in \calP(L)$, let $\Xi_{h}(A) \coloneqq \Xi( h^{-1}[A] ).$ Notice that $\Xi(K) = \Xi_{h}(L)$ and, therefore, if $\Xi$ a probability fam, then $\Xi_{h}$ is also a probability fam.

    More generally, if $\cB$ and $\cC$ are Boolean algebras and $\Xi$ is a (probability) fam on $\cB$, then any homomorphism $f\colon\cC\to\cB$ induces a (probability) fam on $\cC$, namely, $c\mapsto \Xi(f(c))$. 
\end{example}

The following theorem, known as the \emph{compatibility theorem of finitely additive measures}, is one of the fundamental tools for establishing extension theorems of finitely additive measures.

\begin{theorem}\label{m48}
    Let $\cB$ be a Boolean algebra and, for $d \in\{0,1\}$, let $\cB_d$ be a Boolean sub-algebra of $\cB$ with a finitely additive measure $\Xi_d\colon\cB_d\to[0,\infty)$. Then the following statements are equivalent.
    
    \begin{enumerate}[label=\rm(\alph*)]
        \item\label{m48a} There is a finitely additive measure $\Xi$ on the Boolean sub-algebra generated by $\cB_0\cup\cB_1$ extending $\Xi_d$ for $d\in\{0,1\}$.
        
        \item\label{m48b} $\Xi_0(1_\cB)=\Xi_1(1_\cB)$ and, for any $a\in\cB_0$ and $a'\in\cB_{1}$, if $a\leq a'$ then $\Xi_0(a)\leq\Xi_1(a')$.
        
        \item\label{m48c} For any $d,d'\in\{0,1\}$, $a\in\cB_d$ and $a'\in\cB_{d'}$, if $a\leq a'$ then $\Xi_d(a)\leq\Xi_{d'}(a')$.
    \end{enumerate}
\end{theorem}

\begin{corollary}\label{m52}
    Let $\cB$ be a Boolean algebra, $\cC \subseteq\cB$ a sub-algebra with a finitely additive measure $\Xi\colon \cC \to [0,\infty)$, and let $b\in\cB$. If $z\in[0,\infty)$ is between $\sup\{\Xi(a) \colon  a\leq b,\ a\in \cC \}$ and $\inf\{\Xi(a) \colon  b\leq a,\ a \in \cC\},$ then there is a finitely additive measure $\Xi'$ on the Boolean algebra generated by $\cC \cup\{b\}$, extending $\Xi$, such that $\Xi'(b)=z$.

    As a consequence of Zorn's lemma, there is some fam $\Xi'$ on $\cB$ extending $\Xi$.
\end{corollary}

We also use the following particular case of \autoref{m48}.

\begin{theorem}\label{m69}
   Let $\Xi_0$ be a fam on a Boolean sub-algebra of $\cB$ and $\langle  b_i \colon i\in I\rangle\subseteq\cB$. Assume that $0 < \delta\coloneqq \Xi_0(1_\cB)< \infty$ and, for every finite $J\subseteq I$ and $b \in \dom(\Xi_0)$, if $\Xi_0(b)>0$ then $b \wedge\bigwedge_{i\in J} b_i\neq0_\cB$. Then, there exists a finitely additive measure $\Xi$ on $\cB$ extending $\Xi_0$ such that $\Xi(b_i)=\delta$ for every $i\in I$.
\end{theorem}

Recall that a \emph{field of sets over $X$} is a sub-algebra of $\pts(X)$ under the set operations. We review integration over a field of sets. Fix a non-empty set $X$, a field of sets $\cB$ over $X$, and a finitely additive measure $\Xi \colon \cB \to [0, \infty).$ Motivated by the definition of Riemann's integral, if  $f \colon X \to \mathbb{R}$ is a bounded function we can naturally define $\int_{X} f d \Xi$, if it exists, by approximating with upper and lower sums over finite partitions of $X$ from $\cB$. In this case, we say that $f$ is \emph{$\Xi$-integrable}  (see \cite[Def.~3.5.3]{uribethesis}). For example, any bounded function is $\Xi$-integrable when $\dom(\Xi) = \calP(X)$ (see \cite[Thm.~3.5.10]{uribethesis}). Fundamentally, the integral with respect to finitely additive measures behaves similarly to the Riemann integral, that is, we have available the basic properties of the integral such as finite additivity, linearity, and monotonicity (see \cite[Sec.~3.5]{uribethesis}). We also can integrate over subsets of $X$: for a bounded function $f\colon X \to\bbR$, if $E \subseteq X$ and $\chi_{E} f$ is $\Xi$-integrable, we define $\displaystyle\int_{E} f d \Xi \coloneqq \int_{X} \chi_{E} f d \Xi.$ 

If $E \in \cB$  and $f$ is $\Xi$-integrable, then 
$\chi_{E} f$ is $\Xi$-integrable. In general, this is the context in which we will use integration over subsets.

For example, we can calculate the integral with respect to $\Xi^{u}$: 

\begin{example}\label{t42}
    Let $X$ be a non-empty set. If $u\subseteq X$ is finite and non-empty, then $$\int_{X} fd\Xi^u=\frac{1}{|u|}\sum_{k\in u} f(k).$$
\end{example}


In general, we will not be interested in finitely additive measures assigning a positive measure to finite sets. For this reason, we introduce the notion of \emph{free finitely additive measure}. 

\begin{definition}\label{m12}
    \index{finitely additive measure!free finitely additive measure}

    Let $X$ be a non-empty set, $\cB$ a field of sets over $X$ 
    and $\Xi$ is a finitely additive measure on $\cB$.
    \begin{enumerate}[label = \normalfont (\arabic*)]
        \item $\Pbf^\Xi=\Pbf^\cB$ denotes the set of finite partitions of $X$ into sets in $\dom(\Xi)=\cB$.
        \item We say that $\Xi$ is a \emph{free finitely additive measure} if, for any $x\in X$, $\{x\}\in\cB$ and $\Xi(\{x\}) = 0$.
    \end{enumerate}
\end{definition}

Free finitely additive measures are characterized as follows. 

\begin{lemma}\label{m420}
    Let $\cB$ be a field of sets over $X$ and let $\Xi$ be a finite fam on $\cB$ with $\delta \coloneqq \Xi(X)$. Then the following statements are equivalent.
    
    \begin{enumerate}[label = \rm (\roman*)]
        \item\label{m420i} For any $\varp>0$, any finite $F\subseteq X$ and any 
        $P\in\Pbf^\Xi$, 
        there is some non-empty finite $u\subseteq X\menos F$  such that, for all $b\in P$,
        $\left|\delta \Xi^{u}(b) - \Xi(b)\right| < \varp.$
        \item\label{m420ii} All finite sets in $\cB$ have measure zero.
    \end{enumerate}
    Even more, the $u$ in~\ref{m420i} can be found disjoint with all $b\in P$ of measure zero.
\end{lemma}

Inspired by the previous characterization, we define the following type of fams.

\begin{definition}\label{m392}
    Let $\Xi$ be a finite fam on a field of sets over $X$ and $\delta \coloneqq \Xi(X)$. We say that $\Xi$ has the \emph{uniform approximation property (uap)} if, for any $\varp>0$ and any $P\in \Pbf^\Xi$, there is some non-empty finite $u\subseteq X$  such that, for all $b\in P$,
    $\left|\delta \Xi^{u}(b) - \Xi(b)\right| < \varp.$
\end{definition}

The integral can be approximated by finitely supported fams as follows.

\begin{theorem}[{\cite{CMUP}}]\label{t85}
   Let $\Xi_{0}$ be a finitely additive measure on a field of sets $\cB$ over $X$, and let $\delta \coloneqq \Xi_0(X) < \infty$. 
   Let $I$ be an index set and, for each $i\in I$, let $K_i$ be a closed subset of $\R$ and $f_i\colon X\to\R$ bounded.  
   Then the following statements are equivalent.
   \begin{enumerate}[label = \rm (\Roman*)]
   \item\label{t85.0I} For any $P\in\Pbf^{\Xi_0}$, $\varepsilon>0$, any finite set $J\subseteq I$, and any open $G_i\subseteq \R$ containing $K_i$ for $i\in J$, there is some non-empty finite $u\subseteq X$ and a probability measure $\Xi^-$ on $\pts(u)$ such that: 
   \begin{enumerate}[label=\normalfont(\roman*)]
       \item\label{t85.0a} $\displaystyle \left |\Xi_0(b)-\delta\Xi^-(b\cap u) \right | <\varepsilon$ for any $b\in P$, and
       
       \item\label{t85.0b} $\displaystyle\delta\int_u f_i d\Xi^-=\frac{\delta}{|u|}\sum_{k\in u}f_i(k)\Xi^-(\{k\}) \in  G_i$ for any $i\in J$.
   \end{enumerate}
   \item\label{t85II} There is some fam $\Xi$ on $\pts(X)$ extending $\Xi_0$ such that, for any $i \in I$, 
   $\displaystyle\int_{X} f_id \Xi \in K_i.$
   \end{enumerate}
   When $\Xi_0$ has the uap, the measure  $\Xi^-$ in~\ref{t85.0I} can be found uniform, concretely, we can add~\ref{t85I} below to the list of equivalences. In this case, $\Xi$ can be found with the uap in~\ref{t85II}
   \begin{enumerate}[resume*]
   \item\label{t85I} For any $P\in\Pbf^{\Xi_0}$, $\varepsilon>0$, any finite set $J\subseteq I$, and any open $G_i\subseteq \R$ containing $K_i$ for $i\in J$, there is some non-empty finite $u\subseteq X$ such that:
   \begin{enumerate}[label=\normalfont(\roman*)]
       \item\label{t85a} $\left |\Xi_0(b)-\delta \Xi^{u}(b) \right | <\varepsilon$ for any $b\in P$, and
       
       \item\label{t85b} $\displaystyle\frac{\delta}{|u|}\sum_{k\in u}f_i(k) \in  G_i$ for any $i\in J$.
   \end{enumerate}
   \end{enumerate}
\end{theorem}

Finitely additive measures with the uap containing finite sets of positive measure have an interesting characterization.

\begin{theorem}[{\cite{CMUP}}]\label{t94}
    Let $\cB$ be a field of sets over $X$ and $\Xi$ a finite fam on $\cB$ with $\delta\coloneqq \Xi(X)$. Then $\Xi$ has the uap iff either all finite sets in $\cB$ have measure zero, or there are $d < \omega$ and $P_*\in \Pbf^\Xi$ (called the ``\emph{frame of $\Xi$}'') such that, for any $b\in P_*$:
    \begin{enumerate}[label = \normalfont (\roman*)]
        \item $\displaystyle\Xi(b)=\delta \frac{k_b}{d}$ for some $0\leq k_b \leq d$.
        \item Any $b'\subseteq b$ in $\cB$ has either measure zero or $\Xi(b)$.
        \item If $b$ is finite then $k_b \leq |b|$.
    \end{enumerate}
\end{theorem}

\begin{corollary}\label{t96}
    In \autoref{m69}, if $\cB=\pts(X)$ and $\Xi_0$ has the uap, then $\Xi$ can be found with the uap.
\end{corollary}

Integrability and the integrals are related when we integrate with respect to larger finitely additive measures: 

\begin{theorem}[{\cite[Thm.~3.5.27]{uribethesis}}]\label{t93}
    Let $\cB_{0}$ and $\cB_{1}$ be fields of sets over $X$ such that $\cB_{0} \subseteq \cB_{1},$ and let $\Xi_{0}$ and $\Xi_{1}$ be finitely additive measures on $\cB_{0}$ and $\cB_{1}$, respectively, such that $\Xi_{0} \subseteq \Xi_{1}.$ Let $f \colon X \to \RRR$ be a bounded function. Then, whenever $f$ is $\Xi_0$-integrable, $f$ is $\Xi$-integrable, in which case $\displaystyle \int_{X} f \, d \Xi_{0} = \int_{X} f d  \Xi_{1}.$ 
\end{theorem}

We are going to extend finitely additive measures along forcing iterations while preserving the values of the integrals. 
For example, if $M\subseteq N$ are transitive models of ZFC and we have a finitely additive measure over $\calP(\omega)^{M},$ there can be new reals in $N$ and $\calP(\omega)^{M} \subsetneq \calP(\omega)^{N}$, but we are interested in keeping the value of an integral on a fam over $\pts(\omega)^M$ when extending it to $\pts(\omega)^N$. So we will use the following result, which is obtained by simply relativizing \autoref{t93}. 

\begin{corollary}[{\cite[Thm.~3.5.29]{uribethesis}}]\label{t102}
    Let $M\subseteq N$ be transitive models of $\ZFC$ such that $\Xi_{0},\cB_{0} \in M$, and $\Xi_{1},\cB_{1} \in N.$ Assume that $\Xi_{0}$ and $\Xi_{1}$ are finitely additive measures on the fields of sets $\cB_{0}$ and $\cB_{1}$ over some $X\in M$, respectively, such that $\Xi_{0} \subseteq \Xi_{1}$ (and $\cB_{0} \subseteq \cB_{1}$). Then $N\models$``$f$ is $\Xi_1$-integrable'' whenever $M\models$``$f$ is $\Xi_0$-integrable'', in which case $\displaystyle \left(  \int_{X} f d\Xi_{0}  \right)^{M} = \left( \int_{X} f d \Xi_{1} \right)^{N}. $
\end{corollary}


\section{The intersection number for forcing notions}\label{4.1}

The first part of this section is based on \cite[Sec.~3]{IntNumU24}, where the reader can find a more detailed presentation of the intersection number in the context of forcing theory. Afterward, we prove \autoref{Main:inter} (\autoref{s14}) and propose a more general version of the intersection number.

\begin{definition}\label{s2}
    Let $\bbP$ be a forcing notion and $Q \subseteq \bbP.$ 

    \begin{enumerate}[label=\rm(\arabic*)]
        \item For a finite sequence $\bar{q} = \langle q_{i} \colon i < n \rangle \in {}^n \bbP,$ we define $$ i_{\ast}^{\bbP}(\bar{q}) \coloneqq \max \{ \vert F \vert \colon F \subseteq n \conj  \{ q_{i} \colon i \in F\} \ \text{has a lower bound in} \ \bbP  \}. $$ 

        \item  The \emph{intersection number of $Q$ in $\bbP$}, denoted by $\rm{int}^{\bbP}(Q),$ is defined by $$\Int^{\bbP}(Q) \coloneqq \inf \left \{ \frac{i_{\ast}^{\bbP}(\bar{q})}{n} \colon \bar{q} \in {}^n Q  \conj   n \in \omega \smallsetminus \{ 0 \}  \right\}.$$
    \end{enumerate}
    We stipulate $\Int^\bbP(\emptyset)=1$. We omit the upper index $\Por$ when the context is clear. 
\end{definition}

It is clear that if $\bar{q} \in {}^{n}Q$ then $1 \leq i_{\ast}(\bar{q}) \leq n.$  As a consequence, $\Int(Q)$ is a real number in $[0, 1].$  Also, 

\begin{lemma}[{\cite[Lem.~3.3]{IntNumU24}}]\label{s8}
    Let $\bbP$ be a forcing notion and $Q \subseteq \bbP$. Then,

    \begin{enumerate}[label=\rm(\arabic*)]
        \item\label{s8d} For any $p \in \bbP,$ $\Int(\{ p \}) = 1.$ Moreover, $\Int(Q)=1$ iff $Q$ is centered in $\Por$.

        \item\label{s8.5} If $Q$ is finite, then $\Int(Q) \geq \frac{1}{\vert Q \vert} > 0.$ 

        \item\label{s8.6} Let $Q$ be an anti-chain in $\bbP.$ Then $\Int(Q) = \frac{1}{\vert Q \vert}$ if $Q$ is finite, and $\Int(Q) = 0$ when $Q$ is infinite.  

        \item\label{s8g} If $m \in (1, \omega)$ and $\Int(Q) \geq 1 - \frac{1}{m + 1},$ then $Q$ is $m$-linked.
    \end{enumerate}
\end{lemma}

Regarding the behavior of the intersection number under complete embeddings, we have the following result: 

\begin{lemma}[{\cite[Lem.~3.6]{IntNumU24}}]\label{s6}
    Let $\bbP, \bbQ$ be forcing notions, $\iota \colon \bbP \to \bbQ$ a complete embedding and $Q \subseteq \bbP.$ Then $\Int^{\bbP}(Q) = \Int^{\bbQ}(\iota[Q]).$ As a consequence, for $R \subseteq \bbQ$, $\Int^{\bbQ}(R) \leq \Int^{\bbP}(\iota^{-1}[R]).$  
\end{lemma}

Kelley proved that finitely additive measures can be used to define subsets of Boolean algebras whose intersection number is bounded by a given value:

\begin{lemma}[{\cite[Prop.~1]{Kelley59}}]\label{s9}
    Let $\cB$ be a Boolean algebra and   $\Xi \colon \cB \to [0, 1]$ a finitely additive measure. Let $\bbP \coloneqq \cB \menos \{ 0_{\cB} \}$ and $\delta \in [0, 1].$ If $Q \coloneqq \{ p \in \bbP \colon \Xi(p) \geq \delta \},$ then $\Int^{\bbP}(Q) \geq \delta.$\footnote{A recent proof of this result can be found in \cite[Thm.~3.7]{IntNumU24}.} 
\end{lemma}

Recall that $\bbB$ denotes the random forcing notion.

\begin{corollary}\label{s10}
      Let $C \subseteq {}^{\omega}2$ be non-empty and $\delta \in [0, 1].$ If  $Q \coloneqq \{ p \in \bbB \colon \Leb_{C}(p) \geq \delta \},$ then $\Int(Q) \geq \delta.$
\end{corollary}

As an attempt to get the converse of \autoref{s9}, the so-called \emph{Kelley's theorem} arose:

\begin{theorem}[{\cite[Thm.~2]{Kelley59}}]\label{s12}
    Let $\cB$ be a Boolean algebra. If $Q \subseteq \cB$ then there exists a probability finitely additive measure $\Xi \colon \cB \to [0, 1]$ such that $\inf \{ \Xi(b) \colon b \in Q \} = \Int(Q).$
\end{theorem}

With \autoref{s12}, we now have what is necessary to prove \autoref{Main:inter}, one of the most important tools of this paper:  

\begin{mainlemma}\label{s14}
    Let $\bbP$ be a forcing notion, $\delta \in[0, 1]$ and $Q \subseteq \bbP$ such that $\Int^{\bbP}(Q) \geq \delta.$ Let $\varp > 0$,  $n < \omega$ and $\bar{p} = \langle  p_{i} \colon i < n \rangle \in Q^{n}.$ Define 
    $$ \Sigma \coloneqq \{ \sigma \in {}^{n}2 \colon \exists q \in \bbP  \forall i < n\ [ (\sigma(p_{i}) = 0 \Rightarrow q \leq p_{i}) \conj  (\sigma(p_{i}) = 1 \Rightarrow p_{i}   \perp  q)] \}.$$ 
    Then, there exists a function $f \colon \Sigma \to [0, 1]_{\QQQ}$ such that: $$\forall i < n\ \left( \sum \{ f(\sigma) \colon \sigma \in \Sigma \conj \sigma(i) = 0 \} > \delta - \varp \right) \ \text{and} \ \sum_{\sigma \in \Sigma} f(\sigma) = 1.$$
\end{mainlemma}

\begin{PROOF}[\textbf{Proof}]{\ref{s14}}
    Let $(\cB, \iota)$ be the completion of $\Por$. Define $A \coloneqq \iota[Q] \subseteq \cB \setminus \{ 0_{\cB} \}.$ By \autoref{s6}, we have that $\Int^{\cB}(A) = \Int^{\bbP}(Q) \geq \delta.$ On the other hand, by \autoref{s12} there exists a finitely additive measure $\Xi \colon \cB \to [0, 1]$ such that, for any $a \in A, \, \Xi(a) \geq \delta.$ 

    For $\sigma \in {}^{n} 2$ define $ b_{\sigma} \coloneqq \bigwedge_{i < n} \iota(p_{i})^{\sigma(i)}$. 
    It is easy to check that $b_{\sigma} \neq 0_{\cB}$ iff $\sigma \in \Sigma.$ 
    Now, we can find a sequence of rational numbers $\langle f(\sigma) \colon \sigma \in \Sigma \rangle$ such that, for each $\sigma \in \Sigma,$ $\vert \Xi(b_{\sigma}) - f(\sigma) \vert < \frac{\varp}{2^{n}}$ and $\sum_{\sigma \in \Sigma} f(\sigma) = 1.$ For $i < n,$ we have,    
    \begin{equation*}
        \begin{split}
            \sum \{ f(\sigma) \colon \sigma \in \Sigma \conj  \sigma(i) = 0 \} & > \sum \left \{ \Xi(b_{\sigma}) - \frac{\varp}{2^{n}} \colon \sigma \in \Sigma \conj  \sigma(i) = 0 \right \} \\ 
            & = \Xi(\iota(p_{i})) - \sum  \left \{ \frac{\varp}{2^{n}} \colon \sigma \in \Sigma \conj  \sigma(i) = 0 \right \} \\
            & \geq \Xi(\iota(p_{i})) - \varp \geq \delta - \varp,
        \end{split}
    \end{equation*}
which proves the result. 
\end{PROOF}

As the last part of this section, we propose a generalization of the intersection number in terms of probability fams on finite sets. They have deep connections with the intersection number (see \autoref{i20.7.1}) and characterize some particular cases of the linkedness notion introduced in the next section (see \autoref{i2} and \autoref{i11.05.1}).

\begin{definition}\label{i20}
    For any non-empty set $K$, $\bfI_{K}^{\fin}$ is the class of all sequences $\bar{I} = \langle I_{k} \colon k \in K \rangle$ of finite non-empty pairwise disjoint sets. Given some $\bar{I} \in \bfI_{K}^{\fin}$, define $W_{\bar{I}} \coloneqq \bigcup_{k \in K} I_{k}$. When the context is clear, we just write ``$W$'', omitting the index. 

    Fix a forcing notion $\bbP$, a finite set $K$, $\bar I \coloneqq \la I_k \colon k\in K\ra \in \bfI^\fin_K$, and a probability fam $\Xi$ on $\pts(K)$. 
    \begin{enumerate}[label = \normalfont (\arabic*)]
        \item For any $\bar q \in {}^W \Por$, define
        \[ \mu^{\Xi, \bbP}_{\bar I}(\bar q)  \coloneqq \max \left \{  \int_{K} \frac{ \vert \{ \ell \in I_{k} \colon q \leq q_{\ell} \}  \vert}{\vert I_{k} \vert} d \Xi(k) \colon q \in \bbP  \right  \} .\footnote{Recall that 
        $\int_{K} \frac{ \vert \{ \ell \in I_{k} \colon q \leq q_{\ell} \}  \vert}{\vert I_{k} \vert} d \Xi(k) = \sum_{k\in K}\frac{|\{\ell\in I_k\colon q\leq q_\ell\}|}{|I_k|}\Xi(\{k\})$ in this case when $K$ is finite.} \]

        Notice that this set is finite, so the maximum exists.
       
        \item For any $Q \subseteq \Por$, define $M^{\Xi, \bbP}_{\bar{I}}(Q) \coloneqq \min \{\mu^\Xi_{\bar I}(\bar q) \colon \bar q \in {}^W Q\}$.

        The minimum exists because this set is finite (even when $Q$ is infinite). In fact, its size is below $ \prod_{k \in K} ( \vert I_{k} \vert + 1 ) \leq ( \max\{ \vert I_{k} \vert \colon k \in K \} + 1)^{\vert K \vert}$. 
    \end{enumerate}

    We omit the upper index $\Por$ in $\mu^{\Xi, \bbP}_{\bar I}$ and $M^{\Xi, \bbP}_{\bar{I}}$ when the context is clear. Notice that $\Int^\Por(Q) = \inf\{M^{\Xi_0}_{\bar I}(Q)\colon \bar I\in \bfI^{\fin}_{K_0}\}$ when $K_0$ is a singleton and $\Xi_0$ is the unique probability fam on $\pts(K_0)$.
\end{definition}

By the definition, it is clear that $\max \left \{ \frac{\Xi(\{ k \})}{\vert I_{k} \vert} \colon k \in K  \right\}$ is a lower bound of  $\mu_{\bar{I}}^{\Xi}(\bar{q})$, and this value is reached when $\ran(\bar{q})$ is an antichain in $\bbP$ and $\bar{q}$ has no repetitions (see \autoref{i29.3}). On the other hand, $\mu_{\bar{I}}^{\Xi}(\bar{q}) \leq 1$ for any $\bar{q} \in \power{W}{Q}$. As a consequence, $ M_{\bar{I}}^{\Xi}(Q) \in \left [\max \left \{ \frac{\Xi(\{ k \})}{\vert I_{k} \vert} \colon k \in K  \right \} , 1 \right].$ In the following results, we are going to study some combinatorial properties of $Q$ to reach the extreme values of $M_{\bar{I}}^{\Xi}(Q)$. Let us start with the lower bound.

\begin{lemma}\label{i29.3}
    Within the context of \autoref{i20}, if $Q$ has at least $\vert W \vert$-many incompatible conditions then $M_{\bar{I}}^{\Xi}(Q) =  \max \left\{  \frac{\Xi(\{ k \})}{\vert I_{k} \vert} \colon k \in K \right \}$. The converse  holds when the fam $\Xi$ is uniform, i.e.\ $\Xi(\{k\})=\frac{1}{|K|}$ for all $k\in K$. 
\end{lemma}

\begin{PROOF}[\textbf{Proof}]{\ref{i29.3}}    
    Assume that $Q$ has at least $\vert W \vert$-many incompatible conditions, so we can find some $\bar{q} \in {}^{W} Q$ without repetitions such that $\ran(\bar{q})$ is an antichain in $\bbP$. It is clear that $\mu_{\bar{I}}^{\Xi}(\bar{q}) = \max \left\{  \frac{\Xi(\{ k \})}{\vert I_{k} \vert} \colon k \in K \right \}$, and since this is a lower bound of $M_{\bar{I}}^{\Xi}(Q)$, we are done. 
    
    To show the converse, assume that $\Xi$ is uniform and  $Q$ has ${<} \vert W \vert$ incompatible elements and let $\bar{q} \in \power{W}{Q}$. Then, there are $\ell_0\neq\ell_1$ in $W$ such that $q_{\ell_{0}}$ and $q_{\ell_{1}}$ are compatible, hence we can find $k_0,k_1\in K$ such that $\ell_e\in I_{k_e}$ for $e \in \{0, 1\}$ and $\mu^\Xi_{\bar I}(\bar q) \geq \frac{\Xi(\{k_0\})}{|I_{k_0}|} + \frac{\Xi(\{k_1\})}{|I_{k_1}|} > \max \left \{ \frac{\Xi(\{ k \})}{\vert I_{k} \vert} \colon k \in K \right \}$, where the last inequality holds because $\Xi$ is uniform. Thus, since $\bar{q}$ is arbitrary, $M_{\bar{I}}^{\Xi}(Q) > \max \left\{  \frac{\Xi(\{ k \})}{\vert I_{k} \vert} \colon k \in K \right \}.$
\end{PROOF}

Now, we deal with the upper bound. 

\begin{lemma}\label{i20.7}
     In the context of \autoref{i20}, 
    $M_{\bar{I}}^{\Xi}(Q) = 1$ iff $Q$ is $\vert W_{0} \vert$-linked, where $W_{0} \coloneqq \bigcup \{ I_k \colon k\in K,\ \Xi(\{k\})\neq 0\}$. 
\end{lemma}

\begin{PROOF}[\textbf{Proof}]{\ref{i20.7}}
    Let $Q \subseteq \bbP$ and assume that $M_{\bar{I}}^{\Xi}(Q) = 1$. Let $Q' \subseteq Q$ be of size ${\leq}  \vert W_{0} \vert$. Fix $\bar{q} \in {}^{W} Q$ such that $\bar{q} \rest W_{0}$ lists $Q'$. Since $M_{\bar{I}}^{\Xi}(Q) = 1,$ in particular we have that $\mu_{\bar{I}}^{\Xi}(\bar{q}) = 1$, hence there exists a $q \in \bbP$ such that, for any $k \in  K$ with $\Xi(\{ k \}) \neq 0$ and any $\ell \in I_{k},$ $q \leq q_{\ell}$. Consequently, for every $\ell \in W_{0},$ $q \leq q_{\ell.}$ Thus, $Q$ is $\vert W_{0} \vert$-linked.  
    
    To prove the converse, assume that $Q$ is $\vert W_{0} \vert$-linked. Setting $K_{0} \coloneqq \{ k \in K \colon \Xi(\{ k \}) \neq 0 \},$ for any $\bar{q} \in \power{W}{Q}$, $\mu_{\bar{I}}^{\Xi}(\bar{q}) = \mu_{\bar{I} \rest K_{0}}^{\Xi \rest \pts(K_{0})}(\bar{q} \rest W_{0})$. Therefore, it is enough to prove that, for all $\bar{q} \in {}^{W_0}{Q}$, $\mu_{\bar{I \rest K_{0}}}^{\Xi \rest \pts(K_{0})}(\bar{q}) =\allowbreak 1$. Since $Q$ is $\vert W_{0} \vert$-linked and $\vert \ran(\bar{q}) \vert \leq \vert W_{0} \vert$, there exists some $q \in \bbP$ such that, for any $k \in K_{0}$ and any $\ell \in I_{k}$, $q \leq q_{\ell}$. As a consequence, 
    $$ \mu_{\bar{I} \rest K_{0}}^{\Xi \rest \pts(K_{0})}(\bar{q}) \geq  \sum_{k \in K_{0}} \frac{ \vert \{ \ell \in I_{k} \colon q \leq q_{\ell} \} \vert }{\vert I_{k} \vert} \Xi( \{ k \} ) = \sum_{k \in K_{0}} \Xi(\{  k \}) = 1.$$ 
    Thus, $M_{\bar{I}}^{\Xi}(Q) = 1$. 
\end{PROOF}

Fixing $K$ and $\Xi$, and varying on all possible partitions in $\bfI^\fin_K$, we get a combinatorial property to characterize when the upper bound of $M_{\bar{I}}^{\Xi}(Q)$ is reached. 

\begin{corollary}\label{i20.5}
    Let $\bbP$ be a forcing notion and $Q \subseteq \bbP$. Then,  the following statements are equivalent:
    \begin{enumerate}[label = \normalfont (\roman*)]
        \item\label{i20.5.a} $Q$ is centered. 
        
        \item\label{i20.5.b} For any  $\bar{I} \in \bfI_{K}^{\fin}$, $M_{\bar{I}}^{\Xi}(Q) = 1$.

        \item\label{i20.5.c} $\Int^{\bbP}(Q) = 1$. 
    \end{enumerate}
\end{corollary}

The following is a very interesting characterization of the intersection number.

\begin{theorem}\label{i20.7.1}
    Let $\bbP$ be a forcing notion and $Q \subseteq \bbP$. Then, $\Int^\Por(Q) = \inf_{\bar{I} \in \bfI_{K}^{\fin}} M_{\bar{I}}^{\Xi, \bbP}(Q)$.
\end{theorem}

\begin{PROOF}[\textbf{Proof}]{\ref{i20.7}}        
    Let $\bar q = \la q_\ell\colon \ell\in W\ra \in \power{W}{Q}$ and $\bar I\in\bfI^\fin_K$. Let $\Sigma$ as in \autoref{s14} and, for $\varp>0$, find $f\colon \Sigma\to[0,1]_\QQQ$ as the same result applied to $\delta \coloneqq \Int(Q)$. For any $q\in \Por$ there is some $q'\leq q$ in $\Por$ such that, for all $\ell\in W_{\bar I}$, either $q'\leq q_\ell$ or $q'\perp q_\ell$. Then, there is some $\sigma'\in \Sigma$ such that $\sigma'(\ell)=0$ iff $q'\leq q_\ell$. Hence, by using $c_k\coloneqq \frac{\Xi(\{k\})}{|I_k|}$,
    $$\sum_{k\in K}c_k \vert \{\ell\in I_k\colon q\leq q_\ell\} \vert \leq \sum_{k\in K}c_k \vert \{\ell\in I_k\colon q'\leq q_\ell\} \vert = \sum_{k\in K} c_k\vert\{\ell\in I_k\colon \sigma'(\ell) = 0\}\vert.$$
    Therefore,
    \begin{align*}
        \mu_{\bar{I}}^{\Xi}(\bar q) & = \max\set{\sum_{k\in K} c_{k} |\set{\ell \in I_k}{\sigma(\ell)=0}|}{\sigma\in \Sigma}\\ 
         & \geq \sum_{\sigma\in\Sigma} \left( \sum_{k\in K} c_{k} |\set{\ell \in I_k}{\sigma(\ell)=0}| \right) f(\sigma) =  \sum_{\sigma\in\Sigma}\sum_{k\in K}\sum_{\substack{\ell\in I_k\\ \sigma(\ell)=0}} c_{k} f(\sigma) \\
         & = \sum_{k\in K} \sum_{\ell\in I_k} c_{k}  \sum_{\substack{\sigma\in\Sigma\\ \sigma(\ell) = 0}} f(\sigma)  \geq [\Int(Q) - \varp] \sum_{k \in K} \sum_{\ell \in I_{k}} c_{k}  = \Int(Q) - \varp.
    \end{align*}
    Since $\varp$ is arbitrary, $\mu_{\bar{I}}^{\Xi}(\bar{q}) \geq \Int(Q)$. Therefore, since $\bar{q}$ is arbitrary, it follows that $M_{\bar{I}}^{\Xi}(Q) \geq \Int(Q)$. Thus, $\Int(Q) \leq \inf_{\bar{I} \in \bfI_{K}^{\fin}} M_{\bar{I}}^{\Xi}(Q)$. 

    On the other hand, to show the converse inequality, let $n \in \omega \setminus \{ 0 \}$ and $\bar{p} \in {}^{n} Q$. Consider $\bar{J} \in \bfI_{K}^{\fin}$ such that, for any $k \in K,$ $\vert J_{k} \vert = n$ and let $f_{k} \colon J_{k} \to n$ a bijective function. Set $W \coloneqq W_{\bar J}$. Let us define $\bar{q} = \langle q_{\ell} \colon \ell \in W \rangle \in \power{W}{Q}$ such that for $k \in K$ and $\ell \in I_{k}$, $q_{\ell} \coloneqq p_{f_{k}(\ell)}$.  Then, we can find some $q \in \bbP$ such that 
    \begin{equation*}
        \begin{split}
            \mu^{\Xi,\Por}_{\bar J}(\bar q) & = \int_{K} 
         \frac{ \vert \{ \ell \in J_{k} \colon q \leq q_{\ell} \} \vert }{\vert J_{k} \vert} d \Xi(k) 
          = \frac{1}{n} \int_{K} \vert \{ i < n \colon  q \leq p_{i} \} \vert d \Xi(k)\\
          & = \frac{1}{n} \vert \{ i < n \colon q \leq p_{i} \} \vert \int_{K} d \Xi(k)  = \frac{1}{n} \vert \{ \ell < n \colon q \leq p_{i} \} \vert \leq \frac{i_{\ast}^{\bbP}(\bar{p})}{n}. 
        \end{split}
    \end{equation*}

    Therefore, $ \inf_{\bar{I} \in \bfI_{K}^{\fin}}  M_{\bar{I}}^{\Xi, \bbP}(Q) \leq M_{\bar{J}}^{\Xi, \bbP}(Q) \leq \mu^{\Xi,\Por}_{\bar J}(\bar q) \leq  \frac{i_{\ast}^{\bbP}(\bar{p})}{n}.$ Finally, since $n$ is arbitrary, it follows that $ \inf_{\bar{I} \in \bfI_{K}^{\fin}}  M_{\bar{I}}^{\Xi, \bbP}(Q) \leq \Int^{\bbP}(Q).$ 
\end{PROOF}

As a consequence, we can generalize~\autoref{s9} in the context of $M_{\bar{I}}^{\Xi}(Q)$.  

\begin{corollary}\label{i9.m}
    Let $\cB$ be a Boolean algebra, $\Xi_{0}$ a probability fam on it, $\bar{I} = \langle I_{k} \colon k \in K \rangle \in \bfI_{K}^{\fin}$, 
    and $\delta \in [0, 1].$ If $Q \coloneqq \{ b \in \cB \colon \Xi_{0}(b) \geq \delta \},$ then $M_{\bar{I}}^{\Xi, \cB^+}(Q) \geq \delta.$
\end{corollary}

We close this section by showing a relation between $M_{\bar{I}}^{\Xi}(Q)$ and $M_{\bar{J}}^{\Xi_{h}}(Q)$ where $\Xi_{h}$ is as in \autoref{m37.1}. To fix some notation, if $h \colon K \to L$ is a function, for any $l \in L$, denote $K_{l} \coloneqq h^{-1}[\{ l \}]$.

\begin{lemma}\label{i20.f}
    Let $K$ and $L$ be finite sets. Assume that $h \colon K \to L$ is a surjective finite-to-one function.  Let $\bar{I} = \langle I_{k} \colon k \in K \rangle \in \bfI_{K}^{\fin}$ and for any $l \in L$, define $J_{l} \coloneqq \bigcup \{ I_{k} \colon k \in K_{l} \}$ and $\bar{J} \coloneqq \langle J_{l} \colon l \in L \rangle \in \bfI_{L}^{\fin}$. 
    If, for any $k, k' \in K$, $h(k) = h(k')$ implies $\Xi(\{ k \}) = \Xi( \{ k' \})$ and $\vert I_{k} \vert = \vert I_{k'} \vert $, then, for every $\bar{q} \in \power{W}{Q}$, $\mu_{\bar{I}}^{\Xi}(\bar{q}) \geq \mu_{\bar{J}}^{\Xi_{h}}(\bar{q})$. As a consequence, $M_{\bar{I}}^{\Xi}(Q) \geq M_{\bar{J}}^{\Xi_{h}}(Q)$. In particular, if $h$ is a bijective function, then $M_{\bar{I}}^{\Xi}(Q) = M_{\bar{J}}^{\Xi_{h}}(Q)$.
\end{lemma}

\begin{PROOF}[\textbf{Proof}]{\ref{i20.f}}
    Let $\bar{q} = \langle q_{\ell} \colon \ell \in W \rangle \in \power{W}{Q}$. By \autoref{i20}, we can find a $q \in \bbP$ such that $\mu_{\bar{J}}^{\Xi_{h}}(\bar{q}) = \sum_{l \in L} \frac{\vert \{  \ell \in J_{l} \colon q \leq q_{\ell} \} \vert}{\vert J_{l} \vert} \Xi_{h}(\{ l \})$. Since $h(k) = h(k')$ implies $\Xi(\{ k \}) = \Xi(\{ k' \})$ and $\vert I_{k} \vert = \vert I_{k'} \vert$, we have that $\frac{\Xi(K_{h(k)})}{\vert J_{h(k)} \vert} \leq \frac{\Xi(\{k \})}{\vert I_{k} \vert}$ for any $k \in K$.  As a consequence, 
    \begin{equation*}
        \begin{split}
            \mu_{\bar{J}}^{\Xi_{h}}(\bar{q}) & = \sum_{l \in L} \frac{\vert \{  \ell \in J_{l} \colon q \leq q_{\ell} \} \vert}{\vert J_{l} \vert} \Xi_{h}(\{ l \})  = \sum_{l \in L} \left( \sum_{k \in K_{l}} \frac{ \vert \{ \ell \in I_{k} \colon q \leq q_{\ell} \} \vert}{\vert J_{l} \vert}  \right) \Xi(K_{l}) \\
            & = \sum_{k\in K}\frac{ \vert \{ \ell \in I_{k} \colon q \leq q_{\ell} \} \vert}{\vert J_{h(k)} \vert} \Xi(K_{h(k)})  \leq \sum_{k \in K} \frac{\vert \{  \ell \in I_{k} \colon q \leq q_{\ell} \}  \vert}{\vert I_{k} \vert} \Xi(\{ k \}) \leq \mu_{\bar{I}}^{\Xi}(\bar{q}). 
        \end{split}
    \end{equation*}

    Thus, since $\bar{q}$ is arbitrary, it follows that $M_{\bar{J}}^{\Xi_{h}}(Q) \leq M_{\bar{I}}^{\Xi}(Q)$. Finally, when $h$ is bijective, to get the converse inequality it is enough to use $h^{-1}$. 
\end{PROOF}

\section{FAM-linked sets}\label{sec:famlim}

In this section, we introduce the notion of $(\Xi, \bar{I}, \varp)$-linkedness and investigate their basic properties and characterizations. 
For instance, $(\Xi, \bar{I}, \varp)$-linkedness has connections with the intersection number and with forcing with ultrafilter limits from~\cite{GMS}, as well as with the notion of $\Fr$-linkedness (\autoref{n73}), which is a general notion for not adding dominating reals (see~\autoref{ExmPrs}~\ref{ExmPrsb}). This new linkedness notion is the main ingredient to iterate with fams. 

\subsection{Strong fam limits as a linkedness property}\label{i0}
\ 

Below, we define the notion of \emph{$(\Xi, \bar{I}, \varp)$-linkedness}, which 
generalizes the notion of \emph{strong fam limits for intervals} from~\cite[Def.~1.7]{KST} and \cite[Def.~1.10]{KST}. The cited reference only considers free probability fams on $\pts(\omega)$ and the characterization presented in \autoref{i6}, but we allow probability fams on any $\pts(K)$.

\begin{definition}\label{i2}
    Let $K\neq \emptyset$, $\Xi$ a probability fam on $\Pwf(K)$, 
    $\bar{I} = \langle  I_{k} \colon k \in K \rangle \in \bfI^{\fin}_K$, $W\coloneqq W_{\bar I} = \bigcup_{k\in K}I_k$, $\varp_{0} \in [0, 1)$ and let $\bbP$ be a forcing notion. 
    Say that $Q$ is \emph{$(\Xi, \bar{I}, \varp_{0})$-linked (in $\Por$)} if $Q \subseteq \bbP$ and there is a function $\lim^{Q,\Xi, \bar{I}, \varp_{0}} \colon \power{W}{Q} \to \bbP$, and a $\bbP$-name $\dot{\Xi}^{\ast}$ of a probability fam on $\calP(K)$ extending $\Xi$ such that, for any $\bar{q} \in \power{W}{Q},$  
    $$ \limit^{Q,\Xi, \bar{I}, \varp_{0}}( \bar{q}) \Vdash \text{``} \int_{K} \frac{\vert \{ \ell \in I_{k} \colon q_{\ell} \in \dot{G}_{\bbP} \} \vert}{\vert I_{k} \vert}   d \dot{\Xi}^{\ast}(k) \geq 1 - \varp_{0}\text{''}.$$ 
    We omit some or all of the superscripts ``$Q,\Xi, \bar{I}, \varp_{0}$'' in ``$\limit^{Q,\Xi, \bar{I}, \varp_{0}}$'' when the context is clear.
\end{definition}

\begin{remark}\label{i2.2}
In \autoref{i2}:
\begin{enumerate}[label=\normalfont (\arabic*)]
    \item\label{i2.2-1} In general, $K$ is determined by $\Xi$ as the top element of $\dom\Xi$, and $W = W_{\bar I}$ is determined by $\bar{I}$.

    \item\label{i2.2-2} If $\bar{J} = \langle J_{k} \colon k \in K \rangle \in \bfI_{K}^{\fin}$ and $|I_{k}|=|J_k|$ for all $k \in K$, then $(\Xi,\bar{I},\varp_0)$-linkedness is equivalent to $(\Xi,\bar{J}, \varp_0)$-linkedness (see \autoref{i24.01}). For this reason, we could have used a function $\bfi \colon K \to \omega\menos\{0\}$ instead of $\bar{I}$ (and $W$) in the definition. 

    \item The first assumptions imply that there is a finite-to-one map from $W$ onto $K$. Hence, $K$ is infinite iff $W$ infinite, in which case $|K|=|W|$. 

    \item The value $\varp_0 = 0$ is not considered in~\cite{KST}. We realized that it can be included in this framework and that it helps to describe sets with ultrafilter limits. See more details in \autoref{i24},~\autoref{i25} and\ref{i26}.
\end{enumerate}
\end{remark}

\begin{remark}\label{i2.4}
    In our applications (\autoref{k62}, \autoref{presanti} and~\ref{appl}), $K=W=\omega$ and $\Xi$ is free. If convenient, there is no problem if the reader only considers this case throughout the paper (except when we deal with finite $K$).
\end{remark}

Unless otherwise specified, until the end of this section, we fix a non-empty set $K$, $\bar{I} \in \bfI_{K}^{\fin}$, $W \coloneqq W_{\bar{I}}$, a probability fam $\Xi$ on $\calP(K)$, a forcing notion $\bbP$, $Q \subseteq \bbP,$ and $\varp_{0} \in [0, 1).$ 

A first easy example of a $(\Xi, \bar{I}, \varp_{0})$-linked subset is a singleton:  

\begin{example}\label{i4}
    For any $p \in \bbP$, $\{ p \}$ is $(\Xi, \bar{I}, \varp_{0})$-linked.
    
    Indeed, let $Q_0 \coloneqq \{ p \}$.  For $\bar{q} \in {}^{W} \! Q_0$ define $\limit \colon {}^{W} \! Q_0 \to \bbP$ in the natural way: $\limit(\bar{q}) \coloneqq p$. Also, consider $\dot{\Xi}^{\ast}$ as a $\bbP$-name of a fam extending $\Xi$, which is possible because, in any generic extension $V[G]$, we can extend $\Xi$ to a fam with domain $\calP(K) \cap V[G]$. It is clear that, if $G$ contains $p$, then $\displaystyle\int_{K} \frac{\vert \{ \ell \in I_{k} \colon q_{\ell} \in G_{\bbP} \} \vert}{\vert I_{k} \vert}   d \dot{\Xi}^{\ast}(k) = 1 \geq 1 - \varp_{0}.$
\end{example}

\begin{example}\label{i4.2}
    Let $\la\cB,\mu\ra$ be a measure algebra, i.e.\ $\mu$ is a probability measure on $\cB$ such that the only measure zero point is $0_\cB$. If $s\in\cB^+=  \cB\menos\{0_\cB\}$ then $\{ b\in\cB \colon \mu(b\cap s)\geq (1-\varp_0)\mu(s)\}$ is $(\Xi,\bar I,\varp_0)$-linked when $\Xi$ is free. The case $\varp_0=0$ is trivial (and does not need that $\Xi$ is free) since it is enough to set any limit as $s$. The case $\varp_0\in(0,1)$ is more complicated. It is solved in~\cite{MUR23} (using the characterization in \autoref{i6}) by extending ideas from Shelah~\cite{Sh00} in the case of random forcing. See more in~\autoref{i15.2}.
\end{example}

\subsection{Basic properties}\label{sub:fambasic}
\ 

We present below many basic properties of $(\Xi,\bar I, \varp_0)$-linked sets. 
First note that the upwards closure of a linked set is also linked.

\begin{lemma}\label{i6.1}
    $Q$ is $(\Xi, \bar{I}, \varp_0)$-linked iff  $Q^{\uparrow}$ is. 
\end{lemma}

\begin{PROOF}[\textbf{Proof}]{\ref{i6.1}}
    The direction from right to left is clear. Conversely, assume that $Q$ is $(\Xi, \bar{I}, \varp_0)$-linked. Define $\lim^{Q^{\uparrow}} \colon  {}^{W} \! (Q^{\uparrow}) \to \bbP$ such that, for $\bar{q} = \langle q_{\ell} \colon \ell \in W \rangle \in  {}^{W} \!(Q^{\uparrow}),$  $ \lim^{Q^{\uparrow}}(\bar{q}) \coloneqq \lim^{Q}(\bar{p}^{\bar{q}})$ where $\bar{p}^{\bar{q}} = \langle p_{\ell} \colon \ell \in W \rangle \in  {}^{W} \!  Q$ is chosen such that, for any $\ell \in W,$ $p_{\ell} \leq q_{\ell}.$ Straightforward calculations show that this limit is as required. 
\end{PROOF}

In the following results, we deal with some combinatorial facts of $(\Xi, \bar{I}, \varp_{0})$-linked sets. 

\begin{lemma}\label{i6.1.1}
    If $Q$ is $(\Xi, \bar{I}, \varp_{0})$-linked then, for any $\bar q \in {}^W Q,$  $\mathrm{ran}(\bar{q})$ is predense below $\lim^{\Xi}(\bar q)$.
\end{lemma}

\begin{PROOF}[\textbf{Proof}]{\ref{i6.1}}
    Fix $\dot \Xi^{\ast}$ as in \autoref{i2}. Let $\bar{q} = \langle q_{\ell} \colon \ell \in W \rangle \in \power{W}{Q}$ and let $p \leq \lim^{\Xi}(\bar{q}).$ Then,  $$ p \forces \text{``} \int_{K} \frac{\vert \{ \ell \in I_{k} \colon q_{\ell} \in \dot{G}_{\bbP} \} \vert}{\vert I_{k} \vert} d  \dot{\Xi}^{\ast}(k) \geq 1 - \varp_{0} \text{''}. $$  As a consequence, there are $r \in \bbP,$ $k \in K$ and $\ell \in I_{k}$ such that $r \leq p$ and $r \forces$``$q_{\ell} \in \dot{G}_{\bbP}$'', that is, $r \leq^{\bullet} q_{\ell}$, hence $p \parallel q_{\ell}.$ Thus, $\ran(\bar{q})$ is predense below $\lim^{\Xi}(\bar{q}).$  
\end{PROOF}

If $\Xi$ is suitable, a $(\Xi, \bar{I}, \varp_{0})$-linked set cannot contain large antichains. 

\begin{lemma}\label{i7}
    If $Q$ is $(\Xi, \bar{I}, \varp_{0})$-linked and, for any $k \in K,$ $\frac{\Xi(\{ k \})}{|I_k|} < 1 - \varp_{0},$ then $Q$ cannot contain antichains of size $\vert W \vert.$  
\end{lemma}

\begin{PROOF}[\textbf{Proof}]{\ref{i7}}
    Assume that $Q$ is  $(\Xi, \bar{I}, \varp_{0})$-linked. Towards contradiction, suppose that $A \subseteq Q$ is an antichain in $\bbP$ such that $\vert A \vert = \vert W \vert$. So we can define $\bar{q} = \langle q_{\ell} \colon \ell \in W \rangle \in \power{W}{Q}$ such that $\ran(\bar{q}) = A$ and, for any $i \neq j$ in $W$, $q_{i}  \perp  q_{j}.$
    
    Fix $\dot{\Xi}^{\ast}$ as in \autoref{i2}. Work in $V[G]$ for some generic filter $G$ containing $\lim^{\Xi}(\bar{q}).$ By \autoref{i6.1.1}, $A$ is a maximal antichain below $\lim(\bar q)$, so $\vert \{ \ell \in I_{k_{0}} \colon q_{\ell} \in  G \} \vert = 1$ for some $k_{0} \in K$ and, for any $k \in K \menos \{ k_{0} \},$ $\vert \{ \ell \in I_{k} \colon q_{\ell} \in G \} \vert = 0.$ As a consequence, we have that:  $$1 - \varp_{0} \leq \int_{K} \frac{ \vert \{ \ell \in I_{k} \colon q_{\ell} \in G \} \vert }{\vert I_{k} \vert} d \Xi^{\ast}(k) = \int_{\{ k_{0} \}} \frac{1}{\vert I_{k_{0}} \vert} d \Xi^{\ast} (k) =  \frac{\Xi( \{ k_{0} \})}{\vert I_{k_{0}} \vert} < 1 - \varp_{0},$$ which is a contradiction.
\end{PROOF}

Measure zero sets do not affect the linkedness property.

\begin{lemma}\label{i6.2}
    Assume that $K_0\subseteq K$ and $\Xi(K_0)=0$. Then the set $Q$ is $(\Xi,\bar{I},\varp_0)$-linked in $\Por$ iff it is $(\Xi\frestr\pts(K\menos K_0),\bar{I}\frestr (K\menos K_0),\varp_0)$-linked.
\end{lemma}

\begin{PROOF}[\textbf{Proof}]{\ref{i6.2}}
    We can assume that $Q\neq\emptyset$,\footnote{According to \autoref{i2}, the empty set is $(\Xi,\bar I,\varp_0)$-linked.} so choose some $q'\in Q$. Let $\Xi'\coloneqq\Xi\frestr\pts(K\menos K_0)$ and $\bar J\coloneqq \bar{I}\frestr (K\menos K_0)$. 
    If $Q$ is $(\Xi,\bar{I},\varp_0)$-linked witnessed by $\dot\Xi^*$ and $\lim^\Xi$, then $Q$ is $(\Xi',\bar{J},\varp_0)$-linked witnessed by a $\Por$-name of $\dot\Xi^*\frestr \pts(K\menos K_0)$ and the limit function $\limit^{\Xi'}\colon {}^{W_{\bar{J}}} Q\to \Por$ defined by $\limit^{\Xi'}(\bar q) \coloneqq \limit^\Xi(\bar q')$ where $q'_\ell\coloneqq q_\ell$ for $\ell\in W_{\bar J}$, and $q'_\ell\coloneqq q'$ otherwise.

    Conversely, if $Q$ is $(\Xi',\bar{J},\varp_0)$-linked witnessed by $\dot \Xi^\bullet$ and $\limit^{\Xi'}$,then $Q$ is $(\Xi,\bar{I},\varp_0)$-linked witnessed by a $\Por$-name $\dot \Xi^*$ of the fam on $\pts(K)$ defined by $\dot \Xi^*(B)\coloneqq \dot \Xi^\bullet (B\menos K_0)$, and the limit function defined by $\lim^\Xi(\bar q)\coloneqq \lim^{\Xi'}(\bar q\frestr (K\menos K_0))$.
\end{PROOF}

The notion of \emph{Fr\'echet-linkedness} was introduced by the second author in~\cite{mejiavert}, who proved that no $\sigma$-Fr\'echet-linked forcing notions add dominating reals. Afterward, this result was enhanced in~\cite{BCM}. 

\begin{definition}[{\cite[Def.~3.24]{mejiavert}}]\label{n73}
    Let  $\mu$ an infinite cardinal.
    
    \begin{enumerate}[label=\rm(\arabic*)]
        \item A set $Q_0\subseteq\Por$ is \emph{Fr\'echet-linked in $\bbP$}, abbreviated by ``\emph{$\Fr$-linked}'', if for any sequence $\bar{p} = \langle p_{n} \colon n < \omega \rangle \in {}^{\omega} Q_0$ there exists some $q \in \bbP$ such that $q \forces_{\bbP} \text{``} \vert \{ n < \omega \colon p_{n} \in \dot{G} \} \vert = \aleph_{0}\text{''}.$

        \item $\bbP$ is \emph{$\mu$-$\rm{Frechet}$-linked}, abbreviated by ``\emph{$\mu$-$\Fr$-linked}'', if $\Por = \bigcup_{\alpha < \mu} Q_{\alpha}$ for some sequence $\langle Q_{\alpha} \colon \alpha < \mu \rangle$  of $\Fr$-linked subsets of $\bbP.$

        \item $\bbP$ is \emph{$\mu$-$\Fr$-Knaster} if, for any $A \in [\bbP]^{\mu},$ there is some $\Fr$-linked $A' \in [A]^{\mu}.$ \index{$\mu$-$\Fr$-Knaster}
    \end{enumerate}
\end{definition}

This notion behaves well for not adding dominating reals as illustrated in~\autoref{ExmPrs}~\ref{ExmPrsb}. On the other hand, as a consequence of~\autoref{i7}, when $K=W=\omega,$ $\Xi$ is a free fam and $Q$ is $(\Xi, \bar{I}, \varp_{0})$-linked in $\Por$, then $Q$ cannot contain infinite antichains, which is, in particular, a property of Fr\'echet-linked subsets (see~\cite[Ex.~5.3]{mejiavert} and~\cite[Rem.~3.1~(1)]{BCM}). Motivated by this, we present a result that shows that our linkedness property is stronger than Fr\'echet-linkedness.

\begin{lemma}\label{i8} 
    If $\Xi$ is free, $K=\omega$ and $Q$ is $(\Xi, \bar{I}, \varp_{0})$-linked, then $Q$ is $\Fr$-linked. 
\end{lemma}

\begin{PROOF}[\textbf{Proof}]{\ref{i8}}
    Without loss of generality, we can assume that $\bar{I}$ a partition of $\omega$ into finite sets, i.e.\ $W=\omega$. 
    Let  $\bar{q} = \langle q_{\ell} \colon \ell < \omega \rangle \in {}^{\omega}Q.$
    We will prove that there exists some $q \in \bbP$ such that $q \forces_{\bbP} ``\exists^{\infty}\ell < \omega\, (q_{\ell} \in \dot{G}_{\bbP})"$. To see this, it suffices to prove  that $\lim(\bar{q})$ forces this. Let $G$ be a $\bbP$-generic filter over $V$ such that $\lim(\bar{q}) \in G.$  Fix $\Xi^{\ast}$ as in \autoref{i2}.
    
    Working in $V[G]$, define $f \colon \omega \to \RRR$ by $f(k) \coloneqq \frac{\vert \{ \ell \in I_{k} \colon q_{\ell} \in G_{\bbP} \} \vert}{\vert I_{k} \vert}$ for all $k < \omega$. Towards contradiction, suppose that $\{ \ell < \omega \colon q_{\ell} \in G \}$ is finite. Since $\Xi$ is a free fam, this implies that $\Xi(D) = 0$ where $D \coloneqq \{ k < \omega \colon f(k) \neq 0 \}$. Therefore, $\int_{\omega} f d  \Xi^{\ast} = 0$, which contradicts  \autoref{i2}. 
\end{PROOF}

We can generalize \autoref{i8} for arbitrary $\calP(K)$ using more natural hypotheses.

\begin{theorem}\label{i9}
    Assume that $\kappa$ is an infinite cardinal, $|K|=\kappa$, the $\Xi$-measure of any subset of cardinality ${<} \kappa$ is zero, and $\bbP$ forces that $[K]^{<\kappa}\cap V$ is cofinal in $[K]^{<\kappa}$.\footnote{For example, when $\bbP$ is ${<}\kappa$-distributive, or $\kappa$-cc with $\kappa$ regular.}  
    If $Q$ is $(\Xi, \bar{I}, \varp_{0})$-linked, then $Q$ is $\Fr_{\kappa}$-linked, that is, for any sequence $\bar{q} = \langle q_{\ell} \colon \ell < \kappa \rangle \in {}^{\kappa} Q$, there exists some condition $q \in \bbP$ such that $ q \forces \text{``} \,  \vert \{ \ell <\kappa \colon q_{\ell} \in \dot{G}_{\bbP} \} \vert = \kappa\text{''}.$ 
\end{theorem}

In \autoref{i5} and~\ref{i5.1} below, we study the behavior of $(\Xi, \bar{I}, \varp_{0})$-linked subsets under complete and dense embeddings. If $\iota \colon \bbP \to \bbQ$ is a complete embedding, then  $\iota^{\ast} \colon V^{\bbP} \to V^{\bbQ}$ denotes the canonical transformation of names. 

\begin{lemma}\label{i5}
    Let $\bbQ$ be a forcing notion and let $\iota \colon \bbP \to \bbQ$ be a complete embedding. If $Q$ is $(\Xi, \bar{I}, \varp_{0})$-linked in $\Por$, then  so is $\iota [ Q ]$ in $\bbQ.$ 
\end{lemma}

\begin{PROOF}[\textbf{Proof}]{\ref{i5}}
    Assume that $Q \subseteq \bbP$ is $(\Xi, \bar{I}, \varp_{0})$-linked witnessed by  $\lim^{Q}$ and $\dot{\Xi}_{Q}^{\ast}$. Consider $S \coloneqq \iota [ Q ]$ and, for any $\bar{r} = \langle r_{\ell} \colon \ell \in W \rangle \in  \power{W}{S},$ define $\limit^{S}(\bar{r}) \coloneqq \iota \left(\limit^{Q}(\bar{q}_{\bar{r}}) \right)$ where, for any $\ell \in W$, $q_{\bar{r}, \ell} \in Q$ and $\iota(q_{\bar{r}, \ell}) = r_{\ell}$, which is possible to pick by the definition of $S.$ Also define $\dot{\Xi}_{S}^{\ast}$ as a $\bbQ$-name of a fam extending $\iota^{\ast}(\Xi_{Q}^{\ast})$, which is possible because, in any $\Qor$-generic extension $V[G]$, we can extend $\Xi_{Q}^{\ast}$ to fam with domain $\calP(K) \cap V[G]$.\footnote{Notice that, since $\iota$ is a complete embedding, in general we have that $\iota^{\ast}(\Xi_{Q}^{\ast})$ is just a name of the fam $\Xi^*_Q$,  whose domain is $\pts(K)\cap V[\iota^{-1}[G]]$.}
    
    Now, let $\bar{r} = \langle r_{\ell} \colon \ell \in W \rangle \in \power{W}{S}$. By \autoref{i2},  $$ \limit^{Q}(\bar{q}_{\bar{r}}) \Vdash_\bbP \text{``} \int_{K} \frac{ \vert  \{ \ell \in I_{k} \colon q_{\bar{r}, \ell} \in \dot{G}_{\bbP} \}\vert}{\vert I_{k} \vert} d \dot{\Xi}_{Q}^{\ast}(k) \geq 1 - \varp_{0} \text{''}. $$  
    
    Finally, since $\iota_{\ast}(\dot{G}_{\bbP}) = \iota^{-1}[\dot{G}_{\bbQ}]$, we can apply the absoluteness of the integral (\autoref{t102}), to get that
    $$ \limit^{S}(\bar{r}) \Vdash_{\bbQ} \text{``} \int_{K} \frac{ \vert  \{ \ell \in I_{k} \colon r_{\ell} \in \dot{G}_{\bbQ} \}\vert}{\vert I_{k} \vert} d \dot{\Xi}_{S}^{\ast}(k) \geq 1 - \varp_{0} \text{''},$$   
    which proves that $S$ is $(\Xi, \bar{I}, \varp_{0})$-linked.
\end{PROOF}


\begin{lemma}\label{i5.1}
    Let $\bbQ$ be a forcing notion and $\iota \colon \bbP \to \bbQ$ be a dense embedding. If $S$ is $(\Xi, \bar{I}, \varp_{0})$-linked in $\bbQ$, then  so is $\iota^{-1} [ S ]$ in $\bbP.$ 
\end{lemma}

\begin{PROOF}[\textbf{Proof}]{\ref{i5.1}}
   Assume that $S \subseteq \bbQ$ is $(\Xi, \bar{I}, \varp_{0})$-linked witnessed by  $\limit^{S}$ and $\dot{\Xi}_{S}^{\ast}.$ Consider $R \coloneqq \iota^{-1}[S]$ and, for $\bar{r} = \langle r_{\ell} \colon r \in W \rangle \in \power{W}{R}$, define $\limit^{R}(\bar{r}) \in \bbP$ such that  $\iota(\limit^{R}(\bar{r})) \leq \limit^{S}( \langle \iota(r_{\ell}) \colon \ell \in W \rangle),$ which is possible because $\iota$ is a dense embedding. 
   Consider $\dot{\Xi}_{R}^{\ast} \in V^{\bbP}$ such that $\Vdash_{\bbQ}$``$\dot{\Xi}_{S}^{\ast} = \iota^{\ast}(\dot{\Xi}_{R}^{\ast})$'', which is a name of a fam on $\pts(K)$ extending $\Xi$. Let $\bar{r} = \langle r_{\ell} \colon \ell \in W \rangle \in \power{W}{R}$. Hence, $\bar{q} \coloneqq \langle \iota(r_{\ell}) \colon \ell \in W \rangle \in \power{W}{S}$ and, by \autoref{i2}, 
   $$ \limit^{S}(\bar{q}) \Vdash_\bbQ \text{``} \int_{K} \frac{\vert \{ \ell \in I_{k} \colon \iota(r_{\ell}) \in \dot{G}_{\bbQ} \} \vert}{\vert I_{k} \vert}   d \dot{\Xi}_{S}^{\ast}(k) \geq 1 - \varp_{0}\text{''}.$$ 
    Since, $\iota(\lim^{R}(\bar{r})) \leq \limit^{S}(\bar{q}),$ it follows that  
    $$ \limit^{R}(\bar{r}) \Vdash_\bbP \text{``} \int_{K} \frac{\vert \{ \ell \in I_{k} \colon r_{\ell} \in \dot{G}_{\bbP} \} \vert}{\vert I_{k} \vert}   d \dot{\Xi}_{R}^{\ast}(k) \geq 1 - \varp_{0}\text{''}.$$ Thus, $R$ is $(\Xi, \bar{I}, \varp_{0})$-linked.
\end{PROOF}

\subsection{The finite case}\ 

When $K$ is finite, our linkedness notion reduces to a property for $M^\Xi_{\bar I}(Q)$ (see \autoref{i20}), which shows more interesting aspects about the intersection number.

\begin{theorem}\label{i21}
    If $K$ is finite,  
    then $Q$ is $(\Xi,\bar{I},\varp_0)$-linked in $\bbP$ iff $M^\Xi_{\bar{I}}(Q) \geq 1-\varp_0$.
\end{theorem}

\begin{PROOF}[\textbf{Proof}]{\ref{i21}}
    First assume that $Q$ is $(\Xi,\bar{I},\varp_0)$-linked in $\bbP$. Notice that, since $\pts(K)$ is finite, $\dot\Xi^*$ from \autoref{i2}  is forced to be equal to $\Xi$ and,
    for any $\bar q \in {}^W Q$, we have that
    \[{\lim}^{\Xi}(\bar q) \Vdash \text{``} 1-\varp_0 \leq \int_K \frac{|\{\ell\in I_k\colon q_\ell \in \dot G_{\bbP}\}|}{|I_k|} d\Xi(k)  \text{''}.\] 
    Since $W$ is finite, we can find some $q'\leq \lim^\Xi(\bar q)$ forcing that, for $k\in K$ and $\ell\in I_k$, $q_\ell\in \dot G$ iff $q'\leq q_\ell$. This $q'$ forces that the integral is below $\mu^\Xi_{\bar I}(\bar q)$, so $\mu^\Xi_{\bar I}(\bar q) \geq 1-\varp_0$. Since $\bar q$ is arbitrary, we can conclude that $M^\Xi_{\bar I}(Q) \geq 1-\varp_0$.

    Now we show the converse. Assume that $M^\Xi_{\bar{I}}(Q) \geq 1-\varp_0$. For $\bar q \in {}^W Q$, we define $\lim^\Xi(\bar q)$ as a condition that forces the maximal possible value of $\int_{K} \frac{|\{\ell\in I_k\colon q_\ell \in \dot G\}|}{|I_k|} d \Xi$. 
    It is not hard to show that this maximal value is $\mu^\Xi_{\bar{I}}(\bar q)$, which is ${\geq} M^\Xi_{\bar{I}}(Q) \geq 1-\varp_0$. Therefore, $\lim^\Xi(\bar q)$ forces 
    \[ \int_K \frac{|\{\ell\in I_k\colon q_\ell \in \dot G_{\bbP}\}|}{|I_k|} d\check{\Xi}(k) = \mu^\Xi_{\bar I}(\bar q) \geq 1-\varp_0.\qedhere\]
\end{PROOF}

From \autoref{i21} and \autoref{i20.f}, it follows that: 

\begin{corollary}\label{i24-1}
    Let $K$, $L$ be finite sets. Assume that $h \colon K \to L$ is a surjective finite-to-one function.  Let $\bar{I} = \langle I_{k} \colon k \in K \rangle \in \bfI_{K}^{\fin}$ and for any $l \in L$, define $J_{l} \coloneqq \bigcup \{ I_{k} \colon h(k) = l \}$ and $\bar{J} \coloneqq \langle J_{l} \colon l \in L \rangle \in \bfI_{L}^{\fin}$. 
    If for any $k, k' \in K$, $h(k) = h(k')$ implies $\Xi(\{ k \}) = \Xi( \{ k' \})$ and $\vert  I_k\vert = \vert I_{k'}\vert$,  then $(\Xi_{h}, \bar{J}, \varp_{0})$-linkedness  implies  $(\Xi, \bar{I}, \varp_{0})$-linkedness. Furthermore, the equivalence holds when $h$ is a bijective function.
\end{corollary}

As a consequence of \autoref{i21} and \autoref{i20.7}: 

\begin{corollary}\label{i4.0}
     If $K$ is finite and $Q$ is $\vert W \vert$-linked, then it is $(\Xi, \bar{I}, 0)$-linked.
\end{corollary}

Therefore, from \autoref{i4.0} and ~\autoref{s8}~\ref{s8g}, we conclude:  

\begin{corollary}\label{i4.1}
     For any $\varp_{0} \in \left[ 0, \frac{1}{\vert W \vert + 1} \right]$, if $\Int^{\bbP}(Q) \geq 1 - \varp_{0}$ then $Q$ is $(\Xi, \bar{I}, 0)$-linked.
\end{corollary}

Thanks to \autoref{i20.7.1}, we have a strong connection between the intersection number and linkedness for fams on finite sets. Notice that, in the following result, $\Xi$ is fixed.

\begin{theorem}\label{i11.05.1}
   Assume that $K$ is finite. 
   Then $\Int^{\bbP}(Q) \geq 1- \varp_0$ iff, for any $\bar{J} \in \bfI_{K}^{\fin}$, $Q$ is $(\Xi, \bar{J}, \varp_0)$-linked.
\end{theorem}

\begin{PROOF}[\textbf{Proof}]{\ref{i11.05.1}}
    On the one hand, assume that $\Int^{\bbP}(Q) \geq 1 -\varp_{0}$. By \autoref{i20.7.1}, for any $\bar{J} \in \bfI_{K}^{\fin}$, $1 - \varp_{0} \leq \Int^{\bbP}(Q) \leq M_{\bar{J}}^{\Xi}(Q)$, and therefore, $Q$ is $(\Xi, \bar{J}, \varp_{0})$-linked by virtue of \autoref{i21}. On the other hand, if for any $\bar{J} \in \bfI_{K}^{\fin},$ $Q$ is $(\Xi, \bar{J}, \varp)$-linked, then $M_{\bar{J}}^{\Xi}(Q) \geq 1 - \varp_{0}$ by virtue of \autoref{i21}.  As a consequence,  by \autoref{i20.7.1} again, $ 1 - \varp_{0} \leq \inf_{\bar{J} \in \bfI_{K}^{\fin}} M_{\bar{J}}^{\Xi}(Q) = \Int^{\bbP}(Q)$, which proves the result. 
\end{PROOF}

\begin{corollary}\label{i22}
    Let $\Xi^0$ be the unique probability measure on $\calP(\{0\})$. Then $\Int(Q) \geq 1-\varp_0$ iff $Q$ is $(\Xi^0,\bar{J},\varp_0)$-linked for any $\bar{J} = \la J_0\ra$ (with $J_0$ finite non-empty).
\end{corollary}

\subsection{Comparing parameters}\ 

Relations between parameters $(\Xi, \bar{I}, \varp_{0})$ and $(\Xi', \bar{I'}, \varp_{0}')$ entail relations between the corresponding linkedness notions. For instance, if $\varp_{0} < \varp_{0}'$, then $(\Xi, \bar{I}, \varp_{0})$-linkedness implies $(\Xi, \bar{I}, \varp_{0}')$-linkedness. Concerning partitions, we have the following.

\begin{theorem}\label{i24.0}
    Let 
    $ \bar{J} = \langle J_{k} \colon k \in K \rangle \in \bfI_{K}^{\fin}$. 
    Assume that, for any $k \in K,$ $\vert J_{k} \vert = d_{k} \cdot \vert I_{k} \vert + r_{k}$, where $d_k<\omega$ and $0 \leq r_{k} < \vert I_{k} \vert$.  Let $\varp' \geq 0$ and $R \coloneqq \left \{  k \in K \colon \frac{r_{k}}{\vert J_{k} \vert} \leq \varp' \right \}$. If $(1-\varp')\Xi(R) >\varp_0$ and $Q$ is $(\Xi, \bar{J}, \varp_{0})$-linked, then it is $(\Xi, \bar{I}, \varp_R)$-linked in $\Por$, where $\varp_{R} \coloneqq 1 + \varp_{0} - (1- \varp') \Xi(R)$. In particular, it is $(\Xi, \bar{I}, \varp_0 +\varp')$-linked when $\Xi(R)=1$.
\end{theorem}

\begin{PROOF}[\textbf{Proof}]{\ref{i24.0}}
    Assume that $Q$ is $(\Xi, \bar{J}, \varp_{0})$-linked witnessed by $\limit^{\bar{J}}$ and $\dot{\Xi}^{\ast}$. By the hypothesis, for any $k \in K,$ we can partition $J_{k} = \left( \bigcup_{j < d_{k}} D_{j, k} \right) \cup R_{k},$ 
    such that, for any $j < d_{k}, \vert D_{j, k} \vert =  \vert I_{k} \vert,$ and $\vert R_{k} \vert = r_{k} < \vert I_{k} \vert.$ Therefore, for any $k \in K$ and $j < d_{k},$ we can find a bijection $f_{j, k} \colon D_{j, k} \to  I_{k}$ and an injective function $f_{k} \colon R_{n} \to I_{k}.$

    Now, let $\bar{q} = \langle q_{\ell} \colon \ell \in W_{\bar{I}} \rangle \in {}^{W_{\bar{I}}} Q$ and let us define a sequence $\bar{p} = \langle p_{i} \colon i \in W_{\bar{J}} \rangle \in {}^{W_{\bar{J}}} Q$ as follows. Since $\bar{J}$ is a partition of $W_{\bar{J}},$ to define $\bar{p}$ is enough to define it in each set in $\bar{J}$: for $k \in K$ and $i \in J_{k},$
    $$ p_{i} \coloneqq \left\{ 
       \begin{array}{ll} 
         q_{f_{j, k}(i)} & \text{if $i \in D_{j, k}$,} \\[1ex]  
         q_{f_{k}(i)} & \text{if $i \in R_{k}$.} 
       \end{array} \right. $$
    
    It is clear that, for any $k \in K$ and $j < d_{k},$ $\bar{p}  {\rest} D_{j, k} = \bar{q} \circ f_{j, k} $  and   $\bar{p}  {\rest} R_{k} = \bar{q} \circ f_{k}.$
    Hence, by the construction of $\bar{p}$, whenever $d_k\neq 0$, $\Por$ forces
    \begin{equation*}
        \begin{split}
            \frac{\vert \{ \ell \in I_{k} \colon q_{\ell} \in \dot{G}_{\bbP} \} \vert}{\vert I_{k} \vert} & \geq \frac{ \vert \{ i \in J_{k} \colon p_{i} \in \dot{G}_{\bbP}  \} \vert - r_{k} }{\vert I_{k} \vert \cdot d_{k}} = \frac{ \vert \{ i \in J_{k} \colon p_{i} \in \dot{G}_{\bbP}  \} \vert - r_{k} }{\vert J_{k} \vert - r_{k}}\\
            & \geq  \frac{ \vert \{ i \in J_{k} \colon p_{i} \in \dot{G}_{\bbP}  \} \vert - r_{k} }{\vert J_{k} \vert} = \frac{ \vert \{ i \in J_{k} \colon p_{i} \in \dot{G}_{\bbP}  \} \vert}{\vert J_{k} \vert} - \frac{r_{k}}{\vert J_{k} \vert}. 
        \end{split}
    \end{equation*}
    Note that the final inequality still holds when $d_k=0$ because $r_k=|J_k|$ in this case. 
    On the other hand, $\Por$ forces
    \[
        \int_K \frac{r_k}{|J_k|}d\dot\Xi^*(k)  = \int_R \frac{r_k}{|J_k|}d\dot\Xi^*(k) + \int_{K\menos R} \frac{r_k}{|J_k|}d\dot\Xi^*(k) \leq \varp'\Xi(R) + \Xi(K\menos R) = 1-(1-\varp')\Xi(R).
    \]
    As a consequence, $\limit^{\bar J}(\bar{p})$ forces 
    \begin{equation*}
        \begin{split}
            \int_{K} \frac{\vert \{ \ell \in I_{k} \colon q_{\ell} \in \dot{G}_{\bbP} \} \vert}{\vert I_{k} \vert} d \dot{\Xi}^{\ast}(k) 
            & \geq \int_{K} \frac{ \vert \{ i \in J_{k} \colon p_{i} \in \dot{G}_{\bbP}  \} \vert}{\vert J_{k} \vert} d \dot{\Xi}^{\ast}(k) - \int_{K} \frac{r_{k}}{\vert J_{k} \vert} d \dot{\Xi}^{\ast}(k)\\ 
            & \geq 1-\varp_0 - (1-(1-\varp')\Xi(R)) = (1-\varp')\Xi(R) -\varp_0.
        \end{split}
    \end{equation*}
    To conclude, it is enough to define $\limit^{\bar{I}}(\bar{q}) \coloneqq \limit^{\bar{J}}(\bar{p})$ and use the same $\dot{\Xi}^{\ast}$. 
\end{PROOF}

As a result of~\autoref{i24.0}, we infer:

\begin{corollary}\label{i24.01}
    In the context of \autoref{i24.0}, if for any $k \in K,$ $\vert I_{k} \vert$ divides $\vert J_{k} \vert,$ then any $(\Xi, \bar{J}, \varp_{0})$-linked set in $\Por$ is $(\Xi, \bar{I}, \varp_{0})$-linked. In particular, if for any $k \in K$, $\vert I_{k} \vert = \vert J_{k} \vert$, then $(\Xi, \bar{I}, \varp_{0})$-linkedness is equivalent to $(\Xi, \bar{J}, \varp_{0})$-likedness.
\end{corollary}

In the following, we show that our linkedness notion affects the intersection number of a set even when $K$ is infinite. 

\begin{lemma}\label{i4.1.67}
    Let $0<n<\omega$ and assume that, for any $k \in K$, $|I_k|=n$. If $Q$ is $(\Xi, \bar{I}, \varp_0)$-linked then, for all $\bar q=\la q_i\colon i<n\ra\in{}^{n}Q$,  $\frac{i_*(\bar q)}{n}\geq1-\varepsilon_0$. 
\end{lemma}
\begin{PROOF}[\textbf{Proof}]{\ref{i4.1.67}}
    Assume that $Q$ is $(\Xi, \bar{I}, \varp_{0})$-linked witnessed by $\limit^{\Xi}$ and $\dot{\Xi}^{\ast}$. Since each $I_{k}$ has cardinality $n,$ we can enumerate it as $\{\ell_{k,i} \colon i<n\} \subseteq W$. Hence, for any $\ell \in W$, there are unique $k_{\ell} \in K$ and $i_{\ell} < n$ such that $\ell = \ell_{k_{\ell}, i_{\ell}}.$ Let $\bar{q} = \langle q_{i} \colon i < n \rangle \in {}^{n} Q$. Define $\bar{r} \in \power{W}{Q}$ such that, for any $\ell \in W$, $r_{\ell} \coloneqq q_{i_{\ell}}$. Then,  
    \begin{align*}
        \limit^{\Xi}(\bar r) \Vdash `` 1 - \varp_{0} \leq \int_{K} \frac{\vert \{ \ell \in I_{k} \colon r_{\ell} \in \dot{G}_{\bbP} \} \vert}{ n }   d \dot{\Xi}^{\ast}(k) &=\int_{K} \frac{\vert \{ i<n \colon q_{i} \in \dot{G}_{\bbP} \} \vert}{ n }   d \dot{\Xi}^{\ast}(k)\text{''}
    \end{align*}
    Since the integral is on a constant term over a probability fam, $\limit^{\Xi}(\bar r) \Vdash \text{``}  \frac{\vert \{ i<n \colon q_{i} \in \dot{G}_{\bbP} \} \vert}{n }\geq 1 - \varp_{0} \text{''}$, hence there is some $r'\leq\limit^{\Xi}(\bar r)$ such that $\frac{\vert \{ i<n \colon r'\leq q_i \} \vert}{ n }\geq 1 - \varp_{0}$. Thus, $\frac{i_*(\bar q)}{n}\geq1-\varepsilon_0$. 
\end{PROOF}

\begin{corollary}\label{i4.1.69}
    If for all $n<\omega$ there is an $\bar{I}^n = \langle  I_{k}^n \colon k \in K \rangle \in \bfI_{K}^{\fin}$ such that $\forall k\in K\ (|I_k^n|=n)$ and $Q$ is $(\Xi,\bar{I}^n,\varp_0)$-linked, then $\Int^{\bbP}(Q) \geq 1- \varp_0$.
\end{corollary}
 
\begin{theorem}\label{i11.0}
    Assume that 
    $\Xi(\{k\in K \colon |I_k| = n\}) = 0$ for all $n<\omega$ (this implies that $\Xi$ is free). If $Q$ is $(\Xi, \bar{I}, \varp_0)$-linked, then $\Int^{\bbP}(Q) \geq 1- \varp_0$.

    As a consequence, 
    $M^{\Xi'}_{\bar{I}'}(Q)\geq 1-\varp_0$ for any finite $K'$, $\bar I'\in\bfI_{K'}$ and any probability fam $\Xi'$ on $K'$, i.e.\ $Q$ is $(\Xi',\bar{I}',\varp_0)$-linked.
\end{theorem}
\begin{PROOF}[\textbf{Proof}]{\ref{i11.0}}
Let $\varepsilon>0$ and $0<n<\omega$ be arbitrary. Clearly, $\Xi(\{k\in K \colon 
\varepsilon>\frac{n}{|I_k|}\}) = 1$. Pick a partition $\bar{I}^n = \langle  I_{k}^n \colon k \in K \rangle \in \bfI_{K}^{\fin}$ such that $\forall k\in K\ (|I_k^n|=n)$. By using the division algorithm, we can express $|I_k|$ as $d_k\cdot n+r_k$ with $0 \leq r_{k}<n$. By~\autoref{i24.0}, we obtain that $Q$ is $(\Xi,\bar{I}^n,\varp_0+\varepsilon)$-linked. Unfixing $n$, we have $\Int^{\bbP}(Q) \geq 1-(\varp_0+\varepsilon)$ by using~\autoref{i4.1.69}. Since $\varepsilon$ was taken arbitrary, $\Int^{\bbP}(Q) \geq 1-\varp_0$. 
\end{PROOF}

Similar to \autoref{i20.f}, we have a connection between $(\Xi, \bar{I}, \varp)$-linkedness and $(\Xi_{h}, \bar{I}', \varp)$-linked\-ness for $\Xi_{h}$ as defined in \autoref{m37.1} and a suitable $\bar{I}'$:  

\begin{lemma}\label{e102}
    Let $K'$ be a non-empty set and $h\colon K\to K'$. Consider the probability fam $\Xi_{h}$ on $\pts(K')$ defined as in \autoref{m37.1}. Assume that $|I_k|=|I_j|$ whenever $j,k\in K$ and $h(j)=h(k)$, and let $\bar I'\in\Ibf^\fin_{K'}$ such that $|I'_{h(k)}|=|I_k|$ for all $k\in K$. Then, any $(\Xi,\bar I,\varp_0)$-linked set is $(\Xi_{h},\bar I',\varp_0)$-linked. Furthermore, if $\Xi(\{ k \in K  \colon \vert I_{k} \vert = n\}) = 0$ for all $n < \omega,$ then for any $m < \omega,$ $\Xi_{h}(\{ k' \in K' \colon \vert I_{k}' \vert = m \}) = 0.$
\end{lemma}
\begin{PROOF}[\textbf{Proof}]{\ref{e102}}
    Using that $K'\sminus \ran h\in\Ncal(\Xi_{h})$, we get that, for any bounded $f\colon K'\to \R$,
    \[\int_{K'}f d\Xi_{h} = \int_K f\circ h\, d\Xi.\]
    For each $k\in K$, pick some bijection $g_k\colon I_k\to I'_{h(k)}$. 
    Assume that $Q$ is $(\Xi,\bar I,\varp_0)$-linked in $\Por$ witnessed by the limit function $\limit^\Xi$ and the $\Por$-name $\dot\Xi^*$ of a fam extending $\Xi$. Let $W'\coloneqq W_{\bar I'}$. Given $\bar q\in {}^{W'}Q$, define $\bar p\in {}^WQ$ such that $p_\ell \coloneqq  q_{g_k(\ell)}$ for all $\ell\in I_k$ and $k\in K$, and set $\limit^{\Xi_{h}}\bar q \coloneqq  \limit^{\Xi}\bar p$. On the other hand, let $\dot\Xi^{**}$ be a $\Por$-name of the probability fam such that $\dot\Xi^{**}(b)=\dot\Xi^*(h^{-1}[b])$ for all $b\subseteq K'$, which clearly extends $\Xi_{h}$. 

    We prove that $\limit^{\Xi_{h}}$ and $\dot\Xi^{**}$ witness that $Q$ is $(\Xi_{h},\bar I',\varp_0)$-linked. Let $\bar q$ and $\bar p$ as above. Then, $\Por$ forces
    \begin{align*}
        \int_{K'}\frac{\vert \setcolon{ \ell' \in I'_{k'} }{ q_{\ell'} \in \dot G } \vert}{\vert I'_{k'} \vert} d \dot{\Xi}^{**}(k') & = \int_K \frac{\vert \setcolon{ \ell' \in I'_{h(k)} }{ q_{\ell'} \in \dot G } \vert}{\vert I'_{h(k)} \vert} d \dot{\Xi}^{*}(k)\\ 
        & = \int_K \frac{\vert \setcolon{ \ell \in I_{k} }{ p_{\ell} \in \dot G } \vert}{\vert I_{k} \vert} d \dot{\Xi}^{*}(k).
    \end{align*}
    Hence, $\limit^{\Xi_{h}}\bar q = \limit^\Xi \bar p$ forces that this integral is ${\geq}\, 1-\varp_0$.

    Now fix $n<\omega$. If $k\in K$ and $|I'_{h(k)}|=n$, then $|I_k|=|I'_{h(k)}|=n$. Hence, 
    \[h^{-1}[\setcolon{k'\in K'}{|I'_{k'}|=n}]\subseteq \setcolon{k\in K}{|I_k|=n}.\]
    This shows the ``furthermore'' in the statement of the Lemma. 
\end{PROOF}

\subsection{$\omega$-regularity}

Motivated by the notion of \emph{regularity} in the context of ultra-filters (see e.g. \cite{KeislerU}), we introduce a similar notion for finitely additive measures. 

\begin{definition}\label{e101}
    We say that $\Xi$ is \emph{$\omega$-regular} if $K$ can be partitioned into countably many $\Xi$-measure zero sets.
\end{definition}

Due to its relevance in developing the iterations theory with fams, we isolate the assumption of \autoref{i24.0} below.

\begin{definition}\label{i11.00}
     A pair $(\Xi,\bar I)$ is \emph{free} if $\Xi(\setcolon{k\in K }{ |I_k| = n}) = 0$ for all $n<\omega$.
\end{definition}

If $(\Xi,\bar I)$ is a free pair then $\Xi$ is $\omega$-regular, and any $\omega$-regular fam is free. In the case that $K$ is countable, any free fam on $\pts(K)$ is $\omega$-regular.



We show that the existence of a non-$\omega$-regular free probability fam is a large cardinal assumption.

\begin{lemma}\label{e101-1}
    Assume that $\Xi$ is a non-$\omega$-regular fam. Then, there is some $A\in\pts(K)\sminus\Ncal(\Xi)$ such that $\Ncal(\Xi\frestr A)$ is a $\sigma$-ideal in $\pts(A)$.
\end{lemma}
\begin{PROOF}[\textbf{Proof}]{\ref{e101-1}}
    Assume that no such $A$ can exist. Let $\Qor\coloneqq \pts(K)\sminus\Ncal(\Xi)$. Then, 
    \[D\coloneqq \setcolon{a\in\Qor}{a \text{ can be partitioned into countably many $\Xi$-measure zero sets}}\] 
    is open dense in $\Qor$. On the other hand, $\Qor$ has the ccc, 
    so $D$ contains a maximal antichain in $\Qor$, moreover, we can partition $K$ into countably many sets in $D$. This implies that $K$ can be partitioned into countably many $\Xi$-measure zero sets, i.e.\ $\Xi$ is $\omega$-regular. 
\end{PROOF}

\begin{remark}\label{e101-2}
    The fam produced in \autoref{e101-1}, after dividing with $\Xi(A)$, gives an example of a (possibly free) probability fam on $\pts(A)$ whose measure zero ideal is a $\sigma$-ideal. This does not necessarily imply that the fam is $\sigma$-additive: Kumar and Kunen~\cite{KuKu}, assuming the existence of a measurable cardinal, constructed a FS (finite-support) iteration forcing that there is no real-valued cardinal ${\leq}\, \cfrak$ and that there is a probability fam on $\pts([0,1])$, extending the Lebesgue measure, whose measure zero ideal is $\cfrak$-complete. In this generic extension, the fact that there is no real-valued measurable cardinal ${\leq}\, \cfrak$ implies that this fam is not $\sigma$-additive, even more, $\sigma$-additivity fails when restricted to any set of positive measure.
\end{remark}

\begin{theorem}\label{e101-3}
    The following theories are equiconsistent:
    \begin{enumerate}[label =\normalfont (\arabic*)]
        \item\label{e101-31} $\thzfc+{}$ there is a measurable cardinal.
        \item\label{e101-32} $\thzfc+{}$ there is a non-$\omega$-regular free probability fam on some $\pts(A)$.
    \end{enumerate}
\end{theorem}
\begin{PROOF}[\textbf{Proof}]{\ref{e101-3}}
    We show the non-trivial direction, i.e.\ that the consistency of~\ref{e101-32} implies the consistency of~\ref{e101-31}. By \autoref{e101-1}, we can prove in~\ref{e101-32} that there is a free probability fam $\Xi'$ in some $\pts(A)$ such that $\Ncal(\Xi')$ is a $\sigma$-ideal. This ideal is also $\sigma$-saturated, i.e.\ the poset $\pts(A)\sminus \Ncal(\Xi')$ has the ccc. We let $\kappa$ be the smallest cardinal such that $\pts(\kappa)$ contains a $\sigma$-saturated $\sigma$-ideal $I$ such that $\bigcup I=\kappa$. By \cite[Lem.~10.5]{Je2}, $I$ must be $\kappa$-complete. Following Solovay~\cite{Solmeas}, if $\kappa>\cfrak$ then $\kappa$ is measurable, otherwise there is an inner model where $\kappa$ is measurable. 
\end{PROOF}

Notice that as a consequence of \autoref{e102}, under the same conditions there, we have that $(\Xi_{h},\bar I')$ is a free pair whenever $(\Xi,\bar I)$ is. Using this fact we analyze the role of free pairs in our linkedness property.

\begin{lemma}\label{e103}
    Assume that $\Xi$ is $\omega$-regular. Then there is some $h\colon K\to\omega$ satisfying:
    \begin{enumerate}[label = \normalfont (\alph*)]
        \item For $j,k\in K$, $h(j)=h(k)$ implies $|I_j|=|I_k|$. 
        
        \item $\Xi_{h}$ on $\omega$ is free.
    \end{enumerate}
    Moreover, whenever $(\Xi,\bar I)$ is a free pair, $h$ can be defined by $h(k)\coloneqq |I_k|$.
\end{lemma}
\begin{PROOF}[\textbf{Proof}]{\ref{e103}}
    We can partition $K$ into countably many sets $A_n\in\Ncal(\Xi)$ ($n<\omega$) such that $|I_j|=|I_k|$ whenever $j,k\in A_n$. 
    In the case that $(\Xi,\bar I)$ is a free pair, we can set $A_n\coloneqq \setcolon{k\in K}{|I_k|=n}$. 
    Define $h\colon K\to \omega$ by $h(k)\coloneqq n$ if $k\in A_n$. This map is as required.
\end{PROOF}

As a direct consequence of \autoref{e102} and \autoref{e103}:

\begin{corollary}\label{e104}
    If $(\Xi,\bar I)$ is a free pair then there is some free fam $\Xi'$ on $\pts(\omega)$ such that, for any $\varp\in[0,1)$ and $\bar I'\in\Ibf^\fin_\omega$ satisfying $|I'_n|=n$ for all $0<n<\omega$, in any forcing notion, any $(\Xi,\bar I,\varp)$-linked set is $(\Xi',\bar I',\varp)$-linked. In particular, $(\Xi',\bar I')$ is a free pair.
\end{corollary}

The previous can be strengthened as follows.

\begin{lemma}\label{e106}
    If $(\Xi,\bar I)$ is a free pair, $\varp'>0$ and $g\colon \omega\to\omega\menos\{0\}$ then there is some free fam $\Xi'$ on $\pts(\omega)$ such that, for any $0\leq\varp<1-\varp'$ and $\bar I'\in\Ibf^\fin_\omega$ satisfying $|I'_n|=g(n)$ for all $n<\omega$, in any forcing notion, any $(\Xi,\bar I,\varp)$-linked set is $(\Xi',\bar I',\varp+\varp')$-linked.
\end{lemma}
\begin{PROOF}[\textbf{Proof}]{\ref{e106}}
    By \autoref{e104}, wlog we may assume that $K=\omega$ and $|I_k|=k$ for all $0<k<\omega$. Find an increasing function $f\colon \omega\to\omega$ with $f(0)=0$ such that $\frac{g(n)}{f(n)}<\varp'$ for $0<n<\omega$, and pick a $\bar J\in \Ibf^\fin_\omega$ such that $|J_k|=g(n)$ whenever $f(n)\leq k <f(n+1)$. 
    Thus, by \autoref{i24.0}, for $0\leq\varp<1-\varp'$, any $(\Xi,\bar I,\varp)$-linked set in a forcing notion is $(\Xi,\bar J,\varp+\varp')$-linked. 
    
    Now, define $h\colon \omega\to\omega$ by $h(k)\coloneqq n$ whenever $f(n)\leq k<f(n+1)$. Apply \autoref{e102} to $\bar J$ to find the required $\Xi'$.
\end{PROOF}

The linkedness notion of this chapter is stronger than $\Fr$-linked for $\omega$-regular fams.

\begin{theorem}\label{e085}
    Assume that $\Xi$ is $\omega$-regular. If $Q$ is $(\Xi,\bar I,\varp_0)$-linked in $\Por$ then it is $\Fr$-linked. 
\end{theorem}
\begin{PROOF}[\textbf{Proof}]{\ref{e085}}
    By \autoref{e102} and \autoref{e103}, we can assume wlog that $K=\omega$. 
    Let $\dot\Xi^*$ and $\limit^\Xi$ witness the linkedness of $Q$. 
    The result follows by \autoref{i8}. 
\end{PROOF}

\subsection{Characterizations}\ 

The following is a combinatorial characterization of $(\Xi, \bar{I}, \varp_{0})$-linkedness, which includes a special characterization when $\Xi$ satisfies the uap (e.g.\ free).

\begin{theorem}\label{i6}    
    Let $\limit^{\Xi} \colon \power{W}{Q} \to \bbP$ be a function. Then, the following statements are equivalent: 
    \begin{enumerate}[label=\rm(\alph*)]
        \item\label{i6a} $Q$ is $(\Xi, \bar{I}, \varp_{0})$-linked witnessed by $\lim^{\Xi}$.

        \item\label{i6c} Given  
        \begin{itemize}
            \item $i^{\ast} < \omega$, $ \bar{q}^{i} = \langle q_{\ell}^{i} \colon \ell \in W \rangle \in \power{W}{Q}$ for each $i < i^{\ast},$
            
            \item $P\in \Pbf^\Xi$,  
            
            \item $\varp > 0$ and $q \in \bbP$ such that for all $i < i^{\ast}$, $q \leq \limit^{\Xi}(\bar{q}^{i}),$ 
        \end{itemize} 
        there are a non-empty finite set $u \subseteq K$, a probability measure $\mu$ on $\pts(u)$ and some $q' \leq q$ in $\Por$ such that
        \begin{enumerate}[label=\rm(\arabic*)]
            \item\label{i2ac} $\displaystyle\left| \mu(u\cap b) - \Xi(b) \right| < \varp$ for all $b\in P$ and
    
            \item\label{i2bc} $\displaystyle\int_u \frac{\vert \{ \ell \in I_{k} \colon  q' \leq q_{\ell}^{i} \} \vert}{\vert I_{k} \vert} d\mu(k) > 1 - \varp_{0} - \varp$ for all $i < i^{\ast}.$
        \end{enumerate}
        \end{enumerate}
        When $\Xi$ has the uap, the measure $\mu$ in~\ref{i6c} can be found uniform, concretely,~\ref{i6b} below is equivalent to the above. In this case, $\dot \Xi^*$ (in~\ref{i6a}) can be forced with the uap.
        \begin{enumerate}[resume*]
        \item\label{i6b} Given  
        \begin{itemize}
            \item $i^{\ast} < \omega$, $ \bar{q}^{i} = \langle q_{\ell}^{i} \colon \ell \in W \rangle \in \power{W}{Q}$ for each $i < i^{\ast},$
            
            \item $P\in \Pbf^\Xi$,  
            
            \item $\varp > 0$ and $q \in \bbP$ such that for all $i < i^{\ast}$, $q \leq \limit^{\Xi}(\bar{q}^{i}),$ 
        \end{itemize} 
        there are a non-empty finite set $u \subseteq K$  and some $q' \leq q$ in $\Por$ such that
        \begin{enumerate}[label=\rm(\arabic*)]
            \item\label{i2a} $\displaystyle\left| \frac{\vert u \cap b \vert}{\vert u \vert} - \Xi(b) \right| < \varp$ for all $b\in P$ and 
    
            \item\label{i2b} $\displaystyle\frac{1}{\vert u \vert} \sum_{k \in u} \frac{\vert \{ \ell \in I_{k} \colon  q' \leq q_{\ell}^{i} \} \vert}{\vert I_{k} \vert} > 1 - \varp_{0} - \varp,$ for all $i < i^{\ast}.$
        \end{enumerate}
        \end{enumerate}
\end{theorem}

\begin{PROOF}[\textbf{Proof}]{\ref{i6}}
We present the proof of \ref{i6a}${}\Leftrightarrow{}$\ref{i6b}, i.e.\ we assume that $\Xi$ has the uap. A similar argument can prove the equivalence with~\ref{i6c} (without assuming uap).

\ref{i6a}${}\Rightarrow{}$\ref{i6b}: Assume that $Q$ is $(\Xi, \bar{I}, \varp)$-linked, that is, there exists a $\bbP$-name $\dot{\Xi}^{\ast}$ of a probability fam on $\calP(K)$ extending $\Xi$ such that, for any $\bar{q} \in \power{W}{Q},$ 
    \begin{equation}\label{i6-p1-e3}
        \text{$\limit^{\Xi}(\bar{q}) \Vdash \text{``} \int_{K} \frac{\vert \{ \ell \in I_{k} \colon q_{\ell} \in \dot{G}_{\bbP} \} \vert}{\vert I_{k} \vert}   d  \dot{\Xi}^{\ast}(k) \geq 1 - \varp_{0}\text{''}$.}
        \tag{$\oplus$}
    \end{equation}
    Let $i^{\ast} < \omega,$ consider $\bar{q}^{i} \coloneqq \langle q_{\ell}^{i} \colon \ell \in W \rangle \in \power{W}{Q}$ for any $i < i^{\ast},$ and let $P\in \Pbf^\Xi$. Let $\varp > 0$ and $q \in \bbP$ such that, for any $i < i^{\ast},$ $q \leq \lim^{\Xi}(\bar{q}^{i}).$  Working in the generic extension $V[G]$ for some generic filter $G$ containing $q$ define, for any $k \in K$ and $i < i^{\ast},$ $A_{k}^{i} \coloneqq \{ \ell \in I_{k} \colon q_{\ell}^{i} \in G_{\bbP} \}$ and, for any $i < i^{\ast}$, consider the function $f_{i} \colon K \to [0, 1]_{\QQQ}$ such that, for any $k \in K,$ $f_{i}(k) \coloneqq \frac{\vert A_{k}^{i} \vert}{\vert I_{k} \vert}.$ 
    By \eqref{i6-p1-e3} we have that, for any $i < i^{\ast},$ $ \int_{K} f_{i} d \dot{\Xi}^{\ast} \geq 1 - \varp_{0}$.

    The functions $f_{i}$ are bounded and $\Xi$ is a probability fam with the uap, so we can apply \autoref{t85} to get a finite set $u \subseteq K$ such that: 
    \begin{enumerate}[label = \normalfont (\roman*)]
        \item for any $b\in P,$ $\left \vert  \Xi^{u}(b) - \Xi(b)  \right \vert < \varp,$ and

        \item\label{i6ii} for any $i < i^{\ast},$ $ \displaystyle\left \vert \frac{1}{\vert u \vert}  \sum_{k \in u} f_{i}(k) - \int_{K} f_{i} d \Xi^{\ast}    \right \vert <  \varp. $
    \end{enumerate}

    By~\ref{i6ii} we can conclude that, for any $i < i^{\ast},$ $ \frac{1}{\vert u \vert}  \sum_{k \in u} f_{i}(k) > 1- \varp_{0} - \varp$. 
    Now, we can find $\bbP$-names $\dot{u},$ $\{ \dot{A}^{i}_{k} \}_{k < \omega}$ for $i<i^*$ such that, $\dot{u}[G] = u$ and, for any $k \in K,$ $\dot{A}^i_{k}[G] = A_{k}^i.$ As a consequence, in the ground model, we can find a condition $q' \leq q$ deciding $u = \dot{u}$  and $A^i_{k} = \dot{A}^i_{k}$ for any $k \in u$ and $i<i^*$. Notice that $\ell\in A^i_k$ implies $q'\leq^* q^i_\ell$, and $\ell\in I_k\menos A^i_k$ implies $q'\perp q^i_\ell$. Then, find $q''\leq q'$ such that $q''\leq q^i_\ell$ for all $i<i^*$, $k\in u$ and $\ell\in A^i_k$.   
    As a consequence, we have that: 
    \begin{enumerate}
        \item for any $b\in P,$ $\left \vert  \Xi^{u}(b) - \Xi(b)  \right \vert < \varp,$

        \item for any $i < i^{\ast},$ $ \displaystyle\frac{1}{\vert u \vert}  \sum_{k \in u} \frac{\vert \{ \ell \in I_{k} \colon q'' \leq q_{\ell}^{i} \} \vert}{\vert I_{k} \vert} > 1- \varp_{0} - \varp$.
    \end{enumerate}
    \ref{i6b}${}\Rightarrow{}$\ref{i6a}: We plan to use \autoref{t85} to show the existence of $\dot \Xi^*$ (in the generic extension). Let  $G \subseteq \bbP$ be a generic filter over $V.$ Working in $V[G]$, define $\calL \coloneqq \{ \bar{q} \in \power{W}{Q} \cap V \colon \limit^{\Xi}(\bar{q}) \in G \}$. 
    
    Fix $i^{\ast} < \omega$ and, for each $i < i^{\ast}$, $\bar{q}^{i} = \langle q_{\ell}^{i} \colon \ell \in W \rangle \in \calL.$ 
    Working in the ground model, let $\varp > 0$ and $P\in \Pbf^\Xi$. Define $D_{1}$ as the set of conditions $r \in \bbP$ such that there exists a finite set $u \subseteq K$ satisfying:
    \begin{enumerate}[label=\rm(\arabic*$_{r}$)]
        \item $\left| \Xi^{u}(b) - \Xi(b) \right| < \varp$ for all $b\in P,$ and

        \item $\displaystyle\frac{1}{\vert u \vert} \sum_{k \in u} \frac{\vert \{ \ell \in I_{k} \colon  r \leq q_{\ell}^{i} \} \vert}{\vert I_{k} \vert} > 1 - \varp_{0} - \varp,$ for all $i < i^{\ast}.$
    \end{enumerate}

    On the other hand, let $D_{2} \coloneqq \{ r \in \bbP \colon \exists i  < i^{\ast}\ (r \perp  \limit^{\Xi}(\bar{q}^{i})) \}$ and $D \coloneqq D_{1} \cup D_{2}.$  Notice that $D$ is dense in $\bbP,$ and it is clear that $D \in V$ because its definition does not depend on the generic filter $G.$ 
        
    Working on the generic extension again, since $D \in V$ is dense, $G \cap D \neq \emptyset.$ Choose $q' \in D$ such that $q' \in G.$ Since $\limit^{\Xi}(\bar q^i) \in  G$ for any $i<i^*$, necessarily $q' \in D_{1},$ hence there is a finite set $u \subseteq K$ such that: 
    
    \begin{enumerate}
        \item $\left| \Xi^{u}(b) - \Xi(b) \right| < \varp,$ for all $b\in P$, and
    
        \item $\displaystyle \frac{1}{\vert u \vert} \sum_{k \in u} \frac{ \vert \{ \ell \in I_{k} \colon q_{\ell}^{i} \in G \} }{ \vert I_{k} \vert}  \geq \frac{1}{\vert u \vert} \sum_{k \in u} \frac{\vert \{ \ell \in I_{k} \colon  q' \leq q_{\ell}^{i} \} \vert}{\vert I_{k} \vert} > 1 - \varp_{0} - \varp,$ for all $i < i^{\ast}.$
    \end{enumerate}   
    
    By \autoref{t85}, there exists a fam  $\Xi^{\ast}$ on $\calP(K)$ with the uap extending $\Xi$ such that, for all $\bar{q} = \langle q_{\ell} \colon \ell \in W \rangle \in \calL,$  
    $\displaystyle \int_{K} \frac{\vert \{ \ell \in I_{k} \colon q_{\ell} \in G_{\bbP}  \} \vert}{\vert I_{k} \vert}  d  \Xi^{\ast}(k) \geq 1 - \varp_{0}.$
    
    Back in $V,$ let $\dot{\Xi}^{\ast}$ a $\bbP$-name of $\Xi^{\ast}.$ It is clear that $\dot{\Xi}^{\ast}$ is as required.
\end{PROOF}

\begin{remark}\label{i6r}\

    \begin{enumerate}[label=\rm(\arabic*)]
        \item\label{i6r.1} In particular cases,\footnote{See ~\cite[Def.~1.7]{KST} and \cite[Def.~1.10]{KST}.} a function $\limit^{\Xi}$ satisfying the requirement in \autoref{i6}~\ref{i6b} is called a \emph{strong fam limit for intervals} in \cite{KST}. This notion was generalized in \cite{uribethesis}, where the property of $(\Xi, \bar{I}, \varp_{0})$-linkedness was defined according to \autoref{i6}~\ref{i6b}.  However, this automatically implies that $\Xi$ has the uap, which is too restrictive. For this reason, in this paper, we decided to weak the original definition of $(\Xi, \bar{I}, \varp_{0})$-linked from \cite{uribethesis}, which is easily recovered demanding ``uap'' in \autoref{i2}. 

        \item In the case of random forcing, \ref{i6b}${}\Rightarrow{}$\ref{i6a} in \autoref{i6} appeared implicitly in \cite{Sh00}. The converse implication is a new result, which allowed us to present a more friendly definition of $(\Xi, \bar{I}, \varp_{0})$-linkedness.   

        \item In the version of~\cite{KST}, where $K=\omega$, an arbitrary $k^*<\omega$ is fixed (along with $i^*$, $\bar q^i$, etc.), and it is required that $u$ is found disjoint with $k^*$. This is implied by \autoref{i6}~\ref{i6b} when $\Xi$ is a free fam (see \autoref{i2.3}).

        \item The characterizations presented in this subsection are not needed to develop our theory of iterations with fams, nor in the applications within this paper, but are essential to check \autoref{i4},~\ref{i15.2} and~\ref{i11.000}, whose proofs are sophisticated and contained in other articles. We present these characterizations for the completeness of our theory.
    \end{enumerate}
\end{remark}

The following result is a consequence of~\autoref{i6.2}.

\begin{corollary}\label{i2.3}
    If\/ $\Xi$ is free then $(\Xi,\bar{I},\varp_0)$-linkedness is equivalent to \autoref{i6}~\ref{i6b} \emph{evading finite sets}, that is, there we can find $u \subseteq K \menos F$ for any fixed finite set $F$.
\end{corollary}

\begin{remark}
    As a converse of \autoref{i2.3}, if there is some forcing notion containing some non-empty $(\Xi,\bar{I},\varp_0)$-linked subset evading finite sets, then $\Xi$ must be free. This is a direct consequence of \autoref{m420}. 
\end{remark}

The case $\varp_0=0$ has special relevance.

\begin{theorem}\label{i24}
    Let $\lim^\Xi\colon {}^WQ \to \bbP$. Then, the following statements are equivalent. Moreover, when $\Xi$ has the uap, $\dot\Xi^*$ can be found with the uap in~\ref{i24ii}.
    \begin{enumerate}[label= \normalfont (\roman*)]
        \item\label{i24i} $Q$ is $(\Xi,\bar{I},0)$-linked witnessed by $\lim^\Xi$.
        
        \item\label{i24ii} There is a $\bbP$-name $\dot\Xi^*$ of a fam extending $\Xi$ such that, for any $\bar q\in {}^W Q$, $\lim^\Xi (\bar q)$ forces that, for any $0\leq \delta<1$,
        \[\dot \Xi^*\left(\set{k\in K}{\frac{|\{\ell \in I_k \colon q_\ell \in \dot{G}_{\bbP} \}|}{|I_k|} > \delta}\right) =1.\]
        
        \item\label{i24iii} Given
        \begin{itemize}
        \item $i^{\ast} < \omega$ and $ \bar{q}^{i} = \langle q_{\ell}^{i} \colon \ell \in W \rangle \in \power{W}{Q}$ for each $i < i^{\ast},$

        \item $A\subseteq K$ with $\Xi(A)>0$, and 

        \item $\varp > 0$ and $q \in \bbP$ such that, for all $i < i^{\ast}$, $q \leq \limit^{\Xi}(\bar{q}^{i}),$
        \end{itemize}
        there are $q'\leq q$ in $\bbP$ and $k\in A$ such that
        \[\frac{|\{\ell \in I_k \colon q'\leq q^i_\ell\}|}{|I_k|} > 1-\varp \text{ for all $i<i^*$.}\]
    \end{enumerate}
\end{theorem}
\begin{PROOF}[\textbf{Proof}]{\ref{i24}}
    Notice that, for any probability fam on $\pts(K)$ and any function $f\colon K\to [0,1]$, the integral of $f$ is $1$ iff, for any $0\leq \delta <1$, the set $\{k\in K\colon f(k)>\delta\}$ has measure $1$. For this reason, \ref{i24i}${}\Leftrightarrow{}$\ref{i24ii} is immediate from \autoref{i2}.

    \ref{i24ii}${}\Rightarrow{}$\ref{i24iii}: Fix the given objects in~\ref{i24iii}. Then, $q$ forces that the set below has the same $\dot\Xi^*$-measure as $A$:
    \[\set{k\in A}{\frac{|\{\ell \in I_k \colon q_\ell^i \in \dot{G}_{\bbP} \}|}{|I_k|} > 1-\varp \text{ for all $i<i^*$}}.\]
    In particular, $q$ forces that this set has positive measure, so it is not empty. Then, there are $q'\leq q$ in $\bbP$ and some $k\in K$ forced by $q'$ in the set above, also satisfying that, for any $i<i^*$ and $\ell\in I_k$, $q'\leq q^i_\ell$ iff $q'$ forces $q^i_\ell \in \dot{G}_{\bbP}$. Therefore, $k\in A$ and 
    $\displaystyle\frac{|\{\ell \in I_k \colon q'\leq q^i_\ell\}|}{|I_k|} > 1-\varp$ for all $i<i^*$.

    \ref{i24iii}${}\Rightarrow{}$\ref{i24ii}: For $\bar q = \la q_\ell \colon \ell \in W\ra\in {}^W Q$ and $0\leq \delta <1$, let $\dot w(\bar q,\delta)$ be a $\Por$-name of the set $\{k\in K \colon \frac{|\{\ell \in I_k \colon q_\ell \in \dot{G}_{\bbP} \}|}{|I_k|}>\delta\}$. In virtue of \autoref{m69}, it is enough to show that $\Por$ forces that, for any finite $A\in\pts(K)\cap V$ and a finite $J\subseteq (({}^W Q )\times [0,1))\cap V$ such that $\Xi(A)>0$ and $\lim^\Xi(\bar q)\in \dot G$ for all $(\bar q,\delta)\in J$, $A\cap \bigcap_{(\bar q,\delta)\in J}\dot w(\bar q,\delta) \neq \emptyset$. Moreover, if $\Xi$ has the uap then the desired $\dot\Xi^*$ can be found with the uap by \autoref{t96}.

    Indeed, let $A\subseteq K$ with $\Xi(A)>0$, $i^*<\omega$ and, for $i<i^*$, let $\delta_i\in[0,1)$ and $\bar q^i\in {}^WQ$, and assume that $p\in \Por$ forces $\lim^\Xi(\bar q^i)\in \dot G$ for all $i<i^*$. Then, there is some $q\leq p$ in $\Por$ stronger than all $\lim^\Xi(\bar q^i)$. By~\ref{i24iii}, there are $k\in A$ and $q'\leq q$ such that $\frac{|\{\ell \in I_k \colon q'\leq q^i_\ell\}|}{|I_k|} > \delta$ for all $i<i^*$, where $\delta \coloneqq \max_{i<i^*} \delta_i$. Hence, $q'$ forces that $k\in A\cap \bigcap_{i<i^*} \dot w(\bar q^i,\delta) \subseteq A\cap \bigcap_{i<i^*} \dot w(\bar q^i,\delta_i)$.
\end{PROOF}

Recall that probability fams into $\{0, 1\}$ represent ultrafilters (and have the uap).\footnote{See, e.g. \cite[Sec.~3.2]{uribethesis}.} When $\Xi$ is such a fam, the previous result ensures extensions to ultrafilters.

\begin{corollary}\label{i25}
    Let $D\subseteq\pts(K)$ be an ultrafilter, $\bbP$ a forcing notion, 
    and $\lim^D\colon {}^WQ \to \bbP$. Then the following statements are equivalent.
    \begin{enumerate}[label= \normalfont (\roman*)]
        \item $Q$ is $(D,\bar{I},0)$-linked witnessed by $\lim^D$.
        \item There is a $\bbP$-name $\dot D^*$ of an ultrafilter on $\pts(K)$ extending $D$ such that, for any $\bar q\in {}^W Q$, $\lim^D (\bar q)$ forces that, for any $0\leq \delta<1$,
        \[\set{k\in K}{\frac{|\{\ell \in I_k \colon q_\ell \in \dot G\}|}{|I_k|} > \delta} \in \dot D^*.\]
        \item Given
        \begin{itemize}
        \item $i^{\ast} < \omega$ and $ \bar{q}^{i} = \langle q_{\ell}^{i} \colon \ell \in W \rangle \in \power{W}{Q}$ for each $i < i^{\ast},$

        \item $A\in D$, and 

        \item $\varp > 0$ and $q \in \bbP$ such that for all $i < i^{\ast}$, $q \leq \lim^{D}(\bar{q}^{i}),$
        \end{itemize}
        there are $q'\leq q$ in $\bbP$ and $k\in A$ such that
        \[\frac{|\{\ell \in I_k \colon q'\leq q^i_\ell\}|}{|I_k|} > 1-\varp \text{ for all $i<i^*$.}\]
    \end{enumerate}
\end{corollary}

\begin{corollary}\label{i26}
    Along with the hypothesis of the previous result, assume that $\bar{I}^1 = \la I^1_k\colon k\in K\ra$, $I^1_k \coloneqq\{k\}$. Then, the following statements are equivalent.
    \begin{enumerate}[label= \normalfont (\roman*)]
        \item $Q$ is $(D,\bar{I}^1,0)$-linked witnessed by $\lim^D$.
        \item There is a $\bbP$-name $\dot D^*$ of an ultrafilter on $\pts(K)$ extending $D$ such that, for any $\bar q\in {}^K Q$, $\lim^D (\bar q) \Vdash$``$\{k\in K\colon q_k \in \dot G\}\in \dot D^*$''.
        
        \item Given
        \begin{itemize}
        \item $i^{\ast} < \omega$ and $ \bar{q}^{i} = \langle q_{k}^{i} \colon k \in K \rangle \in \power{K}{Q}$ for each $i < i^{\ast},$

        \item $A\in D$, and 

        \item $q \in \bbP$ such that for all $i < i^{\ast}$, $q \leq \lim^{D}(\bar{q}^{i}),$
        \end{itemize}
        there are $q'\leq q$ in $\bbP$ and $k\in A$ such that $q'\leq q^i_k$ for all $i<i^*$.
    \end{enumerate}
\end{corollary}

This last result corresponds to the notion of \emph{ultrafilter limits} for a subset of a poset, which appears implicitly in \cite{GMS}, but developed in more detail in~\cite{Mejmini,CMR2,Ye}. 

\begin{remark}\label{rem:ufi}
    The second author proposed in~\cite{M24Anatomy} a notion of \emph{ultrafilter limits for intervals}. Concretely, if $D$ is an ultrafilter on $\pts(K)$, we can say that a subset of a poset is \emph{$(D,\bar I,\varp_0)$-uf-linked} if it is $(D,\bar I,\varp_0)$-linked (interpreting $D$ as a fam into $\{0,1\}$) such that the fam in the extension witnessing the linkedness property is an unltrafilter (i.e.\ a fam into $\{0,1\}$). 
\end{remark}

Another version of ultrafilter limits with partitions is presented in.

\section{\texorpdfstring{$\mu$}{}-\texorpdfstring{$\FAM$}{}-linked forcing notions}\label{i:subs}

In this section, we define a new linkedness property, namely \emph{$\mu$-$\calY$-linkedness}, which will be the key notion to iterating using fams. 

\begin{definition}\label{i10}
    Denote by $\calY_*$ be the class of all pairs $(\Xi,\bar{I})$ such that $\Xi$ is a probability fam on some $\pts(K)$ and $\bar{I} = \la I_k: k\in K\ra \in \bfI^\fin_K$.\footnote{These pairs $(\Xi,\bar I)$ may be associated to different sets $K$.} 
    
    Fix a forcing notion $\bbP$, an infinite cardinal $\mu$, and a class $\calY \subseteq \calY_*$.
    \begin{enumerate}[label=\rm(\arabic*)]
    \item\label{i10a} Let $\calY \subseteq \calY_*$ be a class. The forcing notion $\bbP$ is \emph{$\mu$-$\calY$-linked}, if there exists a sequence $\langle Q_{\alpha, \varp} \colon \alpha < \mu \conj  \varp \in (0, 1)_{\QQQ} \rangle$ of subsets of $\bbP$ such that: 
    \begin{enumerate}[label=\rm(\roman*)]
        \item\label{i10aa}  $Q_{\alpha, \varp}$ is $(\Xi, \bar{I}, \varp)$-linked for all $(\Xi,\bar{I})\in \calY$,

        \item\label{i10ab} for every $\varp \in (0, 1)_{\QQQ}$, $\bigcup_{\alpha < \mu} Q_{\alpha, \varp}$ is dense in $\bbP$.
    \end{enumerate}
    
    \item\label{i10b} $\bbP$ is \emph{$\mu$-$\FAM$-linked}, if it is 
    $\mu$-$\calY_{\omega}$-linked, where $\calY_{\omega}$ is the set of all pairs $(\Xi, \bar{I})$ such that $\Xi$ is a free fam on $\calP(\omega)$, and $\bar{I}$ is a partition of $\omega$ into finite non-empty sets.  

    \item\label{i10d} The forcing notion $\Por$ is \emph{uniformly $\mu$-$\calY$-linked} if it is $\mu$-$\calY$-linked witnessed by a sequence $\langle Q_{\alpha, \varp} \colon \alpha < \mu \conj  \varp \in (0, 1)_{\QQQ} \rangle$ that also satisfies: for any $(\Xi,\bar I)\in \calY$, the $\Por$-name $\dot\Xi^*$ of the fam extending $\Xi$ witnessing the $(\Xi,\bar I,\varp)$-linkedness of each $Q_{\alpha, \varp}$ do not depend on $Q_{\alpha, \varp}$, concretely,
    there is some $\Por$-name $\dot \Xi^\ast$ of a fam on $\pts(K)$ ($K\coloneqq  \dom \bar I$) extending $\Xi$ such that, for any $\alpha<\mu$, $\varp\in(0,1)_{\QQQ}$ and $\bar q \in {}^{W} Q_{\alpha,\varp}$ ($W\coloneqq  \bigcup_{k\in K}I_k$), there is some $q\in\Por$ (i.e.\ $q=\lim^{Q_{\alpha,\varp},\Xi,\bar I,\varp}(\bar q)$) forcing
    \begin{equation}\label{i10eq}
        \int_K \frac{|\{ \ell\in I_k \colon q_\ell \in \dot G_{\bbP}\}|}{|I_k|}d\dot \Xi^\ast(k) \geq 1-\varp.
    \end{equation}
    In a similar way, we can define \emph{uniformly $\mu$-$\FAM$-linked}.
\end{enumerate}

    When $\mu = \aleph_{0}$, we write (uniformly) $\sigma$-$\calY$-linked, and $\sigma$-$\FAM$-linked, respectively.
\end{definition}

The uniform version of $\mu$-$\calY$-linked is motivated by the following fact. If $\Por$ is $\mu$-$\calY$-linked witnessed by $\langle Q_{\alpha, \varp} \colon \alpha < \mu \conj  \varp \in (0, 1)_{\QQQ} \rangle$ then, for any $(\Xi,\bar I)\in\calY$ and $\alpha<\mu$ and $\varp\in(0,1)_\QQQ$, by \autoref{i2} there is some $\bbP$-name $\dot{\Xi}^{\alpha,\varp}$ of a probability fam on $\calP(K)$ extending $\Xi$ such that, for any $\bar{q} \in \power{W}{Q_{\alpha, \varp}},$ there exists a condition $q\in\Por$ forcing the integral in \eqref{i10eq}. In principle, $\dot\Xi^{\alpha,\varp}$ depends on $(\Xi,\bar I,\alpha,\varp)$, so we define the uniform version where we get the same fam in the generic extension for $(\alpha,\varp)$ (but depending on $(\Xi,\bar I)$). 

The examples we know so far are uniform and, except for \autoref{i17}, extensive to check. We develop the theory a bit more before presenting examples.

In this subsection, unless otherwise stated, we assume that $\calY$ is a set of pairs as described in \autoref{i10} and $\bbP$ is a forcing notion.

Few modifications to \autoref{i10} as below give us equivalent formulations of the $\mu$-$\calY$-linkedness property. These also apply to the uniform version. In particular,~\autoref{i11.1}~\ref{i11.1d} will be useful to prove that (uniform) $\mu$-$\calY$-linkedness is a forcing property. 

\begin{observation}\label{i11.1}
    Let $\mu$ be an infinite cardinal.
    
    \begin{enumerate}[label= \rm (\arabic*)]
        \item If $\bar \varp = \la\varp_n \colon n<\omega\ra$ is a decreasing sequence of rationals in $(0,1)$ converging to $0$ and $\la Q_{\alpha,\varp_n} \colon \alpha<\mu,\ n<\omega\ra$ is a sequence of subsets of $\Por$ satisfying~\ref{i10aa}--\ref{i10ab} of \autoref{i10}, then $\Por$ is $\mu$-$\calY$-linked  witnessed by the sequence of subsets $Q_{\alpha,\varp}\coloneqq Q_{\alpha,\varp_n}$ where $n$ is the largest number such that $\varp_n\leq\varp$.
    
        \item Similar to the above, $\sigma$-$\calY$-linked is equivalent to: There is a sequence $\la Q_n\colon n<\omega\ra$ of subsets of $\Por$ such that, for some sequence $\varp = \la\varp_n \colon n<\omega\ra$ of rationals in $(0,1)$ converging to $0$, 
        \begin{enumerate}[label = \rm (\roman*)]
            \item each $Q_n$ is $(\Xi,\bar I,\varp_n)$-linked for all $(\Xi,\bar I)\in \calY$, and
            \item $\bigcup_{k\geq n}Q_k$ is dense in $\Por$ for all $n<\omega$.
        \end{enumerate}    
    
        \item\label{i11.1d} By \autoref{i6.1}, a forcing notion $\bbP$ is $\mu$-$\calY$-linked witnessed by $\langle Q_{\alpha, \varp} \colon \alpha < \mu \conj \varp \in (0, 1)_{\QQQ} \rangle$ iff it is witnessed by the sequence $\langle Q_{\alpha, \varp}^{\uparrow} \colon \alpha < \mu \conj \varp \in (0, 1)_{\QQQ} \rangle,$ since $\bigcup_{\alpha < \mu} Q_{\alpha, \varp}$ is dense in $\bbP$ iff, $\bigcup_{\alpha < \mu} Q_{\alpha, \varp}^{\uparrow} = \bbP,$ for all $\varp \in (0, 1)_{\QQQ}.$ As a consequence, we get an equivalent formulation of \autoref{i10}: we can replace condition \autoref{i10} \ref{i10a} \ref{i10ab} by ``For every $\varp \in (0, 1)_{\QQQ}, \, \bigcup_{\alpha < \mu} Q_{\alpha, \varp} = \bbP$''.
    
        \item It can also be assumed in \autoref{i10}~\ref{i10a} that, for $\alpha<\mu$ and $\varp'<\varp$, $Q_{\alpha,\varp'}\subseteq Q_{\alpha,\varp}$. Indeed, assume that $\la Q_{\alpha,\varp} \colon \alpha<\mu,\ \varp\in (0,1)_\QQQ\ra$ witnesses that $\Por$ is $\mu$-$\calY$-linked and, by~\ref{i11.1d}, we may assume that $Q^{\uparrow}_{\alpha,\varp}=Q_{\alpha,\varp}$. For $t\in {}^{<\omega}\mu$ and $\varp\in(0,1)_\QQQ$, define
        $Q'_{t,\varp}\coloneqq  \bigcap\{Q_{t(i),2^{-i}} \colon i<|t|\}$  if $2^{-|t|}\leq\frac{\varp}{2}$, otherwise $Q'_{t,\varp}\coloneqq Q'_{t',\varp}$, where $t'\in {}^{<\omega}\mu$ is the sequence of shortest length extending $t$ such that $2^{-|t'|}\leq\frac{\varp}{2}$ and $t'(i)=0$ for all $|t|\leq i< |t'|$. Then, $\la Q'_{t,\varp} \colon t\in{}^{<\omega}\mu,\ \varp\in (0,1)_\QQQ\ra$ witnesses that $\Por$ is $\mu$-$\calY$-linked and $Q'_{t,\varp'}\subseteq Q'_{t,\varp}$ whenever $\varp'<\varp$.
    \end{enumerate}
    
    All these observations above are valid for the uniform version after natural adjustments.
\end{observation}

Notice that $\varp=0$ was not considered in \autoref{i10}. By allowing it, we obtain a stronger notion, in which case using $\varp=0$ is enough because $Q_{\alpha,\varp}\coloneqq Q_\alpha$ can be used for all $\varp\in(0,1)_\Q$.

\begin{definition}\label{i10.2}
    Let $\mu$ be a cardinal. The forcing notion $\Por$ is said to be \emph{strongly $\mu$-$\calY$-linked} if there is a sequence $\Seq{Q_\alpha}{\alpha<\mu}$ of subsets of $\Por$ such that each $Q_\alpha$ is $(\Xi,\bar I,0)$-linked for all $(\Xi,\bar I) \in \calY$ and $\bigcup_{\alpha<\mu}Q_\alpha$ is dense in $\Por$. \autoref{i11.1}~\ref{i11.1d} can be applied here to get an equivalent formulation by replacing the last requirement by $\bigcup_{\alpha<\mu}Q_\alpha = \Por$.

    In analogy with \autoref{i10}, the notions of \emph{strongly $\mu$-$\FAM$-linked} and the uniform version can be defined, also using $\sigma$ when $\mu = \aleph_0$.
\end{definition}

The most basic example we can present is uniform and strong.

\begin{example}\label{i17}
    Any forcing notion $\bbP$ is uniformly and strongly $\vert \bbP \vert$-$\calY_*$-linked witnessed by its singletons. Indeed, for $(\Xi,\bar I)\in\calY_*$, using the fam-limits described in \autoref{i4}, any $\Por$-name $\dot\Xi^*$ of a fam extending $\Xi$ works.
      
    Consequently, any forcing notion $\bbP$ is uniformly and strongly $\vert \bbP \vert$-$\FAM$-linked. In particular, Cohen forcing is uniformly and strongly $\sigma$-$\calY_*$-linked.
\end{example}

The notion of (strongly) $\mu$-$\calY$-linked may affect the intersection number of the components. This is an immediate consequence of \autoref{i21} and~\autoref{i11.0}.

\begin{lemma}\label{i11.01}
    Assume that, for some $(\Xi', \bar{I}') \in \calY$, $\Xi'(\{k\in\dom\bar I' \colon |I'_k|=n\})=0$ for all $n<\omega$. If $\Por$ is (uniformly) $\mu$-$\calY$-linked witnessed by $\Seq{Q_{\alpha,\varp}}{\alpha<\mu,\ \varp\in(0,1)_\QQQ}$, then $\Int^\Por(Q_{\alpha,\varp}) \geq 1-\varp$ for all $\alpha<\mu$ and $\varp\in (0,1)_\QQQ$. In particular, $\Por$ is (uniformly) $\mu$-$\calY'$-linked witnessed by the same sequence, where $\calY'\coloneqq \calY \cup \{(\Xi,\bar I)\in\calY_* \colon \dom\bar I \text{ is finite}\}$. 
    The above is also valid for (uniformly and) strongly $\mu$-$\calY$-linked forcing notions (using $\varp=0$), which means that every $Q_{\alpha,0}$ is centered in $\Por$.
\end{lemma}

\begin{corollary}\label{i15}
    If $\calY$ is as in \autoref{i11.01}, then 
    every $\mu$-$\calY$-linked forcing notion is $\mu$-$m$-linked for all $2\leq m < \omega$,   and every strongly $\mu$-$\calY$-linked forcing notion is $\mu$-centered. 
\end{corollary}

\begin{PROOF}[\textbf{Proof}]{\ref{i15}}
    Let $\bbP$ be a $\mu$-$\calY$-linked forcing notion witnessed by $\langle Q_{\alpha, \varp} \colon \alpha < \mu \conj  \varp \in (0, 1)_{\bbQ} \rangle$. Assume that $m\leq 2$ and let $\varp_m \coloneqq \frac{1}{m+1}$. 
    By  \autoref{i11.01}, for all $\alpha < \mu$ we have that 
    \[ 1 - \frac{1}{m+1} < 1 - \varp_{m} \leq \Int(Q_{\alpha, \varp_{m}}).\]
    By employing~\autoref{s8}~\ref{s8g}, we get that $Q_{\alpha, \varp_{m}}$ is $m$-linked. Next,  for all $\alpha < \mu$ define $Q_{\alpha} \coloneqq \{ p \in \bbP \colon  \exists q \in Q_{\alpha, \varp_{m}}( q \leq p )\}.$ It is clear that each $Q_{\alpha}$ is $m$-linked because $Q_{\alpha, \varp_{m}}$ is, and $\bigcup_{\alpha < \mu} Q_{\alpha} = \bbP.$  
\end{PROOF}

\begin{corollary}\label{i15.1}
    Under the hypothesis of \autoref{i15}, 
    every $\mu$-$\calY$-linked forcing notion is $\theta$-$m$-Knaster for all $0 < m < \omega$ and any cardinal $\theta$ with $\cf(\theta)>\mu$, and every strongly $\mu$-$\calY$-linked forcing notion has precaliber $\theta$.
\end{corollary}

Next, we present two relevant examples of uniformly $\sigma$-$\FAM$-linked forcing notions.

\begin{example}\label{i15.2}
    Shelah~\cite{Sh00} proved, implicitly, that random forcing is $\sigma$-$\FAM$-linked. This was extended in \cite{MUR23}, where the second and the third authors proved that any measure algebra satisfying the $\mu$-density property\footnote{See \cite[Def.~2.5]{MUR23}. This is equivalent to having Maharam type $\mu$.} is uniformly $\mu$-$\FAM$-linked, even more, uniformly $\mu$-$\calY_{\rm free}$-linked where $\calY_{\rm free}\coloneqq \{(\Xi,\bar I)\in \calY_*\colon \Xi \text{ is free}\}$. Consequently, the \emph{measure algebra adding $\mu$-many random reals} 
    is $\mu$-$\calY_{\rm free}$-linked for any infinite cardinal $\mu$. 

    Recall that random forcing is not $\sigma$-centered, so it cannot be strongly $\sigma$-$\FAM$-linked.
\end{example}

\begin{example}\label{i11.000}
    In the construction of the iteration performed by Kellner, Shelah, and T\v anasie in~\cite{KST}, it was expected to use $\Eor$, the standard $\sigma$-centered forcing notion to add an eventually different real in $\baire$, to control $\non(\Mwf)$, but it does not have strong fam-limits (see~\autoref{i45}). For this reason, they used an alternative forcing $\tilde\Eor$, which comes from a creature framework from Shelah and Horowitz~\cite{HSh}, 
    and they proved that it is $\sigma$-$\FAM$-linked. On the other hand, this forcing is not $\sigma$-centered, so it cannot be strongly $\sigma$-$\FAM$-linked.

    Let $\calY_{\rm uap}\coloneqq \{(\Xi,\bar I)\in \calY_*\colon \Xi \text{ has the uap}\}$. 
    Following~\cite{HSh,KST}, the second author~\cite{M24Anatomy} presents a framework of forcing notions, containing $\tilde{\bbE}$, that are uniformly $\sigma$-$\calY_{\rm uap}$-linked. 
\end{example}

Using the connections with Fr\'echet linkedness, as a direct consequence of~\autoref{i9}, we get:

\begin{lemma}\label{i13}
    Let $\kappa$ be an infinite cardinal. 
    Assume that 
    there exists a $(\Xi, \bar{I}) \in \calY$ such that $K\coloneqq \dom \bar I$ has size $\kappa$ and the $\Xi$-measure of any subset of cardinality ${<} \kappa$ is zero. If $\Por$ is a $\mu$-$\calY$-linked forcing notion and $\Vdash$``\/$[K]^{<\kappa}\cap V$ is cofinal in $[K]^{<\kappa}$", then $\Por$ is $\mu$-$\Fr_{\kappa}$-linked. 
\end{lemma}

\begin{corollary}\label{i13.1}
     If there exists $(\Xi, \bar{I}) \in \calY$ such that $\Xi$ is a free fam and $|\dom\bar I| = \aleph_0$, then every $\mu$-$\calY$-linked forcing notion is $\mu$-$\Fr$-linked. 
\end{corollary}

However, the converse is not true: 

\begin{remark}\label{i45}
    The first and second authors proved that $\sigma$-$\FAM$-linked forcing notions do not increase $\non(\Ewf)$, where $\Ewf$ denotes the ideal on ${}^\omega2$ generated by the $F_\sigma$ measure zero sets (see \autoref{e090} and~\ref{a69});
    on the other hand, they also proved that the standard $\sigma$-centered poset $\Eor$ that adds an eventually different real on ${}^\omega \omega$ increases $\non(\Ewf)$ (see~\cite[Thm.~4.10]{CaraboutE}). Hence, $\Eor$ is not $\sigma$-$\FAM$-linked, but it is $\sigma$-Fr-linked (see~\cite[Lem.~3.29]{mejiavert}). 
\end{remark}

\autoref{i11.1}~\ref{i11.1d} is the key point to prove that $\Xi$-$\calY$-linkedness is a forcing property. Recall that $\varphi(x)$ is a \emph{forcing property} iff for any pair of forcing notions $\bbP, \, \bbQ,$ if there is a dense embedding $\iota \colon \bbP \to \bbQ$, then $\varphi(\bbP)$ is equivalent to $\varphi(\bbQ)$; and $\varphi(x)$ is a \emph{hereditary forcing property} if, for any pair of forcing notions $\bbP$ and $\bbQ$, if $\Por\lessdot \Qor$ and
$\varphi(\bbQ)$ holds, then $\varphi(\bbP)$ holds.  

\begin{theorem}\label{i11.2} 
     $\mu$-$\calY$-linkedness and its strong version are forcing properties.
\end{theorem}

\begin{PROOF}[\textbf{Proof}]{\ref{i11.2}}
    Let $\bbP, \bbQ$ be forcing notions and let $\iota \colon \bbP \to \bbQ$ be a dense embedding. 
   
    Assume  $\bbQ$ is $\mu$-$\calY$-linked witnessed by $\langle Q_{\alpha, \varp} \colon \alpha < \mu \conj \varp \in (0, 1)_{\QQQ} \rangle$. By \autoref{i11.1}, we may assume that $\bigcup_{\alpha<\mu} Q_{\alpha,\varp} = \bbQ$ for any $\varp\in(0,1)_\QQQ$.      
    For any $\alpha < \mu$ and $\varp \in (0, 1)_{\QQQ}$ define $P_{\alpha, \varp} \coloneqq \iota^{-1}[Q_{\alpha, \varp}].$ We will prove that  $\bbP$ is $\mu$-$\calY$-linked witnessed by $\langle P_{\alpha, \varp} \colon \alpha < \mu \conj \varp \in (0, 1)_{\QQQ} \rangle$. 
 
    Let $(\Xi, \bar{I}) \in \calY$.  By~\autoref{i5.1}, we have that each $P_{\alpha, \varp}$ is $(\Xi, \bar{I}, \varp)$-linked and, 
    for $\varp \in (0, 1)_{\QQQ}$ and $p \in \bbP$, we can choose an $\alpha_{0} < \mu$ such that $\iota(p) \in Q_{\alpha_{0}, \varp}$,  hence $p \in P_{\alpha_{0}, \varp}$. Thus $\bigcup_{\alpha < \mu} P_{\alpha, \varp}$ is dense in $\bbP,$ moreover it is equal to $\bbP.$ 

    Conversely, assume $\bbP$ is $\mu$-$\calY$-linked witnessed by $\langle P_{\alpha, \varp} \colon \alpha < \mu \conj \varp \in (0, 1)_{\QQQ} \rangle.$  For any $\alpha < \mu$ and $\varp \in (0, 1)_{\QQQ}$ let $Q_{\alpha, \varp} \coloneqq \iota [ P_{\alpha, \varp} ].$ We will prove that $\bbQ$ is $\mu$-$\calY$-linked witnessed by the sequence $\langle Q_{\alpha, \varp} \colon \alpha < \mu \conj \varp \in (0, 1)_{\QQQ} \rangle.$ 
    
    Note that each $Q_{\alpha, \varp}$ is $(\Xi, \bar{I}, \varp)$-linked by~\autoref{i5}.

    Now, let $\varp \in (0, 1)_{\QQQ}$ and $q \in \bbQ.$ By the density of $\iota,$ there exists a $p \in \bbP$ such that $\iota(p) \leq q.$ Then, there are $p'\leq p$ in $\Por$ and $\alpha_{0} < \mu$  such that $p' \in P_{\alpha_{0}, \varp}.$ As a consequence, $\iota(p') \in Q_{\alpha_{0}, \varp}$ and $\iota(p') \leq q$. Hence $\bigcup_{\alpha < \mu} Q_{\alpha, \varp}$ is dense in $\bbQ.$

    The same proof is valid for the strong version, just replace $\varp$ by $0$.
\end{PROOF}

The uniform version is also a forcing property.

\begin{theorem}\label{i11.3}
    Uniformly  $\mu$-$\calY$-linkedness and its strong version are forcing properties.
\end{theorem}
\begin{PROOF}[\textbf{Proof}]{\ref{i11.3}}
    Let $\bbP, \bbQ$ be forcing notions and let $\iota \colon \bbP \to \bbQ$ be a dense embedding. 

    Assume that $\bbP$ is uniformly $\mu$-$\calY$-linked witnessed by $\langle P_{\alpha, \varp} \colon \alpha < \mu \conj \varp \in (0, 1)_{\QQQ} \rangle.$  For any $\alpha < \mu$ and $\varp \in (0, 1)_{\QQQ}$, let $Q_{\alpha, \varp} \coloneqq \iota [ P_{\alpha, \varp} ].$ By \autoref{i11.2}, $\bbQ$ is $\mu$-$\calY$-linked witnessed by the sequence $\langle Q_{\alpha, \varp} \colon \alpha < \mu \conj \varp \in (0, 1)_{\QQQ} \rangle$, so it is enough to deal with the additional condition in~\autoref{i10}~\ref{i10d}, for which it is enough to proceed like in the proof of \autoref{i5.1}. On the other hand, if $\bbQ$ is uniformly $\mu$-$\calY$-linked witnessed by $\langle Q_{\alpha, \varp} \colon \alpha < \mu \conj \varp \in (0, 1)_{\QQQ} \rangle$, we can proceed analogously as in the proof of \autoref{i5}. 
     
    The same argument is valid for the strong version, just replace $\varp$ by $0$.
\end{PROOF}

\begin{question}\label{Q31}
    Is $\mu$-$\calY$-linkedness (and its variations) a hereditary forcing property?
\end{question}

This is partially true: Let $\iota\colon\Por\to\Qor$ be a complete embedding between forcing notions and assume that $\Qor$ is $\mu$-$\calY$-linked witnessed by $\langle Q_{\alpha, \varp} \colon \alpha < \mu \conj \varp \in (0, 1)_{\QQQ} \rangle$. For each $c=(\Xi,\bar I)\in\calY$, $\alpha<\mu$ and $\varp\in(0,1)_\QQQ$, let $\dot \Xi_{c,\alpha,\varp}$ be a $\Por$-name witnessing that $Q_{\alpha,\varp}$ is $(\Xi,\bar I,\varp)$-linked. If $\dot \Xi_{c,\alpha,\varp}\frestr V^\Por$ has a $\Por$-name $\dot \Xi^-_{c,\alpha,\varp}$ for all $c\in \calY$, $\alpha<\mu$ and $\varp\in(0,1)_\QQQ$, this can be used to prove that $\iota^{-1}[Q_{\alpha,\varp}]$ is $(\Xi,\bar I,\varp)$-linked, and then $\Por$ is $\mu$-$\calY$-linked as in the proof of \autoref{i11.2}. This also works for the strong and uniform versions. However, it is unclear how to proceed without the assumption that each $\dot \Xi_{c,\alpha,\varp}\frestr V^\Por$ has a $\Por$-name.

The following is a characterization of ``uniformly $\mu$-$\calY$-linked" without using integration, which can be proved similarly as \autoref{i6}. The proofs of the linkedness property for \autoref{i15.2} and~\ref{i11.000} use this characterization. 

\begin{theorem}[cf.~{\cite[Thm.~3.10]{M24Anatomy}}]\label{i10.1}
    Let $\mu$ be a cardinal and let $\Seq{Q_{\alpha,\varp}}{\alpha<\mu,\ \varp\in(0,1)_\Q}$ a sequence of subsets of $\Por$ such that $\bigcup_{\alpha<\mu}Q_{\alpha,\varp}$ is dense in $\Por$ for all $\varp\in(0,1)_\QQQ$. Then, the following statements are equivalent. 
    \begin{enumerate}[label=\normalfont (\alph*)]
    \item\label{i10.1i} $\Por$ is uniformly $\mu$-$\calY$-linked witnessed by $\Seq{Q_{\alpha,\varp}}{\alpha<\mu,\ \varp\in(0,1)_\Q}$.

    \item\label{i10.1ii'} For each $(\Xi,\bar I)\in \calY$, $\alpha<\mu$ and $\varp\in(0,1)_\Q$, there is some ${\lim}^{\alpha,\varp}\colon {}^WQ_{\alpha,\varp}\to \Por$ where $W\coloneqq W_{\bar I}$ such that, for any
    \begin{multicols}{2} 
    \begin{itemize}
        \item $(\Xi,\bar I)\in\calY$,
        \item $i^*<\omega$,
        \item $(\alpha_i,\varp_i)\in \mu\times (0,1)_\Q$,
        \item $\bar r^i = \Seq{r^i_\ell}{\ell\in W} \in {}^W Q_{\alpha_i,\varp_i}$ for $i<i^*$,
        \item $P\in \Pbf^\Xi$,
        \item $\varp'>0$, and
        \item $q\in \Por$ stronger than $\lim^{\alpha_i,\varp_i} (\bar r^i)$ for all $i<i^*$,
    \end{itemize}
    \end{multicols}
    there are some $q'\leq q$ in $\Por$, a finite non-empty $u\subseteq\dom\bar I$ and a probability measure $\Xi^-$ on $\pts(u)$ such that
    \begin{enumerate}[label = \normalfont (\arabic*)]
        \item\label{-fm1'} $\left| \Xi^-(u\cap b) - \Xi(b)\right| < \varp'$ for all $b\in P$, and
        \item\label{-fm2'} $\displaystyle \int_{u}\frac{|\set{\ell\in I_k}{q' \leq r^i_\ell}|}{|I_k|}d\Xi^-(k) > 1-\varp_i-\varp'$ for all $i<i^*$.
    \end{enumerate}
    \end{enumerate}
    Furthermore, in the case $\Xi$ has the uap for all $(\Xi,\bar I)\in\calY$ then~\ref{i10.1ii} can be added to the list of equivalent statements, even more, in~\ref{i10.1i}, $\dot\Xi^*$ (as in \autoref{i10}~\ref{i10d}) can be found with the uap.
    \begin{enumerate}[resume*]
    \item\label{i10.1ii} For each $(\Xi,\bar I)\in \calY$, $\alpha<\mu$ and $\varp\in(0,1)_\Q$, there is some ${\lim}^{\alpha,\varp}\colon {}^WQ_{\alpha,\varp}\to \Por$ where $W\coloneqq W_{\bar I}$ such that, for any
    \begin{multicols}{2} 
    \begin{itemize}
        \item $(\Xi,\bar I)\in\calY$,
        
        \item $i^*<\omega$,
        
        \item $(\alpha_i,\varp_i)\in \mu\times (0,1)_\Q$,
        
        \item $\bar r^i = \Seq{r^i_\ell}{\ell\in W} \in Q_{\alpha_i,\varp_i}^W$ for $i<i^*$,
        
        \item $P\in \Pbf^\Xi$,
        
        \item $\varp'>0$, and
        
        \item $q\in \Por$ stronger than $\lim^{\alpha_i,\varp_i} (\bar r^i)$ for all $i<i^*$,
    \end{itemize}
    \end{multicols}
    there are some $q'\leq q$ in $\Por$ and $u\subseteq \dom\bar I$ finite non-empty such that
    \begin{enumerate}[label = \normalfont (\arabic*)]
        \item\label{-fm1} $\displaystyle \left| \frac{|u\cap b|}{|u|} - \Xi(b)\right| < \varp'$ for all $b\in P$, and
        \item\label{-fm2} $\displaystyle \frac{1}{|u|}\sum_{k\in u}\frac{|\set{\ell\in I_k}{q' \leq r^i_\ell}|}{|I_k|} > 1-\varp_i-\varp'$ for all $i<i^*$.
    \end{enumerate}
\end{enumerate}
\end{theorem}

A similar characterization of ``uniformly and strongly $\mu$-$\calY$-linked" can be formulated from \autoref{i24}.

\begin{remark}\label{rem:ufi2}
    In connection with \autoref{rem:ufi}, the notions presented in this section can be defined in the context of \emph{ultrafilter limits for intervals}. Concretely, if $\Xi$ is an ``ultrafilter'' for any $(\Xi,\bar I)\in\calY$, we say that a poset is \emph{$\mu$-$\calY$-uf-linked} if it is $\mu$-$\calY$-uf-linked but each $Q_{\alpha,\varp}$ (witnessing the forcing property) is $(\Xi,\bar I,\varp)$-uf-linked for all $(\Xi,\bar I)\in\calY$. All the results of this section can be easily reformulated in this context.
\end{remark}


\section[Iterating with \texorpdfstring{$\mu$}{}-\texorpdfstring{$\FAM$}{}-linked forcing notions]{Iterating with \texorpdfstring{$\mu$}{}-\texorpdfstring{$\FAM$}{}-linked forcing notions
}\label{4.3}

In this section, we will show that the properties needed to be able to generalize the iteration of \cite{Sh00} and establish extension theorems when iterating using fams (see \autoref{k39} and \autoref{k42}), are precisely those that define the $\mu$-$\calY$-linked forcing notions.  We will define an iteration frame based on \cite{Sh00}, but not iterating with random forcing, but in a more general context: using $\mu$-$\calY$-linked forcing notions. 

Our iterations will allow fam-limit functions as for $(\Xi,\bar I,\varp)$-linked sets. This feature will be described in more detail in \autoref{famlimits}.

\subsection{The iteration structure}\ 

First, we introduce a general framework that will allow us to iterate using fams and even ultrafilters (see \autoref{k28}).

\begin{definition}\label{k4}
     We say that $\bft =\langle  \bbP_{\beta}, \dot{\bbQ}_{\alpha},  \bbP_{\alpha}^{-}, \vec{Q}_{\alpha}, \theta_{\alpha} \colon \beta \leq \pi,\ \alpha <  \pi \rangle$ is a \emph{prepared iteration} if it
    satisfies the following conditions: 
    
    \begin{enumerate}[label=(\arabic*)]
        \item\label{k4a} $ \bbP_{\pi} = \langle \bbP_{\alpha},  \dot{\bbQ}_{\alpha} \colon \alpha < \pi \rangle$ is a FS (finite support) iteration;
        
        \item\label{k4b} $\pi = F^{\bft} \cup U^{\bft}$ with $F^{\bft} \cap U^{\bft} = \emptyset$ and, for any $\alpha < \pi,$ $\dot{\bbQ}_{\alpha}$ is a $\bbP_{\alpha}^{-}$-name, 
        
        \begin{itemize}
            \item if $\alpha \in F^{\bft},$ then  $\vec{Q}_{\alpha} = \langle \dot{Q}_{\zeta, \varp}^{\alpha} \colon \zeta < \theta_{\alpha},\ \varp \in (0, 1)_{\QQQ} \rangle$ is a sequence of $\bbP_{\alpha}^{-}$-names of subsets of $\Qnm_\alpha$ such that, for any $\varp \in (0, 1)_{\bbQ}$, $\Vdash_{\bbP^-_\alpha}$``$\bigcup_{\zeta < \theta_{\alpha}} \dot{Q}_{\zeta, \varp}$ is dense in $\dot{\bbQ}_{\alpha}$", and

            \item if $\alpha \in U^{\bft}$, then  $\vec{Q}_{\alpha} = \langle \dot{Q}_{\zeta, 0}^{\alpha} \colon \zeta < \theta_{\alpha} \rangle$ is a sequence of $\bbP_{\alpha}^{-}$-names of subsets of $\Qnm_\alpha$ such that $\Vdash_{\bbP^-_\alpha}$``$\bigcup_{\zeta < \theta_{\alpha}} \dot{Q}_{\zeta, 0}$ is dense in $\dot{\bbQ}_{\alpha}$";
        \end{itemize}

        \item\label{k4c} for any $\alpha < \pi,  \bbP_{\alpha}^{-} \subsetdot \bbP_{\alpha}^{\bullet},$ where  
        $$\bbP_{\alpha}^{\bullet} \coloneqq 
        \left \{ p \in \bbP_{\alpha} \colon \forall \xi \in \dom(p)\ \left[p(\xi) \in V^{\bbP_{\alpha}^{-}}  \conj    \exists \zeta < \theta_{\xi}  \exists \varp  {\in}  [0, 1)_{\QQQ}\ (  \Vdash_{\bbP_{\xi}^{-}} \text{``} p(\xi) \in \dot{Q}_{\zeta, \varp}^{\xi} \text{''} ) \right ] \right \}. $$
    \end{enumerate}
    
        \begin{figure}[ht]
        \centering
        \begin{tikzpicture}

            \node at (5.4, 2) {$\bbP_{\alpha}$}; 
            
            \node at (5.4, 0) {$\bbP_{\alpha}^{-}$};

            \draw[->, color=redun] (5.9, 2) -- (7.7, 2);

            \draw[->, color=redun] (5.9, 0) -- (7.5, 0);

            \draw[->,  color=redun] (5.4, 0.3) -- (5.4, 1.7);

            \draw[->,  color=redun] (8.4, 0.3) -- (8.4, 1.7);

            \node at (8.4, 2) {$\bbP_{\alpha + 1}$} ; 

            \node at (8.4, 0) {$\bbP_{\alpha}^{-} \ast \dot{\bbQ}_{\alpha}$}; 

            \node[greenun] at (6.75, 2.35) { \small $\dot{\bbQ}_{\alpha}$}; 
            
            \node[greenun] at (6.9, -0.35) {\small $\dot{\bbQ}_{\alpha}$}; 
            \end{tikzpicture}

            \caption{Successor step in prepared iterations.}
            \label{f657}
        \end{figure}   

        Also, if $\kappa$ is an infinite cardinal, $\bft$ is said to be a \emph{$\kappa$-prepared iteration} if it is a prepared iteration and, for any $\alpha < \pi$, $\theta_{\alpha} < \kappa$. 
\end{definition}

As with any iteration, a prepared iteration is constructed by recursion on $\beta\leq \pi$. Having reached $\Por_\alpha$ for some $\alpha<\pi$, it has to be guaranteed that $\Por^\bullet_\alpha$ is dense in $\Por_\alpha$, so that the choice of $\Por^-_\alpha\subsetdot \Por^\bullet_\alpha$ and $\Qnm_\alpha$ makes sense in the definition above. It will be clear later in \autoref{k12} that  $\bbP_{\beta}^{\bullet}$ is dense in $\bbP_{\beta}$ for all $\beta \leq \pi.$ It is clear that $\bbP_{\alpha}^{\bullet} \subseteq \bbP_{\beta}^{\bullet},$ whenever $\alpha < \beta \leq \pi.$  
 
Before we advance more, we introduce the following notation.

\begin{definition}\label{k7}
    Let $\bft = \langle \bbP_{\beta}, \dot{\bbQ}_{\alpha}, \bbP_{\alpha}^{-},  \vec{Q}_{\alpha},  \theta_{\alpha} \colon \beta \leq \pi,\ \alpha <  \pi \rangle$ be a prepared iteration and $\gamma < \pi.$ 

    \begin{enumerate}[label=(\arabic*)]
        \item\label{k7a} We say that $\pi$ is the \emph{length} of $\bft.$ 
    
        \item\label{k7b} The \emph{restriction of $\bft$ to $\gamma$}, denoted by $\bft {\rest}  \gamma,$ is defined in the natural way: $$\bft {\rest} \gamma \coloneqq \langle \bbP_{\beta}, \dot{\bbQ}_{\alpha},  \bbP_{\alpha}^{-}, \vec{Q}_{\alpha},  \theta_{\alpha} \colon \beta \leq \gamma,\  \alpha < \gamma \rangle.$$ 
    \end{enumerate} 
\end{definition}

It is clear that for any prepared iteration $\bft$ and any $\gamma \leq \pi, $ $\bft  {\rest}  \gamma$ is a prepared iteration. Similarly for $\kappa$-prepared iterations. 
 
\begin{notation}\label{k8}
    When we say ``$\bft$ is a $\kappa$-prepared iteration'', it is understood that $\bft$ has length $\pi$ partitioned into sets $\{ F^{\bft}, U^{\bft} \}$, and $\bft = \langle \bbP_{\beta}, \dot{\bbQ}_{\alpha}, \bbP_{\alpha}^{-}, \vec Q_{\alpha}, \theta_{\alpha} \colon \beta \leq \pi,\ \alpha < \pi \rangle,$   unless otherwise specified.
\end{notation}

It is not hard to see that $\bbP_{\beta}^{\bullet}$ is dense in $\bbP_{\beta},$ in fact, we can prove a stronger fact. For this, we define the following set. 

\begin{definition}\label{k9}
    Let $\bft$ be a prepared iteration. For any function  $\bar{\varp} \colon \omega \to (0, 1)_{\QQQ}$ define $D_{\bar{\varp}}^{\bft}$ as the set of conditions $p \in \bbP_{\pi}^{\bullet}$ such that:

    \begin{enumerate}[label=(\alph*)]
        \item\label{k9a} $\dom(p) \cap F^{\bft} = \{ \alpha_{n} \colon n < n^{\ast} \}$ where $n^* \coloneqq |\dom(p) \cap F^\bft|$,

        \item\label{k9b} $\alpha_{n} < \alpha_{m}$ whenever $m < n < n^{\ast},$ that is, $\dom(p)\cap F^{\bft}$ is arranged in decreasing order, and 

        \item\label{k9c} for all $n < n^{\ast},$ there exists some $\zeta < \theta_{\alpha_{n}}$ such that $\forces_{\Por^-_\alpha} \text{``}p(\alpha_{n}) \in \dot{Q}_{\zeta,  \bar{\varp}(n)}^{\alpha_{n}}$\text{''}.  
    \end{enumerate}

    When the context is clear, we denote $D_{\bar{\varp}}^{\bft}$ simply as $D_{\bar{\varp}}.$ Notice that $\alpha <\beta \leq \pi$ implies $D_{\bar{\varp}}^{\bft \rest \alpha} \subseteq D_{\bar{\varp}}^{\bft \rest \beta}.$
\end{definition}

The structure of prepared iterations guarantees that $D_{\bar{\varp}}$ is dense in the final step of the iteration: 

\begin{lemma}\label{k12}
    Let $\bft$ be a prepared iteration. Then, for any $\alpha \leq \pi$ and $\bar{\varp} \colon \omega \to (0, 1)_\QQQ,$ the set $ D_{\bar{\varp}}^{\bft {\rest}  \alpha}$ is dense in $\bbP_{\alpha}.$ In particular, $D_{\bar{\varp}}^{\bft}$ is dense in $\bbP_{\pi}.$ 
\end{lemma}

\begin{PROOF}[\textbf{Proof}]{\ref{k12}}
    Proceed by transfinite induction on $\alpha \leq \pi$. When $\alpha = 0$ or $\alpha$ is a limit ordinal, the result is clear. So we only deal with the successor case. Suppose that $\alpha = \xi + 1$ and assume that $D_{\bar{\epsilon}}^{\bft \, {\rest} \, \xi}$ is dense in $\bbP_{\xi}$ for every function $\bar{\epsilon} \colon \omega \to (0, 1)_\QQQ.$ Let $\bar{\varp} \colon \omega \to (0, 1)_\QQQ$ and let $p \in \bbP_{\alpha}.$ Without loss of generality, we can assume that $\xi \in \dom(p)$. By induction hypothesis, $D_{\bar{\varp}}^{\bft  {\rest} \xi}$ is dense in $\bbP_{\xi}$, so we get that $\bbP_{\xi}^{\bullet}$ is dense in $\bbP_{\xi}$  and by~\autoref{k4}~\ref{k4c}, $\bbP_{\xi}^{-} \subsetdot  \bbP_{\xi}.$ As a consequence, we can find a $\zeta < \theta_{\xi}$, a $\bbP_{\xi}^{-}$-name $\dot{q}$ of a member of $\dot{Q}_{\zeta, \bar{\varp}(0)}^{\xi}$ and $r_{0} \in \bbP_{\xi}$ such that $r_{0} \leq p  {\rest}  \xi$ and $r_{0} \forces_{\xi}\text{``} \dot{q} \leq p(\xi) \text{''},$ this in the case $\xi\in F^{\bft}$, otherwise replace $\bar \varp(0)$ by $0$.
    
    Define $\bar{\epsilon} \colon \omega \to (0, 1)_\QQQ$ by $\bar{\epsilon}(n) \coloneqq \bar{\varp}(n +1)$ for any $n < \omega$ in the case $\xi\in F^{\bft}$, otherwise let $\bar \epsilon \coloneqq \bar \varp$.  Since $r_{0} \in \bbP_{\xi},$ by induction hypothesis, there is some $r_{1} \in D_{\bar{\epsilon}}^{\bft  {\rest} \xi}$ such that $r_{1} \leq r_{0}.$ We set $r \coloneqq r_{1} \cup \{ (\xi, \dot{q}) \}.$ It is clear that $r\in\Por^\bullet_\alpha$ and $r \leq_{\alpha} p.$  Now we show that $r \in D_{\bar{\varp}}^{\bft {\rest} \alpha}.$ This is clear in the case $\xi\in U^{\bft}$, so we assume from now on that $\xi\in F^{\bft}$.
    Since $r_{1} \in D_{\bar{\epsilon}}^{\bft  {\rest}  \xi},$ there is $m^{\ast} < \omega$ such that $\dom(r_{1})\cap F^{\bft} = \{ \alpha_{n} \colon n < m^{\ast} \}$ and this enumeration is decreasing. We define $n^{\ast} \coloneqq m^{\ast} + 1$ and 
    $$ \gamma_{n} \coloneqq \left\{ \begin{array}{ll}
            \alpha_{n -1} &   \text{if  $0 < n < n^{\ast}$,} \\[1ex]
             
            \xi &  \text{if $n = 0$,}
             \end{array}
        \right. $$
    hence, $\dom(r)\cap F^{\bft} = \{ \gamma_{n} \colon n < n^{\ast}\}$ and since $\alpha_{0} < \xi,$ it is a decreasing enumeration. Let $0 < n < n^{\ast}.$ Given that $r_{1} \in D_{\bar{\epsilon}}^{\bft {\rest} \xi},$ by \autoref{k9}~\ref{k9c}, there exists $\xi < \theta_{\alpha_{n-1}}$ such that $\forces_{\Por^-_{\alpha_{n-1}}}$``$r_{1}(\alpha_{n-1}) \in \dot{Q}_{\xi,  \bar{\epsilon}(n -1)}^{\alpha_{n-1}}$''. However, as in this case $\gamma_{n} = \alpha_{n-1}$ and $\bar{\epsilon}(n-1) = \bar{\varp}(n),$ we get that $\forces_{\Por^-_{\gamma_{n}}} $``$r(\gamma_{n}) \in \dot{Q}_{\xi,  \bar{\varp}(n)}^{\gamma_{n}}$''.  Finally, if $n = 0,$ then $\gamma_{0} = \xi$ and we already know that $\forces_{\Por^-_\xi}$``$r(\xi) \in \dot{Q}_{\zeta, \bar{\varp}(0)}^{\xi}$''. Thus, $r \in D_{\bar{\varp}}^{\bft  {\rest} \alpha}.$  \qedhere
\end{PROOF}

As a consequence, as mentioned before: 

\begin{corollary}\label{k14}
        If $\bft$ is a prepared iteration then, for any $\alpha \leq \pi, \bbP_{\alpha}^{\bullet}$ is dense in $\bbP_{\alpha}.$  In particular, for any $\alpha < \pi$,  $\bbP_{\alpha}^{-}  \subsetdot  \bbP_{\alpha}.$ 
\end{corollary}
 
One of the fundamental parameters in Shelah's iteration is what he calls \emph{blueprints} (see \cite[Def.~2.9]{Sh00}). We are going to replace this notion with that of \emph{guardrail}  (as in \cite{GMS} and \cite{KST}) since we consider that this facilitates the formalization and compression of the frame for the iterations using fams.

\begin{definition}\label{k16}
    Let  $\gamma, \zeta$ be an ordinals.
    
    \begin{enumerate}
        \item  A function $g$ is said to be a \emph{guardrail for $\gamma$} if $\gamma \subseteq \dom(g)$ and $\ran(g) \subseteq \Ord\times[0,1)_{\QQQ}$. Denote by $g_0$ and $g_1$ the \emph{components of $g$}, i.e.\ 
        $g_{0} \colon \dom(g) \to \Ord$ and $g_{1} \colon \dom(g) \to [0, 1)_{\QQQ}$ such that, for any $\xi \in \dom(g)$, $g(\xi) = (g_{0}(\xi), g_{1}(\xi)).$ 

        \item A set $\calG$ of guardrails is a \emph{complete set of guardrails over $\zeta$} if, for any $\sigma \colon X \to \zeta \times [0, 1)_{\QQQ}$ with $X \in [\gamma]^{< \aleph_{1}}$  there is some $g \in \calG,$ such that $\sigma \subseteq g.$ 
    \end{enumerate} 
\end{definition}

In a prepared iteration $\bft$,  
a guardrail $g$ is in charge of choosing a subset $\dot{Q}^\alpha_{\zeta, \varp}$ of $\Qnm_\alpha$ at each step $\alpha$ of the iteration,\footnote{Which justifies the name ``guardrail''.} as it gives us information about the coordinates along the iteration (see \autoref{f996}). In this sense, it may happen that $g$ makes  \emph{illegal choices}, i.e.\ $g_0(\alpha)\notin\theta_\alpha$, or $g_1(\alpha)\neq 0$ when $\alpha\in U^\bft$, for some $\alpha$, but such illegal choices are harmless in the sense that they happen to just be ignored. 
For a legal choice, whether $g_1(\alpha)$ is $0$ or not, tells us whether the coordinate $\alpha$ is in $U^{\bft}$ or in $F^{\bft}$.

\begin{figure}[ht]
    \centering
    \begin{tikzpicture}            

        \draw[-, color=redun] (-1, 2) -- (12, 2);

        \node[greenun] at (1, 2) {$\bullet$};
            \node[greenun] at (1, 1.5) {$\alpha_{0}$};
            \node at (1, 2.5) {$\dot{q}_{0} \in \dot{Q}_{\greenunq{\zeta_{0}, \varp_{0}}}^{\alpha_{0}}$};

            \node[greenun] at (4, 2) {$\bullet$};
            \node[greenun] at (4, 1.5) {$\alpha_{1}$};
            \node at (4, 2.5) {$\dot{q}_{1} \in \dot{Q}_{\greenunq{\zeta_{1}, \varp_{1}}}^{\alpha_{1}}$};

            \node[greenun] at (7, 1.6) {$\dots$};

            \node[greenun] at (10, 2) {$\bullet$};
            \node[greenun] at (10, 1.5) {$\alpha_{\ell}$};
            \node at (10, 2.5) {$\dot{q}_{\ell} \in \dot{Q}_{\greenunq{\zeta_{\ell}, \varp_{\ell}}}^{\alpha_{\ell}}$};

            \node[greenun] at (12, 1.5) {$\pi$};
            \node[greenun] at (12, 2) {$\bullet$};

            \node[greenun] at (-1, 1.5) {$0$};
            \node[greenun] at (-1, 2) {$\bullet$};    


            \draw[-, color=redun] (-1, 0) -- (12, 0);

            \node[greenun] at (-1, -0.5) {$0$};
            \node[greenun] at (-1, 0) {$\bullet$};

            \node[greenun] at (1, 0) {$\bullet$};
            \node[greenun] at (1, -0.5) {$\alpha_{0}$};
            \node at (1, 0.5) {$\dot{q}_{0} \in \dot{Q}_{\redunq{g(\alpha_{0})}}^{\alpha_{0}}$};

            \node[greenun] at (4, 0) {$\bullet$};
            \node[greenun] at (4, -0.5) {$\alpha_{1}$};
            \node at (4, 0.5) {$\dot{q}_{1} \in \dot{Q}_{\redunq{g(\alpha_{1})}}^{\alpha_{1}}$};

            \node[greenun] at (7, -0.4) {$\dots$};

            \node[greenun] at (10, 0) {$\bullet$};
            \node[greenun] at (10, -0.5) {$\alpha_{\ell}$};
            \node at (10, 0.5) {$\dot{q}_{\ell} \in \dot{Q}_{\redunq{g(\alpha_{\ell})}}^{\alpha_{\ell}}$};

            \node[greenun] at (12, -0.5) {$\pi$};
            \node[greenun] at (12, 0) {$\bullet$};
    \end{tikzpicture}

    \caption{Guardrails will decide coordinates throughout the iteration in the sense that $g(\alpha_{\ell}) = (\zeta_{\ell}, \varp_{\ell}).$}
    \label{f996}
\end{figure}    

The tool we are going to use to obtain complete sets of guardrails is known as the Engelking and Kar\l owicz Theorem (see ~\cite{Engelking}). Here, $\theta^-\coloneqq \sup\set{\mu\in\Card}{\mu<\theta}$.

\begin{theorem}[{\cite{Rinot}}]\label{k20}
    Let $\theta$, $\mu$ and $\chi$ be cardinals such that  $\aleph_0 \leq \theta^- \leq \mu$ and $\chi \leq 2^{\mu}.$ Then there exists $\calG \subseteq {}^{\chi} \mu$ such that $\vert \calG \vert \leq \mu^{< \theta}$ and every function $f \colon X \to \mu$ with $X \in [\chi]^{< \theta}$ can be extended by a function in $\calG.$ 
\end{theorem}

As a consequence of~\autoref{k20}, we get:

\begin{corollary}\label{k23}
    If\/ $\aleph_{0} \leq \mu$ and $\chi \leq 2^{\mu}$ then there exists a complete set of guardrails $\calG$ for $\chi$ over $\mu$ such that $\vert \calG \vert \leq \mu^{\aleph_{0}}.$ 
\end{corollary}

The following definition is a generalization of \cite[Def.~2.11(d)]{Sh00} and it is based on \cite[Def.~2.33]{KST} with some modifications to adapt it to our formalism.

\begin{definition}\label{k24}
    Let $\bft$ be a prepared iteration, $g$ a guardrail for $\pi$, and $\bar I = \la I_k \colon k\in K\ra \in \bfI_{K}^{\fin}$. Consider a sequence of tuples $\tau = \langle (\dot{q}_{\ell}, \alpha_{\ell},  \varp) \colon \ell \in W \rangle.$ We say that \emph{$\tau$ follows $g$} if, for all $\ell\in W_{\bar{I}}$:

    \begin{enumerate}[label=\rm(\alph*)]
        \item $\varp \in [0, 1)_{\QQQ}$,

        \item $g(\alpha_{\ell}) = (\zeta_{\ell},  \varp)$ for some $\zeta_\ell<\theta_{\alpha_{\ell}}$, and 
    
        \item $\dot{q}_{\ell}$ is a $\bbP_{\alpha_{\ell}}^{-}$-name such that $\Vdash_{\alpha_{\ell}} \text{``} \dot{q}_{\ell} \in \dot{Q}_{g(\alpha_\ell)}^{\alpha_{\ell}}$''.\footnote{So, according on whether $\varp>0$ or $\varp=0$, we must have that all $\alpha_\ell$ are in $F^{\bft}$ or all are in $U^{\bft}$.}  
    \end{enumerate}
\end{definition}

We can now define the \emph{generalized iteration with finitely additive measures}. 

\begin{definition}\label{k28}
    We say that $\calC$ is an \emph{(iteration) bedrock} if it is a collection of tuples $\bfc = (\Xi^{\bfc}, \bar{I}^{\bfc}, g^{\bfc})$ such that  $(\Xi^{\bfc}, \bar{I}^{\bfc}) \in \mathcal{Y}_{\ast}$ and $g^{\bfc}$ is a function. According to~\autoref{i2.2}~\ref{i2.2-1}, $(\Xi^{\bfc}, \bar I^\bfc)$ determines sets $K^{\bfc}$ and $W^{\bfc}$.\footnote{As indicated in \autoref{i2.4}, if convenient, there is no problem if the reader restricts to the case $W^\bfc = K^\bfc =\omega$, $\Xi^\bfc$ is free and $\bar I^\bfc$ is an interval partition of $\omega$.}
    
    Given a bedrock $\calC$, we say that $\bft = \langle \bbP_{\beta},  \dot{\bbQ}_{\alpha},  \bbP_{\alpha}^{-}, \vec{Q}_{\alpha}, \theta_{\alpha}, \vec{\Xi}_{\beta} \colon \beta \leq \pi,\ \alpha < \pi \rangle$ is a \emph{$\calC$-iteration} if it satisfies the following conditions:

    \begin{enumerate}[label=\rm(\arabic*)]
        \item\label{k28a} 
        $\tbf$ (excluding $\vec\Xi_\beta$) is a prepared iteration. 

        \item\label{k28b} $g^\bfc$ is a guardrail for $\pi$ for all $\bfc\in\calC$. 

        \item\label{k28c}  For any $\beta \leq \pi,$ $\vec{\Xi}_{\beta}  =  \langle \dot{\Xi}_{\beta}^{\bfc} \colon  \bfc \in \calC \rangle$, $\dot \Xi^\bfc_0 = \Xi^\bfc$ and  $\Vdash_{\bbP_{\beta}}$``$\dot{\Xi}_{\beta}^{\bfc}$ is a  probability fam on $\mathcal{P}(K^{\bfc})$'' for all  $\bfc \in \calC.$

        \item\label{k28d} $\alpha < \beta \leq \pi \Rightarrow \ \Vdash_{\bbP_{\beta}}$`` $\dot{\Xi}_{\alpha}^{\bfc} \subseteq  \dot{\Xi}_{\beta}^{\bfc}$'' for any $\bfc \in \calC.$

        \item\label{k28e} For any $\alpha < \pi$ and $\bfc \in \calC,$ $\Xi_{\alpha}^{\bfc,-} \coloneqq \Xi_{\alpha}^{\bfc}  {\rest}  (\mathcal{P}(K^\bfc) \cap V^{\bbP_{\alpha}^{-}})$ has a $\bbP_{\alpha}^{-}$-name $\dot\Xi_{\alpha}^{\bfc,-}$.

        \item\label{k28ee} $\dot\calY^-_{\alpha}$ is a $\Por^-_\alpha$-name for $\{ (\dot{\Xi}_{\alpha}^{\bfc,-}, \bar{I}^{\bfc}) \colon \bfc \in \calC \}$. 
        If $\alpha \in F^{\bft}$, then 
        $$\Vdash_{\bbP^-_{\alpha}}\text{``$\dot{\bbQ}_{\alpha}$ is $\theta_{\alpha}$-$\dot{\calY}_{\alpha}^-$-linked witnessed by the sequence $\langle \dot{Q}_{\zeta, \varp}^{\alpha} \colon \zeta < \theta_{\alpha} \conj \varp \in (0, 1)_{\QQQ} \rangle$''}$$ 
        and, if $\alpha \in U^{\bft},$ then 
        $$\Vdash_{\bbP_{\alpha}^{-}}\text{``$\dot{\bbQ}_{\alpha}$ is strongly $\theta_{\alpha}$-$\dot{\calY}_{\alpha}^-$-linked witnessed by the sequence $\langle \dot{Q}^\alpha_{\zeta, 0} \colon \zeta < \theta_{\alpha} \rangle$''}.$$

        \item\label{k28f} For any $\bfc \in \calC$ and $\tau = \{ (\dot{q}_{\ell},  \alpha_{\ell},  \varp) \colon \ell \in W^\bfc \}$ following $g^{\bfc}$, it is satisfied that:  
        
        \begin{enumerate}[label=\rm(\arabic{enumi}\alph*)]
            \item\label{k28fa} If the sequence $\langle \alpha_{\ell} \colon \ell \in W^\bfc \rangle$ is constant with value $\alpha$, then:  
            $$ \Vdash_{\bbP_{\alpha}} \text{``} \limit^{\Xi_{\alpha}^{\bfc,-}}_{\ell\in W^\bfc}\dot{q}_{\ell} \Vdash_{\dot{\bbQ}_{\alpha}} \int_{K^\bfc} \frac{\vert \{ \ell \in I^\bfc_{k} \colon \dot{q}_{\ell} \in \dot{G}(\alpha) \} \vert}{\vert I_{k}^\bfc \vert} d  \dot{\Xi}_{\alpha+1}^{\bfc}(k) \geq 1 - \varp \text{''}. $$

            \item\label{k28fb} If the sequence $\bar \alpha \coloneqq \langle  \alpha_{\ell} \colon \ell \in W^\bfc \rangle$ has no maximum and, for all $\alpha<\sup_{\ell\in W^\bfc}\alpha_\ell$, $\Xi^\bfc_0(\{k\in K^\bfc \colon \bar\alpha\rest I_k \text{ is one-to-one and } \forall \ell\in I_k^\bfc\ (\alpha\leq \alpha_\ell)\}) =1$ then, for all $\varp' > 0$, 
            $$ \Vdash_{\bbP_{\pi}} \text{``} \ \dot{\Xi}_{\pi}^{\bfc} \left( \left \{ k\in K^\bfc \colon \frac{\vert \{ \ell \in I_{k}^\bfc \colon \dot{q}_{\ell} \in \dot{G}(\alpha_{\ell}) \} \vert}{\vert I_{k}^\bfc \vert}  \geq (1 - \varp) (1 - \varp') \right \} \right) = 1 \text{''}.$$ 
        \end{enumerate}
    \end{enumerate}
    If, in addition, $\tbf$ is $\kappa$-prepared, we say that $\tbf$ is a \emph{$\kappa$-$\calC$-iteration}.
\end{definition}

Condition \autoref{k28}~\ref{k28fa} requires a justification: why can we take that limit? Since the sequence $\langle \alpha_{\ell} \colon \ell \in W^{\bfc} \rangle$ is constant with value $\alpha$ and $g^\bfc$ is a guardrail, we have that $\langle \dot q_{\ell} \colon \ell \in W^{\bfc} \rangle$ is a sequence of $\Por^-_\alpha$ names of members of $\dot Q^\alpha_{g^\bfc(\alpha)}$ because $\tau$ follows $g^\bfc.$  As a consequence, in $V^{\bbP_{\alpha}^{-}}$ we have that $\langle \dot{q}_{\ell} \colon \ell \in W^{\bfc} \rangle \in {}^{W^{\bfc}}\dot{Q}_{\zeta, \varp}^{\alpha}$ and, as $\Xi_{\alpha}^{\bfc,-}$ is a fam in $V^{\bbP_{\alpha}^{-}}$ which belongs to $\calY^-_\alpha$, and $\dot{Q}_{g^\bfc(\alpha)}^{\alpha}$ is in particular $(\Xi_{\alpha}^{\bfc,-}, \bar{I}^\bfc, \varp)$-linked in that model, we have that $\limit^{\Xi_{\alpha}^{\bfc,-}}_{\ell\in W^\bfc}\dot{q}_{\ell}$ can be calculated.

The additional condition in~\ref{k28fb} about the supremum is stated for the sake of generality. In our applications where $K^\bfc = W^\bfc =\omega$, and $\Xi^\bfc$ is free, $\la\alpha_\ell\colon \ell<\omega\ra$ is increasing, which implies that the sequence has no maximum and the condition of the supremum is satisfied.

In order to simplify the notation, we introduce: 

\begin{notation}\label{k34}    
        Let $\bft$ be a prepared iteration. For any $ \vec{\Xi} = \langle \Xi_{\beta}^{i} \colon \beta \leq \pi \conj i \in I \rangle,$ we denote $ \bft \sqcup \vec{\Xi} \coloneqq  \langle \bbP_{\beta}, \dot{\bbQ}_{\alpha},  \bbP_{\alpha}^{-}, \vec{Q}_{\alpha},  \theta_{\alpha}, \vec{\Xi}_{\beta} \colon \beta \leq \pi,\ \alpha < \pi \rangle$ where $\vec{\Xi}_{\beta} \coloneqq \langle \dot{\Xi}_{\beta}^{i} \colon i \in I \rangle$ for any $\beta \leq \pi$.  
\end{notation} 

Until the end of this section, we fix an iteration bedrock $\calC$.

In order to construct $\calC$-iterations, we need to properly extend the parameters of the iteration at the successor and limit steps so that the iteration still be a $\calC$-iteration. In general, this is going to come down to simply getting the sequences of finitely additive measures properly extended, so that they continue to satisfy the conditions of \autoref{k28}. In the successor steps, there is no difficulty, since we will show that we can amalgamate finitely additive measures throughout the iteration (see \autoref{k38}). However, the limit step in general is not only the most difficult step in the Shelah's construction (see \cite[Lem.~2.13 and~2.14]{Sh00}), but for reasons we will see later, it is the most delicate point in our presentation of a general theory of iterated forcing using fams. We start with the successor step.

\begin{remark}\label{rem:ufi3}
    In the context of ultrafilter limits for intervals as in \autoref{rem:ufi} and~\ref{rem:ufi2}, we can define a \emph{$\calC$-uf-iteration} as in~\autoref{k28}, but demading in addition that all $\Xi^\bfc$ and $\dot\Xi^{\bfc}_\beta$ are ``ultrafilters''. All the results presented from now on (under natural modifications when required) are valid for $\calC$-uf-iterations as well. In this sense, our iteration theory with fams contains the iteration theory with ultrafilters, even the one for ultrafilter limits for intervals.
\end{remark}

\subsection{Extending at successor steps}\

Assuming we have constructed a $\calC$-iteration up to $\Por_\gamma$, we need to pick $\Por^-_\gamma$ and a desired $\Qnm_\gamma$. Much later in \autoref{a70} and~\ref{a70g} we describes ways to construct $\Por^-_\gamma$ such that $\dot \Xi^\bfc_\gamma\frestr V^{\Por^-_\gamma}$ has a $\Por^-_{\gamma}$-name for all $\bfc\in\calC$. Notice that $\Por^-_\gamma = \Por_0 = \{\la\ \ra\}$ satisfies this, which can be used to develop \emph{finite support products with fams} (and ultrafilters).

After picking $\Por_\gamma$ and $\Qnm_\gamma$, we need to find $\dot\Xi^\bfc_{\gamma+1}$ for all $\bfc\in\calC$. For this purpose, we use the following amalgamation result, which is a consequence of \autoref{m48} and extends~\cite[Claim~1.6]{Sh00}. 

\begin{lemma}\label{k38}
    Let $M \subseteq N$ be transitive models of $\ZFC$, $\bbP \in M$ a forcing notion, $G$ a $\bbP$-generic set over $N$, $K \in M$ and $\Xi_{0} \in M$ and $\Xi_{1} \in N$  finite finitely additive measures on $\calP(K) \cap M$ and $\calP(K) \cap N$, respectively, such that $\Xi_{1}$ extends $\Xi_{0}.$ In $M,$ let $\dot{\Xi}_{0}^{\ast}$ be a $\bbP$-name of a finitely additive measure on $\calP(K) \cap M[G]$ extending $\Xi_{0}.$ Then, in $N,$ there is a $\bbP$-name $\dot{\Xi}_{1}^{\ast}$ such that $\bbP$ forces that $\dot{\Xi}_{1}^{\ast}$ is a finitely additive measure on $\calP(K) \cap N[G]$ extending both $\Xi_{1}$ and $\dot{\Xi}_{0}^{\ast}$ (see \autoref{fig:k38}). Even more, if $\Xi_1$ represents an ultrafilter, and $\dot\Xi^*_0$ is forced to represent an ultrafilter, then $\dot\Xi^*_1$ can be found forced to be an ultrafilter, as well. 
    \begin{figure}[ht]
        \begin{tikzcd}
            N \arrow{r}{\bbP}  & N[G] \\
            M \arrow{r}{\bbP} \arrow[swap]{u} & M[G] \arrow{u} 
        \end{tikzcd}
        \hspace{3cm}
        \begin{tikzcd}
            \Xi_{1} \arrow{r}  & \Xi_{1}^{\ast}  \\
            \Xi_{0} \arrow{r} \arrow[swap]{u} & \Xi_{0}^{\ast} \arrow{u}
        \end{tikzcd}   
        \caption{In the diagram on the left, the inclusion relationships between the models are presented. In the diagram on the right, the extension relation between the finitely additive measures is presented.}\label{fig:k38}
    \end{figure}
\end{lemma}

\begin{PROOF}[\textbf{Proof}]{\ref{k38}}
    In view of \autoref{m48}~\ref{m48b} and \autoref{m52}, it is enough to prove that, if $b\in\mathcal{P}(K)\cap N$ and $\dot a\in M$ is a $\bbP$-name of a subset of $K$ then, in $N$, $\bbP$ forces $\dot a \subseteq b \Rightarrow \dot\Xi_{0}^{\ast}(\dot a)\leq \Xi_{1}(b)$ (it is already clear that, in $N$, $\bbP$ forces $\dot\Xi^{\ast}_{0}(K)=\Xi_{1}(K)=\Xi_{0}(K)$). 

    Assume $p\in\bbP$ and $p\Vdash^{N}_\bbP\text{``}\dot a\subseteq b$''. Define $a_0\coloneqq \{k\in K \colon  \exists q\leq p \ ( q\Vdash^{M}_\bbP k\in\dot a) \}$, which is in $M$. Then $a_0\subseteq b$ and $p\Vdash^{M}_\bbP \text{``} \dot a\subseteq a_0 \text{''}$. Therefore, $\Xi_0(a_0)=\Xi_1(a_0)\leq\Xi_1(b)$ and $p\Vdash^{M}_\bbP \text{``} \dot\Xi^{\ast}_{0}(\dot a)\leq \dot\Xi^{\ast}_{0}(a_0)=\Xi_0(a_0)$'', so $p\Vdash^{N}_\bbP \text{``} \dot\Xi^{\ast}_{0}(\dot a)\leq \Xi_{1}(b)$''. 
\end{PROOF}

We now have all the necessary tools to generalize the \emph{extension theorem} at successor steps.

\begin{theorem}\label{k39}
   Let $\bft$ be a prepared iteration of length $\pi = \gamma +1$. Let $\vec{\Xi} = \langle \dot \Xi_{\beta}^{\bfc} \colon \beta \leq \gamma,\ \bfc \in \calC \rangle$ such that $\bft \rest  \gamma \sqcup \vec \Xi$ is a $\calC$-iteration and assume that: 
   \begin{enumerate}[label = \normalfont (\roman*)]
      \item for any $\bfc \in \calC$, $\dot\Xi^\bfc_\gamma\frestr \pts(K^\bfc)\cap V^{\Por^-_\gamma}$ has a $\Por^-_\gamma$-name $\dot{\Xi}_{\gamma}^{\bfc,-}$;
      
      \item $\dot\calY^-_{\gamma}$ is a $\Por^-_\gamma$-name for $\{ (\dot{\Xi}_{\gamma}^{\bfc,-}, \bar{I}^{\bfc}) \colon \bfc \in \calC \}$; if $\gamma \in F^{\bft}$ then 
    $$\Vdash_{\bbP^-_{\gamma}}\text{`` $\dot{\bbQ}_{\gamma}$ is $\theta_{\gamma}$-$\dot{\calY}_{\gamma}^-$-linked witnessed by the sequence $\langle \dot{Q}_{\zeta, \varp}^{\gamma} \colon \zeta < \theta_{\gamma} \conj \varp \in (0, 1)_{\QQQ} \rangle$''};$$ 
        and, if $\gamma \in U^{\bft}$ then 
        $$\Vdash_{\bbP_{\gamma}^{-}}\text{`` $\dot{\bbQ}_{\gamma}$ is strongly $\theta_{\gamma}$-$\dot{\calY}_{\gamma}^-$-linked witnessed by the sequence $\langle \dot{Q}^\gamma_{\zeta, 0} \colon \zeta < \theta_{\gamma} \rangle$''}.$$
   \end{enumerate}
   
   Then, there exists a sequence $\langle \dot{\Xi}_{\pi}^{\bfc} \colon \bfc \in \calC \rangle$ such that 
   $\bft \sqcup \langle \dot{\Xi}_{\beta}^{\bfc} \colon \beta \leq \pi, \bfc \in \calC \rangle$ is a $\calC$-iteration.
   Moreover, for $\bfc\in\calC$, if $\dot\Xi^\bfc_\gamma$ is forced to represent an ultrafilter, then $\dot\Xi^\bfc_\pi$ can be found forced to represent an ultrafilter in one of the following cases:
   \begin{enumerate}[label = \rm (\roman*)]
       \item\label{k39-i} $\gamma\in U^\tbf$, or
       \item\label{k39-ii} $\gamma\in F^\tbf$ and $g^\bfc(\gamma)\notin \theta_\gamma\times(0,1)_\QQQ$, or
       \item\label{k39-iii} $\gamma\in F^\tbf$, $g^\bfc(\gamma)\in \theta_\gamma\times(0,1)_\QQQ$, and the fam extending $\Xi^{\bfc,-}_\gamma$ in the $\Por^-_\gamma\ast\Qnm_\gamma$-extension witnessing that $\dot Q^\xi_{g^\bfc(\xi)}$ is $(\dot \Xi^{\bfc,-}_\gamma,\bar I^\bfc, g^\bfc_1(\gamma))$-linked represents an ultrafilter.
   \end{enumerate}
\end{theorem}

\begin{PROOF}[\textbf{Proof}]{\ref{k39}}
    Let $\bfc\in\calC$. Set $(\zeta^\bfc,\varp^\bfc)\coloneqq g^\bfc(\gamma)$ if both $g_{0}^{\bfc}(\gamma)<\theta_\gamma$, and $g_{1}^{\bfc}(\gamma) = 0$ iff $\gamma\in U^{\bft}$, which we call the \emph{legal case}. 
    Work in $V^-_\gamma\coloneqq V^{\Por^-_\gamma}$. In the legal case, 
    since $Q_{\zeta^\bfc, \varp^\bfc}^{\gamma}$ is $(\Xi_{\gamma}^{\bfc,-}, \bar{I}^\bfc, \varp^\bfc)$-linked, there exists a $\bbQ_{\gamma}$-name $\dot{\Xi}_{\gamma}^{\bfc,\ast}$ of a fam on $\pts(K^\bfc)$ extending $\Xi_{\gamma}^{\bfc,-}$ such that, for any $\bar{q} \in {}^{W^\bfc}\dot{Q}_{\zeta^\bfc, \varp^\bfc}^{\gamma},$  
    \begin{equation}\label{k39-1}
        \limit^{\Xi_{\gamma}^{\bfc,-}}(\bar q) \forces_{\dot{\bbQ}_{\gamma}} \text{``} \int_{W^\bfc} \frac{\vert \{ \ell \in I_{k}^\bfc \colon \dot{q}_{\ell} \in \dot{G}_{\dot{\bbQ}_{\gamma}} \} \vert }{\vert I_{k}^\bfc \vert} d  \dot{\Xi}_{\gamma}^{\bfc,\ast} (k) \geq 1 - \varp^\bfc \text{''.}
    \end{equation}    
    Moreover, since $\bbP_{\gamma}^{-} \subsetdot   \bbP_{\gamma}$ and by integral absoluteness (see \autoref{t102}), \eqref{k39-1} holds in $V_\gamma\coloneqq V^{\Por_\gamma}$. For the no-legal case, let $\dot \Xi^{\bfc,*}_\gamma$ be a $\Qor_\gamma$-name of a fam on $\pts(K^\bfc)$ extending $\Xi^{\bfc,-}_\gamma$, and representing an ultrafilter if $\Xi^{\bfc,-}_\gamma$ does.
    
    Now, by \autoref{k38}, we can amalgamate as in the following diagrams:
                \[\begin{tikzcd}
                    V_{\gamma} \arrow{r}{\dot{\bbQ}_{\gamma}} & V_{\gamma}[G(\gamma)]\\
                    V_{\gamma}^{-} \arrow{r}{\dot{\bbQ}_{\gamma}} \arrow[swap]{u} & V_{\gamma}^{-}[G(\gamma)] \arrow{u}
                \end{tikzcd}
                \hspace{3cm}
                \begin{tikzcd}
                \Xi_{\gamma}^{\bfc} \arrow{r} & \Xi_{\pi}^{\bfc}\\
                \Xi_{\gamma}^{\bfc,-} \arrow{r} \arrow[swap]{u} & \Xi_{\gamma}^{\bfc,\ast} \arrow{u}
                 \end{tikzcd}
                \]
    That is, in $V_\pi=V^{\Por_\pi},$ there exists a fam $\dot{\Xi}_{\pi}^{\bfc}$ on $\calP(K^\bfc)\cap V_\pi$ extending both $\Xi_{\gamma}^{\bfc}$ and $\Xi_{\gamma}^{\bfc,\ast}.$ In the case when $\Vdash^V_\alpha$``$\dot\Xi^\bfc_\gamma$ represents an ultrafilter", in $V^-_\gamma$, $\Xi^{\bfc,-}_\gamma$ also represents an ultrafilter and $\dot \Xi^{\bfc,\ast}_\gamma$ can be found forced to represent an ultrafilter in any of the cases~\ref{k39-i}--\ref{k39-iii} (by \autoref{i25} for~\ref{k39-i}). Therefore, in $V_\pi$, $\Xi^\bfc_\pi$ can be found to represent an ultrafilter, as well.
    
    As a consequence, in $V$, we get a sequence of $\bbP_{\pi}$-names $\langle \dot{\Xi}_{\pi}^{\bfc} \colon \bfc \in \calC \rangle$ of these fams. Now, we must show that $\bft \sqcup \langle \dot{\Xi}_{\beta}^{\bfc} \colon \beta \leq \pi,\ \bfc \in \calC \rangle$ is a $\calC$-iteration. Notice that requirements~\ref{k28a}--\ref{k28ee} of \autoref{k28} are immediate, so it remains to check the condition~\ref{k28f}. 
    
    Suppose that $\tau = \{ (\dot{q}_{\ell}, \alpha_{\ell}, \varp) \colon \ell \in W^\bfc\}$ follows $g^\bfc$. 

    \noindent\ref{k28fa}: Assume the sequence $\langle \alpha_{\ell} \colon \ell \in W^\bfc \rangle$ is constant with value, say, $\alpha.$ Now, we proceed by cases. When $\alpha < \gamma$ then the result follows because $\bft \, {\rest}   \gamma \sqcup \vec \Xi$ is a $\calC$-iteration; when $\alpha = \gamma,$ by~\eqref{k39-1} we get 
    $$ \forces_{\gamma} \text{``} \limit^{\dot{\Xi}_{\gamma}^{\bfc,-}}_{\ell \in W^\bfc} \dot{q}_{\ell} \forces_{\dot{\bbQ}_{\gamma}} \text{``} \int_{W^\bfc} \frac{\vert \{ \ell \in I_{k}^\bfc \colon \dot{q}_{\ell} \in \dot{G}_{\dot{\bbQ}_{\gamma}} \} \vert }{\vert I_{k}^\bfc \vert} d  \dot{\Xi}_{\pi}^{\bfc} (k) \geq 1 - \varp \text{''} \text{''}. $$
    
    \noindent\ref{k28fb}: Let $\varp' > 0.$ If the sequence $\langle \alpha_{\ell} \colon \ell \in W^\bfc \rangle$ has no maximum then $\alpha_{\ell} < \gamma$  for any $\ell \in W^{\bfc}$, therefore, the result is clear because $\bft  {\rest}  \gamma\sqcup \vec \Xi$ is a $\calC$-iteration.\footnote{In this case, we are not using the new sequence of fams, that is, the iteration up to $\gamma$ had already taken care of this condition.}
            
    Finally, we conclude that $\bft \sqcup \langle \dot{\Xi}_{\beta}^{\bfc} \colon \beta \leq \pi, \, \bfc \in \calC \rangle$ is a $\calC$-iteration.
\end{PROOF}

\begin{remark}
    The notion of ``$(\Xi, \bar{I}, \varp_{0})$-linkedness" is not upward absolute for transitive models of $\ZFC$, and for this reason, it was necessary to resort to $\bbP_{\gamma}^{-}$ to be able to extend the finitely additive measures. 
    
    When $\dot\Xi^\bfc_\gamma$ is forced to be free, $\dot\Xi^\bfc_{\gamma+1}$ is also forced free. However, in the case that it has the uap but is not free, it seems that the frame of $\dot\Xi^\bfc_\gamma$ (see \autoref{t94}) must have a $\Por^-_\gamma$-name to guarantee that $\dot\Xi^\bfc_{\gamma+1}$ is forced with the uap. This is reflected in \autoref{k38} as well: when $\Xi_1$ has the uap but is not free, and its frame is in $M$, if $\dot\Xi^*_0$ is forced with the uap then $\dot\Xi^*_1$ can be found (forced) with the uap.
\end{remark}

\subsection{Extending at limit steps}\label{6.3}\

Now, we deal with the problem of extending at limit steps. When analyzing the limit step extension proof for both random forcing and $\tilde{\bbE}$ (see \cite[Lem. 2.14]{Sh00} and \cite[Lem.~2.39]{KST}, respectively), we notice that there are two fundamental points: having available the Lebesgue measure and a Boolean structure. The problem is that, when considering more general forcing notions, we do not necessarily have a Boolean structure or some measure, so we need other alternatives. It is here where the new concept of intersection number for forcing notions plays its stellar role in our theory because, through \autoref{i11.01} and \autoref{s14}, we will have enough structure available in the iterands of our iteration to be able to extend it at limit steps. All this makes sense in the results of this section, the most important technical results of this paper, which give the final shape to our theory. 

\begin{mainlemma}\label{k40}
    Let $\bbP_{\pi} = \langle \bbP_{\alpha}, \dot{\bbQ}_{\alpha} \colon \alpha < \pi \rangle$ be a finite support iteration with $\pi$ limit, let $K$ and $W$ be sets and let $\bar I = \langle I_{k} \colon k \in K \rangle$ be a partition of $W$ into finite sets. Assume that:
    \begin{enumerate}[label=\rm(\arabic*)]
        \item\label{k40a} For any $\alpha < \pi,$ $\dot{\Xi}_{\alpha}$ is a $\bbP_{\alpha}$-name of a probability fam on $\pts(K)$.
        
        \item\label{k40aa} Either $K$ is finite or, for any $n<\omega$, $\Xi_0(\{k\in K\colon |I_k|=n\}) = 0$.   

        \item\label{k40b} $\alpha < \beta$ implies $\Vdash_{\beta}$``\/$\dot{\Xi}_{\alpha} \subseteq \dot{\Xi}_{\beta}$''.

        \item\label{k40c} For any $\alpha < \pi,$ $\dot{Q}_{\alpha}'$ is a $\bbP_{\alpha}$-name of a  subset of $\dot{\bbQ}_{\alpha}$, $\varp_{\alpha}' \in [0, 1)_{\QQQ}$, and $\Por_\alpha$ forces $\Int^{\dot{\bbQ}_{\alpha}}(\dot{Q}_{\alpha}') \geq 1 - \varp_{\alpha}.$
    \end{enumerate}
    
    Then, there is a $\bbP_{\pi}$-name  $\dot{\Xi}_{\pi}$ of a fam on $\pts(K)$ extending $\bigcup_{\alpha < \pi} \dot{\Xi}_{\alpha}$ such that, whenever
    \begin{enumerate}[label = \rm (\roman*)]
        \item\label{k40i} $\bar \beta\coloneqq\langle \beta_{\ell} \colon  \ell \in W \rangle$ is a sequence of ordinals ${<}\pi$ with supremum $\pi$,
        
        \item\label{k40ii} $\Xi_0(\{k\in K \colon \bar\beta\rest I_k \text{ is one-to-one and } \forall \ell\in I_k\ (\beta_\ell\geq\alpha)\}) =1$ for all $\alpha<\pi$,
        
        \item\label{k40iii} $\dot{r}_{\ell}$  is a  $\bbP_{\beta_{\ell}}$-name such that  $ \Vdash_{\beta_{\ell}}$``\/$\dot{r}_{\ell} \in \dot{Q}_{\beta_{\ell}}'$'' for any $\ell \in W$, and
        
        \item\label{k40iv} $\varp > 0$ and $\langle \varp_{\beta{\ell}}' \colon \ell \in W \rangle$ is constant with value $\varp_{0},$
    \end{enumerate}
    we have $\displaystyle\Vdash_{\bbP_{\pi}} \text{``\/} \dot{\Xi}_{\pi} \left( \left \{ k \in K \colon \frac{\vert \{ \ell \in I_{k} \colon \dot{r}_{\ell} \in \dot{G}(\beta_{\ell}) \} \vert}{\vert I_{k} \vert}  \geq (1 - \varp_{0}) (1 - \varp) \right \} \right) = 1 \text{''}.$ 
    
    Even more, if $\Por_\alpha$ forces that $\dot\Xi_\alpha$ represents an ultrafilter for all $\alpha<\pi$, then $\dot\Xi_\pi$ can be found forced to represent an ultrafilter. The same applies for the uap. 
\end{mainlemma}

\begin{PROOF}[\textbf{Proof}]{\ref{k40}}
    The case when $K$ is finite is trivial (all $\dot\Xi_\alpha$ are the same and \ref{k40i}--\ref{k40iv} are never satisfied), so we assume that $K$ is infinite. Then, by~\ref{k40aa} and \ref{k40b}, each $\dot\Xi_\alpha$ is forced to be free.
    
    With the intention of applying \autoref{m69}, we define $\cI$ as the set of pairs $(\bar{c}, \varp)$ such that $\bar c = \langle (\beta_{\ell}, \dot{r}_{\ell}) \colon \ell \in W \rangle$ is as in~\ref{k40i}--\ref{k40iv}, $\langle \varp_{\beta_{\ell}} \colon \ell \in W \rangle$ is constant with value $\varp^{\bar c}_0$ and $\varp>0$. 
  
    For each $(\bar{c}, \varp) \in \cI,$ let $\dot A_{\bar{c}, \varp}$ be a $\Por_\pi$-name of
    $$ \left \{ k \in K \colon \frac{\vert \{ \ell \in I_{k} \colon \dot{r}_{\ell} \in \dot{G}(\beta_{\ell}) \} \vert}{\vert I_{k} \vert} \geq (1-\varp^{\bar c}_0)(1 - \varp) \right \}.$$ 

    Pick a $\bbP_{\pi}$-name $\dot{\Xi}$ for $\bigcup_{\alpha < \pi} \dot{\Xi}_{\alpha}$, which is forced to be a fam  with $\dom(\dot{\Xi}) = \bigcup_{\alpha < \pi} \calP(K) \cap V_{\alpha}.$ 
    
    Let $C \subseteq \cI$ finite, $p \in \bbP$ and $\dot{a} \in V^{\bbP_{\pi}}$ such that $ p \forces_{\bbP_{\pi}} \text{``} \dot{a} \in \dom(\dot{\Xi}) \text{ and } \dot{\Xi}(\dot{a}) > 0\text{''}.$ Our aim will be to find some $q \in \bbP_{\pi}$ such that $q \leq p$ and $q \forces_{\bbP_{\pi}}$``$\dot{a} \cap \bigcap_{(\bar{c}, \varp) \in C} \dot A_{\bar{c}, \varp} \neq \emptyset$'', because this will allow us to apply \autoref{m69} to extend $\dot{\Xi}$ as required.

    Enumerate $C = \{(\bar c^i,\varp^i) \colon i<i^*\}$ with $i^{\ast} < \omega$, $\bar{c}^{i} = \langle (\beta_{\ell}^{i}, \dot{r}_{\ell}^{i}) \colon \ell < \omega \rangle$ and $\varp^i_0 \coloneqq \varp^{\bar c^i}_0$ for each $i<i^*$. Since $\varp' < \varp $ implies that $\dot A_{\bar{c}, \varp'} \subseteq \dot A_{\bar{c}, \varp}$ is forced, by taking $\varp  \coloneqq \min\{ \varp^{i} \colon i < i^{\ast} \}$ without loss of generality we can assume that, for any $i < i^{\ast}, \varp^{i} = \varp.$ Denote $\dot A_{i} \coloneqq \dot A_{\bar{c}^{i},  \varp}.$
    
    By strengthening $p,$ we can find some $\gamma < \pi$ such that $\dot{a}$ is a $\bbP_{\gamma}$-name and $p \in \bbP_{\gamma}.$ 
    
    Recall that $\langle \varp_{\beta_{\ell}^{i}} \colon \ell \in W  \rangle$ is constant with value $\varp_{0}^{i}.$ Since for any $i < i^{\ast}, \varp(1- \varp_{0}^{i}) > 0,$ there exists some $\varp' > 0$ such that, for all $i < i^{\ast},$ $ \varp' < \varp (1 - \varp_{0}^{i})$. 
    By the choice of $\varp',$ we have that $(1-\varp_{0}^{i} - \varp')(\varp_{0}^{i} + \varp') > 0$ for any $i < i^{\ast}.$ Then, we can find some $0<N_0<\omega$ such that
    \begin{equation}\label{k40-5}
        \text{$\frac{1}{N_0} \frac{(1 - \varp_{0}^{i} - \varp' ) ( \varp_{0}^{i} + \varp')}{[\varp  (1 - \varp_{0}^{i}) - \varp']^{2}} < \frac{1}{i^{\ast} +1 }$}.
    \end{equation}
    Now, $p$ forces in $\Por_\gamma$ that $\dot\Xi_\gamma(\dot a)>0$, $\dot\Xi_\gamma(\{k\in K \colon |I_k|>N_0\}) =1$ (by~\ref{k40aa}) and, by~\ref{k40ii}, $\dot\Xi_\gamma(\{k\in K\colon \la\beta^i_\ell \colon \ell\in I_k\ra \text{ is one-to-one and }\forall \ell\in I_k\ (\beta^i_\ell \geq \gamma)\}) = 1$ for all $i<i^*$. Thus, there are some $k\in K$ and $q_0\leq p$ in $\Por_\gamma$ such that $|I_k|>N_0$, $\{\beta^i_\ell\colon \ell\in I_k\}\subseteq \pi\menos \gamma$ is one-to-one for all $i<i^*$, and $q_0\Vdash_\gamma \text{``} k\in \dot a$''. 
    
    It is clear that $\{ \beta_{\ell}^{i} \colon i < i^{\ast},\  \ell \in I_{k} \}$ is finite, so there exists some $m^{\ast} < \omega$ such that the sequence  $\langle \beta_{m} \colon m < m^{\ast} \rangle$ is an increasing enumeration of it. Notice that $\varp_{0}^{i} = \varp_{0}^{i'} = \varp'_{\beta_m}$ whenever  $\beta_{\ell}^{i} = \beta_{\ell'}^{i'} = \beta_m.$  
    Now, we split the rest of this proof into three parts. First, we build a suitable probability tree and then we find a suitable event with high probability.    

    \textbf{Part 1. The tree construction.}

    We will build a tree $T$ of height $m^{\ast} +1,$ a function $\bfp \colon T \to \bbP_{\pi}$ such that $\bfp(\rho) \coloneqq p_{\rho}\in \Por_{\beta_{|\rho|}}$ ($\beta_{m^*}\coloneqq\pi$) for each $\rho \in T,$ and a probability  space on $\suc(\rho)$ for each $\rho \in T,$ by induction on the level $m \leq m^{\ast}$ as follows.  

    In the base step, we define $\Lev_{0}(T) \coloneqq \{ \langle \, \rangle \}$ and $p_{\langle\,  \rangle} \coloneqq q_{0} \in \Por_\gamma\subseteq \Por_{\beta_0}.$ 

    For the successor step, suppose that we have built $T$ up to level $m$.  In order to define $\Lev_{m+1}(T),$ let $J_{m} \coloneqq \{ (i, \ell)\in i^\ast\times I_k \colon  \beta_{\ell}^{i} = \beta_{m} \}$ and let $\rho \in \Lev_{m}(T).$  Working in $V^{\bbP_{\beta_{m}}},$ we have that, by~\ref{k40c}, $\Int^{\bbQ_{\beta_{m}}}(Q_{\beta_{m}}') \geq 1 - \varp_{\beta_{m}}',$ and by virtue of \autoref{s14}, for 
    $$\Sigma_{\rho} \coloneqq \{ \sigma \in {}^{J_{m}} 2 \colon \exists q \in \bbQ_{\beta_{m}} \forall (i, \ell) \in J_{m}\ [(\sigma(i, \ell) = 0 \Rightarrow q \leq \dot{r}_{\ell}^{i}) \conj  (\sigma(i, \ell) = 1 \Rightarrow q  \perp \dot{r}_{\ell}^{i}) ]\},$$
    
    there exists some function $f_{\rho} \colon \Sigma_{\rho} \to [0, 1]_{\QQQ}$ such that $\sum_{\sigma \in \Sigma_{\rho}} f(\sigma) = 1$ and, for any $(i, \ell) \in J_{m},$ we have that: 
    \begin{equation}\label{k40-9}
        \text{$ \sum \{ f_{\rho}(\sigma) \colon \sigma \in \Sigma_{\rho}, \conj \sigma(i, \ell) = 0 \} \geq 1 - \varp_{\beta_{m}}' - \varp' $}
    \end{equation}

    For each $\sigma \in \Sigma_{\rho},$ choose a witness $q_{\sigma} \in \bbQ_{\beta_{m}}$ for ``$\sigma \in \Sigma_{\rho}$''.
    
    Working the ground model again, there is some $q_{\rho}' \leq p_{\rho}$ in $\Por_{\beta_m}$ deciding $\Sigma_{\rho} = \dot{\Sigma}_{\rho}$ and $f_{\rho} = \dot{f}_{\rho}.$ We then define $\suc(\rho) \coloneqq \{ \concat{\rho}{\langle \sigma \rangle} \colon \sigma \in \Sigma_{\rho} \}$ and $p_{\concat{\rho}{\langle \sigma \rangle}} \coloneqq q_{\rho}' \cup \{ (\beta_{m}, \dot{q}_{\sigma}) \},$ where $\dot{q}_{\sigma}$ is a $\bbP_{\beta_{m}}$-name for $q_{\sigma}$ (decided by $q'_{\rho}$).     
    Finally, 
    $(\suc(\rho), \calP(\suc(\rho)),  \p_{\rho})$ is a probability space where, for $\eta = \concat{\rho}{\langle \sigma \rangle},$ $\p_{\rho}(\eta) \coloneqq f_{\rho}(\sigma).$ 
    By the construction, it is clear that if $\eta = \concat{\rho}{\langle \sigma \rangle},$ then  $p_{\eta} \forces_{\beta_{m+1}}$``$\dot{r}_{\ell}^{i} \in \dot{G}(\beta_{m}) \Leftrightarrow \sigma(i, \ell) = 0$'' for any $(i, \ell) \in J_{m}.$ 

    \textbf{Part 2. Comparing a random variable with the binomial distribution.}

    To fix some notation, for any $(i, \ell) \in i^{\ast} \times I_{k},$ define $m_{i, \ell}$ as the unique $m < m^{\ast}$ such that $\beta_{\ell}^{i} = \beta_{m}.$ 
    For any $i < i^{\ast}$ consider the order $\lhd^i$ on $I_k$ such that $\ell \lhd^i \ell'$ iff $\beta^i_\ell < \beta^i_{\ell'}$, which is equivalent to $m_{i,\ell}<m_{i,\ell'}$. For convenience, we extend $\lhd^i$ to $I_k\cup\{n^*_i\}$ by putting a point $n^*_i$ on top of $I_k$. For $\ell\in I_k$, $\ell+_i 1$ denotes the immediate $\lhd^i$-successor of $\ell$.
    
    Let $T^{i}$ be the complete binary tree of height $\vert I_{k} \vert +1.$  Without loss of generality, to ease the notation, we order the levels of $T^{i}$ by $\lhd^i$, that is, the root is at level $\min_{\lhd^i}(I_{k})$ and $\max(T^{i}) = {}^{I_k}2$ is at level $n^*_i$.  The key to concluding the proof is to define a probability space structure on $T^{i}$ and properly transfer information from $T$ into $T^{i}.$ For this, we are going to define a function $\Phi^{i} \colon T \to T^{i}$ such that, for any $\rho \in T,$  
    $\Phi^{i}(\rho) \coloneqq \langle \rho(m_{i, \ell})(i, \ell)) \colon  \ell \in I_{k} \conj m_{i, \ell} < \mathrm{ht}_{T} (\rho) \rangle,$
    that is, $\Phi^{i}(\rho)(\ell) = \rho(m_{i, \ell})(i, \ell)$ whenever $\ell \in I_{k}$ and $m_{i, \ell} < \mathrm{ht}_{T}(\rho).$ Notice that,

    \begin{enumerate}[label=\rm(P\arabic*)]
        \item\label{k40-1} $\Phi^{i}(\rho)$ has domain $\{ \ell \in I_{k} \colon \ell \lhd^i \ell_{\rho}^{i} \},$ where $\ell_{\rho}^{i} \coloneqq \min (\{ \ell \in I_{k} \colon m_{i, \ell} \geq \mathrm{ht}_{T}(\rho) \} \cup \{ n^{\ast}_i \} \}$. As a consequence, $\Phi^{i}[\Lev_{m_{i, \ell}}(T)]  \subseteq \Lev_{\ell}(T^{i})$ for any $\ell \in I_{k}$ and also $\Phi^{i}[\Lev_{m^{\ast}}(T)] \subseteq \max(T^{i}).$  
        
        \item\label{k40-2} If $\rho \subseteq \eta$ in $T,$ then $\Phi^{i}(\rho) \subseteq \Phi^{i}(\eta).$

        \item\label{k40-3-1}  If $\eta \in T$ and $m_{i, \ell} \leq \mathrm{ht}_{T}(\rho),$ then $\Phi^{i}(\rho {\rest} m_{i, \ell}) = \Phi^{i}(\rho) {\rest} \ell.$
    \end{enumerate}

    For any $i < i^{\ast}$ and $\ell \in I_{k},$ define the random variable $X_{\ell}^{i}$ on $\Lev_{m_{i, \ell} +1}(T)$ such that, for any $\eta \in \Lev_{m_{i, \ell}+1}(T),$  $X_{\ell}^{i}(\eta) \coloneqq 1 - \eta(m_{i, \ell})(i, \ell) \in \{ 0, 1\}.$
    
    Now, let us deal with the probability space structure on $T^{i}$. First, we are going to define a probability space on its levels: for any $\ell \in I_{k} \cup \{n^{\ast}_i\}$ and $s \in \Lev_{\ell}(T^{i}),$ define $$ \p_{\ell }^{i}(s) \coloneqq  \p_{\Lev_{m^{\ast}}(T)}[ \forall \ell' \lhd^i \ell \ (1 - X_{\ell'}^ {i} = s(\ell'))]. $$
    
    Notice that,  
    \begin{equation*}
        \begin{split}
            \sum_{s \in \Lev_{\ell}(T^{i})} \p_{\ell}^{i}(s) & = \sum_{s \in \Lev_{\ell}(T^{i})} \p_{\Lev_{m^{\ast}}(T)}[\forall \ell' \lhd^i \ell\  (1 - X_{\ell'}^{i} = s(\ell'))]\\
            & = \sum_{s \in \Lev_{\ell}(T^{i})} \p_{\Lev_{m^{\ast}}(T)}[\Phi^{i}(\rho)  {\rest}  \ell = s]\\
            & = \p_{\Lev_{m^{\ast}}(T)}[\exists s \in \Lev_{\ell}(T^{i})\ (\Phi^{i}(\rho)  {\rest} \ell = s)] = 1,
        \end{split}
    \end{equation*}
    where the last equality is given by~\ref{k40-1}. Thus, we have that $\Lev_{\ell}(T^{i})$ is a probability space with probability function $ \p_{\ell}^{i}$. 
    In this case we have that the levels of $T^{i}$ induces a probability space on $\suc(s)$: denote $p_{i} \coloneqq 1 - \varp_{0}^{i} - \varp'$ and, for any $\ell \in I_{k}$ and $s \in \Lev_{\ell}(T^{i}),$ we set:
    \[\p_{s}^{i}(\concat{s}{\langle d \rangle}) \coloneqq  
        \left\{ \begin{array}{ll}
             \frac{ \p_{\Lev_{\ell +_i 1}(T^{i})}(\concat{s}{\langle d \rangle})}{\p_{\Lev_{\ell}(T^{i})}(s)} &   \text{if } \p_{\Lev_{\ell}(T^{i})}(s) \neq 0, \\[1ex]
             p_{i} &  \text{if } \p_{\Lev_{\ell}(T^{i})}(s) = 0 \conj d = 0,  \\[1ex]
             1 - p_{i} &  \text{if }  \p_{\Lev_{\ell}(T^{i})}(s) = 0 \conj d = 1.
        \end{array}\right.\]
    Similar as the case of $\Lev_\ell(T^i)$ we can verify that, indeed, $(\suc_{T^{i}}(s),  \calP(\suc_{T^{i}}(s),  \p_{s}^{i})$ is a probability space, and for $\ell \in I_{k} \cup \{ n^{\ast}_i \}$ we have that 
    \[ \p_{\Lev_{\ell}(T^{i})}( s) \coloneqq \prod_{\ell' \lhd^i \ell} \p_{s  {\rest} \ell'}(s \, {\rest}  (\ell' +_i 1)).\]
    Now, we must verify the hypotheses of \autoref{p70}: Let $\ell \in I_{k}$ and $s \in \Lev_{\ell}(T^{i}).$ 
    \begin{equation*}
        \begin{split}
            \p_{\ell +_i 1 }^{i}(\concat{s}{\langle 0 \rangle}) & = \p_{\Lev_{m^{\ast}}(T)}[\Phi(\rho)  {\rest}  (\ell +_i 1) = \concat{s}{\langle 0 \rangle}]\\ 
            & = \p_{\Lev_{m^{\ast}}(T)}[\Phi^{i}(\rho) {\rest}  \ell = s \conj \Phi^{i}(\rho)(\ell) = 0] \\
            & = \p_{\Lev_{m^{\ast}}(T)}[\Phi^{i}( \rho  {\rest}  m_{i, \ell}) = s \conj \Phi^{i}(\rho)(\ell) = 0]\\
            & = \p_{\Lev_{m_{i, \ell}+1}(T)}[\Phi^{i}(\rho  {\rest}  m_{i, \ell}) = s \conj \Phi^{i}(\rho)(\ell) = 0] \\
            & = \sum \left \{ \p_{\Lev_{m_{i, \ell}+1}(T)}(\rho) \colon \Phi^{i}(\rho {\rest} m_{i, \ell}) = s,\ \Phi^{i}(\rho)(\ell) = 0 ,\  \rho \in \Lev_{m_{i, \ell}+1}(T)  \right \},
        \end{split}
    \end{equation*}
    where the third equality is given by virtue of \ref{k40-3-1}, and the equalities prior to the sum are given by the definition of $\Phi^{i}$ and because the probability only depends on $\rho {\rest} m_{i, \ell} +1.$ Now, using that $\p_{\Lev_{k +1}(\calS)}(\concat{\rho}{\langle x \rangle}) = \p_ {\Lev_{k}(\calS)}(\rho) \cdot \p_{\rho}(\concat{\rho}{\langle x \rangle})$ holds in any probability tree $\calS,$ we substitute $\rho = \concat{\eta}{\sigma}$ and $ \concat{s}{\langle 0 \rangle} = t$ to get:
    \begin{equation*}
        \begin{split}
            \p_{\ell+_i 1}^{i}(t) & = \sum  \left  \{ \p_{\Lev_{m_{i, \ell}+1}(T)}(\rho) \colon   \Phi^{i}(\rho {\rest}  m_{i, \ell}) = s \conj \Phi^{i}(\rho)(\ell) = 0 \conj \rho \in \Lev_{m_{i, \ell}+1}(T)  \right \} \\
            & = \sum_{ \substack{\eta \in \Lev_{m_{i, \ell}}(T) \\ \Phi^{i}(\eta) = s}} \left( \sum_{\substack{ \sigma \in \Sigma_{\eta} \\ \sigma(i, \ell) = 0}} \p_{\Lev_{m_{i, \ell}}(T)}(\eta)  \cdot \p_{\eta}(\concat{\eta}{\langle \sigma \rangle})  \right)\\
            & = \sum_{ \substack{\eta \in \Lev_{m_{i, \ell}}(T) \\ \Phi^{i}(\eta) = s}} \left( \sum_{\substack{ \sigma \in \Sigma_{\eta} \\ \sigma(i, \ell) = 0}} \p_{\Lev_{m_{i, \ell}}(T)}(\eta)  \cdot f_{\eta}(\sigma)  \right).
        \end{split}
    \end{equation*}

    By \eqref{k40-9}, 
    \begin{equation*}
        \begin{split}
            \sum_{ \substack{\eta \in \Lev_{m_{i, \ell}}(T) \\ \Phi^{i}(\eta) = s}} \left( \sum_{\substack{ \sigma \in \Sigma_{\eta} \\ \sigma(i, \ell) = 0}} \p_{\Lev_{m_{i, \ell}}(T)}(\eta)  \cdot f_{\eta}(\sigma)  \right) & =\sum_{ \substack{\eta \in \Lev_{m_{i, \ell}}(T) \\ \Phi^{i}(\eta) = s}} \left( \p_{\Lev_{m_{i, \ell}}(T)}(\eta)  \cdot \sum_{\substack{\sigma \in \Sigma_{\eta} \\ \sigma(i, \ell) = 0)}} f_{\eta}(\sigma)  \right) \\
            & \geq   \left(  \sum_{ \substack{\eta \in \Lev_{m_{i, \ell}}(T) \\ \Phi^{i}(\eta) = s}}\p_{\Lev_{m_{i, \ell}}(T)}(\eta) \right) \cdot p_{i} \\
            & = \p_{\Lev_{m_{i, \ell}(T)}}[\Phi^{i}(\eta) = s ] \cdot p_{i}   = \p_{\ell}^{i}(s) \cdot p_{i}.
        \end{split}
    \end{equation*}
    
     In conclusion, for any $\ell \in I_{k}$ and $s \in \Lev_{\ell}(T^{i}),$ we have that
    \begin{equation}\label{k40-e11}
        \text{$\p_{\ell +_i 1}^{i}(\concat{s}{\langle 0 \rangle}) \geq p_{i} \cdot \p_{\ell}^{i}(s).$}
    \end{equation}
    
    If $\p_{\ell}^{i}(s) = 0$ then $p_{i} = \p_{s}^{i}(\concat{s}{\langle 0 \rangle}),$ and on the other hand, if $\p_{\ell}^{i}(s) \neq 0$ then we have, by \eqref{k40-e11}, that $p_{i} \leq \p_{s}^{i}(\concat{s}{\langle 0 \rangle}).$ 
    
    For each $i < i^{\ast}$ define the random variable $Y_{i}$ on $\Lev_{n^{\ast}}(T^{i})$ such that, for $s \in \Lev_{n^{\ast}}(T^{i}),$ 
    $$ Y_{i}(s) \coloneqq \vert \{ \ell \in I_{k} \colon s(\ell) = 0 \} \vert.$$
    Since we are under the hypothesis of \autoref{p70}, we have that, for any $z\in \RRR$, 
    \begin{equation}\label{k40-e13}
        \forall z \in \RRR\ ( \p_{\Lev_{n^{\ast}_i}(T^{i})}[Y_{i} \leq  z] \leq \p_{\Omega_{n^{\ast}}}[\B_{n^{\ast},  p_{i}} \leq z] ) \text{ where $n^\ast \coloneqq |I_k|$.}
    \end{equation}
    Finally, using  $\Phi^{i}$ again, since for $\rho \in \Lev_{m^{\ast}}(T),  Y_{i}(\Phi^{i}(\rho)) = \sum_{\ell \in I_{k}} X_{\ell}^{i}(\rho)$ we can conclude that, for any $i < i^{\ast}$ and any $z \in \RRR,$  
    \begin{equation}\label{k40-e15}
        \p_{\Lev_{m^{\ast}}(T)} \left[ \sum_{\ell \in I_{k}} X_{\ell}^{i} \leq z \right] = \p_{\Lev_{n^{\ast}_i}(T^{i})}[Y_{i} \leq z] \leq \p_{\Omega_{n^{\ast}}}[\B_{n^{\ast},  p_{i}} \leq z] ). 
    \end{equation}

    \textbf{Part 3. Find a suitable $\rho \in \Lev_{m^{\ast}}(T)$ with high probability.}

    For any $i < i^{\ast}$, set $z_{i} \coloneqq \vert I_{k} \vert (1 - \varp_{0}^{i})(1- \varp)$. Consider the following event in $\Lev_{m^{\ast}}(T)$: 
    $$ E_{i} \coloneqq \left  \{ \rho \in \Lev_{m^{\ast}}(T) \colon \vert \{ \ell \in I_{k} \colon \rho(m_{i, \ell})(i, \ell) = 0 \} \vert  > z_{i} \right \}, $$ 
    and define $F_{i} \coloneqq \Lev_{m^{\ast}}(T)\menos E_{i}$. Hence it is clear that
    \begin{equation}\label{k40-20}
        \text{$\p_{\Lev_{m^{\ast}}(T)}(F_{i}) = \p_{\Lev_{m^{\ast}}(T)} \left[ \sum_{\ell \in I_{k}} X_{\ell}^{j} \leq z_{i} \right] $}. 
    \end{equation}
    Therefore, using  \eqref{k40-20}, \eqref{k40-e15}, \eqref{k40-e13},  \eqref{k40-5} and Chebyshev’s inequality, since $n^*=|I_k|>N_0$ we get:
    \begin{equation*}
        \begin{split}
            \p_{\Lev_{m^{\ast}}(T)}(F_{i}) &= \p_{\Lev_{m^{\ast}}(T)} \left [ \sum_{\ell \in I_{k}} X_{\ell}^{i} < z_{i} \right] \leq \p_{\Lev_{n^{\ast}}(T^{i})}[Y_{i} \leq z_{i}] \leq \p_{\Omega_{n^{\ast}}}[\B_{n^{\ast},  p_{i}} \leq z_{i}]\\
            & = \p_{\Omega_{n^{\ast}}}\left[\E[\B_{n^{\ast},  p_{i}}] - z_{i} \leq \E[\B_{n^{\ast},  p_{i}}] - \B_{n^{\ast},  p_{i}} \right]\\ 
           & \leq \p_{\Omega_{n^{\ast}}} \left[ \vert \B_{n^{\ast},  p_{i}} - \E[\B_{n^{\ast},  p_{i}}] \vert \geq \E[\B_{n^{\ast}, p_{i}}] - z_{i} \right ] \leq \frac{\Var[\B_{n^{\ast},  p_{i}}]}{(\E[\B_{n^{\ast},  p_{i}}] - z_{i})^{2}}\\
           & = \frac{\vert I_{k} \vert  (1 - \varp_{0}^{i} - \varp' ) ( \varp_{0}^{i} + \varp')}{\vert I_{k} \vert^{2}  [\varp (1 - \varp_{0}^{i}) - \varp']^{2}} = \frac{1}{\vert I_{k} \vert} \frac{(1 - \varp_{0}^{i} - \varp' )  ( \varp_{0}^{i} + \varp')}{[\varp (1 - \varp_{0}^{i}) - \varp']^{2}} < \frac{1}{i^{\ast} + 1}. 
        \end{split}
    \end{equation*}
    As a consequence, in $\Lev_{m^{\ast}}(T)$ we have 
    $$ \p \left( \bigcup_{i < i^{\ast}} F_{i}  \right) \leq \sum_{i < i^{\ast}} \p(F_{i}) < \sum_{i < i^{\ast}} \frac{1}{i^{\ast} + 1} = \frac{i^{\ast}}{i^{\ast} + 1} < 1,  $$
    hence $\p \left( \bigcap_{i < i^{\ast}} E_{i}  \right) > 0$, 
    which implies that $\bigcap_{i < i^{\ast}} E_{i} \neq \emptyset$, i.e.\ there exists some $\eta \in \bigcap_{i < i^{\ast}} E_{i}.$   Since for any $(i, \ell) \in i^{\ast} \times I_{k}$, 
    $r_{\eta} \forces_{\beta_{m^{\ast}}}\text{``}\dot{r}_{\ell}^{i} \in \dot{G}(\beta_{\beta_{m, \ell}}) \Leftrightarrow \eta(m_{i, \ell})(i, \ell) = 0\text{''},$
    it is clear that $r_{\eta} \forces_{\beta_{m^{\ast}}}$``$k \in \dot a\cap \bigcap_{i < i^{\ast}} \dot{A}_{i}$''.
\end{PROOF}

The following result corresponds to \autoref{k40} for the case $\varp_{0} = 0$: 

\begin{mainlemma}\label{k41}
    \autoref{k40} is valid without assuming~\ref{k40aa} and without the ``one-to-one" condition in~\ref{k40ii}, but if restricting~\ref{k40iv} to $\varp_0=\varp=0$. Concretely, assuming~\ref{k40a},~\ref{k40b} and~\ref{k40c} of \autoref{k40}, there is a $\bbP_{\pi}$-name  $\dot{\Xi}_{\pi}$ of a fam on $\pts(K)$ extending $\bigcup_{\alpha < \pi} \dot{\Xi}_{\alpha}$ such that, whenever
    \begin{enumerate}[label = \rm (\roman*)]
        \item\label{k41i} $\bar \beta\coloneqq\langle \beta_{\ell} \colon  \ell \in W \rangle$ is a sequence of ordinals ${<}\pi$ with supremum $\pi$,
        
        \item\label{k41ii} $\Xi_0(\{k\in K \colon  \forall \ell\in I_k\ (\beta_\ell\geq\alpha)\}) =1$ for all $\alpha<\pi$,
        
        \item\label{k41iii} $\dot{r}_{\ell}$  is a  $\bbP_{\beta_{\ell}}$-name such that  $ \Vdash_{\beta_{\ell}}$``\/$\dot{r}_{\ell} \in \dot{Q}_{\beta_{\ell}}'$'' for any $\ell \in W$, and
        
        \item\label{k41iv} $\langle \varp_{\beta{\ell}}' \colon \ell \in W \rangle$ is constant with value $0$ (i.e.\ $\dot Q'_{\beta_\ell}$ is forced centered for each $\ell\in W$),
    \end{enumerate}
    we have $\Vdash_{\bbP_{\pi}} \text{``\/} \dot{\Xi}_{\pi} \left( \setcolon{ k \in K }{ \forall \ell \in I_{k} \ ( \dot{r}_{\ell} \in \dot{G}(\beta_{\ell}) )} \right) = 1 \text{''}.$ 
    
    Even more, if $\Por_\alpha$ forces that $\dot\Xi_\alpha$ represents an ultrafilter for all $\alpha<\pi$, then $\dot\Xi_\pi$ can be found forced to represent an ultrafilter. The same applies for the uap.
\end{mainlemma}

\begin{PROOF}[\textbf{Proof}]{\ref{k41}}
    Let $\cI$ be the set of $\bar c = \langle (\beta_{\ell}, \dot{r}_{\ell}) \colon \ell \in W \rangle$ as in~\ref{k41i}--\ref{k41iv}. Let $\dot{\Xi}$ be a $\Por_\pi$-name for $\bigcup_{\alpha<\pi}\dot\Xi_\alpha$ and, for each $\bar{c} \in \cI,$ let $\dot A_{\bar{c}}$ be a $\bbP_{\pi}$-name of $\setcolon{ k \in K }{ \forall \ell \in I_{k} \ ( \dot{r}_{\ell} \in \dot{G}(\beta_{\ell}) )}$. 
    
    Let $i^*<\omega$, $C=\{\bar c^i\colon i<i^*\} \subseteq \cI$, $\dot A_{i} \coloneqq \dot A_{\bar{c}^{i}}$ for $i<i^*$, $p \in \bbP$ and let $\dot{a}$ be a $\bbP_{\pi}$-name such that $ p \forces_{\bbP_{\pi}} \text{``} \dot{a} \in \dom(\dot{\Xi}) \text{ and } \dot{\Xi}(\dot{a}) > 0\text{''}$. In virtue of \autoref{m69}, it suffices to prove that there is some $q \leq p$  forcing $\dot{a} \cap \bigcap_{i<i^*} \dot A_{i} \neq \emptyset$.
    
    Just like in the proof of~\autoref{k40}, find 
    $\gamma < \pi$, $k\in K$ and $p_0\leq p$ such that $p_0\in \Por_\gamma$, $\dot a$ is a $\Por_\gamma$-name, $\{\beta^i_\ell\colon \ell\in I_k\}\subseteq \pi\menos \gamma$ for all $i<i^*$ and $p_0\Vdash_\gamma$``$k\in \dot a$''. 

    Since $\{ \beta_{\ell}^{i} \colon i < i^{\ast},\ \ell \in I_{k} \}$ is finite, it has an increasing enumeration $\langle \beta_{m} \colon m < m^{\ast} \rangle$ for some $m^{\ast} < \omega$. 
    By recursion, let us build a sequence of conditions $\langle p_{m} \colon m \leq m^{\ast} \rangle$, where $\beta_{m^*}\coloneqq\pi$, such that, for any $m \leq m^{\ast}$, $p_{m} \in \bbP_{\beta_{m}}$, as follows. For $m = 0$, let $p_{0}$ be the condition we already have. Now, assume that we have defined $p_{m} \in \bbP_{\beta_{m}}$ and consider $J_{m} \coloneqq \{ (i, \ell) \in i^{\ast} \times I_{k} \colon \beta_{\ell}^{i} = \beta_{m} \}$. Since $\Vdash_{\beta_{m}}$``$\dot{r}_{\ell}^{i} \in \dot{Q}_{\beta_{m}}'$'' for any $(i, \ell) \in J_{m}$ and, by \autoref{s8}~\ref{s8d}, $\bbP_{\beta_{m}}$ forces that $\dot{Q}_{\beta_{m}}'$ is centered, there exists some $\bbP_{\beta_{m}}$-name $\dot{q}_{m}$ of a condition in $\dot{\bbQ}_{\beta_{m}}$ such that, for any $(i, \ell) \in J_{m}$, $p_{m} \Vdash_{\beta_{m}}$``$\dot{q}_{m} \leq \dot{r}_{\ell}^{i}$''. Set $p_{m+1} \coloneqq q_{m} \cup \{ (\beta_{m}, \dot{q}_{m}) \}  \in \bbP_{\beta_{m+1}}$. Finally, define $q \coloneqq p_{m^{\ast}}.$ 
    Then $q \leq p_{0} \leq p$ and, for any $(i, \ell) \in J_{m}$, we have that $q \rest \beta_{\ell}^{i} \forces_{\beta_{\ell}^i}$``$q(\beta_{\ell}^{i}) \leq \dot{r}_{\ell}^{i}$''. Thus, $q \forces_{\bbP_{\pi}}$``$k\in \dot{a} \cap \bigcap_{(\bar{c}, \varp) \in C} \dot A_{\bar{c}, \varp} \neq \emptyset$''.  
\end{PROOF}

As a consequence of~\autoref{k40} and~\autoref{k41}, we can extend $\calC$-iterations at limit steps, that is, we get an \emph{extension theorem}: 

\begin{theorem}\label{k42}
    Let $\pi$ be a limit ordinal and let $\calC$ be a bedrock such that, for any $\bfc\in\calC$, $g^\bfc$ is a guardrail for $\pi$ and 
    \begin{enumerate}[label = \rm $(\textrm{\faMotorcycle})$]
        \item\label{k42hyp} either $K^\bfc$ is finite or, for any $n<\omega$, $\Xi^\bfc_0(\{k\in K^\bfc\colon |I_k^\bfc|=n\}) = 0$.
    \end{enumerate}    
    
    Let $\bft$ be a prepared iteration of length $\pi$ and let $\langle \dot{\Xi}_{\alpha}^{\bfc} \colon \alpha < \pi, \ \bfc \in \calC \rangle$ be a sequence such that, for any $\alpha < \pi,$ $\bft {\rest} \alpha \sqcup \langle \dot{\Xi}_{\beta}^{\bfc} \colon \beta \leq \alpha, \, \bfc \in \calC \rangle$ is a $\calC$-iteration. Then there exists a sequence $\langle \dot{\Xi}_{\pi}^{\bfc} \colon \bfc \in \calC \rangle$ such that $\bft  \sqcup  \langle \dot{\Xi}_{\alpha}^{\bfc} \colon \alpha \leq \pi, \ \bfc \in \calC \rangle$ is a $\calC$-iteration. Even more, if $\bfc\in\calC$ and, for any $\alpha<\pi$, $\dot\Xi^\bfc_\alpha$ is forced by $\Por_\alpha$ to represent an ultrafilter, then $\dot \Xi^\bfc_\pi$ can be found forced to represent an ultrafilter. The same applies for the uap.

    Moreover, if all $\dot Q^\alpha_{\zeta,\varp}$ are forced centered in $\Qnm_\alpha$ (by $\Por^-_\alpha$) for any $\alpha<\pi$, then~\ref{k42hyp} is not needed for this result.
\end{theorem}

\begin{PROOF}[\textbf{Proof}]{\ref{k42}}
    Fix $\bfc\in\calC$. When  $K^\bfc$ is finite the construction of $\dot\Xi^\bfc_\pi$ is trivial (it must be forced equal to $\Xi^\bfc$), so assume that $K^\bfc$ is infinite. For $\alpha<\pi$, let $\dot Q'_\alpha\coloneqq \dot Q^\alpha_{g^\bfc(\alpha)}$ if $g_{0}^{\bfc}(\alpha)<\theta_\alpha$ and $g_{1}^{\bfc}(\alpha)=0$ iff $\alpha\in U^{\bft}$, otherwise let $\dot Q'_\alpha$ be a $\Por^-_\alpha$-name of some singleton contained in $\Qnm_\alpha$.  We verify the conditions of \autoref{k40}. By hypothesis and \autoref{k28}, \ref{k40a} and \ref{k40b} are clear, \ref{k40aa} follows from~\ref{k42hyp}
    and, for \noindent\ref{k40c}, $\dot Q'_\alpha$ is forced to be either centered or with intersection number ${\geq}1-g^\bfc_1(\alpha)$, the latter by~\ref{k42hyp} and \autoref{i11.0}. Then we can find a $\dot\Xi^\bfc_\pi$ as in \autoref{k40} (or in \autoref{k41} when all $\dot Q^\alpha_{\zeta,\varp}$ are forced centered).

    We now check the conditions in \autoref{k28} to verify that $\bft  \sqcup  \langle \dot{\Xi}_{\alpha}^{\bfc} \colon \alpha \leq \pi, \ \bfc \in \calC \rangle $ is a $\calC$-iteration. Conditions~\ref{k28a}--\ref{k28ee} are clear. To check~\ref{k28f}, assume that $\tau = \{(\dot q_\ell,\alpha_\ell,\varp) \colon \ell\in W^\bfc\}$ follows $g^\bfc$. The case~\ref{k28fa} is already taken care of by the intermediate steps of the iteration.

    To see~\ref{k28fb}, assume that $\bar\alpha\coloneqq \langle \alpha_{\ell} \colon  \ell \in W^\bfc \rangle$ has no maximum, $\beta\coloneqq \sup_{\ell \in W^\bfc} \alpha_{\ell}$ and, for any $\alpha<\beta$, $\Xi^\bfc(\{k\in K^\bfc \colon \bar\alpha\rest I_k \text{ is one-to-one and } \forall \ell\in I^\bfc_k\ (\alpha_\ell\geq \alpha)\}) =1$. When $\beta < \pi,$ $\bft {\rest} \beta \sqcup \langle \dot{\Xi}_{\gamma}^{\bfc} \colon \gamma \leq \beta, \ \bfc \in \calC\rangle$ is a $\calC$-iteration, hence for all $\varp' > 0,$ 
    \[ \Vdash_{\bbP_{\beta}} \text{``} \ \dot{\Xi}_{\beta}^{\bfc} \left( \left \{ k \in K^\bfc \colon  \frac{\vert \{ \ell \in I_{k}^\bfc \colon \dot{q}_{\ell} \in \dot{G}(\alpha_{\ell}) \} \vert}{\vert I_{k}^\bfc \vert}  \geq (1 - \varp) (1 - \varp') \right \} \right) = 1 \text{''}.\]
        
    In the case $\beta=\pi$, the requirements~\ref{k40i}--\ref{k40iv} of \autoref{k40}  are fulfilled, so, for any $\varp'>0$, 
    \begin{equation*}
        \Vdash_{\bbP_{\pi}} \text{``} \ \dot{\Xi}_{\pi}^{\bfc} \left( \left \{ k \in W^\bfc \colon \frac{\vert \{ \ell \in I_{k}^\bfc \colon  \dot{q}_{\ell} \in \dot{G}(\alpha_{\ell}) \} \vert}{\vert I_{k}^\bfc \vert}  \geq (1 - \varp) (1 - \varp') \right \} \right) = 1 \text{''}.\qedhere 
    \end{equation*}
\end{PROOF}


\section{FAM-limits in finite support iterations}\label{famlimits}

Suppose that $\bft$ is a prepared iteration and $\bar{p} = \langle p_{i} \colon i \in L \rangle \subseteq \bbP_{\pi}^{\bullet}.$ Let $i \in L$ and $\xi \in \dom(p_{i}).$ Since $p_{i} \in \bbP_{\pi}^{\bullet},$ there are $\zeta < \theta_{\xi}$ and $\varp \in (0, 1)_{\QQQ}$ such that $\forces_{\xi}$``$p_{i} (\xi) \in \dot{Q}_{\zeta, \varp}^{\xi}$''. If for $\xi$ we could find a \emph{suitable} $W \subseteq L$ such that, for all $\ell \in W$, $\xi \in \dom(p_{\ell})$ and $\forces_{\xi}$``$p_{\ell}(\xi) \in \dot{Q}_{\zeta, \varp}^{\xi}$'', where $\zeta$ and $\varp$ do no depend on $\ell,$ then we could apply a fam limit to $\la p_\xi(\ell) \colon \ell\in W\ra$ (as in \autoref{k28}~\ref{k28fa}). This idea can be used to define a fam-limit to some sequence $\la p_\ell\colon \ell\in W\ra$ in $\Por^\bullet_\pi$. After a complete analysis of the conditions that are required to formalize this idea, we arrive at the notion of \emph{uniform $\Delta$-system}: 

\begin{definition}\label{k50}
    Let $\bft$ be a prepared iteration and let $L$  be a set. 
    \begin{enumerate}[label = \normalfont(\arabic*)]
        \item\label{k50-1}   Say that $\bar{p} = \langle p_{l} \colon  l \in L \rangle \subseteq \bbP_{\pi}$ is a \emph{uniform-$\Delta$-system with parameters $(\Delta, \vec{\alpha}, n^{\ast},  r^{\ast}, \zeta^{\ast},  \varp^{\ast})$} when:        
        \begin{enumerate}[label=\rm(\alph*)]
            \item\label{k50a} For any  $\ell \in L$, $p_{\ell} \in \bbP_{\pi}^{\bullet}.$
    
            \item\label{k50b} $\{ \dom(p_{\ell}) \colon \ell \in L \}$ forms a $\Delta$-system with root $\Delta.$
    
            \item\label{k50c} $ n^{\ast} < \omega$ and, for any $\ell \in L$, $\dom(p_{\ell}) = \{ \alpha_{n,\ell} \colon n < n^{\ast} \}$ is enumerated increasingly, that is, $m < n < n^{\ast} \Rightarrow \, \alpha_{m,\ell} < \alpha_{n,\ell}.$
    
            \item\label{k50d} $r^{\ast} \subseteq n^{\ast}$ and, for any $\ell \in L $ and $n < n^{\ast}$, $n \in r^{\ast} \Leftrightarrow \alpha_{n,\ell} \in \Delta.$ So, whenever $n \in r^{\ast},$ the sequence $\langle \alpha_{n, \ell} \colon \ell\in L \rangle$ is constant with value, say, $\alpha_{n}^{\ast}.$ 
    
            \item\label{k50e} For any $n \in n^{\ast} \menos r^{\ast},$ the sequence $\langle \alpha_{n, \ell} \colon \ell \in L \rangle$ is one-to-one.
    
            \item \label{k50e1}  $\zeta^{\ast} \colon n^{\ast} \times L \to \Ord$  is a function such that $\zeta^{\ast}(n, \ell) < \theta_{\alpha_{n, \ell}}$ and $\langle \zeta^{\ast}(n, \ell) \colon \ell \in L \rangle$ is constant with value $\zeta_{n}^{\ast}$ for all $n \in r^{\ast}.$ 
    
            \item\label{k50g} $\varp^{\ast} \colon n^{\ast} \to [0, 1)_{\QQQ}$ is a function.
    
            \item\label{k50f} For any  $\ell \in L $ and $n < n^{\ast}$, $\Vdash_{\alpha_{n, \ell}}$``$p_{\ell}(\alpha_{n, \ell}) \in \dot{Q}_{\zeta^{\ast}(n, \ell), \varp^{\ast}(n)}^{\alpha_{n, \ell}}$''.
        \end{enumerate}
        Notice that $\varp^*$ determines whether an $\alpha_{n,\ell}$ is in $F^{\bft}$ or $U^{\bft}$, since $\alpha_{n,\ell}\in F^{\bft}$ iff $\varp^*(n) > 0$.

        \item \label{k50-2} We say that an uniform $\Delta$-system with parameters $(\Delta, \vec{\alpha}, n^{\ast}, r^{\ast}, \zeta^{\ast}, \varp^{\ast})$ \emph{follows} a guardrail $g$ for $\pi$ if,  for any $n < n^{\ast}$ and $\ell \in L$, $g(\alpha_{n,\ell}) = (\zeta^{\ast}(n, \ell), \varp^{\ast}(n)).$
    \end{enumerate}
\end{definition}

In general, from a sequence of conditions of regular size $\theta$, we can obtain uniform $\Delta$-systems following guardrails. For this, we will use a refined version of the well-known \emph{$\Delta$-system lemma}:  

\begin{lemma}\label{v7}
    Let $\theta$ be an uncountable regular cardinal and let $\calA = \langle A_{\gamma} \colon \gamma < \theta \rangle$  be a family of finite sets such that $(\bigcup \calA, \lhd)$ is a well-order. Then, there are $E \in [\theta]^{\theta}, \, n^{\ast} < \omega$ and $r^{\ast} \subseteq n^{\ast},$ such that: 

    \begin{enumerate}[label=\rm(\arabic*)]
        \item $\{ A_{\xi} \colon \xi \in E \}$ forms a $\Delta$-system with root $\Delta,$

        \item for any $\xi \in E, \, \dom(A_{\xi}) = \{ a_{n, \gamma} \colon n < n^{\ast} \}$ is arranged in $\lhd$-increasing order,

        \item $a_{n, \xi} \in \Delta$ iff $n \in r^{\ast},$ for any $\xi \in E,$ and

        \item for any $n \in n^{\ast} \menos r^{\ast}$ and $\xi, \zeta \in E,$ if $\xi < \zeta$ then $a_{n, \, \xi} \lhd a_{n \, \zeta}.$  
    \end{enumerate}
\end{lemma}

Although this is not the usual way of stating the $\Delta$-system lemma (see, for example, \cite[Lem.~III.2.6]{Kunen}), it suits our particular needs. The proof is the same as the one presented in \cite[Thm.~1.6]{Kunenv}.

\begin{theorem}\label{k52}
    Let $\theta \geq \kappa$ be uncountable regular cardinals and $\bft$ a $\kappa$-prepared iteration. If $\{ p_{\xi} \colon \xi < \theta \} \subseteq \bbP_{\pi}^{\bullet},$ then there are some $E \in [\theta]^{\theta}$ and some  guardrail $g$ for $\pi$  
    such that $\{p_{\xi} \colon \xi \in E \}$ is a uniform $\Delta$-system following $g$. 
    Moreover, considering the parameters of this uniform $\Delta$-system as in \autoref{k50}~\ref{k50-1}, for $n\in n^*\menos r^*$, $\la \alpha_{n,\xi}\colon \xi\in E\ra$ is increasing.
\end{theorem}

\begin{PROOF}[\textbf{Proof}]{\ref{k52}}
    By \autoref{v7}, we can obtain $E_0 \in [\theta]^{\theta},$ $\Delta$, $n^{\ast} < \omega$ and $r^{\ast} \subseteq n^{\ast}$ such that: 
    \begin{enumerate}[label = \normalfont (\alph*)]
        \item\label{k52pa} $\{ \dom( p_{\xi}) \colon \xi \in E_0 \}$ forms a $\Delta$-system with root $\Delta.$

        \item\label{k52pb} For any $\xi \in E_0$, $\dom(p_{\xi}) = \{ \alpha_{n, \xi} \colon n < n^{\ast} \}$ is arranged in increasing order. 

        \item\label{k52pc} For any $\xi \in E_0$,  $\alpha_{n, \, \xi} \in \Delta$ if, and only if, $n \in r^{\ast}$.

        \item\label{k52pd} For any $n \in n^{\ast} \menos r^{\ast}$ and $\xi, \zeta \in E_0,$ if $\xi < \zeta$ then $\alpha_{n, \, \xi} < \alpha_{n \, \zeta}.$ 
    \end{enumerate}

    By~\ref{k52pb} and~\ref{k52pc} we have that, for any $n \in r^{\ast},$ the sequence $\langle \alpha_{n, \xi} \colon \xi \in E_0 \rangle$ is constant with value, say, $\alpha_{n}^{\ast}.$ 
    Since $p_{\xi} \in \bbP_{\xi}^{\bullet}$ for every $\xi \in E_0,$ for any $n < n^{\ast}$ we can find $\varp_{n, \xi} \in [0, 1)_{\QQQ}$ and $\zeta_{n, \xi} < \theta_{\alpha_{n, \xi}}$ such that $\forces_{\Por^-_{\alpha_{n, \xi}}}$``$p_{\xi}(\alpha_{n, \xi}) \in \dot{Q}_{\zeta_{n, \xi}, \varp_{n, \xi}}^{\alpha_{n, \xi}}$''. 
    
    For any $\bar{\varp} = \langle \varp_{n} \colon n < n^{\ast} \rangle \in {}^{n^*}[0, 1)_{\QQQ},$ define $E_{\bar{\varp}} \coloneqq \{ \xi \in E_0 \colon \forall n < n^{\ast}\ ( \varp_{n, \xi} = \varp_{n} ) \}. $ It is clear that $\{ E_{\bar{\varp}} \colon \bar{\varp} \in {}^{n^{\ast}}[0, 1)_{\QQQ} \}$ is a countable family of pairwise disjoint sets whose union is $E_0.$ So, since $\vert E_0 \vert = \theta$ is regular and uncountable, there exists an $\varp^{\ast} \in {}^{n^{\ast}}[0, 1)_{\QQQ}$ such that $ \vert E_{\varp^{\ast}} \vert = \theta.$  Define $F \coloneqq E_{\varp^{\ast}}$. 
    As a consequence, $\varp_{n, \xi} = \varp^{\ast}(n)$ for all $n < n^{\ast}$ and $\xi \in F.$ 
    
    For any $\bar{\zeta} = \langle \zeta_{n} \colon n \in r^{\ast} \rangle \in \prod_{n \in r^{\ast}} \theta_{\alpha_{n}^{\ast}}$ define $F_{\bar{\zeta}} \coloneqq \{ \xi \in F \colon \forall n \in r^{\ast}\ (\zeta_{n, \xi} = \zeta_{n}) \},$ hence we have that  $\{ F_{\bar{\zeta}} \colon \bar{\zeta} \in \prod_{n \in r^{\ast}} \theta_{\alpha_{n}^{\ast}} \}$ is a family of pairwise disjoint sets whose union is $F.$ Since $\vert F \vert = \theta$ is regular and $ \left\vert \prod_{n \in r^{\ast}} \theta_{\alpha_{n}^{\ast}} \right\vert < \kappa\leq \theta$, there is some $\zeta^{\ast}\in \prod_{n \in r^{\ast}} \theta_{\alpha_{n}^{\ast}}$ such that $|F_{\zeta^*}| = \theta$. Define $E \coloneqq F_{\zeta^{\ast}},$ hence $\zeta_{n, \xi} = \zeta_{n}^{\ast}$ for all $\forall n \in r^{\ast}$ and $\xi \in E.$  
    It is clear that $\{ p_{\xi} \colon \xi \in E \}$ is a uniform $\Delta$-system.

    We define a guardrail $g \colon \pi \to \kappa^-\times[0,1)_\QQQ$ such that, for any $\alpha < \pi,$ 
    $$g(\alpha) \coloneqq \left\{ \begin{array}{ll}
        (\zeta_{n, \xi}, \varp_{n, \xi}) &   \text{if } \exists n < n^{\ast} \exists \xi \in E\ (\alpha = \alpha_{n, \xi}), \\[1ex]
        \left( 2022, \frac{1}{6} \right) &  \text{if }  \alpha \notin \bigcup_{\xi \in E} \dom(p_{\xi}).
    \end{array}
   \right.$$ 

    Finally, it is clear that $\{ p_{\xi} \colon \xi \in E \}$ follows $g$.
\end{PROOF}

As far as guardrails are concerned, uniformity \emph{\textbf{only}} depends on the restriction of the guardrail to the parameter $\vec{\alpha}$. Hence, it is always possible to reduce uniform $\Delta$-systems to a countable uniform $\Delta$-systems following a guardrail in some complete set of guardrails. Formally,

\begin{theorem}\label{k54}
    Let $\bft$ be $\kappa$-prepared iteration, $\calG$ a complete set of guardrails for $\pi$ over $\kappa^-,$ and let $\langle p_{\ell} \colon \ell\in L \rangle$ be a uniform $\Delta$-system with $L$ countable.  
    Then there exists a guardrail $g \in \calG$ 
    such that  $\bar{p} =  \langle p_{\ell} \colon \ell \in L \rangle$ follows $g$ (with the same parameters). 
\end{theorem}

\begin{PROOF}[\textbf{Proof}]{\ref{k54}}
    Let $(\Delta, \vec{\alpha}, n^{\ast},  r^{\ast}, \zeta^{\ast}, \varp^{\ast})$ be the parameters of the $\Delta$-system and 
    let $X \coloneqq \{ \alpha_{n, \ell} \colon \ell \in L \conj n < n^{\ast} \}$, which is countable. 
    Define $\sigma\colon X\to \kappa^-\times [0,1)_\QQQ$ by $\sigma(\alpha_{n,\ell})\coloneqq(\zeta^*(n,\ell),\varp^*(n))$. By virtue of the completeness of $\calG$, there exists some $g \in \calG$ such that $\sigma \subseteq g.$ It is clear that the uniform $\Delta$-system follows $g$.
\end{PROOF}

Notice that \autoref{k54} is a good motivation for the definition of a complete set of guardrails (see \autoref{k16}).

Fix an iteration bedrock $\calC$. 
In the following, we show that our definition of uniform $\Delta$-system allows us to define the desired limit that we mentioned at the beginning of this subsection. 

\begin{definition}\label{k56}    
     Let $\bft$ be a $\calC$-iteration, $\bfc \in \calC $ and $\bar{p} = \langle p_{\ell} \colon \ell \in W^\bfc \rangle \subseteq \bbP_{\pi}$ be a uniform $\Delta$-system with parameters $(\Delta, \,  n^{\ast}, \vec{\alpha}, \, \, r^{\ast}, \, \zeta^{\ast}, \, \varp^{\ast})$ following $g^\bfc$.  We define the \emph{limit} $\limit_{\bft}^{\bfc}(\bar{p}),$ as follows:

     \begin{enumerate}[label=\rm(\alph*)]
         \item $\dom(\limit_{\bft}^{\bfc}(\bar{p})) = \Delta,$

         \item For any $n \in r^{\ast}$, $\limit_{\bft}^{\bfc}(\bar{p})(\alpha_{n}^{\ast})$ is a $\Por^-_{\alpha^*_n}$-name of $\limit^{\Xi_{\alpha_{n}^{\ast}}^{\bfc,-}} ( \langle p_{\ell}(\alpha_{\ell}^{\ast}) \colon \ell \in W^\bfc \rangle).$
     \end{enumerate}

     When the context is clear, we simply write  $\limit^{\bfc}(\bar{p})$ to refer to $\limit_{\bft}^{\bfc}(\bar{p}).$ This limit is well-defined: By \autoref{k50}~\ref{k50d}, for any $n \in r^{\ast},$ $\langle \alpha_{n, \ell} \colon \ell \in W^{\bfc} \rangle$ is constant with value $\alpha_{n}^{\ast}.$  Since $\bar{p}$ is a uniform $\Delta$-system following $g$ we have that, for any $\ell \in W^\bfc$ and $n < n^{\ast}$,  $p_{\ell}(\alpha_{n, \ell})$ is a $\bbP_{\alpha_{n,\ell}}^{-}$-name. Also, in $V^{\bbP_{\alpha_{n}^{\ast}}^{-}}$, $\langle p_{\ell}(\alpha_{n}^{\ast}) \colon \ell \in W^\bfc \rangle \in {}^{W^\bfc} \dot{Q}_{g(\alpha_{n}^{\ast})}^{\alpha_{n}^{\ast}}$ and therefore, $\limit^{\Xi_{^{\alpha_{n}^{\ast}}}^{\bfc,-}}(\langle p_{\ell}(\alpha_{n}^{\ast}) \colon \ell \in W^{\bfc} \rangle) \in \dot{\bbQ}_{\alpha_{n}^{\ast}}$ is defined.
\end{definition}

In the following theorem, we present some properties of the limit that we have just defined. Some of them appear implicitly in the proof of \cite[Lem.~3.4]{Sh00}. Property~\ref{k60d} is particularly interesting, since, as in any $(\Xi,\bar I,\varp_0)$-linked set, $\limit^{\bfc}(\bar{p})$ forces ``many'' conditions to fall into the generic filter.

\begin{theorem}\label{k60}
    Let $\bft$ be a $\calC$-iteration,  $\bfc \in \calC$ and let $\bar{p} = \langle p_{\ell} \colon \ell \in W^\bfc \rangle$ be a uniform $\Delta$-system following $g\coloneqq g^\bfc$ with parameters $(\Delta, \vec{\alpha},  n^{\ast},  r^{\ast}, \zeta^{\ast}, \varp^{\ast}).$ 
    Assume that, for any $n\in n^*\menos r^*$, the sequence $\la\alpha_{n,\ell}\colon \ell\in W^\bfc\ra$ has no maximum and, for all $\alpha< \sup_{\ell\in W^\bfc}\alpha_{n,\ell}$, 
    \[\Xi^\bfc_0(\set{k\in K^\bfc}{\forall \ell\in I^\bfc_k\ (\alpha\leq\alpha_{n,\ell})}) =1.\] 
    Then, $\limit^{\bfc}(\bar{p})$ satisfies the following properties:  

    \begin{enumerate}[label=\rm(\arabic*)]
        \item\label{k60a} For $n \in r^{\ast},$ $\displaystyle\limit^{\bfc}(\bar{p}) \Vdash_{\bbP_{\pi}} \text{``} \int_{K^\bfc} \frac{\vert \{ \ell \in I_{k} \colon p_{\ell}(\alpha_{n}^{\ast}) \in  \dot{G}(\alpha_{n}^{\ast}) \} \vert}{\vert I_{k} \vert} d  \Xi_{\pi}^{\bfc}(k) \geq 1 - \varp^{\ast}(n)\text{''}.$

        \item\label{k60b} For $n \in n^{\ast} \menos r^{\ast}$ and $\varp > 0,$ 
        \[\Vdash_{\bbP_{\pi}} \text{``} \ \dot{\Xi}_{\pi}^{\bfc} \left( \left \{ k \in K^\bfc \colon \frac{\vert \{ \ell \in I_{k} \colon p_{\ell}(\alpha_{n, \ell})  \in \dot{G}(\alpha_{n, \ell}) \} \vert}{\vert I_{k} \vert}  \geq (1 - \varp^{\ast}(n)) (1 - \varp) \right \} \right) = 1 \text{''},\]
        in particular, $\displaystyle\Vdash_{\bbP_{\pi}} \text{``} \int_{K^\bfc} \frac{\vert \{ \ell \in I_{k} \colon p_{\ell}(\alpha_{n, \ell}) \in  \dot{G}(\alpha_{n, \ell}) \} \vert}{\vert I_{k} \vert} d  \Xi_{\pi}^{\bfc}(k) \geq 1 - \varp^{\ast}(n)\text{''}.$

        \item\label{k60c} $\displaystyle\limit^{\bfc}(\bar{p}) \Vdash_{\bbP_{\pi}} \text{``} \int_{K^\bfc} \frac{\vert \{ \ell \in I_{k} \colon p_{\ell} \in \dot{G} \} \vert}{\vert I_{k} \vert} d  \Xi_{\pi}^{\bfc}(k) \geq 1 - \sum_{n < n^{\ast}} \varp^{\ast}(n)\text{''}.$

        \item\label{k60d} If $0 < \varp < 1 -\sum_{n < n^{\ast}} \varp^{\ast}(n)$ then 
        $\displaystyle\limit^{\bfc}(\bar{p}) \Vdash_{\bbP_{\pi}} \text{``}\, \Xi_{\pi}^{\bfc}(\dot{A}_{\varp}) \geq 1 - \frac{1}{1-\varp}\sum_{n<n^*}\varp^*(n)$'', where 
        $$\dot{A}_{\varp} \coloneqq \left \{ k \in K^\bfc \colon   \frac{\vert \{ \ell \in I_{k} \colon p_{\ell} \in \dot{G} \} \vert}{\vert I_{k} \vert} > \varp \right \}. $$ 
        As a consequence, $\limit^{\bfc}(\bar{p}) \Vdash_{\bbP_{\pi}} $``$\dot{A}_{\varp}$ is infinite'' whenever $\Xi^\bfc_0$ is free. 
    \end{enumerate}
\end{theorem}

\begin{PROOF}[\textbf{Proof}]{\ref{k60}}
    Let $q \coloneqq \lim^{\bfc}(\bar{p})$ and $\tau_n \coloneqq \langle (p_{\ell}(\alpha_{n, \ell}),\alpha_{n, \ell},\varp^*(n)) \colon \ell \in W^\bfc \rangle$ for $n<n^*$. It is clear that $\tau_n$ follows $g^\bfc$.
    
    \noindent\ref{k60a}: Let $n \in r^{\ast}$. As $\langle \alpha_{n, \ell} \colon \ell \in W^\bfc \rangle$ is constant with value $\alpha^*_n$, we have by \autoref{k28}~\ref{k28fa}, that: 
    $$ \Vdash_{\bbP_{\alpha_{n}^{\ast}}} \text{``} \limit_{\ell \in W^\bfc}^{\Xi_{\alpha_{n}^{\ast}}^{\bfc,-}}( p_{\ell}(\alpha_{n}^{\ast})) \Vdash_{\dot{\bbQ}_{\alpha_{n}^{\ast}}} \int_{K^\bfc} \frac{\vert \{ \ell \in I_{k} \colon p_{\ell}(\alpha_{n}^{\ast}) \in \dot{G}(\alpha_{n}^{\ast}) \} \vert}{\vert I_{k} \vert}  d  \dot{\Xi}_{\alpha_{n}^{\ast}+1}^{\bfc}(k) \geq 1 - \varp^*(n) \text{''}. $$ 

    Since $\bbP_{\alpha_{n}^{\ast}} \subsetdot  \bbP_{\pi},$ by integral absoluteness (\autoref{t102}), it follows that $\bbP_{\pi}$ forces the same with respect to $\Xi_{\pi}^{\bfc}.$ 
    Finally, since $q  {\rest}  (\alpha_{n}^{\ast}+1) \forces_{\alpha_{n}^{\ast}+1}$``$q(\alpha_{n}^{\ast}) \in \dot{G}(\alpha_{n}^{\ast})$'', we can conclude that
    $$ q  {\rest}  (\alpha_{n}^{\ast} +1) \forces_{\pi} \text{``} \int_{K^\bfc} \frac{\vert \{ \ell \in I_{k} \colon p_{\ell}(\alpha_{n}^{\ast}) \in \dot{G}(\alpha_{n}^{\ast}) \} \vert}{\vert I_{k} \vert}  d  \dot{\Xi}_{\pi}^{\bfc}(k) \geq 1 - \varp^{\ast}(n) \text{''}.$$
    
    \noindent\ref{k60b}: Let $n \in n^{\ast} \menos r^{\ast}$ and $\varp > 0.$ Since $\tau_n$ follows $g^{\bfc}$, by hypothesis and \autoref{k28}~\ref{k28fb}, we obtain:
    \[\Vdash_{\bbP_{\pi}} \text{``} \ \dot{\Xi}_{\pi}^{\bfc} \left( \left \{ k \in K^\bfc \colon \frac{\vert \{ \ell \in I_{k} \colon p_{\ell}(\alpha_{n, \ell}) \in \dot{G}(\alpha_{n, \ell}) \} \vert}{\vert I_{k} \vert}  \geq (1 - \varp^{\ast}(n)) (1 - \varp) \right \} \right) = 1 \text{''}.\]
    Let $G \subseteq \bbP_{\pi}$ be a generic set over $V$. Working in $V[G],$  define $\upsilon_{n} \colon K^\bfc \to \RRR$ such that, for $k \in K^\bfc,$ 
    $\upsilon_{n}(k) \coloneqq \frac{\vert \{ \ell \in I_{k} \colon p_{\ell}(\alpha_{n, \ell}) \in G(\alpha_{n, \ell}) \} \vert}{ \vert I_{k} \vert},$
    and consider $K_{\varp, n} \coloneqq \{ k \in K^\bfc \colon \upsilon_{n}(k) \geq (1 - \varp^{\ast}(n))(1 - \varp) \}.$ By the above, $\Xi_{\pi}^{\bfc}(K_{\varp, n}) = 1$. Thus,  
    \begin{equation*}
        \int_{K^\bfc} \upsilon_{n} d  \Xi_{\pi}^{\bfc}  = \int_{K_{\varp, n}} \upsilon_{n} d  \Xi_{\pi}^{\bfc}  \geq \int_{K_{\varp, n}} (1 - \varp^{\ast}(n))(1 - \varp)  d \Xi_{\pi}^{\bfc} = (1 - \varp^{\ast}(n))(1 - \varp). 
    \end{equation*}
   
    \noindent\ref{k60c}: Let $G \subseteq \bbP_{\pi}$ be a generic filter over $V$ such that $q \in G.$ Working in $V[G],$ define the functions $\varrho$,  $\varrho_{n} \colon K^\bfc \to \RRR$ for any $n < n^{\ast}$ such that, for $k \in K^\bfc,$ 
    $$ \varrho(k) \coloneqq \frac{  \vert \{ \ell \in I_{k} \colon p_{\ell} \notin G \} \vert }{\vert I_{k} \vert} \text{ and } \varrho_{n}(k) \coloneqq \frac{ \vert  \{ \ell \in I_{k} \colon p_{\ell}(\alpha_{n, \ell}) \notin G(\alpha_{n, \ell}) \} \vert }{\vert I_{k} \vert}. $$

    For $p \in \bbP_{\pi},$ $p \in G \Leftrightarrow \forall \alpha \in \dom(p)\  (p(\alpha) \in G(\alpha)),$ hence we can relate $\varrho$ with $\varrho_{n}$ by $\varrho(k)\leq \sum_{n < n^{\ast}} \varrho_{n}(k)$. 
    Also, by~\ref{k60a} and~\ref{k60b}, it is clear that $\int_{K^\bfc} \varrho_{n} d \Xi_{\pi}^{\bfc}  \leq \varp^{\ast}(n)$ for any $n < n^{\ast}$ and therefore, by basic properties of integration, 
    $$ \int_{K^\bfc} \varrho  d \Xi_{\pi}^{\bfc} \leq \sum_{n < n^{\ast}} \left( \int_{K^\bfc} \varrho_{n} \, d \Xi_{\pi}^{\bfc} \right) \leq \sum_{n < n^{\ast}} \varp^{\ast}(n),$$ 
    which implies that 
    $$\int_{K^\bfc} (1 - \varrho)  d \Xi_{\pi}^{\bfc} = 1 - \int_{K^\bfc} \varrho  d \Xi_{\pi}^{\bfc} \geq 1 - \sum_{n < n^{\ast}} \varp^{\ast}(n). $$ 

    \noindent\ref{k60d}: Let $G \subseteq \bbP_{\pi}$ be a generic filter over $V$ such that $q \in G.$ Working in $V[G],$ consider $\varrho \colon K^{\bfc} \to \RRR$ as in the proof of \ref{k60c}, let $\upsilon \coloneqq  1- \varrho$ and $0 < \varp < 1 - \sum_{n < n^{\ast}} \varp^{\ast}(n).$  Thus, $A_{\varp} = \{k \in K^{\bfc} \colon \upsilon(k) > \varp \}.$ By \ref{k60c}, we have 
    \begin{equation*}
        \begin{split}
             1 - \sum_{n < n^{\ast}} \varp^{\ast}(n)  & \leq \int_{K^{\bfc}} \upsilon d  \Xi_{\pi}^{\bfc}  \leq \int_{K^{\bfc} \menos A_{\varp}} \varp d \Xi_{\pi}^{\bfc} + \int_{A_{\varp}} \upsilon  d  \Xi_{\pi}^{\bfc}  \leq  \varp + \Xi_{\pi}^{\bfc}(A_{\varp})(1 - \varp).
        \end{split}
    \end{equation*}

    Thus,  
    $\displaystyle \Xi^{\bfc}_\pi(A_\varp) \geq \frac{1-\sum_{n<n^*}\varp^*(n)-\varp}{1-\varp} = 1 - \frac{1}{1-\varp}\sum_{n<n^*}\varp^*(n).$
\end{PROOF}

As a first application, we show that the fam-limit guarantees the $\Fr$-Knaster property of a $\calC$ iteration under certain conditions. 
This means that, morally, dominating reals are not added along the iteration.

\begin{theorem}\label{k62}
    Let $\bft$ be a $\kappa$-$\calC$-iteration. If $\calG \coloneqq \set{g^\bfc}{\Xi^\bfc \text{ is a free fam on }\calP(\omega),\ \bfc\in \calC}$ is a complete set of guardrails for $\pi$ over $\kappa^-$ and $\theta \geq \kappa$ is uncountable regular, then the final step $\bbP_{\pi}$ is $\theta$-$\Fr$-Knaster. 
\end{theorem}

\begin{PROOF}[\textbf{Proof}]{\ref{k62}}
    Pick $\bar{\varp} \colon \omega \to (0, 1)$ such that $\sum_{n < \omega} \bar{\varp}(n) < 1$. 
    Let $\{ p_{\xi} \colon \xi < \theta \} \subseteq \bbP_{\pi}$. 
    By \autoref{k12}, for any $\xi < \kappa,$ there exists some $q_{\xi} \in  D_{\bar{\varp}}$ such that $q_{\xi} \leq p_{\xi}.$ 
    
    Since $\theta$ is uncountable regular, by \autoref{k52} we get 
    an $E \in [\theta]^{\theta}$ such that  $\{ q_{\xi} \colon \xi \in E \}$ forms a uniform $\Delta$-system 
    with parameters $(\Delta, n^{\ast}, \vec{\alpha}, r^{\ast}, \zeta^{\ast},  \varp^{\ast}).$ Since $q_\xi\in D_{\bar \varp}$ for any $\xi \in E$, 
    we can get this uniform $\Delta$-system such that,  for the decreasing enumeration $\dom(q_{\xi}) \cap F^{\bft} = \{ \gamma_{n, \xi} \colon n < n^{-}\}$ (note that $n^-$ does not depend on $\xi$) we have that, whenever $n'<n^-$ and $\alpha_{n,\xi}= \gamma_{n',\xi}$,
    \begin{equation}\label{k62-e5}
        \varp^{\ast}(n) = \bar\varp(n').
    \end{equation} 
    Now, we show that $\{ p_{\xi} \colon \xi \in E \}$ is $\Fr$-linked in $\bbP_{\pi}.$ For this, let $\bar{p} \coloneqq \langle p_{\xi_{\ell}} \colon \ell < \omega \rangle$ with each $\xi_{\ell} \in E$ and consider $\{ \beta_{\ell} \colon \ell < \omega \} \subseteq \{ \xi_{\ell} \colon \ell < \omega \}$ increasing with order type $\omega.$ Define $\bar{q} \coloneqq \langle q_{\beta_{\ell}} \colon \ell < \omega \rangle.$ It is clear that it is a countable uniform $\Delta$-system 
    and therefore, by virtue of the completeness of $\calG$ and \autoref{k54}, we can find $g  \coloneqq g^\bfc \in  \calG$ such that $\bar{q}$ is a uniform $\Delta$-system following $g$ with the same parameters. 
    Without loss of generality, we may assume that $\bar I^\bfc$ is a partition of $\omega$.

    By \eqref{k62-e5}, we have that: 
    $\displaystyle 0<1 - \sum_{n < \omega} \bar{\varp}(n) \leq 1 - \sum_{n < n^{\ast}} \varp^{\ast}(n).$     
    Also, by  \autoref{k60}~\ref{k60c},
    $$\limit^\bfc(\bar q)\Vdash_\pi \text{``}  \int_\omega \frac{\vert \{ \ell \in I_{k} \colon q_{\beta_{\ell}} \in \dot{G} \} \vert}{\vert I_{k} \vert} d\dot\Xi^\bfc_\pi \geq 1 - \sum_{n < n^{\ast}} \varp^{\ast}(n) > 0 \text{''}.$$
    Since $\Xi^\bfc_0 = \Xi^\bfc$ is free, $\dot \Xi^\bfc_\pi$ is forced to be free and hence $\limit^\bfc(\bar q)$ forces that $q_{\beta_\ell} \in \dot G_\pi$ for infinitely many $\ell$, which implies that $p_{\xi_\ell}\in\dot G_\pi$ for infinitely many $\ell$.     
\end{PROOF}

\section{Relational systems and preservation theory}\label{sec:pres}

We revisit the preservation theory of unbounded families presented in~\cite[Sect.~4]{CM}, which is a generalization of Judah's and Shelah's~\cite{JS} and Brendle's~\cite{Br} preservation theory. We also provide a reformulation of \cite[Lem.~2.7]{Sh00} that gives us sufficient conditions to force $\cov(\calN)$ singular (even with countable cofinality) after finite support iterations. 

We start reviewing some basic notation about relational systems.

\begin{definition}\label{def:relsys}
We say that $R=\la X, Y, \sqsubset\ra$ is a \textit{relational system} if it consists of two non-empty sets $X$ and $Y$ and a relation $\sqsubset$. 
\begin{enumerate}[label=(\arabic*)]
    \item A set $F\subseteq X$ is \emph{$R$-bounded} if $\exists y\in Y \forall x\in F\ (x \sqsubset y)$. 
    
    \item A set $E\subseteq Y$ is \emph{$R$-dominating} if $\forall x\in X \exists y \in E\ ( x \sqsubset y)$. 

    \item Let $M$ be a set. An object $x\in X$ is \textit{$R$-unbounded over $M$}  if $x \not \sqsubset y$ for all $y\in Y\cap M$.

    \item Let $\theta$ be cardinal and $I$ a set. A family $\set{x_i}{i\in I}\subseteq X$ is \emph{strongly $\theta$-$R$-unbounded} if $|I|\geq\theta$ and, for any $y\in Y$, $|\{i\in I \colon  x_i\sqsubset y\}|<\theta$.
\end{enumerate}

We associate two cardinal characteristics with this relational system $R$: 
\begin{itemize}
    \item[{}] $\bfrak(R) \coloneqq \min\{|F| \colon  F\subseteq X  \text{ is }R\text{-unbounded}\}$, the \emph{unbounding number of $R$}, and
    
    \item[{}] $\dfrak(R) \coloneqq \min\{|D| \colon  D\subseteq Y \text{ is } R\text{-dominating}\}$, the \emph{dominating number of $R$}.
\end{itemize}
\end{definition}

We also define the \emph{dual} $R^\perp \coloneqq \la Y,X,\sqsubset^\perp\ra$ where $y \sqsubset^\perp x$ means $x \not\sqsubset y$. Note that $\bfrak(R^\perp) = \dfrak(R)$ and $\dfrak(R^\perp) = \bfrak(R)$.

Given another relational system $R'=\la X',Y',\sqsubset'\ra$,  say that a pair $(\Psi_-,\Psi_+)\colon R\to R'$ is a \emph{Tukey connection from $R$ into $R'$} if 
 $\Psi_-\colon X\to X'$ and $\Psi_+\colon Y'\to Y$ are functions such that  $\forall x\in X \forall y'\in Y'\ (\Psi_-(x) \sqsubset' y' \Rightarrow x \sqsubset \Psi_+(y'))$. Say that $R$ is \emph{Tukey below} $R'$, denoted by $R\leqT R'$, if there is a Tukey connection from $R$ to $R'$. 
 Say that $R$ is \emph{Tukey equivalent} to $R'$, denoted by $R\eqT R'$, if $R\leqT R'$ and $R'\leqT R$. It is well-known that $R\leqT R'$ implies $\bfrak(R')\leq\bfrak(R)$ and $\dfrak(R)\leq\dfrak(R')$. Hence, $R\eqT R'$ implies $\bfrak(R')=\bfrak(R)$ and $\dfrak(R)=\dfrak(R')$. 

\begin{example}\label{exm:Iwf}
For $\Iwf\subseteq\pts(X)$, define the relational systems: 
\begin{enumerate}[label=(\arabic*)]
    \item $\Iwf \coloneqq \la\Iwf,\Iwf,\subseteq\ra$, which is a directed preorder when $\Iwf$ is closed under unions (e.g.\ an ideal).
    
    \item $\Cv_\Iwf \coloneqq\la X,\Iwf,\in\ra$.
\end{enumerate}
\end{example}

The cardinal characteristics associated with an ideal can be characterized by the preceding relational systems as follows:

\begin{fact}
    If $\Iwf$ is an ideal on $\pts(X)$ containing $[X]^{<\aleph_0}$, then:
    \begin{enumerate}[label= \rm (\alph*)]
        \item $\bfrak(\Iwf)=\add(\Iwf)=\min\set{|\Jwf|}{\Jwf\subseteq\Iwf,\ \bigcup\Jwf\notin\Iwf}$. 
        
        \item $\dfrak(\Iwf)=\cof(\Iwf)=\min\set{|\Jwf|}{\Jwf\subseteq\Iwf,\ \forall A\in\Iwf \exists B\in \Jwf\  (A\subseteq B)}$.
        
        \item $\dfrak(\Cv_\Iwf)=\cov(\Iwf)=\min\set{|\Jwf|}{\Jwf\subseteq\Iwf,\ \bigcup\Jwf=X}$. 
        
        \item $\bfrak(\Cv_\Iwf)=\non(\Iwf)=\min\set{|A|}{A\subseteq X,\ A\notin\Iwf}$. 
    \end{enumerate}
\end{fact}

The existence of strongly unbounded families is equivalent to a Tukey-inequality.

\begin{lemma}[{\cite[Lem.~1.16]{CM22}}]\label{c65}
    Let $R=\la X,Y,\sqsubset\ra$ be a relational system, $\theta$ be an infinite cardinal, and let $I$ be a set of size ${\geq}\theta$.
    \begin{enumerate}[label = \normalfont (\alph*)]
        \item $\Cv_{[I]^{<\theta}}\leqT R$ iff there exists a strongly $\theta$-$R$-unbounded family $\set{x_i}{i\in I}$.

        \item $\bfrak(R)\geq\theta$ iff $R\leqT \Cv_{[X]^{<\theta}}$.
    \end{enumerate}
\end{lemma}

We focus our attention on the following types of nice relational systems.

\begin{definition}\label{def:Prs}
    We say that $R=\langle X,Y,\sqsubset\rangle$ is a \textit{Polish relational system (Prs)} if
    \begin{enumerate}[label=\rm(\arabic*)]
        \item $X$ is a Perfect Polish space,
        
        \item $Y$ is analytic in some Polish space $Z$, and
        
        \item\label{def:Prsc} $\sqsubset=\bigcup_{n<\omega}\sqsubset_{n}$ where $\langle\sqsubset_{n}: n<\omega\rangle$  is some increasing sequence of closed subsets of $X\times Z$ such that, for any $n<\omega$ and for any $y\in Y$,
        $(\sqsubset_{n})^{y}=\set{x\in X}{x\sqsubset_{n}y }$ is closed nowhere dense in $X$.
    \end{enumerate}
\end{definition}

\begin{remark}\label{Prsremark}
    By~\autoref{def:Prs}~\ref{def:Prsc}, $\la X,\Mcal(X),\in\ra \leqT R$ where $\Mwf(X)$ denotes the $\sigma$-ideal of meager subsets of $X$. Therefore, $\bfrak(R)\leq \non(\Mcal)$ and $\cov(\Mcal)\leq\dfrak(R)$.
\end{remark}

For the rest of this section, fix a Prs $R=\langle X,Y,\sqsubset\rangle$ and an infinite cardinal $\theta$. 

\begin{definition}[Judah and Shelah {\cite{JS}}, Brendle~{\cite{Br}}]\label{def:good}
    A forcing notion $\Por$ is \textit{$\theta$-$R$-good} if, for any $\Por$-name $\dot{h}$ for a member of $Y$, there is a non-empty set $H\subseteq Y$ (in the ground model) of size ${<}\theta$ such that, for any $x\in X$, if $x$ is $R$-unbounded over  $H$ then $\Vdash \text{``}x\nsqsubset \dot{h}$''.
    
    We say that $\Por$ is \textit{$R$-good} if it is $\aleph_1$-$R$-good.
\end{definition}

\begin{remark}
    Note that $\theta<\theta'$ implies that any $\theta$-$R$-good poset is $\theta'$-$R$-good. Also,  
    ``$\theta$-$R$-good'' is a hereditary forcing property.
\end{remark}

The previous property ensures that strongly unbounded families are preserved in generic extensions.  

\begin{lemma}[{\cite[Lemma~4.7~(b)]{CM}}]\label{mainpres}
    Let $\kappa$ and $\theta$ be infinite cardinals such that $\kappa\leq\cf(\theta)$. Then, any $\kappa$-$R$-good forcing notion preserving $\theta$ as a cardinal preserves all the strongly $\theta$-$R$-unbounded families from the ground model.
\end{lemma}

An iteration of good forcing notion ensures the existence of strong unbounded families, i.e.\ Tukey connections that allow to calculated values of cardinal invariants.

\begin{theorem}[{\cite[Thm.~4.11]{BCM2}}]\label{Comgood}
    Let $\theta\geq\aleph_1$ be a regular cardinal and let $\la \Por_\xi,\Qnm_\xi \colon\allowbreak  \xi<\pi\ra$ be a finite support iteration such that  $\Por_\xi$ forces that $\Qnm_\xi$ is a non-trivial $\theta$-cc $\theta$-$R$-good poset for $\xi<\pi$.  
    If 
    $\pi\geq\theta$ then $\Por_\pi$ forces $\Cv_{[\pi]^{<\theta}}\leqT R$, which implies $\bfrak(R)\leq\theta$ and $|\pi|\leq\dfrak(R)$. 
\end{theorem}

We now provide some examples of good forcing notions. In general, ``small'' forcing notions are automatically good.

\begin{lemma}[{\cite[Lemma~4.10]{CM}}]\label{smallgoodness}
    If $\theta$ is a regular cardinal then any poset of size ${<}\theta$
    is $\theta$-$R$-good. In particular, Cohen forcing $\Cor$ is $R$-good.
\end{lemma}

\begin{example}\label{ExmPrs}
    The following are Prs that describe the cardinal characteristics of Cicho\'n's diagram.
    \begin{enumerate}[label= \rm (\arabic*)]
        \item\label{ExmPrsa} Consider the Polish relational system $\Ed \coloneqq \la{}^\omega \omega,{}^\omega \omega,\neq^\infty\ra$ where $x =^\infty y$ means that $x(n)= y(n)$ for infinitely many $n$.  
        By~\cite[Thm.~2.4.1 \& Thm.~2.4.7]{BJ} (see also~\cite[Thm.~5.3]{CMlocalc}),  $\bfrak(\Ed)=\non(\Mwf)$ and $\dfrak(\Ed)=\cov(\Mwf)$.

        \item\label{ExmPrsb}  The relational system ${}^\omega \omega \coloneqq \la{}^\omega \omega,{}^\omega \omega,\leq^*\ra$ is Polish. Any $\mu$-$\mathrm{Fr}$-linked poset is $\mu^+$-${}^\omega \omega$-good (see~\cite[Thm. 3.30]{mejiavert}). In particular, any $\mu$-$\rm{FAM}$-linked poset is also $\mu^+$-\,${}^\omega \omega$-good by~\autoref{i13}.
    
        \item\label{ExmPrsd}   For each $k<\omega$, let $\id^k:\omega\to\omega$ such that $\id^k(i)=i^k$ for all $i<\omega$ and $\Hcal\coloneqq\largeset{\id^{k+1}}{k<\omega}$. Let $\Lc^*\coloneqq\la\baire, \Scal(\omega, \Hcal), \in^*\ra$ be the Polish relational system where 
        \[\Swf(\omega, \Hcal)\coloneqq\set{\varphi\colon \omega\to[\omega]^{<\aleph_0}}{\exists h\in\Hcal \forall i<\omega\ (|\varphi(i)|\leq h(i))},\]
        and $x\in^*\varphi$ iff $\forall^\infty n\  (x(n)\in\varphi(n))$. As a consequence of~\cite[Thm.~2.3.9]{BJ} (see also~\cite[Thm.~4.2]{CMlocalc}), $\bfrak(\Lc^*)=\add(\Nwf)$ and $\dfrak(\Lc^*)=\cof(\Nwf)$.

        Any $\mu$-centered forcing notion is $\mu^+$-$\Lc^*$-good (see~\cite{Br,JS}) so, in particular, $\sigma$-centered posets are $\Lc^*$-good. Besides,  Kamburelis~\cite{Ka} showed that any Boolean algebra with a strictly positive finitely additive measure is $\Lc^*$-good (in particular, any sub-algebra of random forcing).
  \end{enumerate}
\end{example}

Recall that any Lebesgue measure zero set is contained in some $G_\delta$ measure zero set, hence the latter type of set forms a basis of $\Ncal$. For this reason, coding Lebesgue measure zero sets comes from coding $G_\delta$ sets, or equivalently, coding $F_\sigma$ sets. We use the coding from~\cite{Sh00} of $F_\sigma$ measure one sets, which is essential in the proof of the consistency of $\cf(\cov(\Ncal))= \omega$.

\begin{definition}\label{a40}
    Let $\bar\delta = \la \delta_j\colon j<\omega\ra\in{}^\omega(0,1)$. Define $\Omega_{\bar\delta}$ as the set of sequences $\bar{a} \colon \omega \to \omega \times [\finseq]^{< \omega}$  such that, for any $j < \omega,$ $\bar{a}(j) =(n_{j}, a_{j})$ satisfies that $\langle n_{j} \colon j < \omega \rangle$ is increasing and, for any $j < \omega,$  $a_{j} \subseteq {}^{n_{j}}2$ and 
    $\displaystyle\frac{\vert a_{j} \vert}{2^{n_{j}}} \geq 1 - \delta_j.$ 

    For any $m<\omega$, let $\Tree_{m}(\bar{a})$ be the well-pruned subtree of ${}^{<\omega}2$ such that  
    $$[\Tree_m(\bar a)]=\bigcap_{n_j\geq m}\bigcup_{t\in a_j}[t]=\bigcap_{n_{j} \geq m} \{ x \in \cantor \colon x \rest n_{j} \in a_{j} \},$$ 
    which is a closed set of measure ${\geq}1-\sum_{n_j\geq m}\delta_j$.
    \begin{enumerate}[label=\rm(\arabic*)]
       \item For $\bar{a} \in \Omega_{\bar\delta}$ let $N[\bar{a}] \coloneqq \cantor \menos \bigcup_{m < \omega} [\Tree_{m}(\bar{a})]$, which is a $G_\delta$ set in $\cantor$, and it has measure zero when $\sum_{j<\omega}\delta_j<\infty$.
    
        \item Define the Prs $\Cn_{\bar\delta} \coloneqq \langle \Omega_{\bar\delta}, {}^{\omega} 2, \lhd \rangle$ where $\bar{a} \lhd x$ means $x \notin N[\bar{a}].$
    \end{enumerate}
    
    When $\sum_{j<\omega}\delta_j<\infty$, any null set in $\cantor$ is a subset of $N[\bar{a}]$ for some $\bar{a} \in \Omega_{\bar\delta}$ (see e.g.~\cite[Lem.~4.8]{CMlocalc}), which implies $\Cn_{\bar\delta} \eqT \Cv_{\calN}^{\perp}$. Hence, $\bfrak(\Cn_{\bar\delta})=\cov(\Nwf)$ and $\dfrak(\Cn_{\bar\delta})=\non(\Nwf)$.
\end{definition}

This section concludes by showing our main preservation result for the covering of $\Nwf$, which significantly reformulates Shelah's result~\cite[Lem.~2.7]{Sh00}. 

Below, we introduce the notion \emph{$\theta$-anti-Bendixson family}, which is motivated by Shelah’s result. To this end,  we present the following auxiliary relational system used to force $\cov(\Ncal)$ singular after finite support iterations.

\begin{definition}\label{a50}
    Define $\calT\coloneqq\{ T \colon T \text{ is a well-pruned subtree of }{}^{<\omega} 2 \text{ of height } \omega\}$, $\calP$ as the collection of perfect subtrees of ${}^{<\omega} 2$, and the relational system $\Be \coloneqq \langle \calT, \calP, \supseteq \rangle$. 
    Notice that $\set{T_i}{i\in I}$ is strongly $\theta$-$\Be$-unbounded iff, for any $J\in[I]^{\theta}$, $\bigcap_{i\in J}[T_i]$ does not contain a perfect set, i.e.\ $\bigcap_{i\in J}[T_i]$ is a countable closed set (by the Cantor-Bendixson Theorem). For this reason, we use the name \emph{$\theta$-anti-Bendixson} for the strongly $\theta$-$\Be$-unbounded sets.
\end{definition}

Before stating the preservation theorem, we need a result on the absoluteness of trees:

\begin{lemma}\label{a58}
    Let $\bbP$ be a forcing notion and $\langle T_{i} \colon i \in I \rangle$ a sequence of well-pruned subtrees of $\finseqb.$ Let $G$ be a $\bbP$-generic filter over $V$ and assume that $\bigcap_{i \in I} [T_{i}]^{V[G]}$ is countable. Then \[\bigcap_{i \in I} [T_{i}]^{V[G]} = \bigcap_{i \in  I} 
     [T_{i}]^{V}. \]
\end{lemma}

\begin{PROOF}[\textbf{Proof}]{\ref{a58}}
    Since for any $i \in I$,  $[T_{i}]^{V} \subseteq [T_{i}]^{V[G]},$ it is clear that $\bigcap_{i \in I} [T_{i}]^{V}$ is countable, so it has an enumeration $\langle x_{n} \colon n < w \rangle \in V$ with $w \leq \omega.$ We show that there is no $x \in \bigcap_{i \in I} [T_{i}]^{V[G]}$ such that, for any $n < w$, $x \neq x_{n}.$ Suppose the contrary, so assume that there are $p \in G$ and a name $\dot{x} \in V^{\bbP}$ such that, in $V$, 
    $ p \Vdash \text{``} \dot{x} \in \bigcap_{i \in I} [T_{i}] \text{ and } \forall n < w\ (\dot{x} \neq x_{n}) \text{''}. $
    Now, work in $V.$ By induction on $k < \omega,$ define $y=\la y(k)\colon k<\omega\ra\in{}^\omega\omega$ and a decreasing sequence $\langle p_{k} \colon k < \omega \rangle$ in $\bbP$ with $p_{0} \leq p$ such that $p_{k} \forces \text{``} \dot{x}(k) = y(k) \text{''}$ and, whenever $k < w,$ there exists an $\ell_{k} < \omega$ such that, $p_{k} \forces \text{``} \dot{x}(\ell_{k}) = y(\ell_{k}) \neq x_{k}(\ell_{k})$''. Then $p_{k} \forces \text{``} y  {\rest}  (k+1) = \dot{x}  {\rest} (k+1) \in T_{i}$'' for all $i \in I,$ so $y \in \bigcap_{i \in I}[T_{i}]^{V} = \{ x_{n} \colon n < w \}.$ However, for $n < w$,  $p_{n} \forces \text{``} y(\ell_{n}) \neq x_{n}(\ell_{n})$'', hence $y \neq x_{n},$ which is a contradiction.   
\end{PROOF}

\begin{theorem}\label{a60}
    Let $\kappa\leq \lambda$ be uncountable cardinals such that $\kappa$ is regular, $\bar\delta \in{}^\omega(0,1)$ and let $I$ be a set. Let $\bbP_{\pi} = \langle \bbP_{\alpha},   \dot{\bbQ}_{\alpha} \colon \alpha < \pi \rangle$ be a finite support iteration of $\kappa$-$\rm{cc}$ forcing notions. Assume that

    \begin{enumerate}[label=\rm(\arabic*)]
        \item\label{a60a} $\{ \bar{a}^{i} \colon i\in I\}$ is strongly $\lambda$-$\Cn_{\bar\delta}$-unbounded (so $|I|\geq\lambda$),

        \item\label{a60b} $\Vdash_{\bbP_{\pi}}$``$\{ \Tree_m(\bar{a}^{i}) \colon i\in I,\ m<\omega\}$ is 
        $\kappa$-anti-Bendixson'',\footnote{This is condition $(\ast \ast)_{\bbP}$ in \cite[Lem.~2.7]{Sh00}.} and

        \item\label{a60c} for any $ \alpha < \pi$,  $\Vdash_{\alpha}$`` $\dot{\bbQ}_{\alpha}$ contains a dense subset of size ${<} \lambda$''.
    \end{enumerate}

    Then, $\Vdash_{\bbP_{\pi}}$``$\{ \bar{a}^{i} \colon i\in I\rangle$ is strongly $\lambda$-$\Cn_{\bar\delta}$-unbounded''.  As a consequence,  $\forces_{\Por_\pi}``\cov(\calN) \leq \lambda$ and $\cov([I]^{<\lambda})\leq\non(\Ncal)$'' when $\sum_{j<\omega}\delta_j<\infty$. 
\end{theorem}
\begin{PROOF}[\textbf{Proof}]{\ref{a60}}
    By induction on $\gamma\leq \pi,$ we show that 
    $\Por_\gamma$ forces that $\{ \bar{a}^{i} \colon i\in I\}$ is strongly $\lambda$-$\Cn_{\bar\delta}$-unbounded.
    We distinguish three cases: 

\noindent\textbf{Case: $\gamma= 0$.} In this case $\Vdash_{\bbP_{0}}$``$\{ \bar{a}^{i} \colon i\in I\}$ is strongly $\lambda$-$\Cn_{\bar\delta}$-unbounded'' by condition~\ref{a60a}.

\noindent\textbf{Case: $\gamma=\xi+1$.} First, work in $V^{\bbP_{\xi}}.$ By the induction hypothesis, $\{ \bar{a}^{i} \colon i\in I\}$ is strongly $\lambda$-$\Cn_{\bar\delta}$-unbounded. So we must show that $\dot{\bbQ}_{\xi}$ still forces this. Towards a contradiction, suppose that there are $q\in \bbQ_{\xi}$ and a $\bbQ_{\xi}$-name $\dot{x}$ of a real number in $\cantor$ such that, 
        \begin{equation}\label{a60-e3}
            \text{$q \Vdash_{\bbQ_{\xi}}``\vert \{ i\in I\colon \dot{x} \notin N[\bar{a}^{i}] \} \vert \geq \lambda\text{''}.$}
        \end{equation}
        Therefore, by transfinite recursion on $\zeta < \lambda,$ we can build two sequences $\langle i_{\zeta} \colon \zeta < \lambda \rangle$ and $\langle q_{\zeta} \colon \zeta < \lambda \rangle$ such that, for any $\zeta < \lambda,$
        \begin{multicols}{2}
        \begin{itemize}
            \item  $i_{\zeta} \in I,$
            
            \item $ q_{\zeta} \leq q$,
            
            \item $i_{\zeta} \notin \{ i_{\rho} \colon \rho < \zeta \},$
            
            \item $q_{\zeta} \Vdash_{\bbQ_{\xi}}$``$\dot{x} \notin N[\bar{a}^{i_{\zeta}}]$''.
        \end{itemize}
        \end{multicols}
        By the last item, for any $\zeta<\lambda$ 
        there are $m_{\zeta} < \omega$ and $q_{\zeta}' \leq q_{\zeta}$ such that $q_{\zeta}' \Vdash_{\bbQ_{\xi}}$``$\dot{x} \in [\Tree_{m_{\zeta}}(\bar{a}^{i_\zeta})]$''.
        
        Since by condition \ref{a60c} $\bbQ_{\xi}$ contains a dense subset of size ${<} \lambda,$ we can find $C_{0} \in [\lambda]^{\kappa}$ and $q' \in \bbQ_{\xi}$ such that, for any  $ \zeta \in C_{0}$,  $q' \leq q_{\zeta}'.$ 
        On the other hand, as $\kappa$ is a regular cardinal, by shrinking $C_0$ if necessary, we can assume that there is 
        some $m < \omega$ 
        such that, for each $\zeta \in C_{0}$,  $m_{\zeta} = m.$ Therefore, for any $\zeta \in C_{0},$ we have that $q' \Vdash_{\bbQ_{\xi}}$``$\dot{x} \in [\Tree_{m}(\bar{a}^{i_\zeta})]$'', and therefore $$q' \Vdash_{\bbQ_{\xi}}\text{``}\dot{x} \in \bigcap_{\zeta \in C_{0}} [\Tree_{m}(\bar{a}^{i_\zeta})]\text{''}.$$  

        Now, working in $V^{\bbP_{\xi}}[G]$ where $G$ is a $\bbQ_{\xi}$-generic set over $V^{\bbP_{\xi}}$ containing $q',$ we have by condition \ref{a60b} that $F_{1} \coloneqq \bigcap_{\zeta \in C_{0}} [\Tree_{m}(\bar{a}^{i_\zeta})]$ is a closed subset of ${}^{\omega}2$ not containing 
        a perfect set. Therefore, by Cantor-Bendixson Theorem, it is a countable set. Hence, by \autoref{a58}, we have that  $F_{1} \in V^{\bbP_{\xi}}.$ Thus $x \coloneqq \dot{x}[G] \in V^{\bbP_{\xi}}.$ However, by induction hypothesis, $\vert \{ i\in I \colon y \notin N[\bar{a}^{i}] \} \vert < \lambda,$ for any $y \in V^{\bbP_{\xi}}$,    
        which contradicts \eqref{a60-e3}.

\noindent\textbf{Case: $\gamma$ limit.} We split into two cases:

\noindent\textbf{Case 1:} $\cf(\gamma) < \kappa$.  Towards a contradiction, assume that there are $p \in \bbP_{\gamma}$ and a $\bbP_{\gamma}$-name of a real number $\dot{x}$ in ${}^{\omega}2,$ such that: 
            \begin{equation}\label{a60-e5}
                \text{$p \Vdash_{\gamma} ``\vert \{ i\in I \colon \dot{x} \notin N[\bar{a}^{i}] \} \vert \geq \lambda \text{''}.$}
            \end{equation}

            As in the successor step, we can find $C \in [\lambda]^{\kappa}$, $\{ i_\zeta\colon \zeta\in C\}\subseteq I$ one-to-one, $\{ p_{\zeta} \colon \zeta \in C \} \subseteq \bbP_{\gamma}$ and some $m<\omega$ such that, for any $\zeta \in C$, $p_{\zeta} \leq p$ and $p_{\zeta} \Vdash_{\gamma}``\dot{x} \in \Tree_m(\bar a^{i_\zeta})$''.  
            
            Since $\cf(\gamma) < \kappa,$ there exists a set $L \subseteq \gamma$ cofinal in $\gamma$ with $\vert L \vert < \kappa.$ So, for each $\zeta \in C,$ there exists some $  \xi_{\zeta} \in L,$ such that $p_{\zeta} \in \bbP_{\xi_{\zeta}}.$             
            Now $\vert C \vert =\kappa >  \vert L \vert,$ so we can find a set $C_{0} \subseteq C$ with $\vert C_{0} \vert = \kappa$ and some $\xi \in L$ such that, for any $\zeta \in C_{0}, \, p_{\zeta} \in \bbP_{\xi}.$ 

            On the other hand, since $\Por_\xi$ has the $\kappa$-cc, we can find some $p'\leq p$ in $\Por_\xi$ such that 
            \begin{equation*}\label{a60-e7}
               p'\Vdash_\xi |\{\zeta\in C_0\colon p_\zeta\in \dot G\}|=\kappa.   
            \end{equation*}
            Let $G$ be a $\bbP_{\gamma}$-generic filter over $V$ with $p' \in G.$ In $V[G],$ $C_{1} \coloneqq \{\zeta\in C_0\colon p_\zeta\in G\}\in V[\Por_\xi\cap G]$ and $\vert C_{1} \vert = \kappa.$  Therefore, by condition~\ref{a60b}, $F_{2} \coloneqq \bigcap_{\zeta \in C_{1}} [\Tree_{m}(\bar{a}^{i_\zeta})]^{V[G]}$ is a countable closed subset of ${}^{\omega}2.$ Then, by \autoref{a58}, $F_{2} \subseteq V[\bbP_{\xi}\cap G]$. Then $x \coloneqq \dot{x}[G] \in F_{2}.$ However, by induction hypothesis, since $\xi < \gamma,$ and $x\in V[\Por_\xi\cap G]$, $\vert \{ i\in I \colon  x \notin N[\bar{a}^{i}] \} \vert < \lambda,$ which contradicts \eqref{a60-e5}.
    
            \textbf{Case 2:} $\cf(\gamma) \geq \kappa$. Let $\dot{x}$ a nice $\bbP_{\gamma}$-name of a real number in $\cantor.$ Since $\bbP_{\gamma}$ has the $\kappa$-$\rm{cc}$ and $\cf(\gamma) \geq \kappa,$ we can find a $\xi < \gamma$ such that $\dot{x}$ is a $\bbP_{\xi}$-name. Therefore, by induction hypothesis, $\Vdash_{\bbP_{\xi}}``\vert \{ i\in I \colon \dot{x} \notin N[\bar{a}^{i}]\} \vert < \lambda$''. Thus, $\bbP_{\gamma}$ forces the same. 
        
            We proved $\forces_{\bbP_{\pi}}$``$\{ \bar{a}^{i} \colon i\in I \rangle$ is strongly $\lambda$-$\Cn$-unbounded''. Thus, by \autoref{c65},  
            we conclude that $\forces_{\bbP_\pi}$``$\Cv_{[I]^{< \lambda}} \preceq_{\mathrm{T}} \Cn_{\bar\delta}$''.
\end{PROOF}

The first two authors proved that forcing iterations with finitely additive measures preserves $\non(\Ecal)$ small (see \autoref{a69}). We close this section by proving a goodness result related to this cardinal invariant. 

\begin{definition}\label{e086}    
    Define the relational system $\Ce\coloneqq \la {}^{\omega}2,\NE,\sqce\ra$ where $\NE$ is the collection of sequences $\bar T = \langle T_n \colon n<\omega \rangle$ such that each $T_n$ is a subtree of ${}^{<\omega}2$ (not necessarily well-pruned), $T_n\subseteq T_{n+1}$ and $\Leb([ T_n]) =0$, i.e.\ $\displaystyle \lim_{n\to\infty}\frac{|T\cap {}^n2|}{2^n}=0$, and $x\sqce \bar T$ iff $x\in [T_n]$ for some $n<\omega$. 
\end{definition}

Notice that $\Ce$ is a Polish relational system and $\Ce\eqT \Cv_\Ecal$.

In some cases, $(\Xi,\bar I,\varp_0)$-linkedness is good for $\Ce$. Recall that, if $T$ is a tree and $t \in T$, then $T\wedge t \coloneqq \set{s\in T}{s \text{ is comparable with $t$ in $T$}}$.

\begin{lemma}\label{e087}
    Assume that $(\Xi,\bar I)$ is a free pair, $\Por$ is a forcing notion, $\varp_0\in[0,1)$ and $Q\subseteq\Por$ is $(\Xi,\bar I,\varp_0)$-linked. If $\dot T$ is a $\Por$-name of a subtree of ${}^{<\omega}2$ such that $\Vdash\text{``\/} \Leb([\dot T])=0$'', then there is some subtree $T'\subseteq{}^{<\omega}2$ such that $\Leb([T'])=0$ and, for any $x\in {}^\omega2\menos[T']$, no $q\in Q$ forces $x\in [\dot T]$.
\end{lemma}
\begin{PROOF}{\ref{e087}}
    Fix $\varp_1$ such that $\varp_0<\varp_1<1$. Assume that $Q$ is $(\Xi,\bar I,\varp_0)$-linked in $\Por$. By \autoref{e106}, there is some free fam $\Xi'$ on $\pts(\omega)$ 
    such that $Q$ is $(\Xi',\bar I',\varp_1)$-linked, where $\bar I'\in \Ibf_\omega^\fin$ such that $|I'_k|=2^k$ for all $k<\omega$. 
    Pick some $\varp'>0$ such that $\varp_1+\varp'<1$. 
    
    Let $\dot T$ be a $\Por$-name of a subtree of ${}^{<\omega}2$ such that $[\dot T]\in\Ncal$. 
    Define 
    \[T'\coloneqq \setcolon{s\in {}^{<\omega}2}{\exists q\in Q\ (q\Vdash\text{``}s\in\dot T\text{''})},\] 
    which is clearly a subtree of ${}^{<\omega}2$. It is enough to show that $[T']\in\Ncal$. Assume the contrary. Hence, by Lebesgue Density Theorem, there is some $s\in T'$ such that $\Leb([s]\cap [T'])>(1-\varp')\Leb([s])$, that is, for all $k<\omega$,
    \begin{equation}\label{e087eq1}
        \frac{|\Lev_{|s|+k}(T'\wedge s)|}{2^k}>1-\varp'.
    \end{equation}    
    For each $t\in\Lev_{|s|+k}({}^{<\omega}2\wedge s)$, if $t\in T'$ choose $q_t\in Q$ forcing that $t\in \dot T$, otherwise $q_t$ can be any member of $Q$. Thus, since $Q$ is $(\Xi',\bar I',\varp_1)$-linked,
    \[\limit \bar q \Vdash \int_\omega \frac{|\setcolon{t\in\Lev_{|s|+k}({}^{<\omega}2\wedge s)}{q_t\in\dot G}|}{2^k} d\dot \Xi^*(k) \geq 1-\varp_1\]
    where $\limit$ and $\dot\Xi^*$ witness linkedness. Then, by~\eqref{e087eq1}, 
    \[\limit \bar q \Vdash \int_\omega \frac{|\Lev_{|s|+k}(\dot T\wedge s)|}{2^k} d\dot \Xi^*(k) \geq 1-\varp_1-\varp'.\]
    Since the function inside the integral is monotone decreasing and $\Xi^*$ is forced to be free, we get that $\displaystyle\limit \bar q \Vdash\text{``}\frac{|\Lev_{|s|+k}(\dot T\wedge s)|}{2^k} \geq 1-\varp_1-\varp'$'' for all $k<\omega$, which implies that $\limit \bar q\Vdash\text{``} \Leb([s]\cap [\dot T])\geq 1-\varp_1-\varp'>0$'', a contradiction.
\end{PROOF}

\begin{theorem}\label{e090}
    Let $\calY\subseteq\calY_*$ and assume that it contains some free pair. If $\mu$ is an infinite cardinal, then any $\mu$-$\calY$-linked poset is $\mu^+$-$\Ce$-good. 
\end{theorem}
\begin{PROOF}[\textbf{Proof}]{\ref{e090}}
    Assume that $\bbP$ is $\mu$-$\calY$-linked witnessed by $\langle Q_{\alpha, \varp} \colon\alpha < \mu \conj \varp \in (0, 1)_\QQQ \rangle$. 

    Let $\dot{\bar{T}}$ such that $\forces_{\bbP}$``$\dot{\bar{T}} \in \NE$'', that is, $\bbP$ forces that $\dot{\bar{T}} = \langle \dot{T}_{n} \colon n < \omega \rangle$, where each $\dot{T}_{n}$ is a $\bbP$-name of a subtree of $\finseq$ such that $\dot{T_{n}} \subseteq \dot{T}_{n+ 1}$ and, for any $n < \omega$, $\Leb([\dot{T}_{n}]) = 0.$ Given $n < \omega,$ since $\bbP$ forces that $\Leb([\dot{T}_{n}]) = 0,$ by \autoref{e087}, for any $\alpha < \mu$ and $\varp \in (0, 1)_\QQQ$,  there exists some subtree $T^{\alpha, \varp}_{n}$ of ${}^{<\omega}2$ such that $\Leb([T^{\alpha, \varp}_{n}]) = 0$ and, for any $x \in \cantor \setminus [T^{\alpha, \varp}_{n}]$, there is no $q \in Q_{\alpha, \varp}$ forcing $x \in [\dot{T}_{n}].$ Concretely, set $T^{\alpha,\varp}_n\coloneqq\setcolon{s\in {}^{<\omega}2}{\exists q\in Q_{\alpha,\varp}\ (q\Vdash\text{``}s\in\dot T_n\text{''})}$. Notice that $T^{\alpha,\varp}_n\subseteq T^{\alpha,\varp}_{n+1}$.

    Now, consider the set  $H \coloneqq \{ \bar{T}^{\alpha, \varp} \colon \alpha < \mu \conj \varp \in (0, 1)_{\QQQ} \}$, where for any $\alpha < \mu$ and $\varp \in (0, 1)_{\QQQ}$, $\bar{T}^{\alpha, \varp} \coloneqq \langle T^{\alpha, \varp}_{n} \colon n < \omega \rangle$. It is clear that $H \subseteq \NE$ and $\vert H \vert < \mu^{+}$. Let us show that $H$ is a witness of $\mu^{+}$-$\Ce$-goodness for $\bbP$. Assume that $x \in \cantor$ is $\Ce$-unbounded over $H$, that is, for any $\alpha < \mu$ and $\varp \in (0, 1)_{\QQQ}$, $x  \not \sqsubset^{\Ce} \bar{T}^{\alpha, \varp}$. Towards contradiction, suppose that there exists some $p \in \bbP$ such that $p \forces$\text{``}$x \sqsubset \dot{\bar{T}}$''. Therefore, we can find $m < \omega$, $\beta < \mu$, $\epsilon \in (0, 1)_\QQQ$ and $q \in \bbP$ such that $q \leq p$, $q \in Q_{\beta, \epsilon}$ and $q \forces$``$x \in [\dot{T}_{m}]$''. On the other hand, since $x \not \sqsubset^{\Ce} \bar{T}^{\beta, \epsilon}$, in particular we have that $x \in \cantor \setminus [\dot{T}_{m}^{\beta, \epsilon}]$, which contradicts \autoref{e087} because $q \forces$``$x \in [\dot{T}_{m}]$''. 
\end{PROOF}

\section{Preserving anti-Bendixson and strongly unbounded families: the role of the fams}
\label{presanti}

In the previous section, we proved in \autoref{a60} sufficient conditions to keep $\cov(\calN)$ ``small'' and singular. In general, building finite support iterations that satisfy the conditions~\ref{a60a} and~\ref{a60c} from \autoref{a60} is not a problem.\footnote{Because we can iterate Cohen forcing to add an anti-Bendixson family.} However, for the condition~\ref{a60b}, we need to guarantee that we can preserve anti-Bendixson families, which is a non-trivial task. For this purpose, we need to appeal to the structure of the iterations using fams (see \autoref{a64}). On the other hand, as we mentioned in the introduction, to force $\cov(\calN)$ with countable cofinality, we need to keep $\gb$ ``small'', which will also be an effect of the fams along the iterations (see \autoref{a68}). 
In addition, we show how iterations with fams do not increase $\non(\Ecal)$.

In this section, fix a sequence $\bar\delta = \la\delta_j \colon j<\omega\ra\in {}^\omega (0,1)$. 
Since Cohen forcing is forcing-equivalent to any countable atomless poset, we can use the following form of Cohen forcing: $\bbC$ is the poset whose conditions are finite sequences of the form $\langle (n_j, a_j) \colon j < k \rangle$ satisfying the following requirements:
\begin{itemize}
    \item $\langle n_j \colon j < k \rangle$ is an increasing sequence of natural numbers, 

    \item for any $j < k$,  $a_j \subseteq {}^{n_j}2$ and $\frac{\vert a_j \vert}{2^{n_j}} \geq 1 - \delta_j.$ 
\end{itemize}   
We order $\bbC$ by end extension. If $G$ is a $\Cor$-generic, in $V[G]$, define $\bar a_G\coloneqq \bigcup G$, which is in $\Omega_{\bar\delta}$. Notice that $N[\bar a_G]\in\Nwf$ whenever $\sum_{j<\omega}\delta_j<\infty$.

\begin{notation}\label{k40.1}
    For $0<j<\omega$ define the function $h^*_j\colon \omega\to \omega\menos\{0\}$ by
    \[h^*_j(k) \coloneqq \left\{
         \begin{array}{ll}
             \binom{2^k}{2^{k-j}} & \text{ if $k\geq j$,} \\[1ex]
             1 & \text{otherwise.}
         \end{array}
      \right.\]
      Note that $h^*_j$ diverges to infinity.
\end{notation}

For the rest of this section, fix an iteration bedrock $\calC$. 
Now, we show how to preserve anti-Bendixson families in $\calC$-iterations.

\begin{theorem}\label{a64}
    Let $\kappa \leq \theta$ be uncountable cardinals with $\theta$ regular. Assume that $\bft$ is a $\kappa$-$\calC$-iteration of length $\pi$ and, for any $0<j<\omega$,
    \[\calG_{j}\coloneqq\set{g^\bfc}{\bfc\in \calC,\ \Xi^\bfc \text{ is a free fam on $\pts(\omega)$ and }\forall k<\omega\ (|I^\bfc_k|=h^*_j(k))}\]
    is a complete set of guardrails over $\kappa^-$.     
    Further assume that $L\subseteq U^{\bft}$, $|L|\geq\theta$ and, for $\alpha\in L$, $\dot{\bbQ}_{\alpha}$ is a $\Por^-_\alpha$-name of $\bbC$ and $Q_{s, 0} \coloneqq \{ s \}$ for any $s \in \bbC$. Then, 
    $$\forces_{\pi} \text{``} \{ \Tree_{m}(\bar{a}^{\alpha}) \colon \alpha \in L,\ m<\omega \} \text{ is } \theta\text{-anti-Bendixson}\text{''},$$ where each $\bar{a}^{\alpha}$ is the Cohen real added by $\dot{\bbQ}_{\alpha}$.   
\end{theorem}

\begin{PROOF}[\textbf{Proof}]{\ref{a64}}
    Towards contradiction assume that there are a condition $p \in \bbP_{\pi},$ a $\bbP_{\pi}$-name $\dot{T}$ of a perfect subtree of $\finseq$ and $m < \omega$ such that $p \forces_{\pi}$``$ \vert \dot{E} \vert \geq \theta$'', where 
    $$\dot{E} \coloneqq \{ \alpha \in L \colon \dot{T} \subseteq \Tree_{m}(\bar{a}^{\alpha})\}.$$
    
    \textbf{Part 1}: Get some suitable $0<j^*<\omega$, $g\in\calG_{j^*}$ and a countable uniform $\Delta$-system following $g$. 

    Pick some $\bar{\varp} \colon \omega \to (0,1)$ such that $\sum_{n < \omega} \bar{\varp}(n) < 1$. 
     By transfinite recursion, we can build   sequences $\{ p_{\xi} \colon \xi < \theta \} \subseteq D_{\bar{\varp}}$ and $\{ \alpha_{\xi} \colon \xi < \theta \} \subseteq L,$ such that, for $\xi<\theta$

     \begin{multicols}{2}
         \begin{itemize}
             \item $p_{\xi} \leq p,$ 
    
             \item $\langle \alpha_{\xi} \colon \xi < \theta \rangle$ has no repetitions,
    
             \item $\alpha_{\xi} \in \dom(p_{\xi}),$
    
             \item $p_{\xi} \forces_{\pi} \text{``} \alpha_{\xi} \in \dot{E}$''.  
         \end{itemize}
     \end{multicols}
     By \autoref{k52} we can get some $H \in [\theta]^{\theta}$ such that $\{ p_{\xi} \colon \xi \in H \}$ forms a uniform $\Delta$-system with parameters $(\Delta, n^{\ast}, \vec{\alpha}, r^{\ast}, \zeta^{\ast}, \varp^{\ast})$. Since $p_\xi\in D_{\bar \varp}$ for $\xi\in H$, for the decreasing enumeration $\dom(q_{\xi}) \cap F^{\bft} = \{ \gamma_{n, \xi} \colon n < n^{-}\}$ (note that $n^-$ does not depend on $\xi$), we can assume that, whenever $n'<n^-$ and $\alpha_{n,\xi}= \gamma_{n',\xi}$, $\varp^{\ast}(n) = \bar\varp(n').$ 
    In addition, by shrinking $H$ and increasing each $p_\xi(\alpha_\xi)$ if necessary, we can assume that
    \begin{enumerate}
         \item there exists some $s^{\ast} = \langle (n_{i}, a_{i}) \colon \ell < m^{\ast} \rangle \in \bbC$ such that, for any $\xi \in H$,  $p_{\xi}(\alpha_{\xi}) = s^{\ast}$, and

         \item there is some $c^{\ast} < n^{\ast}$ such that, for any $\xi \in H$,  $\alpha_{\xi} = \alpha_{c^{\ast}, \xi},$ that is, all $\alpha_{\xi}$'s come from the same column of the $\Delta$-system. Notice that $c^{\ast} \notin r^{\ast}$.
    \end{enumerate} 
    Pick some $j^*<\omega$ such that $m\leq j^*$, $\sup_{i<m^*} n_i<j^*$ and $2^{-j^*}\leq \delta_{m^*}$, and denote $h\coloneqq h_{j^*}$. 
    Consider the interval partition $I=\la I_k\colon k\in\omega\ra$ of $\omega$ such that $|I_k|=h(k)$ and, for $k\geq j^*$, enumerate the subsets of ${}^{k} 2$ of size $2^{k}(1-2^{-j^{\ast}})$ as $\{ a_{\ell} \colon \ell \in I_{k} \}.$ In this case, for $\ell\in I_k$,
    $$ \frac{\vert a_{\ell} \vert}{2^{k}}  = 1 - 2^{-j^{\ast}} \geq 1 - \delta_{m^*}.$$ 

    Choose some increasing $\langle \beta_{\ell} \colon \ell < \omega \rangle \subseteq H$. For $k < \omega$ and $\ell \in I_{k},$ we define a condition $p_{\beta_{\ell}}' \in \bbP_{\pi}$ such that $p_{\beta_{\ell}}' \leq p_{\beta_{\ell}},$ $\dom(p_{\beta_{\ell}}') = \dom(p_{\beta_{\ell}})$ and,\footnote{In the third case $h(k)=1$, so $I_k$ is a singleton.}
    $$ p_{\beta_{\ell}}'(\gamma) 
        \coloneqq \left\{ \begin{array}{ll}
             p_{\beta_{\ell}}(\gamma) &   \text{if $\gamma \neq \beta_{\ell}$}, \\[1ex]
             \concat{s^{\ast}}{\langle (k, a_{\ell}) \rangle} &  \text{if   $\gamma = \beta_{\ell}$ and $k\geq j^*$,}\\[1ex]
             \concat{s^*}{\la (j^*,{}^{j^*} 2) \ra} & \text{if   $\gamma = \beta_{\ell}$ and $k< j^*$.}
        \end{array}\right. $$

    Notice that $\bar{p}' \coloneqq \langle p_{\beta_{\ell}}' \colon \ell < \omega \rangle$ stills forms a uniform $\Delta$-system. As $\calG_{j^*}$ is a complete set of guardrails for $\pi$ over $\kappa^-$, by virtue of \autoref{k54} there exists some $\bfc \in \calC$ such that $\bar{p}'$ is a uniform $\Delta$-system following $g\coloneqq g^\bfc$, $\Xi^\bfc$ is a free fam on $\pts(\omega)$ and $|I^\bfc_k|=h(k)$ for all $k<\omega$. Without loss of generality, we can assume that $\bar I^\bfc = \bar I$.
   
    \textbf{Part 2}: Get a contradiction. 
    
    Let $G$ be a $\bbP_{\pi}$-generic over $V$ with $\lim^{\bfc}(\bar{p}') \in G,$ $T\coloneqq \dot T[G]$, and work in $V[G]$. 
    Define, for any $k < \omega$,  $b_{k} \coloneqq \{ \ell \in I_{k} \colon p_{\beta_{\ell}}' \in G \}.$ 
    Then, by \autoref{k60}~\ref{k60c}, 
    $$\int_\omega \frac{|b_k|}{|I_k|}d\Xi^\bfc_\pi(k) \geq 1-\sum_{n<n^*}\varp^*(n)\geq 1-\sum_{n<\omega}\bar\varp(n)>0.$$

    Now, for $k\geq j^*$ and $\ell\in b_k$, $\Lev_{k}(T) \subseteq  a_{\ell}.$ Indeed, by the definition of $p_{\beta_{\ell}}'$,  $p_{\beta_{\ell}}'(\beta_{\ell}) = \concat{s^{\ast}}{\langle (k, a_{\ell}) \rangle} \in G(\beta_{\ell}).$ Since $\beta_{\ell} \in H,$ it follows that 
    $T \subseteq \Tree_{m}(\bar{a}^{\beta_{\ell}})$. 
    On the other hand, $ \bar{a}^{\beta_{\ell}} = \langle (n_{i}^{\beta_{\ell}}, a_{i}^{\beta_{\ell}}) \colon i < \omega \rangle \supseteq \concat{s^{\ast}}{\langle (k, a_{\ell}) \rangle}$, so $n^{\beta_{\ell}}_{m^*} = k$ and $a^{\beta_\ell}_{m^*}= a_\ell$. Therefore, 
    $$ [T] \subseteq \{ x \in \cantor \colon x  {\rest}  n_{m^{\ast}}^{\beta_{\ell}} \in a_{m^{\ast}}^{\beta_{\ell}} \} = \{ x \in \cantor \colon x  {\rest}  k \in a_{\ell} \}$$ 
    because $k\geq j^*\geq m$, so $\Lev_{k}(T)  \subseteq a_{\ell}.$ 

    As a consequence,  for $k\geq j^*$: 
    \begin{equation*}
        \begin{split}
            \vert b_{k} \vert  & = \vert \{ a_{\ell} \colon \ell \in b_{k} \} \vert \leq \vert \{ a_{\ell} \colon \ell \in I_{k}, \, \Lev_{k}(T) \subseteq a_{\ell} \} \vert\\
            & = \{ a \subseteq {}^{k}2 \colon  \Lev_{k}(T) \subseteq a,  \vert a \vert = 2^{k} (1-2^{-j^{\ast}})   \}\\
            & = \vert \{ a \subseteq {}^{k}2 \menos (\Lev_{k}(T)) \colon \vert a \vert =  2^{k-j^*}  \} \vert \leq \binom{2^{k} - \bfm_k}{2^{k-j^*}}
        \end{split}
    \end{equation*}
    where 
    $\bfm_{k} \coloneqq \min \{ \vert \Lev_{k}(T) \vert, 2^{k}(1- 2^{j^{\ast}}) \}.$ 
    Hence,  
    \begin{equation*}
        \begin{split}
            \frac{\vert b_{k} \vert}{\vert I_{k} \vert} & \leq 
            \frac{ \binom{ 2^{k} - \bfm_{k} }{ 2^{k-j^{\ast}} }  } {  \binom{ 2^{k} }{ 2^{k-j^{\ast}}}} = \frac{\prod_{i < \bfm_{k}}(2^{k} - 2^{k-j^{\ast}} - i)}{\prod_{i < \bfm_{k}} (2^{k} - i)} 
            = \prod_{i < \bfm_{k}} \left( 1 - \frac{2^{k-j^{\ast}}}{2^{k} - i} \right) \\
            & \leq \prod_{i < \bfm_{k}} \left( 1 - \frac{2^{k-j^{\ast}}}{2^{k}} \right) = (1 - 2^{-j^{\ast}})^{\bfm_{k}}.
        \end{split}
    \end{equation*}
    Therefore, since $\Xi^\bfc_\pi$ is free,
    $$\int_\omega(1-2^{-j^*})^{\bfm_k}d\Xi^\bfc_\pi(k)\geq \int_\omega \frac{|b_k|}{|I_k|}d\Xi^\bfc_\pi(k)\geq 1-\sum_{n<\omega}\bar\varp(n),$$
    so $(1-2^{-j^*})^{\bfm_k}\geq 1-\sum_{n<\omega}\bar\varp(n)$ for all $k\geq j^*$ because the sequence $\la(1-2^{-j^*})^{\bfm_k}\colon k<\omega\ra$ is monotone-decreasing and $\Xi^\bfc_\pi$ is a free fam. Since $1-\sum_{n<\omega}\bar\varp(n)>0$, we have that $\la\mbf_k\colon k<\omega\ra$ is bounded and, as a consequence, $\la|\Lev_{k}(T)|\colon k<\omega\ra$ is bounded, which implies that $[T]$ is finite. This contradicts that $T$ is a perfect tree.
\end{PROOF}

Now, we prove a similar result that will allow us to preserve strongly-$\theta$-$\,{}^{\omega} \omega$-unbounded families (added by iterating using Cohen forcing), which is the key point to keep $\gb$ ``small''. Again, the structure of the iteration using finitely additive measures plays a fundamental role.

\begin{theorem}\label{a68}
    Let $\kappa \leq \theta$ be uncountable cardinals such that $\theta$ is regular. Assume that $\bft$ is a $\kappa$-$\calC$-iteration of length $\pi$ 
    and that $\calG \coloneqq \set{g^\bfc}{\Xi^\bfc \text{ is a free fam on }\calP(\omega),\ \bfc\in \calC}$ is a complete set of guardrails over $\kappa^-$.

    Further assume that $L\subseteq U^{\bft}$, $|L|\geq \theta$ and,
    for $\alpha \in L,$ $\dot{\bbQ}_{\alpha}$ is a $\Por^-_\alpha$-name of $\bbC$ and $Q_{s, 0} \coloneqq \{ s \}$ for all $s \in \bbC$. Then, $\forces_{\pi} \text{``} \langle \bar{n}^{\alpha} \colon \alpha \in L \rangle \text{ is strongly $\theta$-\,$\baire$-unbounded''},$
    where $\bar{a}^{\alpha} = \langle (n_{\ell}^{\alpha},  a_{\ell}^{\alpha}) \colon \ell < \omega \rangle$ is the Cohen real added by $\dot{\bbQ}_{\alpha}$ at the step $\alpha$ of the iteration, and $\bar{n}^{\alpha} \coloneqq \langle n_{\ell}^{\alpha} \colon \ell < \omega \rangle$. As a consequence, $\Por_\pi$ forces $\Cv_{[L]^{<\theta}}\leqT \baire$, i.e.\ $\bfrak\leq\theta$ and $|L|\leq\dfrak$.
\end{theorem}

\begin{PROOF}[\textbf{Proof}]{\ref{a68}}

Assume, towards a contradiction, that there are a nice $\bbP_{\pi}$-name $\dot{x}$ of a real number in $\baire$, $m<\omega$ and a condition $p \in \bbP_{\pi}$  such that $p \forces_{\pi}$``$ \vert \{ \alpha \in L \colon \forall i\geq m\ (n^{\alpha}_i \leq \dot{x}(i)) \} \vert \geq \theta$''.  To get a contradiction, we are going to find a suitable guardrail $g$ for $\pi$ over $\kappa^-$ and a countable uniform $\Delta$-system following $g$.

Let $\bar{\varp} \colon \omega \to (0,1)$ such that $\sum_{n < \omega} \bar{\varp}(n) < 1$. As in the proof of \autoref{a64}, we can find a uniform $\Delta$-system $\{ p_{\xi} \colon \xi < \theta \} \subseteq D_{\bar{\varp}}$ with parameters $(\Delta, n^{\ast}, \vec{\alpha}, r^{\ast}, \zeta^{\ast}, \varp^{\ast})$, $\{ \alpha_{\xi} \colon \xi < \theta \} \subseteq L$, $s^{\ast} = \langle (n_{i}, a_{i}) \colon \ell < m^{\ast} \rangle \in \bbC$ and $c^{\ast} \in n^{\ast}\menos r^\ast$ such that, for any $\xi < \theta$,
\begin{multicols}{2}
     \begin{itemize}
         \item $p_{\xi} \leq p,$ 

         \item $\langle \alpha_{\xi} \colon \xi < \theta \rangle$ has no repetitions,

         \item $\alpha_{\xi} \in \dom(p_{\xi}),$
         
         \item $p_{\xi} \forces_{\pi}$``$\forall i \geq m\  (n_{i}^{\alpha_\xi} \leq \dot{x}(i))$'',

         \item for the decreasing enumeration\\ $\dom(q_{\xi}) \cap F^{\bft} = \{ \gamma_{n, \xi} \colon n < n^{-}\}$, whenever $n'<n^-$ and $\alpha_{n,\xi}= \gamma_{n',\xi}$, $\varp^{\ast}(n) = \bar\varp(n')$,
         
        \item $p_{\xi}(\alpha_{\xi}) = s^{\ast}$ and $\alpha_{\xi} = \alpha_{c^{\ast}, \xi}$.
     \end{itemize}
\end{multicols}
    By extending $s^*$ if necessary, we can assume that $m\leq m^*$. 
    Let $\ell^*\coloneqq\sup_{i<m^*}\{n_i +1\}$ and
    pick some increasing $\langle \beta_{\ell} \colon \ell < \omega \rangle \subseteq \theta$. For any $\ell < \omega,$  we define a condition $p_{\beta_{\ell}}' \in \bbP_{\pi}$ such that $p_{\beta_{\ell}}' \leq p_{\beta_{\ell}}$, $\dom(p_{\beta_{\ell}}') = \dom(p_{\beta_{\ell}})$ and, 
    $$ p_{\beta_{\ell}}'(\gamma) 
        \coloneqq \left\{ \begin{array}{ll}
             p_{\beta_{\ell}}(\gamma) &   \text{if $\gamma \neq \beta_{\ell}$}, \\[1ex]
             \concat{s^{\ast}}{\langle (\ell, {}^\ell 2) \rangle} &  \text{if   $\gamma = \beta_{\ell}$ and $\ell\geq \ell^*$,}\\[1ex]
             \concat{s^*}{\la (\ell^*,{}^{\ell^*} 2) \ra} & \text{if   $\gamma = \beta_{\ell}$ and $\ell< \ell^*$.}
             \end{array}\right. $$
    The sequence $\bar{p}' \coloneqq \langle p_{\beta_{\ell}}' \colon \ell < \omega \rangle$ is a uniform $\Delta$-system, so we can find some $\bfc\in\calC$ such that this $\Delta$-system follows $g\coloneqq g^\bfc$ and $\Xi^\bfc$ is a free fam on $\pts(\omega)$ (wlog, we can assume that $\bar I^\bfc$ is a partition of $\omega$ for all $\bfc\in\calC$ such that $g^\bfc\in\calG$).

    Let $G$ be a $\Por_\pi$-generic set over $V$ containing $\lim^\bfc(\bar p')$. By \autoref{k60}~\ref{k60c},
    \[\int_\omega \frac{|\{\ell \in I^\bfc_k\colon p'_{\beta_\ell}\in G\}|}{|I_k^{\bfc}|}d\Xi^\bfc_\pi(k)\geq 1 - \sum_{n < n^{\ast}} \varp^{\ast}(n) \geq 1 - \sum_{n < \omega} \bar{\varp}(n) >0,\]
    which implies that $L\coloneqq\set{\ell<\omega}{p'_{\beta_\ell}\in G}$ is infinite. For $\ell\in L\menos \ell^*$, $p_{\beta_{\ell}}'(\beta_{\ell}) \in G(\beta)$, so we have that $n_{m^{\ast}}^{\beta_{\ell}} = \ell$. On the other hand,  we know that, for any $j \geq  m^{\ast}$, $n^{\beta_\ell}_{j}\leq x(j)$, hence $\ell \leq x(m^{\ast})$ because $m^\ast\geq m$. However, since $L$ is infinite, it follows that $\{ \ell < \omega \colon \ell \leq x(m^{\ast}) \}$ is also infinite, which is a contradiction.
    
    Now, $\langle \bar{n}^{\alpha} \colon \alpha < \theta \rangle$ is strongly $\theta$-$\baire$-unbounded in the final generic extension, so, by \autoref{c65}, we have that $\calC_{[L]^{< \theta}} \preceq_{\rm{T}} \baire.$
\end{PROOF}

Combining the previous argument with the idea of the proof of \autoref{e087}, we show that iterations with fams do not increase $\non(\Ecal)$.

\begin{theorem}\label{a69}
    Let $\kappa \leq \theta$ be uncountable cardinals such that $\theta$ is regular. Assume that $\bft$ is a $\kappa$-$\calC$-iteration of length $\pi$ 
    and that $\calG \coloneqq \set{g^\bfc}{\Xi^\bfc \text{ is a free fam on }\calP(\omega),\ \forall k<\omega\ (|I_k|=2^k),\ \bfc\in \calC}$ is a complete set of guardrails over $\kappa^-$.

    Further assume that $L\subseteq U^{\bft}$, $|L|\geq \theta$ and,
    for $\alpha \in L,$ $\dot{\bbQ}_{\alpha}$ is a $\Por^-_\alpha$-name of $\bbC={}^{<\omega}2$ and $Q_{s, 0} \coloneqq \{ s \}$ for all $s \in \bbC$. Then, $\forces_{\pi} \text{``} \langle \dot{c}^{\alpha} \colon \alpha \in L \rangle \text{ is strongly $\theta$-\,$\Ce$-unbounded''},$
    where $\dot{c}^{\alpha}$ is the Cohen real added by $\dot{\bbQ}_{\alpha}$ at the step $\alpha$ of the iteration. As a consequence, $\Por_\pi$ forces $\Cv_{[L]^{<\theta}}\leqT \Ce$, i.e.\ $\non(\Ecal)\leq\theta$ and $|L|\leq\cov(\Ecal)$.
\end{theorem}

\begin{PROOF}[\textbf{Proof}]{\ref{a69}}
Wlog we can assume that $I^\bfc_k=I_k\coloneqq [2^k-1,2^{k+1}-1)$ for all $\bfc\in\calC$ such that $g^\bfc\in\calG$.

Assume that there are a sequence of $\bbP_{\pi}$-name $\Seq{\dot{T}_n}{n<\omega}$ forced to form a sequence in $\NE$, $m<\omega$, and a condition $p \in \bbP_{\pi}$  forcing $ \vert \{ \alpha \in L \colon \dot c^{\alpha}\in [\dot T_m] \} \vert \geq \theta$. We imitate the proof of \autoref{e087} to reach a contradiction.

Let $\bar{\varp} \colon \omega \to (0,1)$ such that $\sum_{n < \omega} \bar{\varp}(n) < 1$. As in the previous proofs, we can find a uniform $\Delta$-system $\{ p_{\xi} \colon \xi < \theta \} \subseteq D_{\bar{\varp}}$ with parameters $(\Delta, n^{\ast}, \vec{\alpha}, r^{\ast}, \zeta^{\ast}, \varp^{\ast})$, $\{ \alpha_{\xi} \colon \xi < \theta \} \subseteq L$, $s^{\ast}\in {}^{<\omega}2$ and $c^{\ast} \in n^{\ast}\menos r^\ast$ such that, for any $\xi < \theta$,
\begin{multicols}{2}
     \begin{itemize}
         \item $p_{\xi} \leq p,$ 

         \item $\langle \alpha_{\xi} \colon \xi < \theta \rangle$ has no repetitions,

         \item $\alpha_{\xi} \in \dom(p_{\xi}),$
         
         \item $p_{\xi} \forces_{\pi}$``$\dot c^{\alpha_\xi} \in[\dot T_m]$'',

         \item for the decreasing enumeration\\ $\dom(q_{\xi}) \cap F^{\bft} = \{ \gamma_{n, \xi} \colon n < n^{-}\}$, whenever $n'<n^-$ and $\alpha_{n,\xi}= \gamma_{n',\xi}$, $\varp^{\ast}(n) = \bar\varp(n')$,
         
        \item $p_{\xi}(\alpha_{\xi}) = s^{\ast}$ and $\alpha_{\xi} = \alpha_{c^{\ast}, \xi}$.
     \end{itemize}
\end{multicols}
    Pick some increasing $\langle \beta_{\ell} \colon \ell<\omega \rangle \subseteq \theta$ and some bijection $h\colon \omega\to {}^{<\omega}2$ sending each interval $I_k$ onto ${}^k2$. For each $\ell<\omega$, define a condition $p_{\beta_{\ell}}' \in \bbP_{\pi}$ such that $p_{\beta_{\ell}}' \leq p_{\beta_{\ell}}$, $\dom(p_{\beta_{\ell}}') = \dom(p_{\beta_{\ell}})$ and, 
    $$ p_{\beta_{\ell}}'(\gamma) 
        \coloneqq \left\{ \begin{array}{ll}
             p_{\beta_{\ell}}(\gamma) &   \text{if $\gamma \neq \beta_{\ell}$}, \\[1ex]
             \concat{s^{\ast}}{h(\ell)} &  \text{if   $\gamma = \beta_{\ell}$.}
             \end{array}\right. $$
    The sequence $\bar{p}' \coloneqq \langle p_{\beta_{\ell}}' \colon \ell < \omega \rangle$ is a uniform $\Delta$-system, so we can find some $\bfc\in\calC$ such that this $\Delta$-system follows $g\coloneqq g^\bfc$ and $\Xi^\bfc$ is a free fam on $\pts(\omega)$.

    Let $G$ be a $\Por_\pi$-generic set over $V$ containing $\lim^\bfc(\bar p')$. By \autoref{k60}~\ref{k60c},
    \[\int_\omega \frac{|\{\ell \in I^\bfc_k\colon p'_{\beta_\ell}\in G\}|}{|I_k^{\bfc}|}d\Xi^\bfc_\pi(k)\geq 1 - \sum_{n < n^{\ast}} \varp^{\ast}(n) \geq 1 - \sum_{n < \omega} \bar{\varp}(n) >0.\]
    For each $k<\omega$, $|I^\bfc_k|=2^k$ and $p'_{\beta_\ell}\in G$ implies that $\concat{s^*}{h(\ell)}\in T_m\coloneqq \dot T_m[G]$ (because $c^{\alpha_{\beta_{\ell}}}\in [T_m]$). Therefore,
    $\displaystyle\int_\omega \frac{\vert\Lev_{|s^*|+k}(T\wedge s^*) \vert}{2^k}\geq 1- \sum_{n < n^{\ast}} \varp^{\ast}(n)$,
    which implies that 
    \[\displaystyle\frac{|\Lev_{|s^*|+k}(T\wedge s^*)|}{2^k}\geq 1- \sum_{n < n^{\ast}} \varp^{\ast}(n)\] 
    for all $k<\omega$ because the term inside the integral is decreasing. Thus, $\displaystyle\frac{\Leb([T])}{\Leb([s^*])}\geq 1- \sum_{n < n^{\ast}} \varp^{\ast}(n)>0$, a contradiction.
\end{PROOF}

\begin{remark}
    The further assumption ``$L\subseteq U^{\bft}$ and, for $\alpha\in L$, $\dot{\bbQ}_{\alpha}$ is a $\Por^-_\alpha$-name of $\bbC$ and $Q_{s, 0} \coloneqq \{ s \}$ for any $s \in \bbC$'' is not required in the previous theorems. Since any finite support iteration of non-trivial posets adds a Cohen real at limit stages, for any $\alpha<\pi$ limit or $0$ we can pick a Cohen real $\dot{\bar a}^\alpha$ in $V_{\alpha+\omega}$ over $V_\alpha$, and modify the proof to show that, in the final extension, 
    $\{ \Tree_m(\bar{a}^{\alpha}) \colon \alpha < \pi \text{ limit or $0$},\ m<\omega \}$ is $\theta$-anti-Bendixsion for \autoref{a64}, and similarly for the other theorems. However, the proof of this more general version is a bit longer.

    Moreover, in \autoref{a64}, instead of using the countably many functions $h^*_j$, it is possible to use a single function, e.g.\ $h^*(k)\coloneqq\prod_{j=1}^k\binom{2^k}{2^{k-j}}$, and assume instead that the following is a complete set of guardrails for $\pi$ over $\kappa^-$:
    \[\calG\coloneqq\set{g^\bfc}{\bfc\in \calC,\ \Xi^\bfc \text{ is a fam on $\pts(\omega)$ and }\forall k<\omega\ (|I^\bfc_k|=h^*(k))}.\]
    However, we decided to stick to $\la h^*_j\colon 0<j<\omega\ra$ for a clearer proof.
\end{remark}

\section{Applications: singular values in Cicho\'n's diagram}\label{appl}

This section aims to provide applications of our results described in the previous sections.  In particular, to prove~\autoref{Mthm:sepcd} and~\ref{Mthm:sepcd2}, i.e.\ a new separation of the left-hand side of Cicho\'n's diagram with the covering of null possibly singular,  we construct a $\kappa$-$\Fr$-Knaster ccc forcing notion by iterating with fams.

Before engaging in the applications, we present some legal ways to construct $\Por^-_\alpha$ in an iteration with fams. 

\begin{lemma}\label{a70}
    Let $\bbP$ be a \emph{ccc} forcing notion, $Q$ a subset of $\bbP$ and let $\langle \dot{F}_{i} \colon i \in I \rangle$ be a sequence of $\bbP$-names such that, for any $i \in I$, $\forces_{\bbP}$``$\dot{F}_{i} \colon \calP(\omega) \to \RRR$''. Then there exists a forcing notion $\bbP^{-}$ such that $\bbP^{-} \subsetdot \bbP,$ $Q \subseteq \bbP^{-},$  $\vert \bbP^{-} \vert \leq \max \{  2, \vert Q \vert, \vert I \vert \}^{\aleph_{0}},$ and, for any $i \in I$, there is a $\Por^-$-name for $\dot{F}_{i}  {\rest} \calP(\omega) \cap V^{\bbP^{-}}$.
\end{lemma}
\begin{PROOF}[\textbf{Proof}]{\ref{a70}}
    For any nice $\Por$-name $\dot{x}$ of a member of $\calP(\omega)$ and any $i \in I,$ define $E_{i}(\dot{x})$ as a nice $\bbP$-name of $F_{i}(\dot{x})$. For a large enough regular cardinal $\chi,$ we can use a generalization of \cite[Lem.~III.8.4]{Kunen} to find a model $N \prec \calH(\chi)$ such that ${}^{\omega} N \subseteq N$,  $\bbP$, $Q \in N$, $Q \cup \{ \dot{F}_{i} \colon i \in I \} \subseteq N$, $N$ is closed under $E_{i}$ for any $i \in I$ and $\vert N \vert \leq \{ 2, \vert Q \vert, \vert I \vert \}^{\aleph_{0}}.$ We define $\bbP^{-} \coloneqq \bbP \cap N$ as a suborder of $\bbP$. We show that $\bbP^{-}$ is as required. 

    First, it is clear that, for any $p, p' \in \bbP^{-},$ $p \leq_{\bbP^{-}} p'$ implies $p \leq_{\bbP} p'$, and $p \perp_{\bbP^{-}} p'$ implies $p  \perp_{\bbP} p'.$ Therefore, to prove that $\bbP^{-}$ is a complete subposet of $\bbP$, it is enough to show that any maximal antichain in $\bbP^{-}$ is a maximal antichain in $\bbP$. For this, let $A$ be a maximal antichain in $\bbP^{-}$. It is clear that $A \subseteq \bbP$ is an antichain in $\bbP,$ and since $\bbP$ is $\rm{ccc}$,  $\vert A \vert \leq \aleph_{0}.$ As a consequence, $A \in N$ because $A \subseteq N$ and $N$ is closed under countable sequences. 
    For the maximality, we have,  
        \begin{equation*}
            \begin{split}
                A \text{ is maximal in } \bbP^{-} 
                & \Leftrightarrow \forall p\in \bbP \cap N \exists q\in A\ (q \parallel_{\bbP^{-}} p)\\
                 &\Leftrightarrow N \models \text{``} \forall p\in \bbP  \exists q\in A \ (q \parallel_{\bbP} p)\text{''} \\
                & \Leftrightarrow \calH(\chi) \models \text{``} \forall p \in\bbP \exists q\in A\ (q \parallel_{\bbP} p) \text{''}\\  
                &\Leftrightarrow \forall p\in \bbP \exists q\in A\ (q \parallel_{\bbP} p) \\
                & \Leftrightarrow A \text{ is maximal in } \bbP.
            \end{split}     
        \end{equation*}   
    Hence, $\bbP^{-} \subsetdot\bbP.$ Second, since $Q \subseteq \bbP$ and $Q \subseteq N,$ we have that $Q \subseteq \bbP^{-}.$ Note that $\vert \bbP^{-} \vert = \vert \bbP \cap N \vert \leq \vert N\vert \leq \max\{ 2, |Q|, \vert I \vert \}^{\aleph_{0}}.$ 

    Finally, for any $i \in I,$ let $ \tau_{i}$ be a $\bbP^-$-name 
    of the map $\pts(\omega)\cap V^{\bbP^-}\, {\to}\ \R$ such that $\Vdash_{\bbP^-} \text{``} \tau_i(\dot x) = E_i(\dot x)$'' for any nice $\bbP^-$-name $\dot x$ of a member of $\pts(\omega)$. 
    Indeed, for such an $\dot x$, $\dot{x} \in N,$ and it is a nice $\bbP$-name, so $E_{i}(\dot{x}) \in N$ and, since $N$ is $\sigma$-closed, we have that $E_{i}(\dot{x})$ is a $\bbP^{-}$-name. Consequently, $\tau_{i}$ is also a $\bbP^{-}$-name and $\forces_{\bbP}$``$\tau_{i} = F_{i} {\rest} \calP(\omega) \cap V^{\bbP^{-}}$''.
\end{PROOF}

We will also use the following generalization of \autoref{a70}.

\begin{lemma}\label{a70g}
    Let $\mu$ and $\nu$ be infinite cardinals, $\bbP$ a \emph{ccc} forcing notion, $\bbQ$ a $\sigma$-closed $\nu^+$-cc forcing notion, $\dot Q$ a $\bbQ$-name  such that  $\Vdash_{\bbQ}$``\/$\dot Q\subseteq\bbP\textrm{\ and\ }|\dot Q|\leq\vert \mu\vert$", and let $\langle \dot{F}_{i} \colon i \in I \rangle$ be a sequence of $\bbP\times\bbQ$-names such that, for any $i \in I$, $\forces_{\bbP\times\bbQ}$``\/$\dot{F}_{i} \colon \calP(\omega) \to \RRR\textrm{"}$. Then there exists a forcing notion $\bbP^{-} \subsetdot \bbP$ such that $\vert \bbP^{-} \vert \leq \max \{  \mu, \nu, \vert I \vert \}^{\aleph_{0}},$ $\Vdash_{\bbQ} \text{``} \dot Q \subseteq \bbP^{-}$'' and, for any $i \in I$, $\dot{F}_{i}  {\rest} \calP(\omega) \cap V^{\bbP^{-}\times\bbQ}$ has a $\bbP^{-}\times\bbQ$-name. 
\end{lemma}
\begin{PROOF}[\textbf{Proof}]{\ref{a70g}}
    Since $\bbQ$ is $\sigma$-closed, $\Vdash_{\Qor}$``$\bbP$ is ccc" and $\Vdash_{\Por}$``$\bbQ$ is ${<}\aleph_1$-distributive" by Easton’s Lemma. Therefore, $\Pwf(\omega)\cap V^{\Por\times\Qor}=\Pwf(\omega)\cap V^{\bbP}$. 

    Set $\lambda\coloneqq\max\{\mu,\nu,|I|\}$. 
    By recursion on $\eta<\omega_1$, build a sequence $\la \dot\bbP_{\eta}'\colon\eta<\omega_1\ra$ of $\bbQ$-names and a sequence $\la \bbP_{\eta}^{-}\colon\eta<\omega_1\ra$ of subsets of $\Por$ such that: 
     \begin{enumerate}[label=\rm($\bullet_\arabic*$)]
    \item\label{a70:g:1} $\Vdash_{\bbQ}``\dot\bbP_{\eta}'\subsetdot\bbP$, $Q\subseteq\dot\bbP_{\eta}$, $\dot{F}_{i}  {\rest} \calP(\omega) \cap V^{\bbQ\ast \dot \bbP'_{\eta+1}}$ has a $\dot \bbP'_{\eta+1}$-name for any $i \in I$, and $|\dot\Por_{\eta}^{'}|\leq |\lambda|^{\aleph_0}"$, and 
        \item\label{a70:g:2} $\Vdash_{\bbQ}``\dot\bbP'_{\eta}\subseteq \bbP_{\eta}^-\subseteq \dot\bbP'_{\eta+1}"$ and $|\Por_{\eta}^{-}|\leq \lambda^{\aleph_0}$.
    \end{enumerate}

    These sequences are constructed as follows. Set $\Por^-_{-1}\coloneqq\dot Q$. When $\eta<\omega_1$ is not a limit ordinal,  by employing~\autoref{a70} in $V^\bbQ$, in $V$ there is a $\Qor$-name $\dot\bbP_{\eta}'$ such that $\Vdash_\bbQ``\dot\bbP_{\eta}' \subsetdot \bbP$, $\Por_{\eta-1}^{-} \subseteq \dot\bbP_{\eta}'$, $|\dot\bbP_{\eta}'|\leq|\lambda|^{\aleph_0}$ and, for any $i \in I$, there is a $\dot\bbP'_{\eta}$-name for $\dot{F}_{i}  {\rest} \calP(\omega) \cap V^{\bbQ\ast\dot\bbP'_{\eta}}"$. Since $\bbQ$ has the $\nu^+$-cc, we can find a $\bbP_{\eta}^{-} \subseteq \bbP$ such that $\Vdash_{\bbQ}``\dot\bbP_{\eta}'\subseteq\bbP_{\eta}^{-}"$ and $|\bbP_{\eta}^{-}|\leq \lambda^{\aleph_0}$ (this uses that $\Vdash_{\Por}$``$\bbQ$ is ${<}\aleph_1$-distributive").
    
    When $\eta<\omega_1$ is limit, assume we have built $\bbP_\xi^{-}$ and $\dot\bbP_\xi'$ for all $\xi<\eta$. Let $\dot\Por_\eta'$ be a $\bbQ$-name of $\bigcup_{\xi<\eta}\dot\Por_\xi'$ and let $\Por^-_\eta\coloneqq \bigcup_{\xi<\eta}\Por^-_\xi$. Observe that $\Vdash_\bbQ``\dot\Por_\eta'=\Por_\eta^-"$ because $\Vdash_{\bbQ}``\dot\bbP'_{\xi}\subseteq\dot\bbP^-_{\xi}\subseteq\dot\bbP'_{\xi+1}"$ for all $\xi<\eta$. This finishes the construction. 

    Let $\bbP^- \coloneqq \bigcup_{\eta<\omega_1}\bbP_{\eta}^{-}$. Notice that $|\bbP^-|\leq\lambda^{\aleph_0}$ and $\bbP^-\subsetdot\bbP$, the latter because $\Vdash_\bbQ``\bbP^-=\bigcup_{\eta<\omega_1}\dot\bbP'_{\eta}\subsetdot\bbP"$. It remains to show that $\Vdash_\bbQ ``\dot{F}_{i}  {\rest} \calP(\omega) \cap V^{\bbP^{-}\times\bbQ}$ has a $\bbP^{-}$-name" for $i \in I$. By \ref{a70:g:1}, choose a $\bbQ$-name $\dot F^{\eta}_i$ of a $\dot\Por'_{\eta+1}$-name of $\dot{F}_{i}  {\rest} \calP(\omega) \cap V^{\bbQ\ast\dot\bbP'_{\eta+1}}$. Since $\Qor$ forces $\dot\Por'_{\eta+1}\subsetdot \Por^-$, each $\dot F^{\eta}_i$ is forced to be a $\Por^-$-name. So let $\dot F^-_i$ be a $\bbQ$-name of a $\Por^-$-name of $\bigcup_{\eta<\omega_1}F_i^{\eta}$, which is forced to be $\dot{F}_{i}  {\rest} \calP(\omega) \cap V^{\bbP^{-}\times\bbQ}$ because $\Por^-$ is the direct limit of the $\subsetdot$-increasing sequence of ccc posets $\la\dot\Por'_\eta\colon \eta<\omega_1\ra$.
\end{PROOF}

We are finally ready to demonstrate~\autoref{Mthm:sepcd} and \ref{Mthm:sepcd2}. First, we review the forcing notions we use in our forcing constructions, for the sake of completeness.

\begin{enumerate}[label = \rm (\arabic*)]
    \item The forcing $\Loc$ is the standard $\sigma$-linked poset adding a generic slalom $\varphi_{\gen} \in \Scal(\omega,\id)\coloneqq\prod_{n<\omega}[\omega]^{\leq n}$, i.e.\ $x \in^* \varphi_{\gen}$ for any $x\in \baire$ in the ground model (see e.g.~\cite[Def.~4.1]{CMR2}). This forcing increases $\add(\Nwf)$ because $\Nwf\eqT \Lc(\omega,\id)\coloneqq \la\baire,\Scal(\omega,\id),\in^*\ra$ (Bartoszy\'nski's characterization of measure, see~\cite[Thm.~4.2]{CMlocalc}).
    
    \item 
    $\Dor$ denotes Hechler's poset for adding a dominating real over the ground model, so it increases $\bfrak$. Recall that it is $\sigma$-centered.
\end{enumerate}

\begin{theorem}\label{a95}
    Let $\theta \leq \kappa$ be uncountable regular cardinals and let $\lambda$ and $\chi$ be cardinals such that $\kappa\leq\lambda\leq \chi = \chi^{\aleph_0}$ and $\cof([\chi]^{<\nu}) = \chi$ for $\nu\in\{\theta,\kappa,\lambda\}$. Assume that either $\kappa=\lambda$, or $\lambda$ is $\aleph_{1}$-inaccessible and $\log \chi <\lambda$.
    
    Then, there exists a $\kappa$-$\Fr$-Knaster $\rm{ccc}$  forcing notion that forces: 
    \begin{enumerate}[label=\rm(\alph*)]
        \item\label{a95:a} $\add(\calN) = \theta,$ $\add(\calM) = \gb = \non(\Ecal) = \kappa,$ $\cov(\calN) = \lambda$, $\cov(\calM) =  \gc = \chi,$ and $\non(\calM) \in \{ \lambda,  \lambda^{+} \}$ (see~\autoref{f502}). 
    
        \item\label{a95:b} If $\lambda$ is regular then $\non(\Mwf)=\lambda$. 
    
        \item\label{a95:c} If $\cf(\lambda) < \kappa$ then $\non(\Mwf)=\lambda^+$, as the constellation in~\autoref{f1300}.
    \end{enumerate}
\end{theorem}     
    
\begin{figure}[ht]
    \centering
    
    \begin{tikzpicture}
        \small{
        
        \node (aleph1) at (-2,3) {$\aleph_{1}$};
        
        \node (addn) at (0,3){$\add(\calN)$};
        
        \node (covn) at (0,7){$\cov(\calN)$};
        
        \node (nonn) at (9,3) {$\non(\calN)$} ;
        
        \node (cfn) at (9,7) {$\cof(\calN)$} ;
        
        \node (addm) at (3,3) {$\add(\calM)$} ;
        
        \node (covm) at (6,3) {$\cov(\calM)$} ;
        
        \node (nonm) at (3,7) {$\non(\calM)$} ;
        
        \node (cfm) at (6,7) {$\cof(\calM)$};
        
        \node (b) at (3,5) {$\gb$};
        
        \node (d) at (6,5) {$\gd$};
        
        \node (c) at (11,7) {$\gc$};
        
        \draw (aleph1) edge[->] (addn)
        
              (addn) edge[->] (covn)
        
              (covn) edge [->] (nonm)
        
              (cfm)edge [->] (cfn)
        
              (cfn) edge[->] (c);

        \draw (addm) edge [->]  (covm) 

        (addm) edge [->] (b)

        (nonm)edge [->] (cfm)

        (d)  edge[->] (cfm);

        \draw (addn) edge [->]  (addm)
        
           (covm) edge [->]  (nonn)
        
           (nonn) edge [->]  (cfn);
        
        \draw (b)  edge [->] (nonm);
        
        \draw (covm) edge [->] (d);
        
        \draw (b) edge [->] (d);

        \draw (-3,2) -- (12,2);

        \draw (-3, 2) -- (-3, 8);

        \draw (-3, 8) -- (12, 8);

        \draw (12, 2) -- (12, 8);
        }

        \draw[redun] (4.5, 2) -- (4.5, 8);

        \draw[redun] (-3, 5.6) -- (4.5, 5.6);

        \draw[redun] (-1.4, 2) -- (-1.4, 5.6);

        \draw[redun] (1.5, 2) -- (1.5, 5.6);

        \draw[redun, dashed] (1.5, 5.6) -- (1.5, 8);


        \filldraw[color=greenun, fill=greenun!10, very thick](2.25, 6.25) circle (0.4);   

        \node[greenun] at (2.25, 6.25) {$\lambda^{+}$};      

        \filldraw[color=greenun, fill=greenun!10, very thick](7.5,5) circle (0.4);

        \node[greenun] at (7.5, 5) {$\chi$}; 

        \filldraw[color=greenun, fill=greenun!10, very thick](2.2, 4.5) circle (0.4);

        \node[greenun] at (2.2, 4.5) {$\kappa$}; 

        \filldraw[color=greenun, fill=greenun!10, very thick](-0.7, 4.5) circle (0.4);

        \node[greenun] at (-0.7, 4.5) {$\theta$}; 

        \filldraw[color=greenun, fill=greenun!10, very thick](-1.5, 7) circle (0.4);

        \node[greenun] at (-1.5, 7) {$\lambda$}; 

        \filldraw[color=greenun, fill=greenun!10, very thick](-1.5, 7) circle (0.4);

        \node[greenun] at (-1.5, 7) {$\lambda$}; 
    \end{tikzpicture}
    \caption{Five values in  Cicho\'n's diagram with $\cov(\calN)$ possibly singular.}
    \label{f502}
\end{figure}

\begin{PROOF}{\autoref{a95}}
    First assume that $\kappa<\lambda$, $\lambda$ is $\aleph_{1}$-inaccessible and $\log \chi <\lambda$. Since both $\kappa$ and $\log\chi$ are smaller than $\lambda$, there is some cardinal $\mu<\lambda$ such that $\kappa\leq\mu$ and $\chi\leq 2^\mu$. Set $\pi\coloneqq \chi+\chi$ (ordinal sum) and choose a complete set of guardrails $\calG$ for $\pi$ over $\mu$  such that $\vert \calG \vert \leq \mu^{\aleph_{0}}$, which exists by~\autoref{k23}. Note that $\mu^{\aleph_{0}} <\lambda$ because $\lambda$ is $\aleph_1$-inaccessible.

    For each $0<j<\omega$, let $\bar I^j$ be the interval partition of $\omega$ such that $|I^j_k|=h_j^*(k)$ for all $k<\omega$ (see~\autoref{k40.1}). Fix a free fam $\Xi_0$ on $\pts(\omega)$, and define the bedrock $\calC$ as the set of all tuples of the form $(\Xi_0,\bar I^j,g)$ for $0<j<\omega$ and $g\in\calG$. Then, $|\calC|\leq \mu^{\aleph_0}<\lambda$.

    We construct a $\kappa$-$\calC$-iteration $\bft$ of length $\pi$ forcing our desired statement. 

    Consider a partition $\{S_0,S_1, S_2, S_3\}$ of $\pi$ such that each $S_i$ has size $\chi$ and $S_0=\chi$. For each $1\leq i\leq 3$, fix a bijection $\bfb_{i} \colon S_i \to S_i \times \chi$ such that $\bfb_{i}(\xi) = (\alpha, \beta)$ implies $\alpha \leq \xi$, which we call a \emph{book-keeping function}. 
    Set $U^{\bft}\coloneqq S_0\cup S_1\cup S_2$ and $F^{\bft}\coloneqq S_3$. We also fix a sequence $\bar\delta = \la\delta_j \colon j<\omega\ra \in {}^\omega(0,1)$ such that $\sum_{j<\omega}\delta_j<\infty$.
    
    We now construct $\bft$ by transfinite recursion on $\pi$ by using the method in~\autoref{4.3}, as represented in \autoref{f6578}. The induction basis is clear and the limit step follows by \autoref{k42}, so we only need to describe what we do in the successor step so that \autoref{k39} can be applied. 
    \begin{figure}[ht]
    \begin{center}
        \begin{tikzpicture}
            \draw[line width=0.2mm, color=black] (1, 0) -- (4, 0);

            \draw[line width=0.2mm, color=bluet] (4, 0) -- (7, 0);

            \draw[line width=0.2mm, color=black] (7, 0) -- (10, 0);

            \draw[line width=0.2mm, color=bluet] (10, 0) -- (13, 0);

            \node at (1.3, 0.3) {$S_{0}$};
            \node at (2.6, -0.4) {$\dot{\bbQ}_{\xi} \coloneqq \bbC$};

            \node[bluet] at (4.3, 0.3) {$S_{1}$}; 
            \node[bluet] at (5.5, -0.4) {$\dot{\bbQ}_{\xi} \coloneqq \bbL \mathbb{O} \mathbb{C}^{\dot{N}_{\xi}}$};

            \node at (7.3, 0.3) {$S_{2}$};
            \node at  (8.5, -0.4){$\dot{\bbQ}_{\xi} \coloneqq \bbD^{\dot{N}_{\xi}}$};

            \node[bluet] at (10.3, 0.3) {$S_{3}$};
            \node[bluet] at (11.5, -0.4) {$\dot{\bbQ}_{\xi} \coloneqq \bbB^{V^{\bbP_{\xi}^{-}}}$};

            \node at (13.7, 0) {$\chi + \chi$};
        \end{tikzpicture}   
    \end{center}
    \caption{Graphic representation of the iteration (considering that $S_i$ for $1\leq i \leq 3$ are cofinal in $\chi+\chi$).}
    \label{f6578}
    \end{figure}
    
    Assume we have constructed our $\calC$-iteration 
    $\bft_{\xi} =\langle \bbP_{\beta}, \dot{\bbQ}_{\alpha}, \bbP_{\alpha}^{-}, \vec{Q}_{\alpha}, \theta_{\alpha},  \vec{\Xi}_{\beta} \colon \alpha < \xi,  \beta \leq \xi \rangle$  
    up to $\xi<\pi$. We also assume that $\vert\Por_\xi\vert\leq\chi$, which implies $\Vdash_\xi \text{``} \cfrak\leq \chi$''. It will be clear that $\vert\Por_{\xi+1}\vert\leq\chi$ (and also for the limit step). We split into three possible cases:
    
    \noindent\textbf{Case 1:} $\xi \in S_0$. Set $\bbP_{\xi}^{-} \coloneqq \bbP^\bullet_{\xi}$, $\dot{\bbQ}_{\xi} \coloneqq \bbC$,  $\bbP_{\xi + 1} \coloneqq \bbP_{\xi} \ast \dot{\bbQ}_{\xi}$, $\theta_{\xi} \coloneqq \aleph_{0}$ and, for any $s \in \bbC$, let $\dot{Q}_{s, 0}^{\xi} \coloneqq \{ s \}.$ Note that $\langle \bbP_{\beta},  \dot{\bbQ}_{\alpha},  \bbP_{\alpha}^{-},  \vec{Q}_{\alpha},  \theta_{\alpha} \colon \beta\leq \xi+1,\ \alpha < \xi + 1 \rangle$ is a $\kappa$-prepared iteration, so by \autoref{k39}, we can find a sequence of names of fams $\vec{\Xi}_{\xi+ 1} \coloneqq \langle \dot{\Xi}_{\xi + 1}^{\bfc} \colon \bfc \in \calC \rangle$ such that $ \bft_{\xi + 1} \coloneqq  \langle \bbP_{\beta},  \dot{\bbQ}_{\alpha},  \bbP_{\alpha}^{-},  \vec{Q}_{\alpha},  \theta_{\alpha},   \vec{\Xi}_{\beta} \colon \beta \leq \xi + 1,\ \alpha < \xi +1\rangle$ is a $\calC$-iteration.
    
    \noindent\textbf{Case 2:} $\xi \in S_1 \cup S_2$. 

    Denote $\nu_1\coloneqq \theta$ and $\nu_2\coloneqq \kappa$. 
    Fix $i\in\{1,2\}$ and assume $\xi\in S_i$. Since $\bbP_\xi$ has the ccc, $\Vdash_\xi \text{``} \cfrak =\chi$'' and $\vert (\baire)^{\Por^\bullet_\xi}\vert = \chi^{\aleph_0} =\chi$, where $(\baire)^{\Por^\bullet_\xi}$ denotes the set of nice $\Por^\bullet_\xi$-names of members of $\baire$. Also $\cof([\chi]^{<\nu_i})=\chi$, so pick a cofinal family $\{ E_{\xi, \beta} \colon \beta < \chi \}$ in $\left[(\baire)^{\Por^\bullet_\xi}\right]^{<\nu_i}$. For each $\beta<\chi$, let $\dot F_{\xi,\beta}$ be a $\Por^\bullet_\xi$-name of the set in $[\baire]^{<\nu_i}$ obtained from $E_{\xi,\beta}$. Then, $\Por_\xi$ forces that $\{\dot F_{\xi,\beta}\colon \beta<\chi\}$ is cofinal in $[\baire]^{<\nu_i}$.
    
    This construction guarantees that we already have this enumeration for any $\gamma < \xi$ in $S_i$. Now, define $\dot{F}_{\xi} \coloneqq \dot{F}_{\bfb_{i}(\xi)}$, $\theta_\xi\coloneqq |E_{\bfb_i(\xi)}|<\nu_i$ (so $\Vdash_\xi \text{``} |\dot F_\xi|\leq\theta_\xi$'') and pick a $\bbP^\bullet_{\xi}$-name of a transitive model $\dot{N}_{\xi}$ of (a finite large enough fragment of) $\ZFC$ of size $\theta_{\xi}$  such that $\forces_{\xi}$ ``$\dot{F}_{\xi} \subseteq \dot{N}_{\xi}$''. 
    We then define $\bbP_{\xi}^{-} \coloneqq \bbP_{\xi}^\bullet$ and the next stage $\bbP_{\xi + 1} \coloneqq \bbP_{\xi} \ast \dot{\bbQ}_{\xi}$, where
    \[\dot{\bbQ}_{\xi}\coloneqq\left\{\begin{array}{ll}
           \Loc^{\dot{N}_{\xi}}    & \text{if $i=1$,}\\[1ex]
           \bbD^{\dot{N}_{\xi}}  & \text{if $i=2$.} 
        \end{array}\right.\]
    Each $\dot Q^\xi_{\zeta,0}$ for $\zeta<\theta_\xi$ is defined as a singleton for some enumeration of $\Qnm_\xi$ in $V_\xi$, so $\Vdash_\xi \text{``}\Qnm_\xi = \bigcup_{\zeta<\theta_\xi}Q^\xi_{\zeta,0}$''. We can then find
    $ \vec{\Xi}_{\xi} \coloneqq \langle \dot{\Xi}_{\xi}^{\bfc} \colon \bfc \in \calC \rangle$ as in the conclusion of \autoref{k39}. Consequently, we have that 
    $  \bft_{\xi + 1} \coloneqq  \langle \bbP_{\beta}, \dot{\bbQ}_{\alpha}, \bbP_{\alpha}^{-},  \vec{Q}_{\alpha},  \theta_{\alpha},  \vec{\Xi}_{\beta} \colon \beta \leq \xi + 1,\ \alpha < \xi +1 \rangle$ is a $\calC$-iteration.

    \noindent\textbf{Case 3:} $\xi \in S_3$. Denote by $\Omega_{\bar \delta}^{\Por^\bullet_\xi}$ the collection of nice $\bbP^\bullet_\xi$-names of members of $\Omega_{\bar \delta}$. As in Case 2, find a cofinal family $\{\dot E_{\xi,\beta} \colon \beta < \chi \}$ in $[\Omega_{\bar \delta}^{\Por^\bullet_\xi}]^{<\lambda}$, and let $\dot F_{\xi,\beta}$ be the $\Por^\bullet_\xi$-name of the member of $[\Omega_{\bar\delta}]^{<\lambda}$ produced by $E_{\xi,\beta}$.
    Set $\dot F_\xi \coloneqq \dot F_{\bfb_{3}(\xi)}.$ Since all elements of $E\coloneqq E_{\bfb_{3}(\xi)}$ are nice $\bbP^\bullet_{\xi}$-names, each $\dot{\bar a}\in E$ depends on some countable sequence $\la A_{\dot{\bar a},j}\colon j<\omega\ra$ of (countable) maximal antichains in $\bbP^\bullet_{\xi}$. Set $Q\coloneqq \bigcup_{\dot{\bar a}}\bigcup_j A_{\dot{\bar a},j}$, so $|Q|<\lambda$. 
            
    Next, by applying \autoref{a70} we can find a forcing notion $\bbP_{\xi}^{-} \subsetdot \bbP_{\xi}^\bullet$ 
    containing $Q$ such that $\vert \bbP_{\xi}^{-} \vert \leq \max \{2, \vert Q \vert, \vert \calC\vert \}^{\aleph_{0}}$ and $\dot\Xi^\bfc_\xi\frestr\pts(\omega)\cap V^{\Por^-_\xi}$ has a $\Por^-_\xi$-name $\dot\Xi^{\bfc,-}_\xi$ for each $\bfc\in\calC$. Notice that $\vert \bbP_{\xi}^{-} \vert < \lambda$ 
    because $\lambda$ is $\aleph_{1}$-inaccessible, $ \vert \calC\vert  \leq \mu^{\aleph_{0}} < \lambda$ and $\vert Q \vert < \lambda$.  
    We then define the next step of the iteration by $\dot{\bbQ}_{\xi} \coloneqq \bbB^{V^{\bbP_{\xi}^{-}}}$,  $\theta_{\xi} \coloneqq \aleph_{0}$, $\bbP_{\xi + 1} \coloneqq \bbP_{\xi} \ast \dot{\bbQ}_{\xi}$, and let $ \vec{Q}_{\xi} \coloneqq \la\dot{Q}_{t, \varp}^{\xi} \colon t\in{}^{<\omega}2,\ \varp\in(0,1)_\QQQ\ra$ be a sequence of $\Por^-_\xi$-names witnessing that random forcing is strongly $\calY_{**}$-linked in $V^{\Por^-_\xi}$, which exists by~\autoref{i15.2}. We can find $ \vec{\Xi}_{\xi} \coloneqq \langle \dot{\Xi}_{\xi}^{\bfc} \colon \bfc \in \calG \rangle$ as in the conclusion of \autoref{k39}. Once again, we obtain a $\calC$-iteration $ \bft_{\xi + 1} \coloneqq  \langle \bbP_{\beta}, \, \dot{\bbQ}_{\alpha}, \, \bbP_{\alpha}^{-}, \, \vec{Q}_{\alpha}, \, \theta_{\alpha}, \,  \vec{\Xi}_{\beta} \colon \beta \leq \xi + 1,\ \alpha < \xi +1 \rangle .$

    This completes the construction of $\bft \coloneqq \bft_{\pi}$. Denote $\bbP \coloneqq \bbP_{\pi}$, It is clear that $\bbP$ is $\rm{ccc}$ and, by~\autoref{k62}, $\bbP$ is $\kappa$-$\Fr$-Knaster. Also, $\bbP$ forces $\cfrak = \chi$.    
    
    We start by proving that $\Por$ forces $\gb = \non(\Ecal) =\kappa$. To force $\gb \geq \kappa$, it suffices to prove that, in $V_{\pi},$ any $F \subseteq \baire$ of size ${<} \kappa$ is bounded. Indeed, assume that $F \subseteq \baire$ has size ${<} \kappa$ in $V_{\pi}$. Since $\bbP$ is ccc and $\cf(\chi)\geq\kappa$ (which is implied by $\cof([\chi]^{<\kappa})=\chi$), there exists some $\xi \in S_2$ such that $F \in V_{\xi}.$ As a consequence, by the iteration construction, $F \subseteq F_{\xi, \beta}$ for some $\beta < \lambda,$ hence $F \subseteq F_{\alpha},$ where $\alpha \coloneqq \bfb_{2}^{-1}(\xi, \beta) \geq \xi$ and $\alpha\in S_2$. Then, the generic real added by $\Qor_\alpha$ dominates $F_\alpha$, and also $F$. The inequality 
    inequality $\non(\Ecal) \leq \kappa$ is forced by~\autoref{a69}. Since $\ZFC$ proves $\min\{\bfrak,\non(\calN)\}\leq\non(\Ecal)$ (\cite[proof of Thm.~4.2]{BS1992}) and below we show that $\non(\Mcal)\leq\cov(\Mcal)$ is forced, we force $\bfrak\leq \non(\Ecal)$, so the equality $\gb = \non(\Ecal) =\kappa$ is settled.

    Similarly, it can be proved that $\Por$ forces $\theta\leq \bfrak(\Lc(\omega,\id))=\add(\Nwf)$. On the other hand, since $\Por$ is obtained by a finite support iteration of $\theta$-$\Lc^*$-good forcing notions (see \autoref{smallgoodness} and ~\autoref{ExmPrs}~\ref{ExmPrsd}), $\Por$ forces $\add(\Nwf)\leq\theta$ by using~\autoref{Comgood}.

    To force $\cov(\calN) \geq \lambda$, it is enough to show that, in $V_\pi$, any subset of $\Omega_{\bar\delta}$ of size ${<} \lambda$  is $\Cn_{\bar\delta}$-bounded (see \autoref{a40}). Let $F$ be one of such families. Then, there is some $\xi\in S_3$ such that $F\in V_\xi$, so $F\subseteq F_{\xi,\beta}$ for some $\beta<\xi$. Pick $\alpha \in S_3$ such that $\bfb_{3}(\alpha) = (\xi, \beta)$. Then $F\subseteq F_\alpha=F_{\bfb_{B}(\alpha)}$ and, by the construction, $F_\alpha\in V^{\Por^-_\alpha}$. Let $r_\alpha\in\cantor$ be the random real over $V^{\Por^-_\alpha}$ added by $\Qor_\alpha$. Then, $r_\alpha\notin N[\bar a]$ for all $\bar a\in F_\alpha$, so $F$ is $\Cn_{\bar \delta}$-bounded by $r_\alpha$.

    To prove that $\Por$ forces that $\cov(\calN) \leq \lambda$ it suffices to verify the hypothesis of \autoref{a60}. First work in $V_\chi$. 
    For $\alpha<\chi$ let $\bar a^\alpha\in\Omega_{\bar\delta}$ be a Cohen real added by $\Qor_\alpha = \Cor$. Since $\bbP_{\chi}$ is forcing equivalent with $\Cor_\chi$, $\{\bar a^\alpha\colon \alpha<\chi\}$ is strongly $\aleph_1$-$\Cn_{\bar\delta}$-unbounded (in $V_\chi$), which implies strongly $\lambda$-$\Cn_{\bar\delta}$-unbounded. Hence, \ref{a60a} of \autoref{a60} holds. On the other hand, \ref{a60c} of \autoref{a60} is clear and \ref{a60b} follows by~\autoref{a64}. Since the hypothesis of~\autoref{a60} are fulfilled,  $\Por$ forces that $\{\bar a^\alpha\colon \alpha<\chi\}$ is strongly $\aleph_1$-$\Cn_{\bar\delta}$-unbounded, which implies    
    $\cov(\calN) \leq \lambda$.

    It remains to settle the values of $\non(\Mwf)$ and $\cov(\Mwf)$. Let $\lambda^*$ be the smallest regular cardinal ${\geq}\lambda$. We get that $\lambda\leq\lambda^*\leq\lambda^+$ and $\lambda^*\leq\chi$ (in the case $\lambda=\chi$, since $\cof([\chi]^{<\lambda})=\chi$ implies $\cf(\chi)\geq \lambda$, $\lambda$ is regular and $\lambda^*=\lambda$). By~\autoref{Comgood} applied to $(R,\theta)=(\Ed,\lambda^*)$ (see~\autoref{ExmPrs}~\ref{ExmPrsa}), we obtain that $\Por$ forces $\non(\Mwf)\leq\lambda^*$ and $\cov(\calM) \geq \chi$. Since $\forces_\Por \text{``}\cfrak=\chi$'', it is forced that $\cov(\calM)=\chi=\cfrak$; on the other hand, $\Por$ forces $\lambda =\cov(\calN) \leq \non(\calM) \leq \lambda^*\leq\lambda^+$. 
    When $\lambda$ is regular, $\lambda^*=\lambda$, so $\non(\calM)=\lambda$ is forced. This settles~\ref{a95:a} and~\ref{a95:b}.
    
    For~\ref{a95:c}, assume that $\cf(\lambda) < \kappa$. By~\ref{a95:a} we know that $\lambda\leq\non(\Mwf)\leq\lambda^+$, so $\Por$ forces $\non(\Mwf)=\lambda^+$ because, in the other case, $\kappa = \add(\calM) \leq \cf(\non(\calM)) = \cf(\lambda) < \kappa,$ which is not possible.
    
    Finally, the case $\kappa = \lambda$ is known from \cite[Thm.~2]{MejMatrix}, but the argument above is valid for this case under the following changes: $\calG\coloneqq{}^{\pi} \kappa^-$, $U^\tbf\coloneqq \pi$ (so $F^\tbf = \emptyset$) and proceed in Case 3 ($\xi\in S_3$) exactly as in Case 2 (i.e.\ $\Qnm_\xi$ has the form $\Bor^{\dot N_\xi}$).
\end{PROOF}
    
We can prove that the hypothesis ``$\log\chi <\lambda$'' in the previous theorem can be changed by a requirement more compatible with GCH. The proof is inspired from~\cite[Lem.~1.33]{GKS}.
       
\begin{theorem}\label{a96}
    In \autoref{a95}, the assumption ``\/$\log\chi <\lambda$'' can be replaced by ``either $2^{\aleph_0}< \kappa$ or $2^{<\kappa}<\lambda$''.
\end{theorem}
\begin{PROOF}[\textbf{Proof}]{\ref{a96}}
    Let $\mu$ be a regular uncountable cardinal such that $\mu\leq\chi$. Recall the poset $\Fn_{<\mu}(\chi,2)$ of partial functions $\chi\to 2$ with domain of size ${<}\mu$. Without any requirements (like GCH), this poset is ${<}\mu$-closed, has the $(2^{<\mu})^+$-cc and forces $|{}^{<\mu}2|=\mu$ and $|\chi|\leq 2^\mu$, which means that the cardinals in $(\mu,2^{<\mu}]$ are collapsed to $\mu$ when $V\models \mu<2^{<\mu}$.
    
    Under the assumptions of \autoref{a95} (excluding $\log\chi<\lambda$), consider the case when $\kappa<\lambda$ and $\lambda$ is $\aleph_1$-inaccessible. In this theorem, we have assumed $\log\chi <\lambda$ to get a complete set of guardrails of size ${<}\lambda$ (as a consequence of Engelking's and Kar\l owicz's Theorem), but now we plan to go along with the iteration construction without this hypothesis. 
    The trick consists of constructing a ccc $\kappa$-prepared finite support iteration in the ground model $V$ such that it is a $\calC$-iteration in the $\Fn_{<\mu}(\chi,2)$-extension for some suitable bedrock $\calC$ that lives in this generic extension. This means that the sequences of fams are constructed in the $\Fn_{<\mu}(\chi,2)$-extension. Of course, we will need some requirements for $\mu$ to conclude that the iteration constructed in $V$ forces what we want, but these are given along with the construction. We take advantage of $|\chi|\leq 2^\mu$ in the $\Fn_{<\mu}(\chi,2)$-extension to guarantee that $\calC$ is small enough.
    
    So let $\Por^*\coloneqq\Fn_{<\mu}(\chi,2)$. In $V^{\Por^*}$, we have that $|\chi|\leq 2^\mu$, so by \autoref{k23} there is some complete set $\calG$ of guardrails for $\chi+\chi$ over $\nu\coloneqq \max\{\kappa,\mu\}$ of size ${\leq}|\nu|^{\aleph_0}.$ Fix a free fam $\Xi_0$ on $\pts(\omega)$  (which could come from the ground model) and define $\calC$ as the set of all tuples of the form $(\Xi_0,\bar I^j,g)$ for $0<j<\omega$ and $g\in\calG$, exactly as in the proof of \autoref{a95}. Clearly, $|\calC|\leq|\nu|^{\aleph_0}$.
    
    Back in $V$, let $\dot\calG$ and $\dot\calC$ be $\Por^*$-names of the objects constructed above. We now construct, in $V$, a ccc $\kappa$-prepared iteration $\tbf$ of length $\chi+\chi$ such that $\Por^*$ forces that it is a $\dot\calC$-iteration after plugging names of fams $\la \dot \Xi^\bfc_\alpha\colon \alpha\leq\chi+\chi,\ \bfc\in\dot\calC\ra$, i.e.\ each $\dot\Xi^\bfc_\alpha$ is a $\Por^*\times \Por_\alpha$-name of a fam on $\pts(\omega)$. We partition $\chi+\chi = S_0\cup S_1\cup S_2\cup S_3$, set $S_0\coloneqq\chi$ and $F^\tbf\coloneqq S_3$, and fix book-keeping functions $\bfb_i$ for $1\leq i\leq 3$, exactly as in the proof of \autoref{a95}.
    
    In the case $\xi\in U^\tbf = S_0\cup S_1\cup S_2$ we proceed exactly as in the proof of \autoref{a95}, so $\Por^-_\xi = \Por_\xi^\bullet$. Also, $\Por^*\times\Por_\xi$ forces that $\Qnm_\xi$ is uniformly and strongly $|\theta_\xi|$-$\calY_{*}$-linked witnessed by singletons (recall that, by Easton's Lemma, $\Por_\xi$ is ccc in $V^{\Por^*}$). Then, by \autoref{k39}, there is a $\Por^*\times\Por_{\xi+1}$-name $\dot{\vec{\Xi}}_{\xi+1}$ of a sequence $\Seq{\dot\Xi^\bfc_{\xi+1}}{\bfc\in\dot\calC}$ of fams on $\pts(\omega)$ such that $\Por^*$ forces that the iteration up to $\xi+1$ is a $\dot\calC$-iteration after including them.
    
    The case $\xi\in S_3$ is more delicate. We let $\dot F_\xi$ be a $\Por_\xi$-name (chosen by the book-keeping function $\bfb_3$) and $Q\subseteq \Por_\xi^\bullet$ of size ${<}\lambda$ exactly as in the proof of \autoref{a95}. This time, we apply \autoref{a70g} to find some $\Por^-_\xi\subsetdot \Por^\bullet_\xi$ of size ${\leq}\max\{|Q|,\nu,2^{<\mu}\}^{\aleph_0}$ such that, for each $\Por^*$-name $\dot\bfc$ of a member of $\dot\calC$, $\dot\Xi^{\dot \bfc}_\xi\frestr\pts(\omega)\cap V^{\Por^*\times\Por^-_\xi}$ has a $\Por^*\times\Por^-_\xi$-name. We let $\Qnm_\xi$ be a $\Por^-_\xi$-name of random forcing and let $\la\dot Q_{t,\varp}\colon t\in{}^{<\omega}2,\ \varp\in(0,1)_\QQQ\ra$ be the sequence of $\Por^-_\xi$-names witnessing that random forcing is uniformly $\sigma$-$\calY_{**}$-linked. By Easton's Lemma, $\Por^-_\xi$ forces that $\Por^*$ is ${<}\aleph_1$-distribute, so $\Por^*\times\Por^-_\xi$ forces that $\Qnm_\xi$ is the random forcing from $V^{\Por^*\times\Por^-_\xi}$, and so it is uniformly $\sigma$-$\calY_{**}$-linked witnessed by the same sets. Therefore, \autoref{k39} can be applied to find a $\Por^*\times\Por_{\xi+1}$-name $\dot{\vec{\Xi}}_{\xi+1}$ of a sequence of fams on $\pts(\omega)$ as in the previous case. At the end of the proof, we need that $|\Por^-_\xi|<\lambda$ (for $\xi\in S_3$) to guarantee that $\Por_\xi$ forces $|\Qnm_\xi|<\lambda$ (and be able to apply \autoref{a60} at the very end). Since we have $|\Por^-_\xi|\leq \max\{|Q|,\nu,2^{<\mu}\}^{\aleph_0}$, \emph{we require that $2^{<\mu}<\lambda$} to get $|\Por^-_\xi|<\lambda$ (recall that $\lambda$ is $\aleph_1$-inaccessible).
    
    The limit step $\alpha\leq\chi+\chi$ is clearly defined in $V$. In $V^{\Por^*}$, \autoref{k42} can be easily applied to find the names of sequences of fams that we need to get a $\dot\calC$-iteration up to $\alpha$. This completes the construction.
    
    Under the assumption that $2^{<\mu}<\lambda$, the arguments in the proof of \autoref{a95} that do not use fams can be performed verbatim (in $V$), so it remains to show which requirements must be satisfied by $\mu$ so that $\Por\coloneqq\Por_{\chi+\chi}$ forces (over $V$) that $\cov(\Ncal)\leq\lambda$ and $\non(\Ecal)\leq \kappa$.
    
    Work in $V^{\Por^*}$. If $\kappa$ is not collapsed, i.e.\ $V\models$``$\kappa\leq\mu$ or $2^{<\mu}<\kappa$", then we have a $\kappa$-$\calC$-iteration, so $\Por$ forces that the Cohen reals added by $\Por_\chi$ form a $\kappa$-anti-Bendixson family and a strongly $\kappa$-$\Ce$-unbounded family by \autoref{a64} and~\ref{a69}, respectively. However, this family of Cohen reals live in $V^{\Por_\chi}$, thus, by downward absoluteness between $V^{\Por}$ and $V^{\Por^*\times\Por}$, in $V$, $\Por$ forces the same. Therefore, in $V$, $\Por$ forces $\non(\Ecal)\leq\kappa$ and, by \autoref{a60} (applied in $V^{\Por_\chi}$), $\cov(\Ncal)\leq\lambda$.

    In conclusion, we need that $2^{<\mu}<\lambda$ and that either $\kappa\leq\mu$ or $2^{<\mu}<\kappa$ to succeed with the previous argument. In the first case, $\kappa\leq\mu$ and $2^{<\mu}<\lambda$ imply $2^{<\kappa}<\lambda$, so it is enough to use $\mu=\kappa$; in the second case, $2^{<\mu}<\kappa$ implies that $2^{\aleph_0}<\kappa$, so it is enough to use $\mu=\aleph_1$. For this reason, it is enough to assume that either $2^{<\kappa}<\lambda$ or $2^{\aleph_0}< \kappa$ to succeed with the proof of the theorem.
\end{PROOF}

The following well-known result due to Shelah is a direct consequence of \autoref{a96}.

\begin{corollary}[{\cite[Thm.~3.5]{Sh00}}]\label{a88}
    $\mathrm{Con}(\ZFC) \Rightarrow \mathrm{Con}(\ZFC + \cf(\cov(\calN)) = \aleph_{0}).$ 
\end{corollary}

The advantage of \autoref{a96} is that the finite support iteration of ccc forcings for this consistency result can be constructed in a model of GCH, which is not possible with \autoref{a95} if GCH is added to its assumptions, $\log\chi<\lambda\leq \chi$ implies that $\lambda=\chi$ is a successor cardinal. Under \autoref{a96}, by assuming GCH in the ground model, $\lambda$ can be any uncountable cardinal of countable cofinality, $\aleph_1\leq\theta\leq\kappa<\lambda$ with $\theta$ and $\kappa$ regular, and $\chi>\lambda$ such that $\cf(\chi)>\lambda$ (in particular, $\lambda=\aleph_\omega$ can be used).


\section{Questions}

Apart from Fremlin’s conjecture mentioned at the beginning of~\autoref{SecIntro},  solved by Shelah in~\cite{Sh00}, the following is of interest as well:

\begin{question}
    Can $\cov(\Ewf)$ have countable cofinality?  
\end{question}

Miller (see~\cite[Thm. 5.1.18]{BJ}) proved that $\dfrak\leq\cov(\Ewf)$ implies $\cf(\cov(\Ewf))>\omega$. Hence, to provide a positive answer to this question, it must be forced $\dfrak < \cov(\Ewf)$, which implies $\cov(\Ewf) = \cov(\Nwf)$ (see~\cite[Thm.~2.6.9]{BJ}, which was originally proved by Bartoszy\'nski and Shelah~\cite{BS1992}). Hence, it would be necessary to force $\cf(\cov(\Nwf))=\omega$ without increasing $\dfrak$, e.g.~by using $\baire$-bounding forcing.

In case that $\lambda$ is regular in \autoref{a96}, as in~\cite{KST} it is possible to force $\non(\Mcal)$ to be larger than $\lambda^+$, particularly, any desired $\aleph_1$-inaccessible regular value between $\lambda$ and $\chi$. However, we do not know how to do this when $\lambda$ is singular.

\begin{question}
    Can the left side of Cicho\'n's diagram be separated with the value of $\cov(\Nwf)$ singular and $\cov(\Ncal)^+<\non(\Mcal)$? Can this be done for Cicho\'n's maximum as well?
\end{question}

It is well-known that $\add(\Ncal)$, $\add(\Mcal)$ and $\bfrak$ are regular, and it is known how to force, using large cardinals, Cicho\'n's maximum with the cardinals of the right side singular. It is not known how to force this without using large cardinals.

On the other hand, it is unknown how to force Cicho\'n's maximum with $\cov(\Ncal)$ or $\non(\Mcal)$ singular. Moreover, 
it is not known how to force $\non(\Mwf)$ singular with finite support iteration techniques.

\begin{question}
    Is it consistent with ZFC that $\non(\Mwf)<\cov(\Mwf)$ with $\non(\Mwf)$ singular? More generally: Can the left side of Cicho\'n's diagram be separated with the value of $\non(\Mwf)$ singular?
\end{question}

It is also unknown how to force $\aleph_1<\cov(\Mcal)<\non(\Mcal)$ with $\non(\Mcal)$ singular.

{\small

\bibliography{appl}
\bibliographystyle{alpha}
}

\end{document}